\documentclass[11pt, oneside]{amsart}
\usepackage[a4paper,inner=3.0cm,outer=2.5cm,top=2.8cm,bottom=2.8cm]{geometry}
\usepackage[T1]{fontenc}
\usepackage[english]{babel}
\usepackage{lmodern}
\usepackage{amsmath,amssymb,amsthm,mathtools}
\usepackage{tikz-cd}
\usetikzlibrary{arrows.meta}
\usepackage{enumitem}
\usepackage{needspace}
\usepackage{xcolor}
\usepackage{xurl}
\usepackage[pagebackref]{hyperref}
\usepackage[nameinlink,noabbrev]{cleveref}
\usepackage{aliascnt}
\usepackage[protrusion=true,expansion=false]{microtype}
\hypersetup{
  colorlinks=true,
  linkcolor=blue!45!black,
  citecolor=blue!45!black,
  urlcolor=blue!55!black,
  pdftitle={Gersten injectivity with invertible finite coefficients},
  pdfauthor={Niels Feld},
  bookmarksnumbered=true
}
\renewcommand*{\backrefalt}[4]{%
  \ifcase #1\relax
  \or {\footnotesize\textup{[p.~#2]}}%
  \else {\footnotesize\textup{[pp.~#2]}}%
  \fi}

\newtheorem{theorem}{Theorem}[section]
\newaliascnt{proposition}{theorem}
\newtheorem{proposition}[proposition]{Proposition}
\aliascntresetthe{proposition}
\newaliascnt{lemma}{theorem}
\newtheorem{lemma}[lemma]{Lemma}
\aliascntresetthe{lemma}
\newaliascnt{sublemma}{theorem}
\newtheorem{sublemma}[sublemma]{Sublemma}
\aliascntresetthe{sublemma}
\newaliascnt{corollary}{theorem}
\newtheorem{corollary}[corollary]{Corollary}
\aliascntresetthe{corollary}
\theoremstyle{definition}
\newaliascnt{definition}{theorem}
\newtheorem{definition}[definition]{Definition}
\aliascntresetthe{definition}
\newaliascnt{construction}{theorem}
\newtheorem{construction}[construction]{Construction}
\aliascntresetthe{construction}
\newaliascnt{notation}{theorem}
\newtheorem{notation}[notation]{Notation}
\aliascntresetthe{notation}
\theoremstyle{remark}
\newaliascnt{remark}{theorem}
\newtheorem{remark}[remark]{Remark}
\aliascntresetthe{remark}
\crefname{theorem}{Theorem}{Theorems}
\crefname{proposition}{Proposition}{Propositions}
\crefname{lemma}{Lemma}{Lemmas}
\crefname{corollary}{Corollary}{Corollaries}
\crefname{definition}{Definition}{Definitions}
\crefname{sublemma}{Sublemma}{Sublemmas}
\crefname{construction}{Construction}{Constructions}
\crefname{notation}{Notation}{Notations}
\crefname{remark}{Remark}{Remarks}
\crefname{section}{Section}{Sections}
\crefname{subsection}{Section}{Sections}
\crefname{appendix}{Appendix}{Appendices}
\Crefname{appendix}{Appendix}{Appendices}
\DeclareMathOperator{\Spec}{Spec}
\DeclareMathOperator{\Frac}{Frac}
\DeclareMathOperator{\cofib}{cofib}

\newcommand{\Z}{\mathbf Z}
\newcommand{\stackstag}[1]{\href{https://stacks.math.columbia.edu/tag/#1}{\textup{Tag~#1}}}
\newcommand{\stmtlink}[2]{\href{#1}{#2}}
\newcommand{\arfiv}[3]{\stmtlink{https://ar5iv.labs.arxiv.org/html/#1\##2}{#3}}
\newcommand{\arxivhtml}[3]{\stmtlink{https://arxiv.org/html/#1\##2}{#3}}
\numberwithin{equation}{section}
\setlist[enumerate,1]{label=\textup{(\roman*)},leftmargin=*,itemsep=3pt}
\title{Gersten injectivity with invertible finite coefficients}
\author[Niels Feld]{Niels Feld}
\address{IRMAR, University of Rennes, 263 Avenue G\'en\'eral Leclerc, 35000 Rennes, France}
\email{niels.feld@univ-rennes.fr}
\date{26 November 2026}

\begin{document}
\begin{abstract}
	Let $R$ be a noetherian regular local ring with fraction field $F$.
	For every integer $N\geq 2$ invertible in $R$ and every integer $n$,
	we prove that the restriction map
	\[
	K_n(R;\mathbf{Z}/N)\to K_n(F;\mathbf{Z}/N)
	\]
	is injective, with finite coefficients taken on the nonconnective
	algebraic $K$-theory spectrum.
	We also prove injectivity of restriction to $F$ in \'etale
	cohomology with finite locally constant abelian coefficients of
	invertible order, in every cohomological degree.
	In contrast, for every prime $\ell$, we construct an excellent
	discrete valuation ring containing $\mathbf{Q}$ for which
	generic restriction in continuous $\ell$-adic cohomology
	fails to be injective in degree two.
\end{abstract}
\maketitle
\tableofcontents

\section{Introduction}\label{sec:introduction}

Gersten's conjecture seeks to describe the algebraic $K$-theory of a
regular local ring in terms of the $K$-theory of its residue fields.
For a scheme $Y$, let $K(Y)$ denote the nonconnective algebraic
$K$-theory spectrum of perfect complexes on $Y$, and let
$K_n(Y)=\pi_n K(Y)$. For a ring $A$, write $K(A)=K(\Spec A)$.
Let $R$ be a noetherian regular local ring with fraction field $F$,
put $X=\Spec R$, and write $X^{(p)}$ for the set of points of
codimension $p$ in $X$. The filtration by codimension of support, together with localization,
d\'evissage and resolution, gives the coniveau spectral sequence in
algebraic $K$-theory \cite[\S7]{Quillen}. Its rows, augmented by the
$K$-groups of $R$, are the Gersten complexes. In degree $n$, the
complex is
\begin{multline}\label{eq:intro-gersten-complex}
 0\longrightarrow K_n(R)\longrightarrow K_n(F)
 \longrightarrow\bigoplus_{x\in X^{(1)}}K_{n-1}(\kappa(x))\\
 \longrightarrow\bigoplus_{x\in X^{(2)}}K_{n-2}(\kappa(x))
 \longrightarrow\cdots .
\end{multline}
The terms in codimension greater than $\dim R$ are zero. Gersten's
conjecture predicts exactness of \eqref{eq:intro-gersten-complex} for
every $n\geq0$. In Quillen's formulation
\cite[\S7, Conjecture~5.10]{Quillen}, if $\mathcal M^p(R)$ denotes
the category of finitely generated $R$-modules with support in
codimension at least $p$, the maps
$K_n(\mathcal M^{p+1}(R))\to K_n(\mathcal M^p(R))$ should vanish for
all $p,n\geq0$. Localization and d\'evissage show that these vanishing
statements imply exactness of the Gersten complexes. In particular,
resolution and localization identify the assertion for $p=0$ with
injectivity of the first map $K_n(R)\to K_n(F)$.

\paragraph{\textnormal{\bfseries Finite coefficients.}}
For any scheme $Y$ and any natural number $N\geq2$ put
\[
 K(Y;\Z/N)=\cofib\bigl(N:K(Y)\longrightarrow K(Y)\bigr),
 \qquad K_n(Y;\Z/N)=\pi_n K(Y;\Z/N).
\]
This convention gives the natural exact sequence
\begin{equation}\label{eq:coefficient-sequence}
 0\longrightarrow K_n(Y)/N
   \longrightarrow K_n(Y;\Z/N)
   \longrightarrow K_{n-1}(Y)[N]\longrightarrow0,
\end{equation}
where $M[N]$ denotes the kernel of multiplication by $N$ on
an abelian group $M$. The map to $K_{n-1}(Y)[N]$ in \eqref{eq:coefficient-sequence}
is induced by the connecting map of the cofiber sequence.
We write $K^{\mathrm{cn}}(Y)=\tau_{\geq0}K(Y)$ for the connective
cover and take its finite coefficients by the same construction, with
$K^{\mathrm{cn}}_n(Y;\Z/N)=\pi_n(K^{\mathrm{cn}}(Y)/N)$.
On an affine regular noetherian scheme, resolution and vanishing
of the negative $K$-groups identify $K^{\mathrm{cn}}$ with $K$
\cite[Proposition~6.8(b)]{TT}.
We use $\mathbf N=\{0,1,2,\ldots\}$.

Replacing the $K$-theory spectra in the coniveau construction by their
cofibers of multiplication by $N$ gives the Gersten complex with
$\Z/N$-coefficients. Its terms after the augmentation are
$\bigoplus_{x\in X^{(p)}}K_{n-p}(\kappa(x);\Z/N)$, $p\geq0$, and its
first map is $K_n(R;\Z/N)\to K_n(F;\Z/N)$. By
\eqref{eq:coefficient-sequence}, the coefficient groups need not be
quotients of the integral $K$-groups: this complex is not in general
obtained by reducing \eqref{eq:intro-gersten-complex} modulo $N$.

The distinction between integral and finite coefficients, and between
invertible and noninvertible coefficient orders, is essential. In
mixed characteristic, recent preprints
\cite[Theorem~1.1 in each]{OpenAIGersten3,OpenAIGersten5} report
counterexamples to the unrestricted integral conjecture in degrees
three and five, respectively, on ramified regular local surfaces of
residue characteristic five. For a complete
regular local surface of residue characteristic three the first map
with $\Z/3$-coefficients is not injective
\cite[Theorem~1.1]{FeldCounterexample}. \Cref{thm:gersten-injectivity}
proves injectivity with finite coefficients of invertible order for
arbitrary noetherian regular local rings. We do not assert exactness at the subsequent terms of
the Gersten complex.

\subsection*{Main results}

\begin{theorem}\label{thm:gersten-injectivity}
Let $R$ be a noetherian regular local ring with fraction field $F$.
For every natural number $N\geq2$ invertible in $R$ and every
$n\in\Z$, restriction along $R\longrightarrow F$ induces an injection
\[
                 K_n(R;\Z/N)\longrightarrow K_n(F;\Z/N).
\]
\end{theorem}

By \eqref{eq:coefficient-sequence},
\cref{thm:gersten-injectivity} also
implies that $K_n(R)/N\to K_n(F)/N$ is injective. In particular, every
element of $\ker(K_n(R)\to K_n(F))$ belongs to $N K_n(R)$ for every
$N\geq2$ invertible in $R$. This is divisibility in $K_n(R)$; the
assertion does not imply divisibility inside the kernel or integral
injectivity.

\paragraph{\textnormal{\bfseries \'Etale cohomology.}}
We write $H^q_{\mathrm{\acute et}}(T,M)$ for cohomology on the
small \'etale site of a scheme $T$, and
$H^q_{\mathrm{\acute et}}(B,M)=H^q_{\mathrm{\acute et}}(\Spec B,M)$
for a ring $B$. For a morphism $Y\to T$
and a sheaf $M$ on $T_{\mathrm{\acute et}}$, the inverse image of $M$
on $Y_{\mathrm{\acute et}}$ is denoted by $M|_Y$.
For a finite locally constant abelian sheaf $M$ on $X$ and $q\geq0$,
the first map of the augmented coniveau complex in \'etale cohomology
is $H^q_{\mathrm{\acute et}}(X,M)\to H^q_{\mathrm{\acute et}}(F,M|_F)$.
The subsequent terms in codimension $p\geq1$ are
\[
 \bigoplus_{x\in X^{(p)}}
 H^{q+p}_{\mathrm{\acute et},\{x\}}(\Spec\mathcal O_{X,x},M|_{\Spec\mathcal O_{X,x}}),
\]
with the differentials given by localization (cf.\ \cite[\S3]{BlochOgus}).
Let $A$ be a finite \'etale $R$-algebra, and let $M$ be a finite locally
constant abelian sheaf on $\Spec A$ whose stalks have orders invertible
in $A$. By \cref{thm:etale-injectivity}, for every $q\geq0$ the map
\begin{equation}\label{eq:intro-etale-injection}
 H^q_{\mathrm{\acute et}}(A,M)\longrightarrow
 H^q_{\mathrm{\acute et}}(A\otimes_R F,M|_{A\otimes_R F})
\end{equation}
is injective. The coefficients in \eqref{eq:intro-etale-injection}
need not be Tate twists. In degrees one and two,
\cref{prop:degree-one-residue-sequence,prop:semilocal-degree-two-injectivity}
give generic injectivity more generally for connected noetherian regular
semilocal rings, by a direct argument with height-one residues.

\paragraph{\textnormal{\bfseries Continuous cohomology.}}
For an integer $m\geq2$ invertible on the scheme under consideration,
write $\Z/m(j)=\mu_m^{\otimes j}$ for $j\geq0$, with the constant sheaf
in twist zero; negative tensor powers are defined by duality over
$\Z/m$. We write $H^q_{\mathrm{cont}}(-,\Z_\ell(j))$ for the continuous
cohomology of the inverse system of Tate sheaves $\Z/\ell^n(j)$, as in
\cref{def:continuous-tate-cohomology}.

The injectivity in \cref{thm:etale-injectivity} does not extend
to continuous $\ell$-adic cohomology.
For every prime $\ell$, there exists an excellent discrete valuation ring
$R$ containing $\mathbf Q$, with fraction field $F$ purely transcendental
of degree two over $\mathbf Q$, such that
\begin{equation}\label{eq:intro-continuous-kernel}
 \ker\bigl(H^2_{\mathrm{cont}}(R,\Z_\ell(j))\longrightarrow
 H^2_{\mathrm{cont}}(F,\Z_\ell(j))\bigr)
 \simeq\begin{cases}
 0 & j=1,\\
 \Z_\ell & j\ne1
 \end{cases}
\end{equation}
as $\Z_\ell$-modules (\cref{thm:continuous-counterexample}).
For every $j\in\Z$ and $n\geq1$, the corresponding map in degree two
with $\Z/\ell^n(j)$-coefficients is injective by
\cref{thm:etale-injectivity}. Taking ordinary inverse limits
of these finite-coefficient maps
preserves injectivity, but continuous cohomology involves a derived
inverse limit. The kernel in \eqref{eq:intro-continuous-kernel} is a
first derived-limit obstruction (\cref{cor:continuous-kernel-derived-limit}),
not a failure at any finite coefficient level.

\paragraph{\textnormal{\bfseries Semilocal rings.}}
Let $R$ be a regular semilocal domain with fraction field $F$, essentially
of finite type over a field or over $\Z$, whose closed residue fields have
at most one positive characteristic.
\Cref{thm:finite-type-semilocal-gersten-injectivity}, together with
\cref{sub:moore-primary-decomposition}, gives the injectivity of
$K_n(R;\Z/N)\to K_n(F;\Z/N)$ for all $n\in\Z$ and all $N\geq2$
invertible in $R$.
\Cref{cor:finite-etale-algebra-k-injectivity} treats finite
\'etale algebras over a noetherian regular local ring, and
\cref{cor:semilocal-regular-model-k-injection} treats noetherian semilocal
domains $R$ that receive a regular homomorphism from a regular ring
essentially of finite type over a field or over $\Z$, again assuming
that the closed residue fields of $R$ have at most one positive
characteristic. This includes the unramified mixed-characteristic
case in \cref{sub:unramified-regular-over-zp}.
\Cref{cor:semilocal-regular-model-etale-injectivity} gives \'etale
injectivity for noetherian semilocal domains with a regular model,
under the same restriction on the closed residue characteristics,
and \cref{cor:equicharacteristic-semilocal-etale-injection} treats regular
semilocal domains containing a field. In degree two, the residue calculation of
\cref{prop:semilocal-degree-two-injectivity} proves \'etale
injectivity for every connected noetherian regular semilocal
ring; together with \cref{thm:semilocal-finite-connective-comparison}
it gives the same assertion for $K_2(-;\Z/N)$
(\cref{cor:semilocal-K2-injectivity}).
\Cref{thm:postnikov-band-injectivity} treats each finite
Postnikov band of \'etale $K$-theory over a regular local ring.
For continuous cohomology, \cref{cor:continuous-degree-two-obstruction}
expresses the degree-two kernel on a connected noetherian regular
semilocal ring in which $\ell$ is invertible through compatible
height-one residues.

\subsection*{Relation with earlier results}

\paragraph{\textnormal{\bfseries Algebraic $K$-theory.}}
Quillen proved the integral Gersten conjecture for regular local
rings essentially of finite type over a field
\cite[\S7, Theorem~5.11]{Quillen}, and Panin proved it for
arbitrary equicharacteristic regular local rings
\cite[Theorem~A]{Panin}. In equal characteristic the classical results
also give injectivity for every $N\geq2$, including noninvertible
coefficient orders (\cref{lem:equicharacteristic-injection}).

In mixed characteristic, Gillet proved the Gersten conjecture for
discrete valuation rings with coefficients of order prime to the
residue characteristic \cite{Gillet}. Together with the relative
theorem of Gillet and Levine \cite[Corollary~6]{GilletLevine},
this gives full exactness for regular local rings essentially
smooth over a discrete valuation ring with such coefficients.
Geisser and Levine supplied the residue-characteristic-primary
case in this smooth setting \cite[Theorem~8.2]{GeisserLevine}.
Skalit obtained full finite-coefficient exactness for arbitrary
unramified regular local rings \cite[Corollary~5.7]{Skalit}.
For henselian regular local rings the first map is injective with
invertible coefficients by \cite[Corollary~6.13]{FeldHenselian} and
primary decomposition (\cref{sub:moore-primary-decomposition}).
The proof of \cite[Corollary~6.13]{FeldHenselian} also gives a direct
induction on the dimension using Gabber's rigidity \cite{Gabber92} and Gillet's theorem
\cite{Gillet}, analogous to Sakagaito's argument in \'etale cohomology
\cite[Proposition~5]{Sakagaito20}.
For improved Milnor $K$-theory, Dahlhausen proved integral
Gersten exactness in every positive degree for complete discrete
valuation rings with finite residue field
\cite[\arxivhtml{1509.01087v3}{S4.Thmdefin4}{Corollary~4.4}]{Dahlhausen18};
L\"uders extended this result to excellent henselian discrete
valuation rings with finite residue field
\cite[\arxivhtml{2010.02622v4}{S5.Thmtheorem1}{Theorem~5.1}]{LudersMilnor24}.
For regular local surfaces whose reduced closed curves have finite
normalization, an independent proof of the first map by moving
residues among Cartier curves, and of full exactness in every degree,
with invertible finite coefficients, is given in
\cite[Theorem~1.1]{FeldSurface}. The injectivity part is
a special case of
\cref{thm:gersten-injectivity}.
The proof of \cite[Theorem~1.1]{FeldSurface} uses only local constructions
and the basic operations on $K$-theory with supports. By contrast,
\cref{thm:gersten-injectivity} covers
nonhenselian ramified regular local rings of mixed characteristic in
arbitrary dimension, without finiteness hypotheses on normalizations. The full exactness in \cite[Theorem~1.1]{FeldSurface}
is stronger than injectivity alone.

\paragraph{\textnormal{\bfseries \'Etale cohomology.}}
Sakagaito proved generic injectivity with invertible finite Tate coefficients over henselian regular local
rings of arbitrary dimension \cite[Proposition~5]{Sakagaito20}, and
Gersten exactness in every cohomological degree with invertible finite
Tate coefficients over henselian regular local surfaces, without
assuming excellence or equicharacteristic
\cite[Theorem~9]{Sakagaito20}.
For the henselizations of local rings of smooth schemes over a
henselian discrete valuation ring with infinite residue field,
Schmidt and Strunk proved full Gersten exactness for coefficients
pulled back from bounded complexes on the base with locally constant
constructible cohomology sheaves of $\Z/m$-modules, where $m$ is
invertible on the base
\cite[Theorem~6.8 and Corollary~6.10]{SchmidtStrunk23}.
Gaidai-Turlov proves Gersten exactness for regular local rings
geometrically regular over a discrete valuation ring, with finite
locally constant abelian coefficients pulled back from that base
and annihilated by an integer invertible on the base
\cite[Theorem~1.1]{GaidaiTurlov}.
Panin--Zainoulline prove
Gersten exactness for finite locally constant coefficients on
smooth semilocal schemes over a field
\cite[Corollary~1.3]{PaninZainoulline}.\footnote{In the arXiv version,
the corresponding assertion is Corollary~1.2, stated over an infinite
field; Remark~1.5 explains the passage to finite fields by transfers.}
Let $B$ be a Dedekind ring and let $m\geq2$ be a natural number invertible on $B$.
For smooth $B$-schemes and the sheaves $\mu_m^{\otimes b}$, $b\in\Z$,
Gersten exactness follows from Geisser's theorem
\cite[Theorem~1.2]{Geisser04} and the norm-residue theorem
(\cref{rem:etale-dedekind-literature}). Over an excellent discrete
valuation ring, L\"uders treats these coefficients for prime-power moduli within an
axiomatic Bloch--Ogus theory
\cite[\arxivhtml{2201.11553v4}{S2.Thmtheorem9}{Example~2.9} and \arxivhtml{2201.11553v4}{S3.Thmtheorem6}{Theorem~3.6}]{Luders22}.
In equal characteristic \cref{thm:etale-injectivity} also follows
from Panin--Zainoulline \cite[Corollary~1.3]{PaninZainoulline} by
approximation (\cref{rem:etale-equicharacteristic-literature}).
\Cref{thm:etale-injectivity} does not require a
smooth structure over a field or a discrete valuation ring, or that the
coefficients be pulled back from such a base.

\subsection*{Strategy of the proof}

The geometric presentation arguments of Quillen and Gillet--Levine
\cite[\S7]{Quillen}\cite{GilletLevine} use smoothness over a field or
a discrete valuation ring. A ramified mixed-characteristic regular
local ring need not have such a presentation.
We use \'etale sheaves of spectra to represent a class with zero generic
image by a section equipped with a path to the zero section over the
fraction field. In the finite-type case, \cref{thm:generic-path-return}
constructs a regular finite locally free scheme of rank prime to the
coefficient prime on which the sections are homotopic.
By \cref{lem:postnikov-band-trace}, restriction followed by
transfer on a finite Postnikov band is multiplication by that rank.
Together with the coefficient torsion bound
(\cref{sub:finite-band-continuity}), this yields injectivity for finite
bands by \cref{prop:finite-type-band-injection}. We describe the four
steps in the proof of
\cref{thm:gersten-injectivity} for
$N=\ell^\nu$, where $\ell$ is prime and $\nu\geq1$; this suffices by
\cref{sub:moore-primary-decomposition}.

\paragraph{\textnormal{\bfseries 1. The canonical finite comparison.}}
For an affine scheme $X$ on which $N$ is invertible, let
$\mathcal K^{\mathrm{\acute et}}_{X,N}$ be the hypercomplete \'etale
sheaf of spectra associated with $U\mapsto K^{\mathrm{cn}}(U)/N$
(\cref{con:etale-k-sheaves}); for a ring $A$ we write
$\mathcal K^{\mathrm{\acute et}}_{A,N}$ for the sheaf on $\Spec A$.
By \cref{sub:finite-k-etale-coefficients}, its only nonzero homotopy
sheaves are $\Z/N(j)$ in degrees $2j\geq0$. For a local ring $A$ in
which $N$ is invertible, the canonical map
\begin{equation}\label{eq:intro-comparison}
 K_n^{\mathrm{cn}}(A;\Z/N)\xrightarrow{\ \sim\ }
 \pi_nR\Gamma_{\mathrm{\acute et}}(A,
          \tau_{\leq2n}\mathcal K^{\mathrm{\acute et}}_{A,N})
 \qquad(n\geq0)
\end{equation}
is an isomorphism (\cref{thm:local-finite-connective-comparison}).
Here and below, $\tau_{\leq b}$ on spectra denotes truncation in
homotopy degrees, whereas $\tau^{\leq b}$ on complexes denotes
truncation in cohomological degrees.

We construct a filtered refinement of the canonical comparison map in
\eqref{eq:intro-comparison}. We start over a semilocal Dedekind domain
$\Z_{(\Sigma)}$ in which $N$ is invertible (\cref{def:essentially-smooth-sigma}), with the slice
filtration of $\mathrm{KGL}$ on smooth schemes
(\cref{con:classical-source-filtration,sub:classical-filtration-inputs};
\cite[\arxivhtml{2508.09915v1}{S4.E39}{Theorem~4.39} and
\arxivhtml{2508.09915v1}{S6.E21}{Theorem~6.21}]{BEMMotivic25}).
Its grades identify with shifts of Bloch's cycle complexes by
\cite[Theorem~1.1(1)]{BachmannEffective} and
\cite[Theorem~7.18 and Corollary~7.19]{Spitzweck}
(\cref{lem:classical-slice-source-bounds}). \Cref{def:lisse-mod-N} extends the source filtration by left Kan extension
from smooth algebras. The local
upper bound of Geisser \cite[Corollary~4.4]{Geisser04} and the slice
connectivity bound of
\cite[\arxivhtml{2508.09915v1}{S9.E5}{Corollary~9.5}]{BEMMotivic25}
give $j$-connectivity of its stage $\mathrm{Fil}_{\mathrm{lis}}^jK^{\mathrm{cn}}(A)/N$
(\cref{prop:filtered-etale-comparison}\textup{(i)}).
The rigidity argument of
\cref{lem:finite-ordinary-truncation-smooth-kan}, using
\cite[\arxivhtml{2508.09915v1}{S2.E25}{Lemma~2.25}]{BEMMotivic25},
identifies the lisse motivic grades.

\Cref{prop:filtered-etale-comparison} gives the map of filtrations
\[
 \mathrm{Fil}_{\mathrm{lis}}^{\bullet}K^{\mathrm{cn}}(A)/N\longrightarrow
 R\Gamma_{\mathrm{\acute et}}(A,\tau_{\geq2\bullet}
                             \mathcal K^{\mathrm{\acute et}}_{A,N}).
\]
By \cref{prop:filtered-etale-comparison}\textup{(ii), (iii)},
its stage zero is the canonical \'etale unit, and its grade $j$ is
the truncation map
\begin{equation}\label{eq:intro-graded}
 \bigl(\tau^{\leq j}R\Gamma_{\mathrm{\acute et}}(A,\Z/N(j))\bigr)[2j]
 \longrightarrow R\Gamma_{\mathrm{\acute et}}(A,\Z/N(j))[2j].
\end{equation}
The identification of the morphism in \eqref{eq:intro-graded}
uses the normalization in
\cite[\arxivhtml{2508.09915v1}{S6.E8}{Proposition~6.8}, formula (6.10), and
\arxivhtml{2508.09915v1}{S6.E21}{Theorem~6.21(2)}]{BEMMotivic25}
(\cref{lem:classical-grade-normalization}). The cofiber of
\eqref{eq:intro-graded} has homotopy groups only in degrees at most
$j-1$. The stage $\mathrm{Fil}_{\mathrm{lis}}^{n+1}K^{\mathrm{cn}}(A)/N$ is $(n+1)$-connective,
so quotienting by it does not change $\pi_n$. On the target, truncation
at $2n$ removes the stages of index $>n$. The isomorphism \eqref{eq:intro-comparison} is
therefore obtained by descending induction on $j=n,\ldots,0$
(\cref{sub:finite-filtration-comparison,lem:finite-comparison-from-source-bound}).
For comparison, the identification of the graded pieces in positive
weights of the lisse filtration over \(\Z\) on local rings follows from
Elmanto--Morrow \cite[Theorem~7.7]{ElmantoMorrow} and Bouis
\cite[Lemma~2.9 and Corollary~2.12]{BouisWeibel25}, together with
\cite[Corollary~5.6]{Bouis}. \Cref{rem:comparison-literature} explains
this relation; the present proof does not use those comparisons.

\paragraph{\textnormal{\bfseries 2. Returning a generic path.}}
Let $R$ be a regular semilocal domain essentially of finite type over
a field or over $\Z$, with fraction field $F$, in which $\ell$ is
invertible and whose closed residue fields have at most one positive
characteristic. Let $\mathcal F$ be a bounded hypercomplete \'etale sheaf of
spectra on $\Spec R$ with finite locally constant $\ell$-primary
homotopy sheaves. If two sections $z_0,z_1$ of the space
$\Omega^\infty R\Gamma_{\mathrm{\acute et}}(R,\mathcal F)$ are joined by a path
over $F$, then \cref{thm:generic-path-return} constructs a
regular scheme $T$, finite locally free over $R$ of rank prime to
$\ell$, over which the sections are joined by a path (not necessarily
extending the chosen path over $F$). The construction of $T$ in
\cref{thm:generic-path-return} proceeds in three steps, described in
(a)--(c): we represent a generic path by finite Galois data, construct a
smooth projective $R$-scheme $X$ with an open subscheme $U$ on which the
two sections are homotopic, and construct a regular finite locally free
$R$-scheme of rank prime to $\ell$ with a morphism to $U$.

\emph{(a) Finite Galois data.}
The unit-symbol theorem of L\"uders
\cite[\arxivhtml{2201.11553v4}{S3.Thmtheorem8}{Corollary~3.8}]{Luders22}
and its semilocal extension (\cref{prop:semilocal-unit-symbols})
give effacement of positive-degree cohomology with finite locally
constant $\ell$-primary coefficients by finite \'etale covers
(\cref{prop:finite-etale-effacement}).
The proof of \cref{sub:local-symbols-via-smooth-sources} follows
\cite[proof of \arxivhtml{2201.11553v4}{S3.Thmtheorem8}{Corollary~3.8}]{Luders22},
using Geisser's theorem \cite[Theorem~1.2]{Geisser04} on smooth local
rings together with the finite-coefficient Milnor lifting statement
in \cref{sub:unramified-milnor-smooth-local}, stated in
\cite[\arxivhtml{2010.02622v4}{S4.Thmtheorem3}{Corollary~4.3} and \arxivhtml{2010.02622v4}{S4.Thmtheorem4}{Remark~4.4}]{LudersMilnor24},
the norm-residue theorem of Rost and Voevodsky for fields
\cite{Voevodsky11,HaesemeyerWeibel19}, and Gabber's affine analogue of
proper base change \cite[\stackstag{09ZI}]{Stacks}.
By Postnikov induction
(\cref{prop:bounded-finite-etale-descent}), sections of $\mathcal F$ are
computed on the finite \'etale site. Write $\rho_*$ for restriction to
that site (\cref{not:finite-etale-topoi}). By
\cref{cor:finite-site-local-trivialization}, the sheaf of spaces
$\Omega^\infty\rho_*\mathcal F$ becomes constant on a finite \'etale
cover, with value a $\pi$-finite space $\mathsf A$ (finitely many
components and finite homotopy groups, vanishing in high degrees). For
$r\in\{1,2\}$, put
\[
 \mathcal M_r(\mathsf A)=\mathsf A^r\mathbin{/\mkern-6mu/}\operatorname{Aut}(\mathsf A),
\]
where $/\mkern-6mu/$ denotes the homotopy quotient for the diagonal
action induced by evaluation. The diagonal gives
$\mathcal M_1(\mathsf A)\to\mathcal M_2(\mathsf A)$
(\cref{con:marked-moduli}). By \cref{prop:finite-group-stages,prop:finite-galois-data-two-sections},
the two sections, together with their underlying sheaf of spaces, are
represented by a finite quotient $G$ of the \'etale fundamental group,
its connected right $G$-torsor $P_0$, and a map
$\gamma:BG\to\mathcal M_2(\mathsf A)$. By
\cref{lem:finite-stage-homotopy-lift,prop:finite-galois-data-two-sections},
a generic path is represented by a homotopy lift after passage to a
finite quotient of the absolute Galois group of $F$ refining its map
to $G$. The homotopy lift includes the homotopy from its image in
$\mathcal M_2(\mathsf A)$ to the prescribed map (\cref{def:finite-stage}).

\emph{(b) From homotopy lifts to a smooth projective scheme.}
We first reduce $G$ to a finite $\ell$-group, and then reduce the
absolute Galois group to a pro-$\ell$ group. Let $G_\ell\subseteq G$ be
an $\ell$-Sylow subgroup and let $\Spec R'=P_0/G_\ell$. This is finite
\'etale over $R$ of rank prime to $\ell$ and $P_0\to\Spec R'$ is a
$G_\ell$-torsor (\cref{sub:finite-sylow-base-change}). Write $F'=\Frac R'$, and let $L$ be the fixed field of a pro-$\ell$ Sylow subgroup of the
absolute Galois group of $F'$; every finite subextension of $L/F'$
has degree prime to $\ell$ (\cref{sub:pro-ell-sylow-fixed-field}).
Over $L$ the generic lift may be represented on a finite $\ell$-group
$H$ with surjection $\phi:H\twoheadrightarrow G_\ell$, in the homotopy-commutative diagram
of \cref{lem:finite-stage-homotopy-lift,sub:pro-ell-sylow-fixed-field}
\begin{equation}\label{eq:intro-lifting}
\begin{tikzcd}[column sep=large,row sep=large]
 & & \mathcal M_1(\mathsf A) \arrow[d,"\text{diagonal}"] \\
 BH \arrow[r,"B\phi"'] \arrow[urr,dashed,"\delta"] & BG_\ell
       \arrow[r,"\gamma_\ell"'] & \mathcal M_2(\mathsf A)
\end{tikzcd}
\end{equation}
Here $\gamma_\ell=\gamma\circ B(G_\ell\hookrightarrow G)$.
The lifting datum in \eqref{eq:intro-lifting} includes a specified homotopy
between the two composites from $BH$ to $\mathcal M_2(\mathsf A)$. By
\cref{sub:pro-ell-sylow-fixed-field}, the corresponding $H$-torsor over
$L$, with its identification after extension of structure group, is a
marked lift of $P_0\times_{\Spec R'}\Spec L$ along $\phi$
(\cref{def:marked-lift}). The two Sylow reductions ensure
that $|\ker\phi|$ is a power of $\ell$, hence invertible on the base,
while $[F':F]$ and the degrees of all finite subextensions of
$L/F'$ are prime to $\ell$.

Instead of extending this particular torsor over $R'$, we parametrize
marked lifts of $P_0$. This is related to the lifting varieties of
Demarche--Florence
\cite[\arfiv{1711.06585}{S2.Thmtheorem1}{Proposition~2.1} and
\arfiv{1711.06585}{S3.Thmtheorem3}{Example~3.3}]{DemarcheFlorence}.
Let $V^\circ$ be the open scheme of invertible matrices in the
space of endomorphisms of the regular representation of $H$. By \cref{con:lifting-model}, the free quotient
$U=(P_0\times_{\Spec R'}V^\circ)/H$ carries an $H$-torsor $Q$ and a marking
$Q\times^H G_\ell\simeq P_0\times_{\Spec R'}U$. Every marked lift over a
field equipped with an $R'$-algebra structure
is obtained from an appropriate rational point of $U$
(\cref{con:lifting-model,prop:marked-lifting-model-field-versality}).
The marking transports the homotopy $\delta$ to a path over $U$
(\cref{cor:transport-marked-path-by-torsor}), and the lift over $L$
gives an $L$-point of $U_{F'}$. Regard $U$ as an $R$-scheme via
$\Spec R'\to\Spec R$; then $U_F=U_{F'}$ as an $F$-scheme.
The $L$-point factors through a closed point $u$ of $U_{F'}$ with
residue field finite over $F'$. Its degree over $F$ is
$[\kappa(u):F'][F':F]$, prime to $\ell$ by the two Sylow reductions
(\cref{lem:closed-point-from-sylow-field}). By
\cref{lem:lifting-compactification-tame-fibres}, the stabilizers of the
projective quotient stack lie in $\ker\phi$, so the stack is tame.
Combining tame destackification
\cite[\arxivhtml{1905.00872v1}{theoremx2}{Theorem~B} and Remark~1.3]{BerghRydh}
with \cite[\S3.3]{EdidinRydh} gives a smooth projective $R'$-scheme $X$
containing $U$ as an open dense in every geometric fibre. The
modification leaves $U$ unchanged, so the marked lift and the path on
$U$ are preserved
(\cref{lem:destackification-preserves-scheme-open,prop:projective-marked-lifting-model}). Viewed over $R$,
$X$ is still smooth and projective, with geometric fibres pure of
relative dimension $|H|^2$ (\cref{lem:closed-point-from-sylow-field}).

\emph{(c) Multisections of controlled degree.}
We need a regular finite locally free $R$-scheme $T$ of rank prime to
$\ell$, with a morphism $T\to U$. The closure of the chosen generic
point is not necessarily finite or regular over $R$ and may meet
$X\setminus U$. An ordinary complete intersection on $X$ may have
degree divisible by $\ell$. Our construction combines a closure of
controlled degree with a complete intersection on a blowup. Related
finite quasi-sections and moving arguments are due to
Gabber--Liu--Lorenzini
\cite[Theorem~6.3]{GLLHypersurfaces15}\cite[\arfiv{1209.2828}{S2.Thmthm3}{Theorem~2.3}]{GLLIndex13}.
The additional requirements here are regularity of the total space,
avoidance of the boundary, and control of the generic degree modulo
$\ell$ throughout the successive replacements.

Using Gabber's projective $\ell'$-alteration theorem
\cite[Expos\'e~X, Theorems~2.1 and~2.4]{ILO14}, we start with a regular
projective $R$-scheme $C\to X$ whose generic fibre is finite of degree
prime to $\ell$ (\cref{lem:multisection-arithmetic-alteration}). Embed
$C$ as a closed subscheme of $\mathbb P^a_X$, with $a\geq1$,
and let $c\geq1$ be the relative dimension over $R$ of $\mathbb P^a_X$. Let $\pi:W\to\mathbb P^a_X$
be the blowup along $C$, with exceptional divisor $E$. Choose a very
ample invertible sheaf $\mathcal A$ for $\mathbb P^a_X\to\Spec R$ and a sufficiently
large positive multiple $b$ of $\ell$ so that
$\mathcal L_b(\mathcal A)=\pi^*\mathcal A^{\otimes b}\otimes\mathcal O_W(-E)$ is very ample over
$R$ (\cref{sub:multisection-blowup-very-ample}). After a finite \'etale
base change of rank prime to $\ell$ when needed, choose $c$ sections of
the pullback of $\mathcal L_b(\mathcal A)$ (\cref{con:multisection-replacement}).
By \cref{lem:multisection-geometric-incidence,sub:multisection-simultaneous-parameter-choice},
the sections can be chosen simultaneously so that on each geometric
closed fibre the local equations of each initial segment have
linearly independent cotangent classes at every common closed zero.
By \cref{sub:multisection-absolute-cotangent}, their common zero scheme
$Y$ is regular. Let $F_1$ be the residue field of any generic point of
the resulting base. In the calculation of
\cref{prop:multisection-generic-chow-degree}, the exceptional divisor
contributes the negative of the centre's zero-cycle:
\[
 \pi_{F_1*}[Y_{F_1}]=b^c c_1(\mathcal A_{F_1})^c\cap[\mathbb P^a_{X_{F_1}}]-[C_{F_1}]
 \quad\text{in }\operatorname{CH}_0(\mathbb P^a_{X_{F_1}}).
\]
For a finite scheme $W$ over a field $E$, write
$\deg_E(W)=\dim_E\Gamma(W,\mathcal O_W)$.
Since $\ell\mid b$ and $c\geq1$, taking degrees in
\eqref{eq:multisection-generic-degree} gives
$\deg_{F_1}(Y_{F_1})\equiv-\deg_F(C_F)\pmod\ell$, using invariance of
the degree of a finite scheme under extension of the ground field.
Replacing $C$ by $Y$ lowers the maximum dimension of its closed fibres
whenever that maximum is positive, and otherwise does not increase it.
The maximum dimension of the closed fibres of the inverse image of
$X\setminus U$ strictly decreases until they are empty
(\cref{prop:multisection-decrease-closed-fibres}). The generic degrees
remain nonzero modulo $\ell$, with a common residue on the components
after each finite \'etale change of base, by
\eqref{eq:multisection-iterated-degree}. Iteration in
\cref{thm:regular-multisections} produces a regular $T$ finite over the
base and mapping to $U$. Regularity,
together with the local dimension calculation, gives flatness by miracle
flatness
(\cref{sub:multisection-closed-local-dimensions,sub:multisection-finite-flatness}).
This is the geometric content of \cref{thm:regular-multisections,cor:arithmetic-multisections}.
Pulling back the path over $U$ to $T$ proves
\cref{thm:generic-path-return}.

\paragraph{\textnormal{\bfseries 3. Transfer on the whole Postnikov band.}}
Injectivity in \'etale cohomology for each Tate sheaf separately does
not formally imply injectivity on a finite Postnikov band. In the
descent spectral sequence over $F$, the image of a permanent cycle
from $R$ can be a boundary of a differential whose source is not in
the image of restriction. The
example in \cref{rem:band-trace-pitfalls} illustrates the failure of a
gradewise injectivity argument. This is why \cref{thm:generic-path-return,lem:postnikov-band-trace}
are applied to the entire band, with its extension data.

For $0\leq a\leq b$ and a $\Z[1/N]$-algebra $A$, use the notation
for the homotopy groups of sections of a finite Postnikov band
(\cref{not:finite-postnikov-band}):
\[
 \Phi^{[a,b]}_n(A)=\pi_nR\Gamma_{\mathrm{\acute et}}(A,
               \tau_{[a,b]}\mathcal K^{\mathrm{\acute et}}_{A,N}).
\]
For a finite locally free morphism $f:T\to\Spec R$ of constant rank $e$,
\cref{lem:postnikov-band-trace} constructs pullback and transfer
maps on finite Postnikov bands from extension and restriction of
scalars in $K$-theory, with composite multiplication by $e$. In
\cref{lem:postnikov-band-trace}, the bounded truncation of the
sheafification of $Y\mapsto K^{\mathrm{cn}}(Y\times_{\Spec R}T)/N$ on the small
\'etale site of $\Spec R$ is identified with the direct image of the
corresponding truncation on $T$, before taking global sections.
\Cref{lem:postnikov-band-trace} requires only finite local freeness of
$f$, not regularity of $T$. Even when $T$ is not local, the comparison
of \cref{thm:semilocal-finite-connective-comparison} is not needed.

For a regular semilocal domain $R$ as in step 2, apply
\cref{thm:generic-path-return} to sections of $\Omega^\infty\mathcal F$
representing $x\in\Phi^{[a,b]}_n(R)$ and zero, where
$\mathcal F=\Sigma^{-n}\tau_{[a,b]}\mathcal K^{\mathrm{\acute et}}_{R,N}$.
If $x|_F=0$, then $f^*x=0$ for a finite locally free $f$ of rank $e$
prime to $\ell$. The composite $\operatorname{tr}_f f^*$ gives $ex=0$.
The group $\Phi^{[a,b]}_n(R)$ is killed by $N^2$
(\cref{sub:finite-band-continuity}); hence $x=0$. This proves
\cref{prop:finite-type-band-injection}.

\paragraph{\textnormal{\bfseries 4. Approximation and conclusion.}}
To pass to an arbitrary noetherian regular local ring $R$, use
N\'eron--Popescu desingularization in equal characteristic
\cite[\stackstag{07GC}]{Stacks}, and Popescu's answer to Swan's
question in mixed characteristic \cite[Theorem~3.8]{Popescu19}.
By \cref{sub:finite-band-continuity,lem:finite-type-witness}, a class
with zero generic image vanishes after inverting a nonzero element of
$R$. Both the class and this vanishing after localization descend to a
regular local ring essentially of finite type. The descent concerns a
principal localization, not an induced map of fraction fields: the
approximating maps need not be injective.
It is essential that \cref{prop:finite-type-band-injection} covers ramified rings:
in mixed characteristic $(0,p)$, if
$p\in\mathfrak m_R^2$, a finite expression for $p$ as a sum of
products of elements of $\mathfrak m_R$ descends to a relation in the
square of the maximal ideal at a sufficiently late local stage of
\cref{lem:finite-type-witness}.
Thus \cref{thm:postnikov-band-injectivity} gives injectivity on
finite bands for arbitrary $R$. Together with the natural isomorphism
of \cref{thm:local-finite-connective-comparison} and the agreement of
connective and nonconnective $K$-theory on regular noetherian rings
\cite[Proposition~6.8(b)]{TT}, we obtain, for $n\geq0$, the commutative
square
\begin{equation}\label{eq:intro-final-square}
\begin{tikzcd}[column sep=large,row sep=large]
 K_n(R;\Z/N) \arrow[r] \arrow[d,"\sim"'] &
 K_n(F;\Z/N) \arrow[d,"\sim"] \\
 \Phi^{[0,2n]}_n(R) \arrow[r,hook] & \Phi^{[0,2n]}_n(F).
\end{tikzcd}
\end{equation}
The comparison square \eqref{eq:intro-final-square} proves
\cref{thm:gersten-injectivity};
negative degrees vanish by regularity \cite[Proposition~6.8(b)]{TT}
and \eqref{eq:coefficient-sequence}.

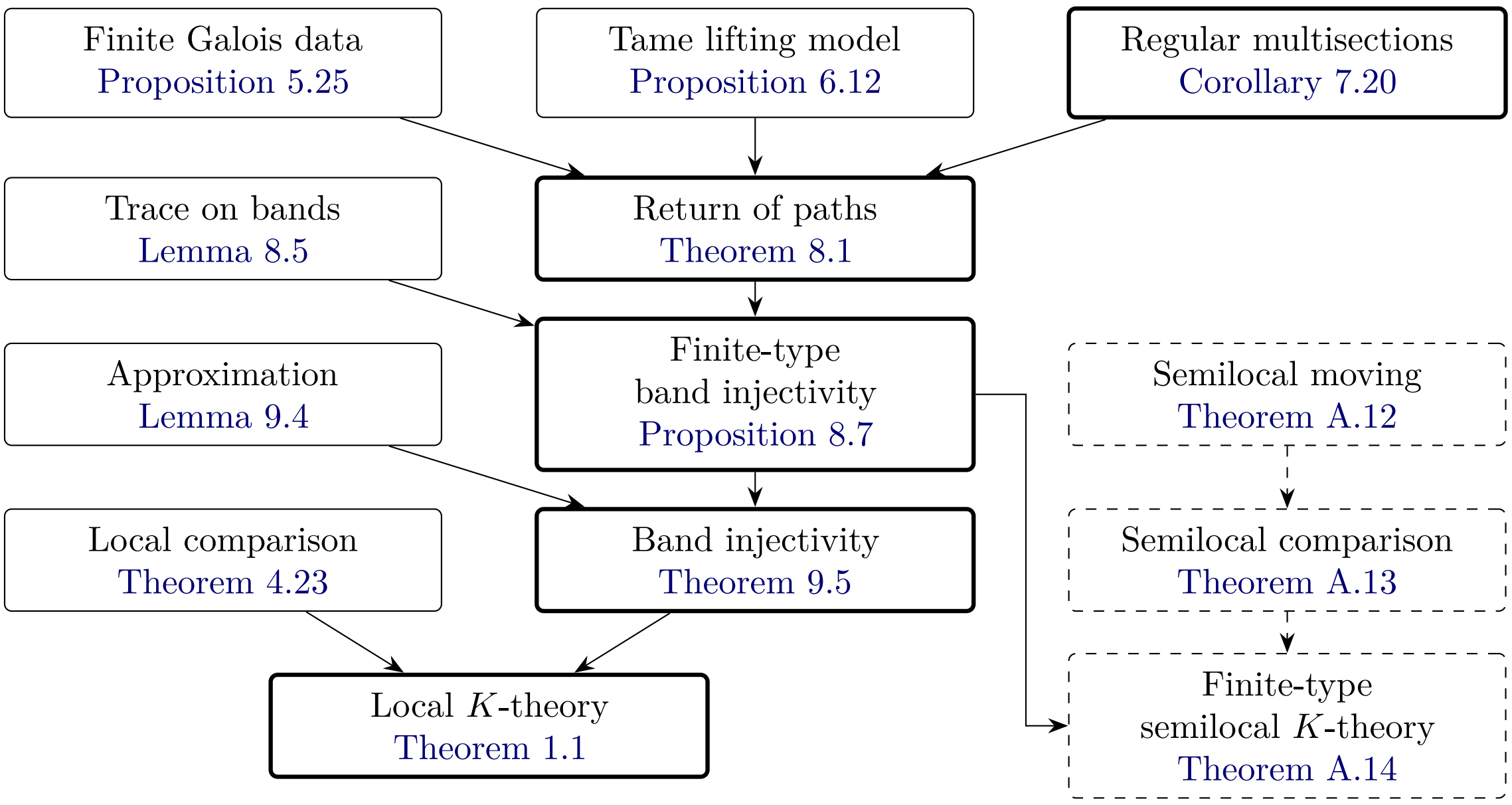
\begin{figure}[!htbp]
\centering
\begin{tikzpicture}[
 >={Stealth[length=2mm]},
 box/.style={draw,rounded corners=2pt,align=center,font=\small,
             text width=4.0cm,inner sep=5pt},
 main/.style={box,very thick},
 aux/.style={box,dashed}]
 \node[box] (data) at (0,0) {Finite Galois data\\\cref{prop:finite-galois-data-two-sections}};
 \node[box] (model) at (5.3,0) {Tame lifting model\\\cref{prop:projective-marked-lifting-model}};
 \node[main] (multi) at (10.6,0) {Regular multisections\\\cref{cor:arithmetic-multisections}};
 \node[main] (ret) at (5.3,-1.65) {Return of paths\\\cref{thm:generic-path-return}};
 \node[box] (trace) at (0,-1.65) {Trace on bands\\\cref{lem:postnikov-band-trace}};
 \node[main] (finite) at (5.3,-3.3) {Finite-type band injectivity\\\cref{prop:finite-type-band-injection}};
 \node[box] (approx) at (0,-3.3) {Approximation\\\cref{lem:finite-type-witness}};
 \node[main] (bands) at (5.3,-4.95) {Band injectivity\\\cref{thm:postnikov-band-injectivity}};
 \node[box] (local) at (0,-4.95) {Local comparison\\\cref{thm:local-finite-connective-comparison}};
 \node[aux] (moving) at (10.6,-3.3) {Semilocal moving\\\cref{thm:semilocal-motivic-vanishing}};
 \node[aux] (semi) at (10.6,-4.95) {Semilocal comparison\\\cref{thm:semilocal-finite-connective-comparison}};
 \node[main] (main) at (2.65,-6.6) {Local $K$-theory\\\cref{thm:gersten-injectivity}};
 \node[aux] (semigi) at (10.6,-6.6) {Finite-type\\semilocal $K$-theory\\\cref{thm:finite-type-semilocal-gersten-injectivity}};
 \draw[->] (data) -- (ret); \draw[->] (model) -- (ret);
 \draw[->] (multi) -- (ret); \draw[->] (ret) -- (finite);
 \draw[->] (trace) -- (finite); \draw[->] (finite) -- (bands);
 \draw[->] (approx) -- (bands); \draw[->] (bands) -- (main);
 \draw[->] (local) -- (main); \draw[->,dashed] (moving) -- (semi);
 \draw[->] (finite.east) -- ++(0.5,0) |- (semigi.west);
 \draw[->,dashed] (semi) -- (semigi);
\end{tikzpicture}
\caption{Principal dependencies of the proofs of
\cref{thm:gersten-injectivity} and
\cref{thm:finite-type-semilocal-gersten-injectivity}. Arrows point from inputs to results
using them; dashed boxes belong to \cref{app:semilocal-moving}.}
\label{fig:intro-dependencies}
\end{figure}

\Cref{fig:intro-dependencies} indicates the principal logical dependencies.

\paragraph{\textnormal{\bfseries Ordinary and continuous \'etale cohomology.}}
For finite locally constant abelian coefficients,
\cref{thm:generic-path-return} applied to Eilenberg--Mac Lane
spectra, together with the ordinary finite locally free trace
(\cref{sub:finite-flat-etale-trace}) and approximation
(\cref{lem:finite-type-witness}), gives
\cref{thm:etale-injectivity} without the comparison in
\eqref{eq:intro-comparison}.
For the discrete valuation ring of \cref{thm:continuous-counterexample},
\cref{lem:finite-level-residue-surjectivity} gives surjectivity of the finite-level
residue maps in degree one. By \cref{lem:rational-function-compatible-classes},
for twists $j\ne1$ the induced map on inverse limits is zero, although
its target is free of rank one over $\Z_\ell$ by
\cref{cor:trait-tate-residues,lem:cyclotomic-valuation}.
\Cref{cor:continuous-degree-two-obstruction}\textup{(i)} identifies the cokernel
of the limit residue map with the kernel in
\eqref{eq:intro-continuous-kernel}; by
\cref{cor:continuous-kernel-derived-limit}, this kernel is an obstruction
in the first derived inverse limits.

\subsection*{Organization}

We fix conventions in \cref{sec:conventions} and recall classical cycle
complexes in \cref{sec:smooth-sources-cycles}.
\Cref{sec:finite-connective-comparison} develops the smooth Kan-extension
argument and proves
\cref{thm:local-finite-connective-comparison}.
\Cref{sec:finite-etale-effacement} gives the finite Galois data
(\cref{prop:finite-galois-data-two-sections});
\cref{sec:tame-lifting-model} constructs the smooth projective models
of marked lifts (\cref{prop:projective-marked-lifting-model}); and
\cref{sec:regular-finite-flat-multisections} proves the existence of
multisections of degree prime to $\ell$
(\cref{cor:arithmetic-multisections}).
\Cref{sec:generic-return} combines these results to prove
\cref{thm:generic-path-return}, and establishes
\cref{lem:postnikov-band-trace}. The proof of
\cref{thm:gersten-injectivity} concludes
in \cref{sec:approximation}.

\Cref{sec:ordinary-etale-injection,sec:continuous-cohomology}
are devoted to ordinary finite and continuous \'etale cohomology,
respectively. \Cref{app:semilocal-moving} establishes the semilocal
comparison, starting from the upper bound for Bloch's complexes in
\cref{thm:semilocal-motivic-vanishing}. The proof of
\cref{thm:semilocal-motivic-vanishing} extends Geisser's argument
\cite[Theorem~4.2 and Corollary~4.4]{Geisser04} using simultaneous moving
(\cref{prop:flat-principal-divisor-effacement}). This gives \cref{thm:semilocal-finite-connective-comparison} and
its application in \cref{thm:finite-type-semilocal-gersten-injectivity}.
The proof of \cref{thm:gersten-injectivity}
uses \cref{thm:local-finite-connective-comparison} but not
\cref{thm:semilocal-finite-connective-comparison}, so
\cref{app:semilocal-moving} can be omitted on a first reading.

We use the hierarchy Sublemma, Lemma, Proposition, Theorem to
indicate the role of a result in the argument. Sublemmas generally
record standard facts or auxiliary verifications, with references
and proofs for convenience. This
level of detail is intended to facilitate independent verification
and future formalization in Lean.

\section{Conventions and elementary reductions}\label{sec:conventions}

Rings are commutative with identity, and ring maps preserve the
identity. We allow the zero ring and write $\Spec(0)=\varnothing$.
A semilocal ring has finitely many maximal ideals (possibly none),
whereas a local ring has exactly one
\cite[\stackstag{07BI}]{Stacks}; a domain is nonzero. The set
$\mathbf N$ of natural numbers includes zero. An empty product of
rings is the zero ring. An empty sum in a ring is zero, and an
empty product of elements is one. A multiplicative subset may contain zero,
in which case the localization is the zero ring
\cite[\stackstag{00CQ}]{Stacks}.

The notation $E/N$ for a spectrum $E$ means the cofiber of
multiplication by $N$. For spectra and sheaves of spectra,
truncation subscripts refer to homotopy degrees:
$\tau_{\leq b}E$ has the homotopy sheaves of $E$ in degrees at
most $b$. For a cohomologically indexed complex $C$,
$\tau^{\leq b}C$ preserves cohomology in degrees at most $b$.
We regard complexes as Eilenberg--Mac Lane spectra, so
$\pi_{-i}C=H^i(C)$ and $\pi_i(C[k])=H^{k-i}(C)$.
For $N\geq2$ invertible on a scheme $X$ and $j\in\Z$, let
$\Z/N(j)=\mu_N^{\otimes j}$ be the Tate sheaf on the small \'etale site,
with $\Z/N(0)=\Z/N$. For $j>0$, put
$\mu_N^{\otimes(-j)}=\underline{\operatorname{Hom}}_{\Z/N}
(\mu_N^{\otimes j},\Z/N)$, the dual of the positive tensor power.

We use the language of $\infty$-categories of \cite{HTT,HA}. The
$\infty$-category of spaces is denoted by $\mathcal S$
\cite[Definition~1.2.16.1]{HTT}, the $\infty$-category of spectra by
$\mathrm{Sp}$ \cite[Definition~1.4.3.1]{HA}, and the derived
$\infty$-category of abelian groups by $\mathcal D(\Z)$
\cite[\S1.3.5]{HA}; a \emph{space} is an object of $\mathcal S$.
The Eilenberg--Mac Lane functor $\mathcal D(\Z)\to\mathrm{Sp}$
preserves limits and colimits and is conservative. For a site
$\mathcal C$, a \emph{sheaf of spaces} (respectively \emph{of
spectra}) on $\mathcal C$ is a functor from
$\mathcal C^{\mathrm{op}}$ to $\mathcal S$ (respectively to
$\mathrm{Sp}$) that satisfies descent for the covering sieves of
$\mathcal C$ \cite[Definition~6.2.2.6]{HTT}. A \emph{geometric
point} of a scheme $X$ is a morphism $\bar x:\Spec\Omega\to X$,
where $\Omega$ is an algebraically closed field
\cite[\stackstag{03PO}]{Stacks}. A connected scheme is understood to
be nonempty \cite[\stackstag{004S}]{Stacks}.

\begin{definition}[Semilocalizations and rings essentially of finite type]
\label{def:arithmetic-semilocal}
Let $A$ be a ring, let $t\in\mathbf N$, and let
$\mathfrak p_1,\ldots,\mathfrak p_t$ be pairwise incomparable prime
ideals of $A$. Put $S=A\setminus\bigcup_{i=1}^t\mathfrak p_i$.
The \emph{semilocalization of $A$ at
$\mathfrak p_1,\ldots,\mathfrak p_t$} is the ring $S^{-1}A$, and
the $\mathfrak p_i$ are the \emph{marked primes} of this
presentation. For $t=0$ the union is empty, $S=A$, and the
semilocalization is the zero ring. For $t\geq1$, by prime avoidance \cite[\stackstag{00DS}]{Stacks},
a prime ideal $\mathfrak q$ of $A$ satisfies
$\mathfrak q\cap S=\varnothing$ if and only if
$\mathfrak q\subseteq\mathfrak p_i$ for some $i$. Hence the maximal
ideals of $S^{-1}A$ are the $t$ pairwise distinct ideals
$\mathfrak p_iS^{-1}A$, and
\[
 S^{-1}A/\mathfrak p_iS^{-1}A=\Frac(A/\mathfrak p_i)
 \qquad(1\leq i\leq t).
\]

Let $R$ be a semilocal ring. The \emph{closed residue fields} of
$R$ are the fields $R/\mathfrak m$, where $\mathfrak m$ runs
through the maximal ideals of $R$.

A ring is \emph{essentially of finite type} over a ring $B$ if it is
isomorphic to $S^{-1}C$ for a finite-type $B$-algebra $C$ and a
multiplicative subset $S\subseteq C$. A finite-type algebra over a
ring essentially of finite type over $B$ is essentially of finite
type over $B$, and so are its localizations and its quotients.
\end{definition}

\begin{sublemma}[Semilocal domains essentially of finite type]
\label{sub:arithmetic-finite-etale}
Let $B_0$ be a field or the ring $\Z$, and let $R$ be a semilocal
domain essentially of finite type over $B_0$. Put $B=\mathbf F_p$ if
$B_0=\Z$ and $R$ has characteristic $p>0$, and $B=B_0$ otherwise.
\begin{enumerate}
\item There are an integral finite-type $B$-algebra $A$, with
$B\to A$ injective, and pairwise incomparable prime ideals
$\mathfrak p_1,\ldots,\mathfrak p_t$ of $A$, $t\geq1$, such that $R$
is isomorphic to the semilocalization of $A$ at
$\mathfrak p_1,\ldots,\mathfrak p_t$. If $R$ is regular, the ring
$A$ can be chosen regular.
\item Suppose that $R$ is regular, and let $R'$ be a nonzero finite
\'etale $R$-algebra. There are a natural number $r\geq1$ and an
isomorphism of $R$-algebras $R'\simeq\prod_{k=1}^rR'_k$ such that
every $R'_k$ is a regular semilocal domain essentially of finite
type over $B_0$. Every closed residue field of $R'_k$ is a finite
extension of a closed residue field of $R$.
\end{enumerate}
\end{sublemma}

\begin{proof}
\textup{(i)} Choose a finite-type $B_0$-algebra $C$, a
multiplicative subset $S\subseteq C$ and an isomorphism
$S^{-1}C\simeq R$, and let $A$ be the image of $C$ in $R$. The ring
$A$ is a domain, $R$ is a localization of $A$, the map $B\to A$ is
injective, and $A$ is of finite type over $B$. Let
$\mathfrak p_1,\ldots,\mathfrak p_t$ be the inverse images in $A$ of
the maximal ideals of $R$; they are pairwise incomparable. Every
element of $A\setminus\bigcup_i\mathfrak p_i$ lies in no maximal
ideal of $R$, so it is a unit of $R$. Hence $R$ is the
semilocalization of $A$ at $\mathfrak p_1,\ldots,\mathfrak p_t$.
Suppose that $R$ is regular. The ring $A$ is excellent
\cite[\stackstag{07QW}]{Stacks}, so its regular locus is the
complement of $V(I)$ for an ideal $I\subseteq A$
\cite[\stackstag{07QT}, \stackstag{07P7}]{Stacks}, and this locus
contains every $\mathfrak p_i$, because $A_{\mathfrak p_i}$ is a
localization of $R$. By prime avoidance
\cite[\stackstag{00DS}]{Stacks}, there is
$h\in I\setminus\bigcup_i\mathfrak p_i$. The ring $A[1/h]$ is
regular, and $R$ is its semilocalization at the ideals
$\mathfrak p_iA[1/h]$.

\smallskip\noindent
\textup{(ii)} The ring $R'$ is noetherian, and it is semilocal
because it is finite over $R$. It is regular
\cite[\stackstag{0AH7}]{Stacks}, hence normal
\cite[\stackstag{0567}]{Stacks}, and therefore a finite product
$\prod_{k=1}^rR'_k$ of normal domains
\cite[\stackstag{030C}]{Stacks}. Each $R'_k$ is a quotient and a
localization of $R'$, so it is a regular semilocal domain, and it is
essentially of finite type over $B_0$ because $R'$ is of finite type
over $R$. A maximal ideal $\mathfrak n$ of $R'_k$ lies over a
maximal ideal $\mathfrak m$ of $R$
\cite[\stackstag{00GR}]{Stacks}, and $R'_k/\mathfrak n$ is a
quotient of the finite $R/\mathfrak m$-algebra $R'/\mathfrak mR'$.
\end{proof}
\begin{sublemma}[Primary decomposition]
\label{sub:moore-primary-decomposition}
Let $N\geq2$ be a natural number and let $N=\prod_{i=1}^rN_i$ be
its factorization into powers of pairwise distinct primes.
There is an equivalence
\[
 E/N\xrightarrow{\ \sim\ }\bigoplus_{i=1}^rE/N_i
\]
in the $\infty$-category $\mathrm{Sp}$, natural in
$E\in\mathrm{Sp}$.
\end{sublemma}

\begin{proof}
Let $\mathbb S$ be the sphere spectrum, the unit of the smash
product $\otimes$ on $\mathrm{Sp}$, which preserves colimits in
each variable \cite[\S4.8.2]{HA}. For an integer $m$,
multiplication by $m$ on $E\simeq E\otimes\mathbb S$ is
$\mathrm{id}_E\otimes m$, so that
$E/m\simeq E\otimes(\mathbb S/m)$ naturally in $E$. It therefore
suffices to construct an equivalence
$\mathbb S/N\simeq\bigoplus_{i=1}^r\mathbb S/N_i$.

Put $d_i=N/N_i$. Since $N_id_i=N$ in $\pi_0\mathbb S=\Z$, the
square
\[
\begin{tikzcd}
 \mathbb S \arrow[r,"N"] \arrow[d,"d_i"'] &
 \mathbb S \arrow[d,equal] \\
 \mathbb S \arrow[r,"N_i"'] & \mathbb S
\end{tikzcd}
\]
commutes up to homotopy. A choice of homotopy induces a map of
cofibres $q_i:\mathbb S/N\to\mathbb S/N_i$. Let
\[
 q=(q_i)_{i=1}^r:\mathbb S/N\longrightarrow
   \prod_{i=1}^r\mathbb S/N_i=\bigoplus_{i=1}^r\mathbb S/N_i .
\]
A homomorphism of abelian groups is an isomorphism if its
localization at every prime is one. It therefore suffices to prove
that $\pi_n(q)\otimes\Z_{(\ell)}$ is an isomorphism for every
$n\in\Z$ and every prime $\ell$. If an
integer $m$ is prime to $\ell$, then $m$ acts invertibly on
$\pi_*(\mathbb S)\otimes\Z_{(\ell)}$, and the long exact sequence
of the cofibre gives $\pi_*(\mathbb S/m)\otimes\Z_{(\ell)}=0$.
If $\ell\nmid N$, the source and target of
$\pi_n(q)\otimes\Z_{(\ell)}$ are thus zero. If $\ell\mid N$, let
$i$ be the unique index such that $N_i$ is a power of $\ell$.
The summands $\mathbb S/N_k$ with $k\neq i$ vanish after
tensoring their homotopy groups with $\Z_{(\ell)}$. The vertical
maps $d_i$ and $\mathrm{id}$ of the square induce isomorphisms
on $\pi_*(\mathbb S)\otimes\Z_{(\ell)}$, since $d_i$ is prime to
$\ell$. The five lemma, applied to the long exact sequences of
the two rows, shows that $\pi_n(q_i)\otimes\Z_{(\ell)}$ is an
isomorphism.
\end{proof}

\begin{sublemma}[Freeness over a semilocal ring
  {\cite[\stackstag{02M9}]{Stacks}}]
\label{sub:semilocal-projective-free}
Let $R$ be a semilocal ring and let $P$ be a finite locally free
$R$-module of constant rank $e\in\mathbf N$. Then $P\simeq R^e$.
\end{sublemma}

\begin{proof}
The statement is \cite[\stackstag{02M9}]{Stacks}, whose proof is
omitted there; it is formalized in Mathlib as
\href{https://leanprover-community.github.io/mathlib4_docs/Mathlib/RingTheory/LocalRing/Module.html\#Module.free_of_flat_of_finrank_eq}{\texttt{Module.free\_of\_flat\_of\_finrank\_eq}}.
We follow the hints of \cite[\stackstag{02M9}]{Stacks}.
For $R=0$ the assertion is immediate. Suppose $R\neq0$, let
$\mathfrak m_1,\ldots,\mathfrak m_s$ be its maximal ideals,
and let $J=\bigcap_{i=1}^s\mathfrak m_i$ be its Jacobson radical.
The Chinese remainder theorem \cite[\stackstag{00DT}]{Stacks}
and tensor product with $P$ give
\[
 R/J\xrightarrow{\sim}\prod_{i=1}^s R/\mathfrak m_i,
 \qquad
 P/JP\xrightarrow{\sim}\prod_{i=1}^s P/\mathfrak m_iP.
\]
Each $P/\mathfrak m_iP$ has dimension $e$ over
$R/\mathfrak m_i$. Choose a basis in each factor. The $e$
tuples of corresponding basis vectors form an $R/J$-basis
of $P/JP$. Lift them to $P$ to obtain an $R$-linear map
$\varphi:R^e\to P$ whose reduction modulo $J$ is an
isomorphism. Since $P$ is finite and $J$ is the Jacobson
radical, Nakayama's lemma
\cite[\stackstag{00DV}, part~(6)]{Stacks}
implies that $\varphi$ is surjective.

The module $P$ is projective by
\cite[\stackstag{00NX}]{Stacks}, so the exact sequence
\[
 0\longrightarrow\ker\varphi\longrightarrow R^e
   \xrightarrow{\ \varphi\ }P\longrightarrow0
\]
splits. Thus $\ker\varphi$ is a direct summand of $R^e$,
and in particular is finitely generated. Reduction of this
split sequence modulo $J$ is exact. Since the reduction of
$\varphi$ is an isomorphism, it gives
$(\ker\varphi)/J(\ker\varphi)=0$. Nakayama's lemma
\cite[\stackstag{00DV}, part~(2)]{Stacks}
therefore gives $\ker\varphi=0$, proving the assertion.
\end{proof}

\section{Smooth sources and classical cycles}
\label{sec:smooth-sources-cycles}

We fix the base and the conventions for classical cycles used in
\cref{lem:classical-slice-source-bounds}. The proof of
\cref{thm:local-finite-connective-comparison} uses Geisser's upper
bound \cite[Corollary~4.4]{Geisser04}. The continuity and descent
statements of \cref{sub:classical-cycle-spread,sub:cycle-descent-localization}
are also used in the proof of the semilocal upper bound in
\cref{thm:semilocal-motivic-vanishing}.

\begin{definition}[The semilocal base and smooth sources]
\label{def:essentially-smooth-sigma}
Let \(\Sigma\) be a nonempty finite set of prime numbers. Define
\[
 \Z_{(\Sigma)}=
   \left(\Z\setminus\bigcup_{p\in\Sigma}p\Z\right)^{-1}\Z,
\]
a principal ideal domain whose maximal ideals are
\(p\Z_{(\Sigma)}\), for \(p\in\Sigma\), and whose fraction field is
\(\mathbf Q\).
An \emph{essentially smooth} \(\Z_{(\Sigma)}\)-algebra is a ring
\(M^{-1}C\), where \(C\) is a smooth finitely presented
\(\Z_{(\Sigma)}\)-algebra and \(M\subseteq C\) is a multiplicative
subset. It is \emph{semilocal essentially smooth} if one can take
\(M=C\setminus\bigcup_{i=1}^s\mathfrak p_i\) for some
\(s\in\mathbf N\) and pairwise incomparable prime ideals
\(\mathfrak p_1,\ldots,\mathfrak p_s\) of \(C\)
(\cref{def:arithmetic-semilocal}).
\end{definition}

\begin{notation}[Codimension, Cartier restriction, and the function \(\delta\)]
\label{not:codimension-delta}
Let \(T\) be a noetherian scheme. For an integral closed subscheme
\(W\subseteq T\) with generic point \(\eta_W\), put
\[
 \operatorname{codim}(W,T)=\dim\mathcal O_{T,\eta_W}.
\]
If \(T=\Spec C\) and \(T'=\Spec M^{-1}C\) for a multiplicative
subset \(M\subseteq C\), and if \(\eta_W\in T'\), then
\(\mathcal O_{T',\eta_W}=\mathcal O_{T,\eta_W}\), so that
\begin{equation}\label{eq:codim-localization}
 \operatorname{codim}(W\cap T',T')=\operatorname{codim}(W,T).
\end{equation}

Let \(D\subseteq T\) be a locally principal closed subscheme
\cite[\stackstag{01WR}]{Stacks}, and let \(W\subseteq T\) be an
integral closed subscheme with \(W\not\subseteq D\). Then
\(W\cap D=W\times_TD\) is an effective Cartier divisor on \(W\).
The \emph{Cartier restriction} of \(W\) to \(D\) is the cycle
\[
 D\cdot[W]=\sum_V
   \operatorname{length}\bigl(\mathcal O_{W\cap D,\eta_V}\bigr)\,[V],
\]
where \(V\) runs through the irreducible components of \(W\cap D\).
We extend \(D\cdot(-)\) by linearity to \(\Z\)-linear combinations
of integral closed subschemes of \(T\) not contained in \(D\).

Let \(T\) be of finite type over \(\Z_{(\Sigma)}\), or a localization
of such an affine scheme. For a point \(x\in T\) with image \(s\in\Spec\Z_{(\Sigma)}\), define
\[
 \delta(x)=\operatorname{trdeg}_{\kappa(s)}\kappa(x)
                     +1-\operatorname{ht}(s).
\]
For an integral closed subscheme \(V\) of a finite-type \(T\),
\(\delta(\eta_V)\) is Levine's dimension \cite[\S1.8]{LevineLocalization},
which need not be the Krull dimension of \(V\).
If \(T\) is smooth of constant relative dimension \(n\) over
\(\Z_{(\Sigma)}\), its connected components are integral by regularity
\cite[\stackstag{07NF},\stackstag{0567},\stackstag{033M}]{Stacks}, dominant by
flatness, and have function fields of transcendence degree \(n\) over
\(\mathbf Q\). The dimension formula over the universally catenary base
\(\Z_{(\Sigma)}\) \cite[\stackstag{02JU}]{Stacks}, applied to the
connected component of \(x\), gives
\begin{equation}\label{eq:codim-delta}
 \dim\mathcal O_{T,x}=n+1-\delta(x)\qquad(x\in T).
\end{equation}
If \(T'\to T\) is a finite morphism and \(x'\in T'\) has image
\(x\in T\), then \(\delta(x')=\delta(x)\), because
\(\kappa(x')/\kappa(x)\) is a finite extension.
Since \(s\mapsto1-\operatorname{ht}(s)\) is a dimension function on the
base, \(\delta\) is a dimension function on finite-type schemes
\cite[\stackstag{02JW}]{Stacks}, and hence on their localizations (chains
between surviving primes are unchanged). For an integral closed subscheme
\(V\subseteq T\) and \(x\in V\), we have
\(\delta(\eta_V)-\delta(x)=\dim\mathcal O_{V,x}\) by \cite[\stackstag{02IA}]{Stacks}.
\end{notation}

\begin{sublemma}[Commuting Cartier restrictions]
\label{sub:cartier-restrictions-commute}
Let \(T\) be a noetherian catenary scheme, let \(D_1,D_2\subseteq T\) be
locally principal closed subschemes, and let \(W\subseteq T\) be an
integral closed subscheme. Suppose that \(W\not\subseteq D_1\),
\(W\not\subseteq D_2\), that no irreducible component of
\(W\cap D_1\) is contained in \(D_2\), and that no irreducible
component of \(W\cap D_2\) is contained in \(D_1\). Then
\[
 D_2\cdot\bigl(D_1\cdot[W]\bigr)=D_1\cdot\bigl(D_2\cdot[W]\bigr)
\]
as cycles on \(W\cap D_1\cap D_2\). Equivalently, on the group
\(Z(T)_{D_1,D_2}\) of \(\Z\)-linear combinations of such \(W\), the
square of Cartier restrictions
\[
\begin{tikzcd}[column sep=large]
 Z(T)_{D_1,D_2} \arrow[r,"D_1\cdot"] \arrow[d,"D_2\cdot"'] &
 Z(D_1)_{D_2} \arrow[d,"D_2\cdot"] \\
 Z(D_2)_{D_1} \arrow[r,"D_1\cdot"'] & Z(D_1\cap D_2)
\end{tikzcd}
\]
commutes, where \(Z(-)\) denotes the group of cycles and
\(Z(D_i)_{D_j}\) its subgroup on integral closed subschemes of
\(D_i\) not contained in \(D_j\).
\end{sublemma}

\begin{proof}
For Chow groups of schemes of finite type over \(\Z_{(\Sigma)}\),
this is \cite[\stackstag{0B73}]{Stacks}. We prove the equality of
cycles. The scheme \(W\) is catenary \cite[\stackstag{02IX}]{Stacks}.
By the principal ideal theorem \cite[\stackstag{0BBZ}]{Stacks}, every
component of \(W\cap D_i\) has codimension one in \(W\), and every
component of its intersection with the other divisor has codimension
one in that component. By catenarity, each such component has
codimension two in \(W\). The whole double intersection is the
finite union of the intersections of the components of \(W\cap D_i\)
with the other divisor, so each of its irreducible components is one
of the components just considered. Since all have codimension two in
the integral scheme \(W\), no one is properly contained in another. Consequently
both iterated restrictions are supported on components of the whole
double intersection. Let \(x\) be a generic point of
\(W\cap D_1\cap D_2\), let \(A=\mathcal O_{W,x}\) with maximal ideal
\(\mathfrak m\), and let \(a,b\in\mathfrak m\) be the images of local
equations of \(D_1\) and \(D_2\). The ring \(A\) is a two-dimensional
noetherian local domain, \(a\) and \(b\) are nonzero, and the ideal
\((a,b)\) is \(\mathfrak m\)-primary. A prime ideal
\(\mathfrak q\neq\mathfrak m\) of \(A\) containing \(a\) has height one and does not contain \(b\). Let \(\mathcal Q(a)\) be
the set of these primes. The coefficient of
\([\overline{\{x\}}]\) in \(D_2\cdot(D_1\cdot[W])\) is
\[
 c(a,b)=\sum_{\mathfrak q\in\mathcal Q(a)}
   \operatorname{length}(A_{\mathfrak q}/aA_{\mathfrak q})\cdot
   \operatorname{length}\bigl(A/(\mathfrak q+bA)\bigr),
\]
and the coefficient in \(D_1\cdot(D_2\cdot[W])\) is \(c(b,a)\).

For a finite \(A\)-module \(N\) such that \(N/bN\) and
\(N[b]=\ker(b:N\to N)\) have finite length, put
\(e_b(N)=\operatorname{length}(N/bN)-\operatorname{length}(N[b])\).
By the snake lemma \cite[\stackstag{07JW}]{Stacks} and additivity
of length \cite[\stackstag{00IV}]{Stacks}, \(e_b\) is additive in
short exact sequences of such modules. Choose a filtration of
\(A/aA\) with successive quotients \(A/\mathfrak p_k\), where
\(\mathfrak p_k\) is a prime ideal containing \(a\)
\cite[\stackstag{00L0}]{Stacks}. Then
\(\mathfrak p_k\in\mathcal Q(a)\cup\{\mathfrak m\}\). For every
subquotient \(N\) in this filtration, \(N/bN\) and \(N[b]\) are
finite modules supported at \(\mathfrak m\), since \(aN=0\) and
\((a,b)\) is \(\mathfrak m\)-primary. They thus have finite length.
We have
\(e_b(A/\mathfrak m)=1-1=0\), and
\(e_b(A/\mathfrak q)=\operatorname{length}(A/(\mathfrak q+bA))\)
for \(\mathfrak q\in\mathcal Q(a)\), since \(b\) is a nonzero
element of the domain \(A/\mathfrak q\). Localizing the filtration
at \(\mathfrak q\) shows that \(\mathfrak q\) occurs
\(\operatorname{length}(A_{\mathfrak q}/aA_{\mathfrak q})\) times.
Therefore
\[
 c(a,b)=e_b(A/aA)
  =\operatorname{length}\bigl(A/(a,b)\bigr)
   -\operatorname{length}\bigl((aA:b)/aA\bigr).
\]
Since \(a\) and \(b\) are nonzerodivisors, the rule \(u\mapsto w\),
where \(bu=aw\), defines an isomorphism
\((aA:b)/aA\xrightarrow{\sim}(bA:a)/bA\), whose inverse is given by
the same rule with \(a\) and \(b\) exchanged. Hence
\(c(a,b)=c(b,a)\).
\end{proof}

\begin{definition}[Bloch's cycle complex]
\label{def:bloch-complex}
Let \(X\) be a scheme of finite type over \(\Z_{(\Sigma)}\), or the
spectrum of a localization of a finite-type
\(\Z_{(\Sigma)}\)-algebra. Assume that \(X\) has a finite open-and-closed
 decomposition on each piece of which \(\dim\mathcal O_{X,x}=c-\delta(x)\)
for a constant integer \(c\). This includes smooth schemes by
\eqref{eq:codim-delta}, their effective Cartier divisors
\cite[\stackstag{00KW}]{Stacks} (the local dimension drops by one),
and localizations of these schemes. For a smooth morphism \(Y\to X\) of constant relative
dimension \(r\), the flat dimension formula \cite[\stackstag{00ON}]{Stacks}
and the dimension formula on the smooth fibre give, at \(y\) over \(x\),
\(\dim\mathcal O_{Y,y}=\dim\mathcal O_{X,x}+r-\operatorname{trdeg}_{\kappa(x)}\kappa(y)\).
Since \(\delta(y)-\delta(x)=\operatorname{trdeg}_{\kappa(x)}\kappa(y)\),
the constant on the source is \(c+r\). In particular it holds
on products with simplex faces, and is preserved by \'etale pullback. Let \(j,m\in\mathbf N\). Put
\[
 \Delta^m_{\Z_{(\Sigma)}}=\operatorname{Spec}
       \Z_{(\Sigma)}[v_0,\ldots,v_m]/(v_0+\cdots+v_m-1).
\]
For a nonempty subset \(I\subseteq\{0,\ldots,m\}\), define the closed face
\[
 \Delta^I_{\Z_{(\Sigma)}}=
 \operatorname{Spec}\!\left(
 \Z_{(\Sigma)}[v_i\mid i\in I]\big/\bigl(\textstyle\sum_{i\in I}v_i-1\bigr)
 \right)\hookrightarrow\Delta^m_{\Z_{(\Sigma)}}.
\]
The closed immersion is induced by the surjective
\(\Z_{(\Sigma)}\)-algebra homomorphism sending the class of
\(v_i\) to the class of \(v_i\) for \(i\in I\), and to zero
for \(i\notin I\). Its ideal is generated by the coordinate
classes \(v_i\), \(i\notin I\). Increasing enumeration of \(I\)
identifies this face with \(\Delta^{|I|-1}_{\Z_{(\Sigma)}}\).

An integral closed subscheme
\(W\subseteq X\times_{\Z_{(\Sigma)}}\Delta^m_{\Z_{(\Sigma)}}\) is
\emph{admissible of codimension \(j\)} if
\(\operatorname{codim}(W,X\times_{\Z_{(\Sigma)}}\Delta^m_{\Z_{(\Sigma)}})=j\)
and if, for every nonempty \(I\subseteq\{0,\ldots,m\}\), every
irreducible component \(V\) of
\[
 W\times_{\Delta^m_{\Z_{(\Sigma)}}}\Delta^I_{\Z_{(\Sigma)}}
 \ \subseteq\ X\times_{\Z_{(\Sigma)}}\Delta^I_{\Z_{(\Sigma)}}
\]
satisfies
\(\operatorname{codim}(V,X\times_{\Z_{(\Sigma)}}\Delta^I_{\Z_{(\Sigma)}})\geq j\).
The group \(z^j(X,m)\) is the free abelian group on the admissible
integral closed subschemes of codimension \(j\) of
\(X\times_{\Z_{(\Sigma)}}\Delta^m_{\Z_{(\Sigma)}}\). A
\(\Z\)-linear combination of integral closed subschemes of
\(X\times_{\Z_{(\Sigma)}}\Delta^m_{\Z_{(\Sigma)}}\) is
\emph{admissible of codimension \(j\)} if it belongs to
\(z^j(X,m)\).

For \(m\geq1\) and \(0\leq i\leq m\), let
\(d_i:\Delta^{m-1}_{\Z_{(\Sigma)}}\hookrightarrow
\Delta^m_{\Z_{(\Sigma)}}\) be the face immersion for
\(I=\{0,\ldots,m\}\setminus\{i\}\), using this increasing enumeration.
Thus its pullback on coordinate classes sends \(v_r\) to
\(v_r\) if \(r<i\), to zero if \(r=i\), and to \(v_{r-1}\)
if \(r>i\). The image of \(\mathrm{id}_X\times d_i\) is the effective
Cartier divisor \(\{v_i=0\}\) of
\(X\times_{\Z_{(\Sigma)}}\Delta^m_{\Z_{(\Sigma)}}\), and
\[
 d_i^*=\{v_i=0\}\cdot(-):z^j(X,m)\longrightarrow z^j(X,m-1)
\]
is the Cartier restriction of \cref{not:codimension-delta}.
Indeed, on a piece with constant \(c\), an admissible \(W\) has
\(\delta(\eta_W)=c+m-j\). A component of its intersection with a face
cut out by \(r\) coordinates has codimension at most \(r\) in \(W\)
by the generalized principal ideal theorem \cite[\stackstag{0BBZ}]{Stacks}, hence has \(\delta\) at
least \(c+m-j-r\). Admissibility gives the opposite inequality,
so its codimension in the face is exactly \(j\). In particular no
component lies in a codimension-one face, and the components obtained
by restriction meet the remaining faces properly, by the original
face conditions. The schemes in question are catenary, being of finite
type over \(\Z_{(\Sigma)}\) or localizations of such schemes
\cite[\stackstag{02JB}, \stackstag{02J0}]{Stacks}. The Cartier restrictions commute on codimension-two faces
by \cref{sub:cartier-restrictions-commute}.
The differential is
\[
 \partial=\sum_{i=0}^m(-1)^i d_i^*:z^j(X,m)\longrightarrow z^j(X,m-1)
 \qquad(m\geq1),
\]
and \(\partial=0\) on \(z^j(X,0)\) \cite[\S3]{Geisser04}. For two
such schemes \(X_1\) and \(X_2\),
\[
 z^j(X_1\sqcup X_2,m)=z^j(X_1,m)\oplus z^j(X_2,m),
\]
compatibly with \(\partial\). We put \(z^j(X,m)=0\) when \(j<0\) or
\(m<0\), and we use the cohomological grading
\(\Z(j)^{\mathrm{cla}}(X)^q=z^j(X,2j-q)\).
\end{definition}

\begin{notation}[Operations on cycles]
\label{not:cycle-operations}
We use the flat pullbacks and proper pushforwards of
\cite[\S3]{Geisser04}. Let \(f:T'\to T\) be a finite morphism and
\(g:T'\to T\) a flat morphism of noetherian schemes. For integral
closed subschemes \(W'\subseteq T'\) and \(W\subseteq T\),
\[
 f_*[W']=[\kappa(\eta_{W'}):\kappa(\eta_{f(W')})]\cdot[f(W')],
 \qquad
 g^*[W]=\sum_V\operatorname{length}
     \bigl(\mathcal O_{W\times_TT',\eta_V}\bigr)\,[V],
\]
where \(V\) runs through the irreducible components of
\(W\times_TT'\). For morphisms of schemes as in
\cref{def:bloch-complex}, we apply these operations to the products
of \(f\) and \(g\) with \(\Delta^m_{\Z_{(\Sigma)}}\), whenever the
codimension and face conditions hold. In the applications to smooth
pullbacks and finite pushforwards these conditions follow from the
relative dimension formula and preservation of \(\delta\) under finite
maps, respectively; the appropriate change of codimension for a
finite map is the difference of the dimension constants.

For the following identities, assume that \(T,T'\) are of finite type
over \(\Z_{(\Sigma)}\), that the morphisms are over this base, and
that \(g\) is flat of fixed relative dimension \(r\). Let
\(D\subseteq T\) be a locally principal closed subscheme, with
inverse images \(f^{-1}(D)\) and \(g^{-1}(D)\). If
\(W'\not\subseteq f^{-1}(D)\) and \(W\not\subseteq D\), then
\begin{equation}\label{eq:cartier-push-pull}
 D\cdot f_*[W']=f_*\bigl(f^{-1}(D)\cdot[W']\bigr),
 \qquad
 g^{-1}(D)\cdot g^*[W]=g^*\bigl(D\cdot[W]\bigr).
\end{equation}
These assertions may be checked locally on \(T\), where \(D\) is the
zero scheme of a section of \(\mathcal O_T\). By flatness, every
component of \(g^{-1}(W)\) dominates \(W\)
\cite[\stackstag{02R8}]{Stacks}, so none is contained in
\(g^{-1}(D)\). Restricting to the finite morphism \(W'\to f(W')\)
and the flat morphism \(g^{-1}(W)\to W\), the two sides of each
identity have the same class in the Chow group of
\(f(W')\cap D\), respectively of \(g^{-1}(W\cap D)\), by
\cite[\stackstag{02TA}]{Stacks} and \cite[\stackstag{0B71}]{Stacks},
where the Gysin class of a cycle with no component in the divisor
is its Cartier restriction \cite[\stackstag{02T8}]{Stacks}. For the
first support all components have \(\delta\)-dimension
\(\delta(\eta_{W'})-1\), by the principal ideal theorem, the dimension
function of \cref{not:codimension-delta}, and finiteness.
For the second they have \(\delta\)-dimension
\(\delta(\eta_W)-1+r\), since the components of \(W\cap D\) have
\(\delta\)-dimension \(\delta(\eta_W)-1\) and flat pullback has fixed
relative dimension \(r\) \cite[\stackstag{02R8}, \stackstag{02RE}]{Stacks}.
The cycles in the respective
identities have exactly these dimensions, and the Chow groups of the
supports in these degrees are groups of cycles, since there are no
integral closed subschemes of dimension one larger
\cite[\stackstag{02RW}]{Stacks}.

The identities also hold for diagrams obtained from such finite-type
models by localizing the target and making the corresponding base
changes of the source and the divisor. Indeed, integral cycles on a
localization are restrictions of their closures in the model, and
restriction discards the components whose generic points do not
survive. The lengths and residue field degrees at the surviving
points are unchanged. The second identity also holds for any ordinary
localization morphism \(g:\Spec M^{-1}C\to\Spec C\) of noetherian
affine schemes: its two sides consist of precisely the same surviving
components of \(W\cap D\), with unchanged lengths. For the same
reason, finite pushforward commutes with such restriction to a
localization of the target and the induced localization of the source.
\end{notation}

\begin{sublemma}[Continuity of the cycle complex]
\label{sub:classical-cycle-spread}
Let \(C\) be a finite-type \(\Z_{(\Sigma)}\)-algebra such that
\(\Spec C\) satisfies the dimension condition of \cref{def:bloch-complex}.
Let \(S\subseteq S'\) be multiplicative subsets of \(C\), and let
\(j\in\mathbf N\). Flat pullback along the localizations
\(S^{-1}C[1/g]\to S'^{-1}C\), for \(g\in S'\), induces an
isomorphism of complexes
\[
 \varinjlim_{g\in S'}z^j\bigl(\Spec S^{-1}C[1/g],\bullet\bigr)
   \xrightarrow{\ \sim\ }z^j\bigl(\Spec S'^{-1}C,\bullet\bigr),
\]
where \(S'\) is preordered by divisibility.
\end{sublemma}

\begin{proof}
Geisser uses this continuity without proof: in the proof of
\cite[Corollary~3.3(a)]{Geisser04} the cycle complex is said to
commute with direct limits of rings, and the proofs of
\cite[Theorem~4.2 and Corollary~4.4]{Geisser04} reduce to local
schemes by the same property. We give the argument.

Since \(S\subseteq S'\), once the case \(S=\{1\}\) is known the
source is
\[
 \varinjlim_{g\in S'}\varinjlim_{s\in S}
     z^j\bigl(\Spec C[1/(sg)],\bullet\bigr)
 \simeq\varinjlim_{g\in S'}z^j\bigl(\Spec C[1/g],\bullet\bigr),
\]
by cofinality of multiplication \(S\times S'\to S'\) on the
localization preorders (for each \(g\), the objects whose product is
divisible by \(g\) form a nonempty filtered preorder).
Assume therefore that \(S=\{1\}\), and write
\(S\) for \(S'\). For \(m\in\mathbf N\) and a nonempty
\(I\subseteq\{0,\ldots,m\}\), let \(C_m\) and \(C_I\) be the
coordinate rings of \(\Spec C\times_{\Z_{(\Sigma)}}\Delta^m_{\Z_{(\Sigma)}}\)
and \(\Spec C\times_{\Z_{(\Sigma)}}\Delta^I_{\Z_{(\Sigma)}}\). We say
that a prime ideal of \(C_m\) or \(C_I\) \emph{meets} \(S\) if it
contains the image of an element of \(S\).

\smallskip\noindent
\emph{The map.}
The integral closed subschemes of \(\Spec S^{-1}C_m\) are the
\(V(\mathfrak qS^{-1}C_m)\), where \(\mathfrak q\) is a prime ideal
of \(C_m\) which does not meet \(S\). Let \(g\in S\), and let
\(\mathfrak q\) be a prime ideal of \(C_m\) with \(g\notin\mathfrak q\).
By \cref{not:cycle-operations}, the flat pullback of
\([V(\mathfrak qC_m[1/g])]\) to \(\Spec S^{-1}C_m\) is
\([V(\mathfrak qS^{-1}C_m)]\) if \(\mathfrak q\) does not meet
\(S\), and zero otherwise. The irreducible components of
\(V(\mathfrak qS^{-1}C_m)\times_{\Delta^m_{\Z_{(\Sigma)}}}\Delta^I_{\Z_{(\Sigma)}}\)
are the \(V(\mathfrak rS^{-1}C_I)\), where \(\mathfrak r\) runs
through the prime ideals of \(C_I\) minimal over \(\mathfrak qC_I\)
which do not meet \(S\). By \eqref{eq:codim-localization}, flat pullback
therefore preserves admissibility of codimension \(j\). It
commutes with the face maps \(d_i^*\) by
\eqref{eq:cartier-push-pull}.

\smallskip\noindent
\emph{Injectivity.}
Let \(g\in S\), and let
\(\alpha=\sum_kn_k[V(\mathfrak q_kC_m[1/g])]\in z^j(\Spec C[1/g],m)\)
have image zero, where the \(\mathfrak q_k\) are pairwise distinct
and \(n_k\neq0\). Then every \(\mathfrak q_k\) contains the image of
some \(s_k\in S\), and \(\alpha\) restricts to zero in
\(z^j(\Spec C[1/(g\prod_ks_k)],m)\).

\smallskip\noindent
\emph{Surjectivity.}
Let \(\mathfrak q\) be a prime ideal of \(C_m\) which does not meet
\(S\), such that \(V(\mathfrak qS^{-1}C_m)\) is admissible of
codimension \(j\). For each nonempty \(I\), let
\(\mathfrak r_{I,1},\ldots,\mathfrak r_{I,n_I}\) be the prime ideals
of \(C_I\) minimal over \(\mathfrak qC_I\) and of height less than \(j\). By admissibility and
\eqref{eq:codim-localization}, each \(\mathfrak r_{I,k}\) contains
the image of some \(s_{I,k}\in S\). Put \(g=\prod_{I,k}s_{I,k}\).
The irreducible components of
\(V(\mathfrak qC_m[1/g])\times_{\Delta^m_{\Z_{(\Sigma)}}}\Delta^I_{\Z_{(\Sigma)}}\)
correspond to the prime ideals of \(C_I\) minimal over
\(\mathfrak qC_I\) not containing \(g\), all of height at least \(j\). Thus
\([V(\mathfrak qC_m[1/g])]\in z^j(\Spec C[1/g],m)\), and its image
is \([V(\mathfrak qS^{-1}C_m)]\).
\end{proof}

\begin{sublemma}[Descent and localization for the cycle complex]
\label{sub:cycle-descent-localization}
Let \(j\in\mathbf N\).
\begin{enumerate}
\item Let \(X\) be a smooth scheme of finite type over \(\Z_{(\Sigma)}\).
Let \(\Z(j)^{\mathrm{cla}}_X\) be the complex of sheaves
\(U\mapsto\Z(j)^{\mathrm{cla}}(U)\) on the small Zariski site and on
the finite-type basis of the small Nisnevich site of \(X\), with
flat pullbacks as restriction maps. The degreewise sheaf property
follows from \cite[Lemma~3.1]{Geisser04}, applied on each piece of
constant relative dimension. The canonical maps
\[
 \Z(j)^{\mathrm{cla}}(X)\longrightarrow
 R\Gamma_{\mathrm{Zar}}\bigl(X,\Z(j)^{\mathrm{cla}}_X\bigr)
 \longrightarrow
 R\Gamma_{\mathrm{Nis}}\bigl(X,\Z(j)^{\mathrm{cla}}_X\bigr)
\]
are equivalences in the derived \(\infty\)-category
\(\mathcal D(\Z)\).
\item Let \(A\) be an essentially smooth \(\Z_{(\Sigma)}\)-algebra
and let \(h\in A\) be a nonzerodivisor. Let
\(i:\Spec(A/hA)\hookrightarrow\Spec A\) and
\(u:\Spec A[1/h]\hookrightarrow\Spec A\) be the closed and open
immersions. For \(j\geq1\) there is a long exact sequence
\begin{multline}\label{eq:cycle-localization-sequence}
 \cdots\longrightarrow
 H_m\bigl(z^{j-1}(\Spec(A/hA),\bullet)\bigr)
 \xrightarrow{\ i_*\ }
 H_m\bigl(z^j(\Spec A,\bullet)\bigr)\\
 \xrightarrow{\ u^*\ }
 H_m\bigl(z^j(\Spec A[1/h],\bullet)\bigr)
 \longrightarrow
 H_{m-1}\bigl(z^{j-1}(\Spec(A/hA),\bullet)\bigr)
 \longrightarrow\cdots .
\end{multline}
In cohomological grading, the first term is
\(H^{2j-m-2}\bigl(\Z(j-1)^{\mathrm{cla}}(\Spec(A/hA))\bigr)\) and the
second is \(H^{2j-m}\bigl(\Z(j)^{\mathrm{cla}}(\Spec A)\bigr)\).
\end{enumerate}
\end{sublemma}

\begin{proof}
\emph{Proof of \textup{(i)}.}
For an open covering \(X=U\cup V\), put \(Z=X\setminus U\),
with the reduced closed subscheme structure. The subcomplex of
\(z^j(X,\bullet)\) generated by admissible cycles supported on
\(Z\times\Delta^\bullet\) identifies by restriction with the
corresponding subcomplex of \(z^j(V,\bullet)\), since \(Z\subseteq V\).
On each connected component of \(X\), of relative dimension \(n\),
these are the dimension-indexed complexes \(z_{n+1-j}(Z,\bullet)\) on the
corresponding part of \(Z\), by \eqref{eq:codim-delta} and
\cite[\S1.8]{LevineLocalization}.
Levine's localization theorem over the semilocal base \(\Z_{(\Sigma)}\)
\cite[Theorem~1.7 and the identity following it]{LevineLocalization}
identifies the respective cofibers with \(z^j(U,\bullet)\) and
\(z^j(U\cap V,\bullet)\), by taking the direct sum over the connected
components of \(X\). Thus the commutative square of restrictions
\[
\begin{tikzcd}
 z^j(X,\bullet)\arrow[r]\arrow[d] & z^j(V,\bullet)\arrow[d] \\
 z^j(U,\bullet)\arrow[r] & z^j(U\cap V,\bullet)
\end{tikzcd}
\]
is homotopy cartesian, since the induced map between the fibres of
its vertical arrows identifies both fibres with the same complex
of cycles supported on \(Z\times\Delta^\bullet\). The Brown--Gersten criterion, as applied in
\cite[proof of Theorem~3.2(b)]{Geisser04} (see also \cite[Corollary~3.3(b)]{Geisser04}), gives
Zariski descent: the canonical map to the derived sections of the
associated complex of sheaves induces an equivalence on stalks, and
both presheaves of complexes satisfy Mayer--Vietoris on the
finite-dimensional noetherian Zariski site. This proves the first
equivalence. The second is \cite[Proposition~3.6]{Geisser04}.

\smallskip\noindent
\emph{Proof of \textup{(ii)}.}
Write \(A=S^{-1}C\), where \(C\) is a smooth finitely presented
\(\Z_{(\Sigma)}\)-algebra and \(S\subseteq C\) is a multiplicative
subset, and write \(h=h_0/s_0\) with \(h_0\in C\) and \(s_0\in S\).
The ideal \(\ker(h_0:C\to C)\) is finitely generated and its
localization at \(S\) is zero. After replacing \(C\) by \(C[1/s]\)
for a suitable \(s\in S\), we may assume that \(h_0\) is a
nonzerodivisor of \(C\). For \(s\in S\), put
\[
 X_s=\Spec C[1/s],\qquad D_s=\Spec\bigl(C[1/s]/(h_0)\bigr),
 \qquad U_s=\Spec C[1/(sh_0)] .
\]
For every integer \(q\), Levine's theorem
\cite[Theorem~1.7 and \S1.8]{LevineLocalization} (see also
\cite[Theorem~3.2(a)]{Geisser04} for the discrete valuation ring case)
gives a distinguished triangle
\[
 z_q(D_s,\bullet)\xrightarrow{\ i_*\ }z_q(X_s,\bullet)
   \xrightarrow{\ u^*\ }z_q(U_s,\bullet)
\]
of cycle complexes indexed by the dimension \(\delta\) of
\cref{not:codimension-delta}: the group \(z_q(T,m)\) is generated
by the integral closed subschemes of
\(T\times_{\Z_{(\Sigma)}}\Delta^m_{\Z_{(\Sigma)}}\) with
\(\delta=q+m\) at the generic point, subject to the corresponding
condition on faces. We may take the direct sum of these triangles
for the connected components of \(X_s\) and their intersections with
\(D_s\) and \(U_s\), choosing \(q\) separately on each component.
On such a component of relative dimension \(n\) over \(\Z_{(\Sigma)}\),
\eqref{eq:codim-delta} shows that an integral closed subscheme of
\(X_s\times_{\Z_{(\Sigma)}}\Delta^I_{\Z_{(\Sigma)}}\) has codimension
at least \(j\) if and only if \(\delta\leq n+|I|-j\) at its generic
point. Since \(h_0\) is a nonzerodivisor,
\(\dim\mathcal O_{D_s\times\Delta^I,x}=\dim\mathcal O_{X_s\times\Delta^I,x}-1\)
for every point \(x\) of
\(D_s\times_{\Z_{(\Sigma)}}\Delta^I_{\Z_{(\Sigma)}}\)
\cite[\stackstag{00KW}]{Stacks}. An integral closed subscheme of
\(D_s\times_{\Z_{(\Sigma)}}\Delta^I_{\Z_{(\Sigma)}}\) therefore has
codimension at least \(j-1\) if and only if \(\delta\leq n+|I|-j\)
at its generic point. Consequently, on this component and for
\(q=n+1-j\), the triangle is
\[
 z^{j-1}(D_s,\bullet)\xrightarrow{\ i_*\ }z^j(X_s,\bullet)
   \xrightarrow{\ u^*\ }z^j(U_s,\bullet).
\]
Its long exact sequences of homology groups are compatible with
the restrictions from \(s\) to a multiple of \(s\). Their filtered
colimit over \(s\in S\) is exact, and
\cref{sub:classical-cycle-spread}, applied to the algebras
\(C/(h_0)\), \(C\) and \(C[1/h_0]\) and to the images of \(S\),
identifies it with
\eqref{eq:cycle-localization-sequence}.
\end{proof}

\section{A finite \texorpdfstring{\'etale}{etale} comparison for connective \texorpdfstring{$K$}{K}-theory}
\label{sec:finite-connective-comparison}

The comparisons in \cref{thm:local-finite-connective-comparison,thm:semilocal-finite-connective-comparison}
recover each nonnegative connective \(K\)-group with invertible
finite coefficients from the corresponding finite Postnikov
truncation on the small \'etale site.
\Cref{lem:classical-grade-normalization} identifies the graded
unit maps with the canonical motivic-to-\'etale comparison.
On semilocal rings, these are the cohomological truncation maps
of \cref{prop:filtered-etale-comparison}.
Their degree bounds then give the comparison by the
finite-filtration argument of
\cref{sub:finite-filtration-comparison}.

\subsection{\texorpdfstring{\'Etale $K$}{Etale K}-theory and base change}

\begin{notation}[\'Etale sites and sheaves of spectra]
\label{not:etale-sites-spectra}
Let \(S\) be a quasi-compact quasi-separated scheme. Fix a big \'etale site
\((\mathit{Sch}/S)_{\mathrm{\acute et}}\) of \(S\), chosen large
enough to contain the schemes under consideration
\cite[\stackstag{021A}]{Stacks}. For an object \(X\) of this site,
\((\mathit{Sch}/X)_{\mathrm{\acute et}}\) is its big \'etale site,
that is, the localization of
\((\mathit{Sch}/S)_{\mathrm{\acute et}}\) at \(X\), and
\(X_{\mathrm{\acute et}}\) is its small \'etale site
\cite[\stackstag{021B}]{Stacks}. For a ring \(A\), we also write
\(A_{\mathrm{\acute et}}\) for \((\Spec A)_{\mathrm{\acute et}}\).
In this section we use the full
subsites of quasi-compact quasi-separated schemes, with the induced
\'etale topology, for both big and small sites. These subsites contain
the affine bases: restriction gives the same categories of sheaves
(by the comparison lemma for a basis, cf.\
\cite[\stackstag{021E}, \stackstag{04HR}]{Stacks}; for sheaves of spaces
one uses descent for covering sieves as in
\cite[Proposition~6.2.2.7]{HTT}). The same holds for sheaves of
spectra, by considering their constituent pointed spaces, and for
hypercomplete sheaves. In particular the small site of \(X\) still
contains its final object \(X\).
Let \((\mathit{Sm}/S)_{\mathrm{\acute et}}\subseteq
(\mathit{Sch}/S)_{\mathrm{\acute et}}\) be the full subcategory of
smooth \(S\)-schemes with the induced topology. Thus its objects
are quasi-compact and quasi-separated.

For a site \(\mathcal C\), let
\[
 \mathrm{Sh}^{\wedge}(\mathcal C;\mathrm{Sp})\subseteq
 \mathrm{Sh}(\mathcal C;\mathrm{Sp})\subseteq
 \mathrm{PSh}(\mathcal C;\mathrm{Sp})
 =\operatorname{Fun}(\mathcal C^{\mathrm{op}},\mathrm{Sp})
\]
be the \(\infty\)-categories of hypercomplete sheaves, of sheaves,
and of presheaves of spectra on \(\mathcal C\). For a presheaf
\(\mathcal P\) and \(n\in\Z\), let \(\pi_n\mathcal P\) be the sheaf of
abelian groups on \(\mathcal C\) associated with
\(U\mapsto\pi_n(\mathcal P(U))\). A morphism \(\alpha\) of presheaves
is \emph{\(\infty\)-connective} if \(\pi_n(\alpha)\) is an isomorphism
for every \(n\in\Z\), and a sheaf \(\mathcal F\) is
\emph{hypercomplete} if
\(\operatorname{Map}(\beta,\mathcal F)\) is an equivalence for every
\(\infty\)-connective morphism \(\beta\) of sheaves; compare
\cite[\S6.5.2]{HTT} for sheaves of spaces. The two inclusions into
\(\mathrm{PSh}(\mathcal C;\mathrm{Sp})\) have left adjoints
\[
 L_{\mathcal C}:\mathrm{PSh}(\mathcal C;\mathrm{Sp})\longrightarrow
     \mathrm{Sh}(\mathcal C;\mathrm{Sp}),\qquad
 L^{\wedge}_{\mathcal C}:\mathrm{PSh}(\mathcal C;\mathrm{Sp})
     \longrightarrow\mathrm{Sh}^{\wedge}(\mathcal C;\mathrm{Sp}),
\]
called \emph{sheafification} and \emph{hypercomplete sheafification};
their units \(\mathcal P\to L_{\mathcal C}\mathcal P\) and
\(\mathcal P\to L^{\wedge}_{\mathcal C}\mathcal P\) are
\(\infty\)-connective. By definition, an \(\infty\)-connective
morphism of hypercomplete sheaves is an equivalence. Hence a morphism
\(\alpha:\mathcal P\to\mathcal F\) from a presheaf to a hypercomplete
sheaf is \(\infty\)-connective if and only if the induced morphism
\(L^{\wedge}_{\mathcal C}\mathcal P\to\mathcal F\) is an equivalence.

For a functor \(u:\mathcal C'\to\mathcal C\), let
\(u^*:\mathrm{PSh}(\mathcal C;\mathrm{Sp})\to
\mathrm{PSh}(\mathcal C';\mathrm{Sp})\) be
\(\mathcal P\mapsto\mathcal P\circ u^{\mathrm{op}}\). If \(\mathcal C\)
is \((\mathit{Sch}/S)_{\mathrm{\acute et}}\) or
\((\mathit{Sm}/S)_{\mathrm{\acute et}}\), \(X\) is an object of \(\mathcal C\),
and \(u:X_{\mathrm{\acute et}}\to\mathcal C\) is
\((U\to X)\mapsto(U\to X\to S)\), we write
\(\mathcal P|_{X_{\mathrm{\acute et}}}\) for \(u^*\mathcal P\). For a
sheaf of spectra \(\mathcal F\) on \(X_{\mathrm{\acute et}}\), put
\[
 R\Gamma_{\mathrm{\acute et}}(X,\mathcal F)=\mathcal F(X),
\]
and write \(R\Gamma_{\mathrm{\acute et}}(A,\mathcal F)\) when
\(X=\Spec A\). For a bounded-below complex of abelian sheaves this is
the usual complex of derived global sections, under the
Eilenberg--Mac Lane construction of \cref{sec:conventions}.
\end{notation}

\begin{sublemma}[Restriction between \'etale sites]
\label{sub:big-site-restriction}
Let \(S\) be a quasi-compact quasi-separated scheme, and use the sites
of \cref{not:etale-sites-spectra}. Let \(u:\mathcal C'\to\mathcal C\) be one
of the four functors
\[
\begin{aligned}
 (\mathit{Sm}/S)_{\mathrm{\acute et}}&\hookrightarrow
 (\mathit{Sch}/S)_{\mathrm{\acute et}},&\qquad
 (\mathit{Sch}/X)_{\mathrm{\acute et}}&\longrightarrow
 (\mathit{Sch}/S)_{\mathrm{\acute et}},\\
 X_{\mathrm{\acute et}}&\longrightarrow
 (\mathit{Sch}/S)_{\mathrm{\acute et}},&\qquad
 X_{\mathrm{\acute et}}&\longrightarrow
 (\mathit{Sm}/S)_{\mathrm{\acute et}},
\end{aligned}
\]
where \(X\) is an object of the target site,
\((\mathit{Sch}/X)_{\mathrm{\acute et}}\) is the big \'etale site of
\(X\), and the last three functors are
\((U\to X)\mapsto(U\to X\to S)\).
\begin{enumerate}
\item The functor \(u^*\) carries sheaves to sheaves and hypercomplete
sheaves to hypercomplete sheaves, and the squares
\[
\begin{tikzcd}[column sep=large]
 \mathrm{PSh}(\mathcal C;\mathrm{Sp})
   \arrow[r,"L_{\mathcal C}"]\arrow[d,"u^*"'] &
 \mathrm{Sh}(\mathcal C;\mathrm{Sp})\arrow[d,"u^*"] \\
 \mathrm{PSh}(\mathcal C';\mathrm{Sp})\arrow[r,"L_{\mathcal C'}"] &
 \mathrm{Sh}(\mathcal C';\mathrm{Sp})
\end{tikzcd}
\qquad
\begin{tikzcd}[column sep=large]
 \mathrm{PSh}(\mathcal C;\mathrm{Sp})
   \arrow[r,"L^{\wedge}_{\mathcal C}"]\arrow[d,"u^*"'] &
 \mathrm{Sh}^{\wedge}(\mathcal C;\mathrm{Sp})\arrow[d,"u^*"] \\
 \mathrm{PSh}(\mathcal C';\mathrm{Sp})
   \arrow[r,"L^{\wedge}_{\mathcal C'}"] &
 \mathrm{Sh}^{\wedge}(\mathcal C';\mathrm{Sp})
\end{tikzcd}
\]
commute: for every presheaf \(\mathcal P\) on \(\mathcal C\), the maps
\(u^*\mathcal P\to u^*L_{\mathcal C}\mathcal P\) and
\(u^*\mathcal P\to u^*L^{\wedge}_{\mathcal C}\mathcal P\) induce
equivalences
\(L_{\mathcal C'}u^*\mathcal P\simeq u^*L_{\mathcal C}\mathcal P\) and
\(L^{\wedge}_{\mathcal C'}u^*\mathcal P\simeq
u^*L^{\wedge}_{\mathcal C}\mathcal P\).
\item For every sheaf \(\mathcal F\) on \(\mathcal C\) and all
\(a,b,n\in\Z\), there are canonical identifications
\[
 u^*\tau_{\geq a}\mathcal F\simeq\tau_{\geq a}u^*\mathcal F,\qquad
 u^*\tau_{\leq b}\mathcal F\simeq\tau_{\leq b}u^*\mathcal F,\qquad
 u^*\pi_n\mathcal F\cong\pi_nu^*\mathcal F .
\]
\item Let \(\mathcal C\) be \((\mathit{Sch}/S)_{\mathrm{\acute et}}\)
or \((\mathit{Sm}/S)_{\mathrm{\acute et}}\). For every sheaf
\(\mathcal F\) on \(\mathcal C\) and every object \(X\) of
\(\mathcal C\),
\begin{equation}\label{eq:big-site-evaluation}
 \mathcal F(X)=R\Gamma_{\mathrm{\acute et}}
       \bigl(X,\mathcal F|_{X_{\mathrm{\acute et}}}\bigr).
\end{equation}
A morphism \(\alpha\) of sheaves on \(\mathcal C\) is an equivalence
if and only if \(\alpha|_{X_{\mathrm{\acute et}}}\) is an equivalence
for every object \(X\) of \(\mathcal C\). A sheaf \(\mathcal F\) on
\(\mathcal C\) is hypercomplete if and only if
\(\mathcal F|_{X_{\mathrm{\acute et}}}\) is hypercomplete for every
object \(X\) of \(\mathcal C\).
\item Let \(\mathcal C\) be one of the sites
\((\mathit{Sch}/S)_{\mathrm{\acute et}}\),
\((\mathit{Sm}/S)_{\mathrm{\acute et}}\), \(X_{\mathrm{\acute et}}\),
and let \(\alpha:\mathcal P\to\mathcal Q\) be a morphism of presheaves
of spectra on \(\mathcal C\) such that \(\alpha(U)\) is an equivalence
for every affine object \(U\) of \(\mathcal C\). Then
\(L^{\wedge}_{\mathcal C}(\alpha)\) is an equivalence.
\end{enumerate}
\end{sublemma}

\begin{proof}
\emph{Sheaves of sets and of abelian groups.}
The functor \(X_{\mathrm{\acute et}}\to
(\mathit{Sch}/X)_{\mathrm{\acute et}}\) is continuous and
cocontinuous, and on sheaves of sets restriction along it is both the
inverse image \(i_X^{-1}\) and the direct image \(\pi_{X,*}\) of
morphisms of topoi \cite[\stackstag{021G}]{Stacks}. The functor
\((\mathit{Sch}/X)_{\mathrm{\acute et}}\to
(\mathit{Sch}/S)_{\mathrm{\acute et}}\) is the localization at the
object \(X\); on sheaves of sets restriction is the inverse image of
the localization morphism, and its left adjoint forgets the map to
(the sheaf represented by) \(X\) \cite[\stackstag{00XZ}]{Stacks}.
An \'etale covering of a smooth \(S\)-scheme consists of smooth
\(S\)-schemes; for the quasi-compact quasi-separated sites every
covering admits a refinement by affine objects. Hence each of the
four functors \(u\) is continuous and cocontinuous, and for every
presheaf of abelian groups
\(\mathcal A\) on \(\mathcal C\) the canonical map
\begin{equation}\label{eq:restriction-sheafification-abelian}
 (\mathcal A\circ u^{\mathrm{op}})^{\#}\longrightarrow
 \mathcal A^{\#}\circ u^{\mathrm{op}}
\end{equation}
is an isomorphism: by \cite[\stackstag{00XM}]{Stacks} the two
presheaves have isomorphic sheafifications, and the right-hand side
is already a sheaf by \cite[\stackstag{00WW}]{Stacks}.

\smallskip\noindent
\emph{Proof of \textup{(i)} and \textup{(ii)}, except for
hypercompleteness on the smooth site.}
The functor \(u^*\) carries sheaves of spaces to sheaves, by
\v{C}ech descent for \'etale coverings. It also commutes with their
sheafification: the restriction of a covering-sieve inclusion in
\(\mathcal C\) becomes an equivalence after sheafification in
\(\mathcal C'\), since its pullback to an object of \(\mathcal C'\)
admits a covering refinement there. Since restriction preserves
colimits, it sends the local equivalences generated by these sieve
inclusions to local equivalences \cite[Proposition~6.2.2.7]{HTT}. The same assertions
hold for spectra, by stabilization. By
\eqref{eq:restriction-sheafification-abelian}, applied to
\(\mathcal A=(U\mapsto\pi_n\mathcal P(U))\), it commutes with
homotopy sheaves, hence preserves \(\infty\)-connective morphisms.
The truncation \(\tau_{\leq b}\mathcal F\) is the sheafification of
\(U\mapsto\tau_{\leq b}(\mathcal F(U))\), and
\(\tau_{\geq a}\mathcal F\) is the fibre of
\(\mathcal F\to\tau_{\leq a-1}\mathcal F\); this gives
\textup{(ii)}.

For the last three functors \(u\), we describe the left adjoint
on spectral sheaves. Write \(u_!\) for left Kan extension along
\(u^{\mathrm{op}}\) on presheaves; the left adjoint on sheaves is
\(L_{\mathcal C}u_!\). In the case of the localization at \(X\),
\((u_!\mathcal P)(T)=\bigoplus_{h:T\to X}\mathcal P(T\xrightarrow{h}X)\),
where \(h\) runs through the maps over \(S\). For either functor
from the small site of \(X\), the pointwise formula for left Kan
extension \cite[Lemma~4.3.2.13]{HTT} gives
\[
 (u_!\mathcal P)(T)=\bigoplus_{h:T\to X}
 \operatorname*{colim}_{(T\to V\to X)^{\mathrm{op}}}\mathcal P(V\to X),
\]
where the factorizations of \(h\) are through objects \(V\) of
\(X_{\mathrm{\acute et}}\). For each \(h\), their category has
finite limits. The final object is the factorization through \(X\)
itself, and fibre products are computed in schemes over \(X\)
\cite[\stackstag{021D}]{Stacks}. The equalizer of two maps
\(V\rightrightarrows W\) over \(X\) is the pullback of the diagonal
\(W\to W\times_XW\), which is an open immersion since \(W\to X\)
is unramified \cite[\stackstag{02GE}]{Stacks}. These limits remain
quasi-compact and quasi-separated: morphisms between quasi-compact
quasi-separated schemes are quasi-compact and quasi-separated
\cite[\stackstag{01KU}, \stackstag{01KV}, \stackstag{03GI}]{Stacks},
and the equalizer is a quasi-compact open subscheme of \(V\).
The opposite of the category of factorizations is therefore filtered.

Let \(\mathcal P\) have zero homotopy sheaves. In either formula,
each element of \(\pi_i(u_!\mathcal P)(T)\) comes from finitely many
classes in \(\pi_i\mathcal P(V\to X)\) (with \(V=T\) in the slice
case), since homotopy groups of spectra commute with sums and
filtered colimits. Each such class vanishes on an \'etale covering
of \(V\). Pulling back these coverings along \(T\to V\), taking a
common refinement and, if necessary, affine open coverings, kills
the element locally on \(T\). This proves that \(L_{\mathcal C}u_!\)
sends presheaves with zero homotopy sheaves to sheaves with zero
homotopy sheaves. Applied to cofibers, it follows that this left
adjoint on sheaves preserves \(\infty\)-connective morphisms. Its
right adjoint \(u^*\) thus preserves hypercomplete sheaves by
\cref{not:etale-sites-spectra}.

If \(u^*\) preserves hypercomplete sheaves, then for a presheaf
\(\mathcal P\) the map
\(u^*\mathcal P\to u^*L^{\wedge}_{\mathcal C}\mathcal P\) is
\(\infty\)-connective with hypercomplete target, so
\(L^{\wedge}_{\mathcal C'}u^*\mathcal P\simeq
u^*L^{\wedge}_{\mathcal C}\mathcal P\) by
\cref{not:etale-sites-spectra}. This proves \textup{(i)} for the last
three functors.

\smallskip\noindent
\emph{Proof of \textup{(iii)}.}
Formula \eqref{eq:big-site-evaluation} is the definition of
\(R\Gamma_{\mathrm{\acute et}}\) in \cref{not:etale-sites-spectra}.
It shows that a morphism \(\alpha\) of sheaves on \(\mathcal C\) is an
equivalence as soon as all the \(\alpha|_{X_{\mathrm{\acute et}}}\)
are. Let \(\mathcal F\) be a sheaf on \(\mathcal C\) all of whose
restrictions \(\mathcal F|_{X_{\mathrm{\acute et}}}\) are
hypercomplete, and let
\(\eta:\mathcal F\to L^{\wedge}_{\mathcal C}\mathcal F\) be the unit.
For every object \(X\), the morphism \(\eta|_{X_{\mathrm{\acute et}}}\)
is \(\infty\)-connective by \textup{(ii)}, and its target is
hypercomplete by \textup{(i)} for the restriction to the small
site of \(X\). It is therefore an equivalence, so \(\eta\) is an
equivalence and \(\mathcal F\) is hypercomplete. The converse
follows from the preservation of hypercompleteness in \textup{(i)}
for the functors from the small site.

\smallskip\noindent
\emph{Hypercompleteness on the smooth site.}
Let \(\mathcal F\) be a hypercomplete sheaf on
\((\mathit{Sch}/S)_{\mathrm{\acute et}}\), and let \(\mathcal F'\) be
its restriction to \((\mathit{Sm}/S)_{\mathrm{\acute et}}\). For a
smooth \(S\)-scheme \(X\), the sheaf
\(\mathcal F'|_{X_{\mathrm{\acute et}}}=
\mathcal F|_{X_{\mathrm{\acute et}}}\) is hypercomplete. By
\textup{(iii)}, \(\mathcal F'\) is hypercomplete. The last paragraph
of the proof of \textup{(i)} now applies to the inclusion of the
smooth site.

\smallskip\noindent
\emph{Proof of \textup{(iv)}.}
For \(n\in\Z\), let \(\mathcal A\to\mathcal B\) be the map
\(U\mapsto\pi_n(\alpha(U))\) of presheaves of abelian groups. Its
kernel and cokernel presheaves vanish on affine objects. Every object
of \(\mathcal C\) has a Zariski covering by affine objects of
\(\mathcal C\), so their sheafifications vanish, and
\(\pi_n(\alpha)\) is an isomorphism because sheafification is exact.
Thus \(\alpha\) is \(\infty\)-connective, and so is
\(L^{\wedge}_{\mathcal C}(\alpha)\), which is therefore an
equivalence by \cref{not:etale-sites-spectra}.
\end{proof}
\begin{construction}[\'Etale \(K\)-theory and base change]
\label{con:etale-k-sheaves}
Let \(N\geq2\) be a natural number, and use
\cref{not:etale-sites-spectra} with \(S=\Spec\Z[1/N]\). Let
\(K^{\mathrm{cn}}/N\) be the presheaf of spectra
\(U\mapsto K^{\mathrm{cn}}(U)/N\) on
\((\mathit{Sch}/\Z[1/N])_{\mathrm{\acute et}}\). For an object \(X\)
of this site, that is, a quasi-compact quasi-separated scheme on which
\(N\) is invertible, put
\[
 \mathcal K^{\mathrm{big}}_N
   =L^{\wedge}_{(\mathit{Sch}/\Z[1/N])_{\mathrm{\acute et}}}
       \bigl(K^{\mathrm{cn}}/N\bigr),
 \qquad
 \mathcal K^{\mathrm{\acute et}}_{X,N}
   =L^{\wedge}_{X_{\mathrm{\acute et}}}
       \bigl((K^{\mathrm{cn}}/N)|_{X_{\mathrm{\acute et}}}\bigr).
\]
We write \(\mathcal K^{\mathrm{\acute et}}_{A,N}\) when
\(X=\Spec A\). By \cref{sub:big-site-restriction}\textup{(iv)},
these sheaves depend only on the values of \(K^{\mathrm{cn}}/N\) on
affine schemes. In formulas for sheaves of spectra, an abelian
sheaf denotes its Eilenberg--Mac Lane sheaf.

By \cref{sub:big-site-restriction}\textup{(i)}, the unit maps induce
a canonical equivalence
\begin{equation}\label{eq:big-small}
 \mathcal K^{\mathrm{\acute et}}_{X,N}\xrightarrow{\ \sim\ }
 \mathcal K^{\mathrm{big}}_N\big|_{X_{\mathrm{\acute et}}}.
\end{equation}
For \(-\infty\leq a\leq b\leq+\infty\),
\cref{sub:big-site-restriction}\textup{(ii)} and
\eqref{eq:big-site-evaluation} then give
\begin{equation}\label{eq:big-small-truncated}
 R\Gamma_{\mathrm{\acute et}}\bigl(X,
      \tau_{\geq a}\tau_{\leq b}
      \mathcal K^{\mathrm{\acute et}}_{X,N}\bigr)
 \simeq
 \bigl(\tau_{\geq a}\tau_{\leq b}
             \mathcal K^{\mathrm{big}}_N\bigr)(X).
\end{equation}
Thus the left-hand side of \eqref{eq:big-small-truncated} is
contravariantly functorial in all morphisms of such schemes. The presheaf \(K^{\mathrm{cn}}/N\) is connective, and so is its
hypercomplete sheafification. For \(a\leq0\leq b\), there is therefore
a canonical equivalence
\(\tau_{\geq a}\tau_{\leq b}\mathcal K^{\mathrm{big}}_N
\simeq\tau_{\leq b}\mathcal K^{\mathrm{big}}_N\).
The composite
\[
 K^{\mathrm{cn}}/N\longrightarrow\mathcal K^{\mathrm{big}}_N
   \longrightarrow\tau_{\leq b}\mathcal K^{\mathrm{big}}_N
\]
and this equivalence, evaluated using
\eqref{eq:big-small-truncated}, give a natural transformation
\begin{equation}\label{eq:truncated-unit}
 K^{\mathrm{cn}}(X)/N\longrightarrow
 R\Gamma_{\mathrm{\acute et}}\bigl(X,
      \tau_{\geq a}\tau_{\leq b}
      \mathcal K^{\mathrm{\acute et}}_{X,N}\bigr)
 \qquad(a\leq0\leq b).
\end{equation}
For \(n\in\mathbf N\), the \emph{finite comparison map} is the case
\(a=-\infty\), \(b=2n\) of \eqref{eq:truncated-unit}:
\begin{equation}\label{eq:finite-comparison-map}
 K^{\mathrm{cn}}(X)/N\longrightarrow
 R\Gamma_{\mathrm{\acute et}}\bigl(X,
      \tau_{\leq2n}\mathcal K^{\mathrm{\acute et}}_{X,N}\bigr).
\end{equation}

Let \(f:Y\to X\) be a morphism of quasi-compact quasi-separated
schemes on which \(N\) is invertible, let \(f^*\) and \(f_*\) be the inverse and direct image
functors between sheaves of spectra on \(X_{\mathrm{\acute et}}\) and
\(Y_{\mathrm{\acute et}}\), and put
\(f^{*,h}=L^{\wedge}_{Y_{\mathrm{\acute et}}}\circ f^*\). For an
object \(U\) of \(X_{\mathrm{\acute et}}\), the restriction maps
\[
 \mathcal K^{\mathrm{big}}_N(U)\longrightarrow
 \mathcal K^{\mathrm{big}}_N(U\times_XY)
\]
of the big sheaf define through \eqref{eq:big-small} a morphism
\(\mathcal K^{\mathrm{\acute et}}_{X,N}\to
f_*\mathcal K^{\mathrm{\acute et}}_{Y,N}\) of sheaves on
\(X_{\mathrm{\acute et}}\). The direct image \(f_*\) preserves
hypercomplete sheaves, because its left adjoint \(f^*\) preserves
\(\infty\)-connective morphisms. By adjunction and the universal
property of \(L^{\wedge}_{Y_{\mathrm{\acute et}}}\), this morphism
corresponds to a morphism
\begin{equation}\label{eq:k-sheaf-base-change}
 f^{*,h}\mathcal K^{\mathrm{\acute et}}_{X,N}\longrightarrow
 \mathcal K^{\mathrm{\acute et}}_{Y,N}.
\end{equation}
Since the restriction maps are those of the single functor
\(\mathcal K^{\mathrm{big}}_N\), the morphisms
\eqref{eq:k-sheaf-base-change} for a composite \(Y'\to Y\to X\) are
coherently identified with the composites of the morphisms for
\(Y'\to Y\) and \(Y\to X\). The unit
\(K^{\mathrm{cn}}/N\to\mathcal K^{\mathrm{big}}_N\) is a natural
transformation of presheaves on the big site; its naturality for
\(U\times_XY\to U\) uses the pullback of perfect complexes along
that morphism. Evaluating at \(U=X\), and using
\cref{sub:big-site-restriction}\textup{(ii)} for the truncations, shows that
\eqref{eq:k-sheaf-base-change} is compatible with the maps
\eqref{eq:truncated-unit} for \(X\) and \(Y\).
\end{construction}
\begin{sublemma}[The homotopy sheaves of finite \'etale \(K\)-theory]
\label{sub:finite-k-etale-coefficients}
Let \(X\) be an affine scheme and let \(N\geq2\) be a natural number
invertible on \(X\).
\begin{enumerate}
\item The sheaf \(\mathcal K^{\mathrm{\acute et}}_{X,N}\) of
\cref{con:etale-k-sheaves} is connective:
\(\pi_n\mathcal K^{\mathrm{\acute et}}_{X,N}=0\) for \(n<0\).
\item Suppose that \(N=\ell^\nu\), where \(\ell\) is prime and
\(\nu\geq1\). For every \(j\in\mathbf N\) there is a canonical
isomorphism of sheaves on \(X_{\mathrm{\acute et}}\)
\[
 \lambda_j:\Z/N(j)\xrightarrow{\ \sim\ }
 \pi_{2j}\mathcal K^{\mathrm{\acute et}}_{X,N},
 \qquad\text{and}\qquad
 \pi_{2j+1}\mathcal K^{\mathrm{\acute et}}_{X,N}=0 .
\]
\item In the situation of \textup{(ii)}, let \(f:Y\to X\) be a
morphism of affine schemes. For every \(j\in\mathbf N\), the square
\[
\begin{tikzcd}[column sep=huge]
 f^*\Z/N(j)\arrow[r,"f^*\lambda_j"]\arrow[d,"\simeq"'] &
 f^*\pi_{2j}\mathcal K^{\mathrm{\acute et}}_{X,N}
   \arrow[d,"\pi_{2j}\eqref{eq:k-sheaf-base-change}"] \\
 \Z/N(j)\arrow[r,"\lambda_j"] &
 \pi_{2j}\mathcal K^{\mathrm{\acute et}}_{Y,N}
\end{tikzcd}
\]
commutes, where the left vertical map is the canonical isomorphism
of Tate sheaves.
\end{enumerate}
\end{sublemma}

\begin{proof}
\emph{Stalks.}
Let \(\bar x\) be a geometric point of \(X\), let
\(\mathcal O^{\mathrm{sh}}_{X,\bar x}\) be the strict henselian local
ring at \(\bar x\), and let \(\kappa\) be its residue field. The
ring \(\mathcal O^{\mathrm{sh}}_{X,\bar x}\) is the filtered colimit
of the rings \(\mathcal O(U)\) over the affine \'etale neighborhoods
\((U,\bar u)\) of \(\bar x\) \cite[\stackstag{04HX}]{Stacks}, and
\(K^{\mathrm{cn}}(-)/N\) commutes with filtered colimits of rings.
Since the stalk of a presheaf agrees with the stalk of its
sheafification \cite[\stackstag{03PT}]{Stacks}, for every \(n\in\Z\)
\begin{equation}\label{eq:k-sheaf-stalk}
 \bigl(\pi_n\mathcal K^{\mathrm{\acute et}}_{X,N}\bigr)_{\bar x}
 =\varinjlim_{(U,\bar u)}\pi_n\bigl(K^{\mathrm{cn}}(U)/N\bigr)
 =\pi_n\bigl(K^{\mathrm{cn}}
     (\mathcal O^{\mathrm{sh}}_{X,\bar x})/N\bigr)
 \xrightarrow{\ \sim\ }\pi_n\bigl(K^{\mathrm{cn}}(\kappa)/N\bigr).
\end{equation}
The last map is induced by
\(\mathcal O^{\mathrm{sh}}_{X,\bar x}\to\kappa\); it is an
isomorphism by Gabber's rigidity theorem, in the connective form for
a henselian pair on which \(N\) is invertible
\cite{Gabber92} (see \cite[Theorem~1.4]{CMM21}). The groups in \eqref{eq:k-sheaf-stalk}
vanish for \(n<0\), and a sheaf of abelian groups on
\(X_{\mathrm{\acute et}}\) with zero stalks is zero
\cite[\stackstag{03PU}]{Stacks}. This proves \textup{(i)}.

\smallskip\noindent
\emph{Degree zero.}
The rank of finitely generated projective modules defines a map of
presheaves \(K_0^{\mathrm{cn}}(U)/N\to\Gamma(U,\Z/N)\) on the affine
objects \(U\) of \(X_{\mathrm{\acute et}}\). Every finitely
generated projective module over the local ring
\(\mathcal O^{\mathrm{sh}}_{X,\bar x}\) is free
(\cref{sub:semilocal-projective-free}), so this map is an
isomorphism on the stalk at \(\bar x\). Its inverse after
sheafification is \(\lambda_0\). On regular schemes it agrees with
the isomorphism \(\lambda_0\) of \cite[\arxivhtml{2508.09915v1}{S6.E7}{Lemma~6.7}]{BEMMotivic25},
as seen from the normalization on the residue fields.

\smallskip\noindent
\emph{Positive degrees.}
Let \(N=\ell^\nu\). For a ring \(A\) and \(i\geq1\), the map
\(\pi_i(K^{\mathrm{cn}}(A)/N)\to\pi_i(K(A)/N)\) is an isomorphism:
in the coefficient sequences \eqref{eq:coefficient-sequence} of the
two theories, the outer groups \(K_i(A)/N\) and \(K_{i-1}(A)[N]\)
agree, because \(K^{\mathrm{cn}}_m(A)=K_m(A)\) for \(m\geq0\). Hence
\(\pi_i\mathcal K^{\mathrm{\acute et}}_{X,N}\) is the \'etale
sheaf associated with \(U\mapsto\pi_i(K(U)/N)\) for \(i\geq1\). For
these sheaves, \cite[\arxivhtml{2508.09915v1}{S6.E7}{Lemma~6.7}]{BEMMotivic25} provides the natural
isomorphisms \(\lambda_j\) for \(j\geq1\) and the vanishing in odd
degrees, on the \'etale site of every quasi-compact quasi-separated
scheme on which \(\ell\) is invertible. This proves \textup{(ii)}.

\smallskip\noindent
\emph{Base change.}
The rank map and the isomorphisms \(\lambda_j\) of
\cite[\arxivhtml{2508.09915v1}{S6.E7}{Lemma~6.7}]{BEMMotivic25} are natural in the scheme, and the
morphism \eqref{eq:k-sheaf-base-change} is induced by pullback of
perfect complexes (\cref{con:etale-k-sheaves}). The sheaf \(\mu_N\)
on \(Y_{\mathrm{\acute et}}\) is the pullback of the sheaf
\(\mu_N\) on \(X_{\mathrm{\acute et}}\), which gives the left
vertical isomorphism in \textup{(iii)} and the commutativity of the
square.
\end{proof}
\subsection{Smooth left Kan extensions}

\begin{notation}[Categories of smooth sources]
\label{not:smooth-source-index}
Let \(\mathrm{CAlg}_{\Z_{(\Sigma)}}\) be the category of commutative
\(\Z_{(\Sigma)}\)-algebras, and let
\[
 \mathrm{Sm}_\Sigma\subseteq\mathrm{EssSm}_\Sigma
   \subseteq\mathrm{CAlg}_{\Z_{(\Sigma)}}
\]
be the full subcategories of smooth finitely presented
\(\Z_{(\Sigma)}\)-algebras and of essentially smooth
\(\Z_{(\Sigma)}\)-algebras (\cref{def:essentially-smooth-sigma}).
For a \(\Z_{(\Sigma)}\)-algebra \(R\), let
\[
 \mathrm{Sm}_{\Sigma/R}\subseteq\mathrm{EssSm}_{\Sigma/R}
\]
be the categories of pairs \((D,\varphi)\), where \(D\) is an object
of \(\mathrm{Sm}_\Sigma\), respectively of \(\mathrm{EssSm}_\Sigma\),
and \(\varphi:D\to R\) is a homomorphism of
\(\Z_{(\Sigma)}\)-algebras. A morphism
\((D,\varphi)\to(D',\varphi')\) is a homomorphism of
\(\Z_{(\Sigma)}\)-algebras \(\psi:D\to D'\) with
\(\varphi'\circ\psi=\varphi\). We use these ordinary categories as
index categories through their nerves.

Suppose that \(R\) is nonzero and semilocal, with maximal ideals
\(\mathfrak m_1,\ldots,\mathfrak m_s\). Let
\[
 \mathrm{EssSm}^{\mathrm{sl}}_{\Sigma/R}\subseteq
 \mathrm{EssSm}^{\mathrm{inc}}_{\Sigma/R}\subseteq
 \mathrm{EssSm}_{\Sigma/R}
\]
be the following full subcategories. A pair \((D,\varphi)\) belongs
to \(\mathrm{EssSm}^{\mathrm{inc}}_{\Sigma/R}\) if the prime ideals
\(\varphi^{-1}(\mathfrak m_1),\ldots,\varphi^{-1}(\mathfrak m_s)\) of
\(D\) are pairwise incomparable. It belongs to
\(\mathrm{EssSm}^{\mathrm{sl}}_{\Sigma/R}\) if \(D\) is semilocal
and its maximal ideals are exactly the \(s\) pairwise distinct
ideals \(\varphi^{-1}(\mathfrak m_i)\).
\end{notation}

\begin{lemma}[A final index consisting of semilocal rings]
\label{lem:semilocal-smooth-source-finality}
Let \(R\) be a nonzero semilocal \(\Z_{(\Sigma)}\)-algebra, let
\(\mathcal V\) be a presentable \(\infty\)-category, and let
\(G:\mathrm{Sm}_\Sigma\to\mathcal V\) be a functor. Let
\(\widetilde G:\mathrm{EssSm}_\Sigma\to\mathcal V\) be the left Kan
extension of \(G\) along the inclusion
\(\iota:\mathrm{Sm}_\Sigma\hookrightarrow\mathrm{EssSm}_\Sigma\):
\[
\begin{tikzcd}
 \mathrm{Sm}_\Sigma \arrow[r,"G"] \arrow[d,hook,"\iota"'] &
 \mathcal V \\
 \mathrm{EssSm}_\Sigma. \arrow[ur,"\widetilde G"'] &
\end{tikzcd}
\]
\begin{enumerate}
\item For a smooth finitely presented \(\Z_{(\Sigma)}\)-algebra
\(C\) and a multiplicative subset \(M\subseteq C\),
\[
 \widetilde G(M^{-1}C)\simeq\varinjlim_{f\in M}G(C[1/f]),
\]
where \(M\) is preordered by divisibility.
\item The canonical maps
\[
 \operatorname*{colim}_{(D,\varphi)\in\mathrm{EssSm}^{\mathrm{sl}}_{\Sigma/R}}
     \widetilde G(D)
 \longrightarrow
 \operatorname*{colim}_{(D,\varphi)\in\mathrm{EssSm}_{\Sigma/R}}
     \widetilde G(D)
 \longleftarrow
 \operatorname*{colim}_{(C,\varphi)\in\mathrm{Sm}_{\Sigma/R}}G(C)
\]
are equivalences in \(\mathcal V\).
\item For every object \((D,\varphi)\) of
\(\mathrm{EssSm}^{\mathrm{sl}}_{\Sigma/R}\), the ring \(D\) is a
nonzero semilocal essentially smooth \(\Z_{(\Sigma)}\)-algebra in
the sense of \cref{def:essentially-smooth-sigma}.
\end{enumerate}
\end{lemma}

\begin{proof}
Let \(\mathfrak m_1,\ldots,\mathfrak m_s\) be the maximal ideals
of \(R\). We use the criterion of
\cite[Theorem~4.1.3.1]{HTT}: a functor \(u:\mathcal A\to\mathcal B\)
of ordinary categories is cofinal if, for every object \(b\) of
\(\mathcal B\), the comma category \(b/u\) of pairs
\((a,b\to u(a))\) has weakly contractible nerve.

\smallskip\noindent
\emph{Proof of \textup{(i)}.}
Put \(D=M^{-1}C\). The left Kan extension is given by
\(\widetilde G(D)=\operatorname*{colim}_{(C',\psi)\in\mathrm{Sm}_{\Sigma/D}}G(C')\)
\cite[Lemma~4.3.2.13]{HTT}. Consider the functor
\[
 u:M\longrightarrow\mathrm{Sm}_{\Sigma/D},\qquad
 f\longmapsto\bigl(C[1/f],\,C[1/f]\to D\bigr).
\]
Let \((C',\psi)\) be an object of \(\mathrm{Sm}_{\Sigma/D}\). An
object of \((C',\psi)/u\) is a pair \((f,\chi)\), where \(f\in M\)
and \(\chi:C'\to C[1/f]\) is a homomorphism whose composite with
\(C[1/f]\to D\) is \(\psi\). Since \(C'\) is finitely presented and
\(D=\varinjlim_{f\in M}C[1/f]\), such a pair exists, and two such
pairs first admit a common localization (invert the product of their denominators). At that
stage, equality in \(D\) implies equality after inverting one more
element of \(M\) \cite[\stackstag{00QO}]{Stacks}. Thus
\((C',\psi)/u\) is a nonempty filtered category, whose nerve is
weakly contractible \cite[Lemma~5.3.1.20]{HTT}. Hence \(u\) is
cofinal, which gives \textup{(i)}.

\smallskip\noindent
\emph{The second map in \textup{(ii)}.}
Let \(G'\) and \(\widetilde G'\) be the left Kan extensions of
\(G\) and \(\widetilde G\) to \(\mathrm{CAlg}_{\Z_{(\Sigma)}}\). The
two colimits are \(G'(R)\) and \(\widetilde G'(R)\)
\cite[Lemma~4.3.2.13]{HTT}, and \(G'\simeq\widetilde G'\) by
transitivity of left Kan extensions
\cite[Proposition~4.3.2.8]{HTT}.

\smallskip\noindent
\emph{Cofinality of
\(\mathrm{EssSm}^{\mathrm{inc}}_{\Sigma/R}\subseteq\mathrm{EssSm}_{\Sigma/R}\).}
Choose \(r_1,\ldots,r_s\in R\) with
\(r_i\bmod\mathfrak m_k=\delta_{ik}\), using the Chinese
remainder theorem \cite[\stackstag{00DT}]{Stacks}. For an object
\((D,\varphi)\) of \(\mathrm{EssSm}_{\Sigma/R}\), put
\[
 T(D,\varphi)=\bigl(D[t_1,\ldots,t_s],\,\varphi_T\bigr),\qquad
 \varphi_T|_D=\varphi,\quad\varphi_T(t_i)=r_i .
\]
The ring \(D[t_1,\ldots,t_s]\) is essentially smooth: if
\(D=M^{-1}C\), then \(D[t_1,\ldots,t_s]=M^{-1}(C[t_1,\ldots,t_s])\).
For \(i\neq k\), the element \(t_i\) belongs to
\(\varphi_T^{-1}(\mathfrak m_k)\) and not to
\(\varphi_T^{-1}(\mathfrak m_i)\), which contains \(t_i-1\). Hence
the prime ideals \(\varphi_T^{-1}(\mathfrak m_i)\) are pairwise
incomparable, and \(T\) is a functor
\(\mathrm{EssSm}_{\Sigma/R}\to\mathrm{EssSm}^{\mathrm{inc}}_{\Sigma/R}\).
The inclusions \(D\hookrightarrow D[t_1,\ldots,t_s]\) form a natural
transformation \(\epsilon:\mathrm{id}\to T\) of endofunctors of
\(\mathrm{EssSm}_{\Sigma/R}\).

Let \((D_0,\varphi_0)\) be an object of \(\mathrm{EssSm}_{\Sigma/R}\),
and let \(\mathcal J\) be the comma category of morphisms
\(\psi:(D_0,\varphi_0)\to(D,\varphi)\) with \((D,\varphi)\) in
\(\mathrm{EssSm}^{\mathrm{inc}}_{\Sigma/R}\). The assignment
\(\psi\mapsto\epsilon_D\circ\psi\) is an endofunctor \(\Theta\) of
\(\mathcal J\). The maps \(\epsilon_D\) define a natural
transformation \(\mathrm{id}_{\mathcal J}\to\Theta\), and the maps
\(T(\psi):T(D_0,\varphi_0)\to T(D,\varphi)\) define a natural
transformation to \(\Theta\) from the constant functor with value
\(\epsilon_{D_0}\), because
\(\epsilon_D\circ\psi=T(\psi)\circ\epsilon_{D_0}\). On nerves, these
two transformations give a homotopy between the identity and a
constant map. Thus \(\mathcal J\) has weakly contractible nerve, and
the inclusion is cofinal.

\smallskip\noindent
\emph{Cofinality of
\(\mathrm{EssSm}^{\mathrm{sl}}_{\Sigma/R}\subseteq\mathrm{EssSm}^{\mathrm{inc}}_{\Sigma/R}\),
and \textup{(iii)}.}
Let \((D,\varphi)\) be an object of
\(\mathrm{EssSm}^{\mathrm{inc}}_{\Sigma/R}\), and let
\(D_{\mathrm{sl}}\) be the semilocalization of \(D\) at the prime
ideals \(\varphi^{-1}(\mathfrak m_1),\ldots,\varphi^{-1}(\mathfrak m_s)\)
(\cref{def:arithmetic-semilocal}). The elements of
\(D\setminus\bigcup_i\varphi^{-1}(\mathfrak m_i)\) are mapped by
\(\varphi\) to units of \(R\), so \(\varphi\) extends uniquely to
\(\varphi_{\mathrm{sl}}:D_{\mathrm{sl}}\to R\). By
\cref{def:arithmetic-semilocal}, the maximal ideals of
\(D_{\mathrm{sl}}\) are the \(s\) pairwise distinct ideals
\(\varphi^{-1}(\mathfrak m_i)D_{\mathrm{sl}}=\varphi_{\mathrm{sl}}^{-1}(\mathfrak m_i)\).
If \(D=M^{-1}C\) with \(C\) smooth and finitely presented, then
\(D_{\mathrm{sl}}\) is the semilocalization of \(C\) at the inverse
images of the \(\mathfrak m_i\). Hence
\((D_{\mathrm{sl}},\varphi_{\mathrm{sl}})\) belongs to
\(\mathrm{EssSm}^{\mathrm{sl}}_{\Sigma/R}\), and \(D_{\mathrm{sl}}\)
is a nonzero semilocal essentially smooth
\(\Z_{(\Sigma)}\)-algebra. Applied to an object \((D,\varphi)\) of
\(\mathrm{EssSm}^{\mathrm{sl}}_{\Sigma/R}\), for which
\(D=D_{\mathrm{sl}}\), this proves \textup{(iii)}.

Let \((P,\varphi_P)\) be an object of
\(\mathrm{EssSm}^{\mathrm{sl}}_{\Sigma/R}\) and let
\(\psi:(D,\varphi)\to(P,\varphi_P)\) be a morphism. For
\(a\in D\setminus\bigcup_i\varphi^{-1}(\mathfrak m_i)\), the element
\(\psi(a)\) lies in no maximal ideal
\(\varphi_P^{-1}(\mathfrak m_i)\) of \(P\), so it is a unit. Hence
\(\psi\) factors uniquely through \(D\to D_{\mathrm{sl}}\). Thus
\(D\to D_{\mathrm{sl}}\) is an initial object of the comma category
of morphisms from \((D,\varphi)\) to objects of
\(\mathrm{EssSm}^{\mathrm{sl}}_{\Sigma/R}\), whose nerve is therefore
weakly contractible. Composing the two cofinal inclusions gives the
first equivalence in \textup{(ii)}.
\end{proof}

\begin{sublemma}[Fibre products of henselian surjections]
\label{sub:henselian-fibre-product}
Let \(A_1\to A_0\leftarrow A_2\) be surjective ring homomorphisms
with kernels \(I_1\) and \(I_2\), such that \((A_1,I_1)\) and
\((A_2,I_2)\) are henselian pairs. Put \(D=A_1\times_{A_0}A_2\), with
projections \(p_i:D\to A_i\). Then \(p_1\) and \(p_2\) are
surjective, and \((D,\ker p_1)\) and \((D,\ker p_2)\) are henselian
pairs.
\end{sublemma}

\begin{proof}
By symmetry it suffices to treat \(p_1\). It is surjective because
\(A_2\to A_0\) is. Its kernel is \(J=\{(0,b):b\in I_2\}\). We verify
Gabber's criterion \cite[\stackstag{09XI}\,(5)]{Stacks} for
\((D,J)\). Let \(x=(0,b)\in J\) and \(y=(y_1,y_2)\in D\). Then
\(1-yx=(1,1-y_2b)\), and \(1-y_2b\) is a unit of \(A_2\) with image
\(1\) in \(A_0\), since \(I_2\) is contained in the Jacobson radical
of \(A_2\). Hence \((1,(1-y_2b)^{-1})\) is an element of \(D\)
inverse to \(1-yx\), and \(J\) is contained in the Jacobson radical
of \(D\). Let \(n\geq1\) and
\[
 f(T)=T^n(T-1)+a_nT^n+\cdots+a_1T+a_0,\qquad a_k=(0,b_k)\in J .
\]
By Gabber's criterion \cite[\stackstag{09XI}\,(5)]{Stacks} for \((A_2,I_2)\), the polynomial
\(T^n(T-1)+b_nT^n+\cdots+b_0\) has a root \(\alpha_2\in1+I_2\). The
pair \(\alpha=(1,\alpha_2)\) lies in \(D\) and in \(1+J\), and
\(f(\alpha)=0\), because the first component of \(f\) is
\(T^n(T-1)\).
\end{proof}

\begin{sublemma}[Continuity and rigidity of truncated \'etale cohomology]
\label{sub:truncated-etale-finitary-rigid}
Let \(B\) be a ring in which the integer \(N\geq2\) is invertible,
and let \(j\geq0\). Consider the functor from \(B\)-algebras to
\(\mathcal D(\Z)\) given by
\begin{equation}\label{eq:truncated-etale-functor}
 A\longmapsto\tau^{\leq j}R\Gamma_{\mathrm{\acute et}}(A,\Z/N(j)).
\end{equation}
\begin{enumerate}
\item The functor \eqref{eq:truncated-etale-functor} commutes with
filtered colimits.
\item For every henselian pair \((A,I)\) in which \(A\) is a
\(B\)-algebra, the functor \eqref{eq:truncated-etale-functor} sends
\(A\to A/I\) to an equivalence.
\end{enumerate}
\end{sublemma}

\begin{proof}
For every \(B\)-algebra \(A\), the sheaf \(\Z/N(j)\) on
\((\Spec A)_{\mathrm{\acute et}}\) is the pullback of the sheaf
\(\Z/N(j)\) on \((\Spec B)_{\mathrm{\acute et}}\).

\textup{(i)} We may assume that \(\Lambda\) is a directed set
\cite[\stackstag{0032}]{Stacks}. By
\cite[\stackstag{03Q6}]{Stacks}, applied to the identity of
\(\Spec B\) and the sheaf \(\Z/N(j)\),
\[
 \varinjlim_{\lambda}H^q_{\mathrm{\acute et}}
      (\Spec A_\lambda,\Z/N(j))
 \xrightarrow{\ \sim\ }
 H^q_{\mathrm{\acute et}}(\Spec A,\Z/N(j))
 \qquad(q\geq0).
\]
Filtered colimits are exact in \(\mathcal D(\Z)\), so they commute
with cohomology and with \(\tau^{\leq j}\).

\textup{(ii)} The sheaf \(\Z/N(j)\) is torsion, so Gabber's affine
analogue of proper base change \cite[\stackstag{09ZI}]{Stacks} gives
\(H^q_{\mathrm{\acute et}}(\Spec A,\Z/N(j))\xrightarrow{\sim}
H^q_{\mathrm{\acute et}}(\Spec A/I,\Z/N(j))\) for every \(q\).
\end{proof}

\begin{lemma}[Smooth left Kan extension of ordinary \texorpdfstring{\'etale}{etale} truncations]
\label{lem:finite-ordinary-truncation-smooth-kan}
Let \(B\) be a ring in which the integer \(N\geq 2\) is
invertible, let \(j\geq 0\), and let \(A\) be a \(B\)-algebra. Let
\(\mathrm{Sm}_{B/A}\) be the category of pairs \((P,\varphi)\), where
\(P\) is a smooth finitely presented \(B\)-algebra and
\(\varphi:P\to A\) is a homomorphism of \(B\)-algebras; a morphism
\((P,\varphi)\to(P',\varphi')\) is a homomorphism of \(B\)-algebras
\(\psi:P\to P'\) with \(\varphi'\circ\psi=\varphi\).
\begin{enumerate}
\item The category \(\mathrm{Sm}_{B/A}\) is sifted.
\item The map
\begin{equation}\label{eq:truncated-etale-kan}
 \operatorname*{colim}_{(P,\varphi)\in\mathrm{Sm}_{B/A}}
 \tau^{\leq j}R\Gamma_{\mathrm{\acute et}}(P,\Z/N(j))
 \longrightarrow\tau^{\leq j}R\Gamma_{\mathrm{\acute et}}(A,\Z/N(j))
\end{equation}
induced by the pullbacks along
\(\Spec\varphi:\Spec A\to\Spec P\) is an equivalence in
\(\mathcal D(\Z)\).
\end{enumerate}
For \(B=\Z_{(\Sigma)}\), the category \(\mathrm{Sm}_{B/A}\) is the
category \(\mathrm{Sm}_{\Sigma/A}\) of
\cref{not:smooth-source-index}.
\end{lemma}

\begin{proof}
\textup{(i)} The category \(\mathrm{Sm}_{B/A}\) contains
\((B,B\to A)\). Let \((P,\varphi)\) and \((Q,\psi)\) be two objects.
The \(B\)-algebra \(P\otimes_BQ\) is smooth and finitely presented
\cite[\stackstag{00T4}, \stackstag{00TD}]{Stacks}, and by the
universal property of the tensor product
\((P\otimes_BQ,\varphi\otimes\psi)\) is an initial object of the
category of objects of \(\mathrm{Sm}_{B/A}\) equipped with morphisms
from \((P,\varphi)\) and from \((Q,\psi)\). A category with an
initial object has weakly contractible nerve, so the diagonal
\(\mathrm{Sm}_{B/A}\to\mathrm{Sm}_{B/A}\times\mathrm{Sm}_{B/A}\) is
cofinal \cite[Theorem~4.1.3.1]{HTT}. This is the definition of a
sifted category \cite[Definition~5.5.8.1]{HTT}.

\textup{(ii)} The spectrum underlying
\((\tau^{\leq j}R\Gamma_{\mathrm{\acute et}}(A,\Z/N(j)))[j]\)
is connective.
Consider the functor
\[
 \Phi:\{\text{\(B\)-algebras}\}\longrightarrow\mathcal S,\qquad
 A\longmapsto\Omega^\infty\bigl((\tau^{\leq j}
 R\Gamma_{\mathrm{\acute et}}(A,\Z/N(j)))[j]\bigr).
\]
The functor \(\Omega^\infty\) preserves sifted colimits of connective
spectra and is conservative on connective spectra
\cite[\S1.4.3]{HA}. By
\cref{sub:truncated-etale-finitary-rigid}\textup{(i)}, \(\Phi\)
preserves filtered colimits. By
\cref{sub:truncated-etale-finitary-rigid}\textup{(ii)}, \(\Phi\)
sends \(A\to A/I\) to an equivalence for every henselian pair
\((A,I)\). Let \(A_1\to A_0\leftarrow A_2\) be surjective maps of
\(B\)-algebras forming henselian pairs with their respective kernels.
By \cref{sub:henselian-fibre-product}, the projections of the fibre
product also form henselian pairs with their kernels. Thus \(\Phi\) sends
the four maps of the cartesian square
\[
\begin{tikzcd}
 A_1\times_{A_0}A_2\arrow[r,"p_2"]\arrow[d,"p_1"'] & A_2\arrow[d] \\
 A_1\arrow[r] & A_0
\end{tikzcd}
\]
to equivalences, hence the square to a cartesian square of spaces.
These are the hypotheses of Mathew's criterion
\cite[\arxivhtml{1908.02162v2}{A1.SS0.SSS1}{Proposition~A.0.1}]{EHK};
its specialization to finitary rigid functors is
\cite[\arxivhtml{2508.09915v1}{S2.E25}{Lemma~2.25}]{BEMMotivic25}.
It follows that
\[
 \operatorname*{colim}_{(P,\varphi)\in\mathrm{Sm}_{B/A}}\Phi(P)
 \longrightarrow\Phi(A)
\]
is an equivalence of spaces. By \textup{(i)}, this map is
\(\Omega^\infty\) of the map of connective spectra underlying the
\(j\)-fold shift of \eqref{eq:truncated-etale-kan}, which is
therefore an equivalence. The functor
\(\mathcal D(\Z)\to\mathrm{Sp}\) preserves colimits and is
conservative, so \eqref{eq:truncated-etale-kan} is an equivalence.
\end{proof}
\subsection{The classical source filtration}

\begin{construction}[The classical source filtration]
\label{con:classical-source-filtration}
In this section, a \emph{decreasing filtration} is a functor
\(F_\bullet:(\mathbf N,\leq)^{\mathrm{op}}\to\mathrm{Sp}\),
\(j\mapsto F_j\); a map of filtrations is a natural transformation.
The spectrum \(F_j\) is the \emph{stage} \(j\) of \(F_\bullet\), and
\[
 \operatorname{gr}_jF=\cofib\bigl(F_{j+1}\to F_j\bigr)
\]
is its \emph{grade} \(j\).

Use the base \(\Z_{(\Sigma)}\) of
\cref{def:essentially-smooth-sigma}, and let \(\mathrm{KGL}\) be the
motivic spectrum over \(\Z_{(\Sigma)}\) representing homotopy
algebraic \(K\)-theory \(KH\). For \(j\in\mathbf N\), let
\(f_j\mathrm{KGL}\) be its \(j\)-effective cover, with canonical
maps \(f_{j+1}\mathrm{KGL}\to f_j\mathrm{KGL}\). We use the notation of
\cite[\arxivhtml{2508.09915v1}{S4.E15}{Definition~4.15},
\arxivhtml{2508.09915v1}{S4.E37}{Remarks~4.37--4.38} and
\arxivhtml{2508.09915v1}{S4.E39}{Theorem~4.39(1)}]{BEMMotivic25}:
for a smooth quasi-compact quasi-separated \(\Z_{(\Sigma)}\)-scheme \(P\),
\[
 \mathrm{Fil}_{\mathbb A}^j KH(P)=(f_j\mathrm{KGL})(P)\in\mathrm{Sp},
 \qquad
 \mathrm{gr}_{\mathbb A}^j KH(P)\simeq \Z(j)^{\mathbb A}(P)[2j],
\]
where the first term is the spectrum of sections of the effective
cover over \(P\) and \(\Z(j)^{\mathbb A}(P)\in\mathcal D(\Z)\) is the
\(\mathbb A^1\)-motivic complex. These are the stages and grades of the
slice filtration of \(KH(P)\). In the notation of \cite{BEMMotivic25},
\(f_j\mathrm{KGL}=\mathrm{Fil}^j_{\mathrm{slice}}\mathrm{KGL}\) and the
cofiber of \(f_{j+1}\mathrm{KGL}\to f_j\mathrm{KGL}\) is
\(s^j\mathrm{KGL}\). We use the same notation on essentially smooth
\(\Z_{(\Sigma)}\)-algebras \(D\) by evaluating the intrinsic filtration
and complexes of \cite{BEMMotivic25} on \(\Spec D\); these agree with
the filtered colimits of their values on smooth localizations
by \cite[\arxivhtml{2508.09915v1}{S4.E39}{Theorem~4.39(2)}]{BEMMotivic25}.
The base-change equivalences in
\cite[\arxivhtml{2508.09915v1}{S9.E3}{Theorem~9.3}]{BEMMotivic25}
identify the covers and slices with those pulled back from \(\Spec\Z\).
\end{construction}

\begin{sublemma}[Functoriality, continuity, grades, and stage zero]
\label{sub:classical-filtration-inputs}
The filtration \(\mathrm{Fil}_{\mathbb A}^{\bullet}KH\) of
\cref{con:classical-source-filtration} has the following properties.
\begin{enumerate}
\item Let \(\mathit{Sm}/\Z_{(\Sigma)}\) be the category of smooth
      quasi-compact quasi-separated \(\Z_{(\Sigma)}\)-schemes and all
      morphisms of \(\Z_{(\Sigma)}\)-schemes between them. The assignments
      \(P\mapsto \mathrm{Fil}_{\mathbb A}^{\bullet}KH(P)\) and \(P\mapsto \Z(j)^{\mathbb A}(P)\) are functors
      \[
       \mathrm{Fil}_{\mathbb A}^{\bullet}KH:(\mathit{Sm}/\Z_{(\Sigma)})^{\mathrm{op}}
         \longrightarrow
         \operatorname{Fun}(\mathbf N^{\mathrm{op}},\mathrm{Sp}),
       \qquad
       \Z(j)^{\mathbb A}:(\mathit{Sm}/\Z_{(\Sigma)})^{\mathrm{op}}
         \longrightarrow\mathcal D(\Z).
      \]
\item For a smooth finitely presented \(\Z_{(\Sigma)}\)-algebra
      \(C\), a multiplicative subset \(W\subseteq C\), and
      \(j\in\mathbf N\), the restriction maps induce equivalences
      \[
       \varinjlim_{f\in W}\mathrm{Fil}_{\mathbb A}^jKH(C[1/f])\xrightarrow{\ \sim\ }
           \mathrm{Fil}_{\mathbb A}^jKH(W^{-1}C),
       \qquad
       \varinjlim_{f\in W}\Z(j)^{\mathbb A}(C[1/f])\xrightarrow{\ \sim\ }
           \Z(j)^{\mathbb A}(W^{-1}C),
      \]
      compatible with the maps \(\mathrm{Fil}_{\mathbb A}^{j+1}KH\to \mathrm{Fil}_{\mathbb A}^jKH\). Equivalently, by
      \cref{lem:semilocal-smooth-source-finality}\textup{(i)}, the
      functors \(\mathrm{Fil}_{\mathbb A}^jKH\) and \(\Z(j)^{\mathbb A}\) on \(\mathrm{EssSm}_\Sigma\) are
      the left Kan extensions of their restrictions to
      \(\mathrm{Sm}_\Sigma\).
\item For every smooth quasi-compact quasi-separated
      \(\Z_{(\Sigma)}\)-scheme or essentially smooth
      \(\Z_{(\Sigma)}\)-algebra \(P\) and every
      \(j\in\mathbf N\), there is a natural cofiber sequence
      \[
      \mathrm{Fil}_{\mathbb A}^{j+1}KH(P)\longrightarrow \mathrm{Fil}_{\mathbb A}^jKH(P)
                       \longrightarrow \Z(j)^{\mathbb A}(P)[2j].
      \]
\item For every smooth affine \(\Z_{(\Sigma)}\)-scheme or
      essentially smooth \(\Z_{(\Sigma)}\)-algebra \(P\), the
      stage-zero map and the canonical comparison maps are natural
      equivalences
      \[
      \mathrm{Fil}_{\mathbb A}^0KH(P)\xrightarrow{\sim}KH(P)
            \xleftarrow{\sim}K(P)
            \xleftarrow{\sim}K^{\mathrm{cn}}(P).
      \]
\end{enumerate}
\end{sublemma}

\begin{proof}
The functors in \textup{(i)} and the cofiber sequences in
\textup{(iii)} are those of
\cite[\arxivhtml{2508.09915v1}{S4.E15}{Definition~4.15}, \arxivhtml{2508.09915v1}{S4.E21}{Remark~4.21} and \arxivhtml{2508.09915v1}{S4.E39}{Theorem~4.39(1)--(2)}]{BEMMotivic25}.
Their evaluation over the fixed base \(\Z_{(\Sigma)}\) gives
the same filtration, because smooth pullback commutes with
effective covers and slices
\cite[\arxivhtml{2508.09915v1}{S3.E41}{Remark~3.41} and \arxivhtml{2508.09915v1}{S4.E21}{Remark~4.21}]{BEMMotivic25}. The equivalences
in \textup{(ii)}, including the transition maps for the open
immersions \(\Spec C[1/ff']\to\Spec C[1/f]\), are
\cite[\arxivhtml{2508.09915v1}{S4.E39}{Theorem~4.39(2)}]{BEMMotivic25}.

The filtration is exhaustive, and its grades in negative indices
vanish, by \cite[\arxivhtml{2508.09915v1}{S4.E20}{Remark~4.20} and \arxivhtml{2508.09915v1}{S4.E39}{Theorem~4.39(1)}]{BEMMotivic25}.
Hence \(\mathrm{Fil}_{\mathbb A}^0KH(P)\to KH(P)\) is an equivalence. The schemes \(P\) in
\textup{(iv)} are regular, noetherian, and affine, so
\(K(P)\to KH(P)\) is an equivalence by Quillen's homotopy
invariance of the \(K\)-theory of regular noetherian rings
\cite[\S6]{Quillen}, and \(K^{\mathrm{cn}}(P)\to K(P)\) is an
equivalence because the negative \(K\)-groups of a regular
noetherian ring vanish \cite[Proposition~6.8(b)]{TT}. This proves
\textup{(iv)}.
\end{proof}
\begin{lemma}[Integral bounds on semilocal sources]
\label{lem:classical-slice-source-bounds}
Let \(P\) be a nonzero semilocal essentially smooth
\(\Z_{(\Sigma)}\)-algebra (\cref{def:essentially-smooth-sigma}).
Assume that \(H^q(\Z(j)^{\mathrm{cla}}(\Spec P))=0\) for every
\(q>j\) and \(j\in\mathbf N\). For every \(j\in\mathbf N\),
\[
 H^q\bigl(\Z(j)^{\mathbb A}(P)\bigr)=0\quad(q>j),\qquad
 \pi_i\mathrm{Fil}_{\mathbb A}^jKH(P)=0\quad(i<j).
\]
\end{lemma}

\begin{proof}
\emph{The grades.}
Let \(C\) be a smooth finitely presented \(\Z_{(\Sigma)}\)-algebra.
The base $\Spec\Z_{(\Sigma)}$ is a Dedekind scheme. Thus
$f_1(H\Z)=0$ for Spitzweck's motivic Eilenberg--Mac Lane spectrum on
this base, as recalled in \cite[Remark~1.2]{BachmannEffective}.
By Bott periodicity for \(\mathrm{KGL}\) and the identification of
its zero slice with that spectrum \cite[Theorem~1.1(1), Example~1.3 and
Lemma~2.1]{BachmannEffective}, which does not require $2$ to be
invertible, the complex \(\Z(j)^{\mathbb A}(C)\) computes the cohomology of weight \(j\)
represented by that spectrum on \(\Spec C\). By
\cite[Theorem~7.18 and Corollary~7.19]{Spitzweck}, this is the
Nisnevich hypercohomology of the complex of sheaves
\(\Z(j)^{\mathrm{cla}}_{\Spec C}\). With
\cref{sub:cycle-descent-localization}\textup{(i)} we obtain
isomorphisms
\begin{equation}\label{eq:slice-grade-cycle}
 H^q\bigl(\Z(j)^{\mathbb A}(C)\bigr)\cong
 H^q\bigl(\Z(j)^{\mathrm{cla}}(\Spec C)\bigr)
 \qquad(q\in\Z),
\end{equation}
compatible with restriction along the open immersions
\(\Spec C[1/f]\to\Spec C\), which acts on the right-hand side by
flat pullback of cycles.

Write \(P=W^{-1}C\) as in \cref{def:essentially-smooth-sigma}, and
apply \eqref{eq:slice-grade-cycle} to the algebras \(C[1/f]\),
\(f\in W\). By \cref{sub:classical-filtration-inputs}\textup{(ii)}
and by \cref{sub:classical-cycle-spread}, applied with \(S=\{1\}\)
and \(S'=W\), the colimit over \(f\in W\) gives
\[
 H^q\bigl(\Z(j)^{\mathbb A}(P)\bigr)\cong
 H^q\bigl(\Z(j)^{\mathrm{cla}}(\Spec P)\bigr)
 \qquad(q\in\Z).
\]
The right-hand side vanishes for \(q>j\) by hypothesis.

\smallskip\noindent
\emph{The stages.}
Let \(d=\dim P\), which is also its valuative dimension
\cite[Proposition~2.3.2(9)]{EHIK}. The absolute base-change equivalences of
\cite[\arxivhtml{2508.09915v1}{S9.E3}{Theorem~9.3(1)}]{BEMMotivic25}
identify the $k$-effective covers of $\mathrm{KGL}$ pulled back from
$\Spec\Z$ to the base $\Spec\Z_{(\Sigma)}$, and then to $\Spec P$,
with the respective intrinsic $k$-effective covers. Thus the connectivity bound of
\cite[\arxivhtml{2508.09915v1}{S9.E5}{Corollary~9.5}]{BEMMotivic25}
applies to \(\mathrm{Fil}_{\mathbb A}^kKH(P)\): it is
\((k-d)\)-connective, also on essentially smooth localizations by
\cite[\arxivhtml{2508.09915v1}{S4.E39}{Theorem~4.39(2)}]{BEMMotivic25}.
Fix \(j\) and choose \(k\geq j+d\), so that \(\mathrm{Fil}_{\mathbb A}^kKH(P)\) is
\(j\)-connective. For \(j\leq r<k\), the vanishing already proved
gives
\[
 \pi_i\bigl(\Z(r)^{\mathbb A}(P)[2r]\bigr)=H^{2r-i}\bigl(\Z(r)^{\mathbb A}(P)\bigr)=0
 \qquad(i<r),
\]
so \(\Z(r)^{\mathbb A}(P)[2r]\) is \(j\)-connective. By descending induction on
\(r\) from \(k-1\) to \(j\), the cofiber sequences
\[
 \mathrm{Fil}_{\mathbb A}^{r+1}KH(P)\longrightarrow \mathrm{Fil}_{\mathbb A}^rKH(P)
       \longrightarrow \Z(r)^{\mathbb A}(P)[2r]
\]
of \cref{sub:classical-filtration-inputs}\textup{(iii)} show that
\(\mathrm{Fil}_{\mathbb A}^rKH(P)\) is \(j\)-connective. The case \(r=j\) is the assertion.
\end{proof}
\subsection{The sheafified filtration and its graded maps}

For the rest of this section, fix \(N=\ell^\nu\), where \(\ell\) is a
prime number invertible in \(\Z_{(\Sigma)}\) and \(\nu\geq1\).

\begin{construction}[The sheafified classical filtration]
\label{con:classical-etale-filtration}
Let \((\mathit{Sch}/\Z_{(\Sigma)})_{\mathrm{\acute et}}\) be the
localization of the site
\((\mathit{Sch}/\Z[1/N])_{\mathrm{\acute et}}\) of
\cref{con:etale-k-sheaves} at the object \(\Spec\Z_{(\Sigma)}\), and
put \(\mathcal C_\Sigma=(\mathit{Sm}/\Z_{(\Sigma)})_{\mathrm{\acute et}}\)
(\cref{not:etale-sites-spectra}).
By \cref{sub:classical-filtration-inputs}\textup{(i)},
\(P\mapsto \mathrm{Fil}_{\mathbb A}^{\bullet}KH(P)/N\) is a decreasing filtration of presheaves
of spectra on \(\mathcal C_\Sigma\). Applying
\(L^{\wedge}_{\mathcal C_\Sigma}\) stagewise gives the filtration
\(L^{\wedge}_{\mathcal C_\Sigma}(\mathrm{Fil}_{\mathbb A}^{\bullet}KH/N)\)
of hypercomplete sheaves. Its grades are the cofibers of successive
stage maps. Since \(L^{\wedge}_{\mathcal C_\Sigma}\) preserves cofiber
sequences,
\cref{sub:classical-filtration-inputs}\textup{(iii)} gives
\begin{equation}\label{eq:sheafified-grade}
 \operatorname{gr}_j\bigl(L^{\wedge}_{\mathcal C_\Sigma}
     (\mathrm{Fil}_{\mathbb A}^{\bullet}KH/N)\bigr)
 \simeq L^{\wedge}_{\mathcal C_\Sigma}\bigl((\Z(j)^{\mathbb A}/N)[2j]\bigr).
\end{equation}
The maps of presheaves
\(\mathrm{Fil}_{\mathbb A}^0KH/N\to KH/N\leftarrow K/N\leftarrow K^{\mathrm{cn}}/N\) on
\(\mathcal C_\Sigma\) are equivalences on affine objects, by
\cref{sub:classical-filtration-inputs}\textup{(iv)}. Hence
\cref{sub:big-site-restriction}\textup{(iv)} and \textup{(i)} give
\begin{equation}\label{eq:sheafified-stage-zero}
 L^{\wedge}_{\mathcal C_\Sigma}(\mathrm{Fil}_{\mathbb A}^0KH/N)\simeq
 L^{\wedge}_{\mathcal C_\Sigma}\bigl(K^{\mathrm{cn}}/N\bigr)\simeq
 \mathcal K^{\mathrm{big}}_N\big|_{\mathcal C_\Sigma},
\end{equation}
where the sheaf \(\mathcal K^{\mathrm{big}}_N\) of
\cref{con:etale-k-sheaves} is restricted along
\(\mathcal C_\Sigma\hookrightarrow
(\mathit{Sch}/\Z_{(\Sigma)})_{\mathrm{\acute et}}\to
(\mathit{Sch}/\Z[1/N])_{\mathrm{\acute et}}\). For a smooth
\(\Z_{(\Sigma)}\)-scheme \(P\),
\cref{sub:big-site-restriction}\textup{(i)} gives
\begin{equation}\label{eq:sheafified-stage-small}
 L^{\wedge}_{\mathcal C_\Sigma}(\mathrm{Fil}_{\mathbb A}^jKH/N)|_{P_{\mathrm{\acute et}}}\simeq
 L^{\wedge}_{P_{\mathrm{\acute et}}}
   \bigl((\mathrm{Fil}_{\mathbb A}^jKH/N)|_{P_{\mathrm{\acute et}}}\bigr).
\end{equation}
When \(P\) is a ring, \(P_{\mathrm{\acute et}}\) denotes the small
\'etale site of \(\Spec P\).
\end{construction}

\begin{lemma}[The classical filtration at geometric stalks]
\label{lem:classical-filtration-stalks}
Let \(P\) be a smooth finitely presented \(\Z_{(\Sigma)}\)-algebra,
let \(\bar x\) be a geometric point of \(\Spec P\), and let
\(r,j\in\mathbf N\). With the notation of
\cref{con:classical-etale-filtration}:
\begin{enumerate}
\item the stalk at \(\bar x\) of
\(\pi_i L^{\wedge}_{P_{\mathrm{\acute et}}}
\bigl((\Z(r)^{\mathbb A}/N)|_{P_{\mathrm{\acute et}}}\bigr)\) is zero for every
\(i\neq0\);
\item the stalk at \(\bar x\) of
\(\pi_i((\mathrm{Fil}_{\mathbb A}^jKH/N)|_{P_{\mathrm{\acute et}}})\) is zero for
\(i<2j\), and the map
\[
 \bigl(\pi_i((\mathrm{Fil}_{\mathbb A}^jKH/N)|_{P_{\mathrm{\acute et}}})\bigr)_{\bar x}
 \longrightarrow
 \bigl(\pi_i\mathcal K^{\mathrm{\acute et}}_{P,N}\bigr)_{\bar x},
\]
induced by the stage-zero map and sheafification, is an isomorphism for every \(i\geq2j\).
\end{enumerate}
\end{lemma}

\begin{proof}
Let \(\mathcal O^{\mathrm{sh}}\) be the strict henselian local ring
of \(\Spec P\) at \(\bar x\). We evaluate the intrinsic filtration and
motivic complexes of \cite{BEMMotivic25} on
\(\Spec\mathcal O^{\mathrm{sh}}\). The stalk of a presheaf of abelian groups agrees with the stalk of
its sheafification \cite[\stackstag{03PT}]{Stacks}; the unit of
hypercomplete sheafification induces isomorphisms on homotopy sheaves
(\cref{not:etale-sites-spectra}). The cofiber of multiplication by \(N\)
commutes with filtered colimits. By
\eqref{eq:sheafified-stage-small} and
\cite[\arxivhtml{2508.09915v1}{S6.E1}{Remark~6.1} and \arxivhtml{2508.09915v1}{S4.E39}{Theorem~4.39(2)}]{BEMMotivic25}, for every
\(i\in\Z\),
\begin{equation}\label{eq:classical-stage-stalk}
 \begin{aligned}
 \bigl(\pi_i((\mathrm{Fil}_{\mathbb A}^jKH/N)|_{P_{\mathrm{\acute et}}})\bigr)_{\bar x}
 &=\varinjlim_{(U,\bar u)}\pi_i\bigl(\mathrm{Fil}_{\mathbb A}^jKH(U)/N\bigr)\\
 &=\pi_i\bigl(\mathrm{Fil}_{\mathbb A}^jKH(\mathcal O^{\mathrm{sh}})/N\bigr),
 \end{aligned}
\end{equation}
where \((U,\bar u)\) runs through the affine \'etale neighborhoods
of \(\bar x\). The same argument identifies the stalk in
\textup{(i)} with \(H^{-i}(\Z(r)^{\mathbb A}(\mathcal O^{\mathrm{sh}})/N)\).

By \cite[\arxivhtml{2508.09915v1}{S6.E21}{Theorem~6.21(2)}]{BEMMotivic25}, the smooth schemes under consideration over the
mixed-characteristic Dedekind domain \(\Z_{(\Sigma)}\) satisfy the
conditions \((\mathrm A)\)--\((\mathrm H)\) of the class
\(\mathfrak{BL}^\ell\) of that reference, and by
\cite[\arxivhtml{2508.09915v1}{S6.E2}{Lemma~6.2(3)}]{BEMMotivic25} so do their strict localizations.
In the notation of \cite{BEMMotivic25} the coefficient prime is
\(p\) and the modulus is \(p^m\); here they are \(\ell\) and
\(N=\ell^\nu\), and \(\ell\) is invertible on
\(\Spec\mathcal O^{\mathrm{sh}}\). Condition \((\mathrm B)\) gives
\[
 H^q\bigl(\Z(r)^{\mathbb A}(\mathcal O^{\mathrm{sh}})/N\bigr)=0\qquad(q\neq0),
\]
which is \textup{(i)}. By \cite[\arxivhtml{2508.09915v1}{S9.E5}{Corollary~9.5}]{BEMMotivic25},
\(\mathrm{Fil}_{\mathbb A}^kKH(\mathcal O^{\mathrm{sh}})/N\) is
\((k-\dim\mathcal O^{\mathrm{sh}})\)-connective.
Choose \(k\geq2j+\dim\mathcal O^{\mathrm{sh}}\), so that
\(\mathrm{Fil}_{\mathbb A}^kKH(\mathcal O^{\mathrm{sh}})/N\) is \(2j\)-connective. For \(j\leq r<k\), the spectrum
\(\mathrm{gr}_{\mathbb A}^rKH(\mathcal O^{\mathrm{sh}})/N\simeq
(\Z(r)^{\mathbb A}(\mathcal O^{\mathrm{sh}})/N)[2r]\) has homotopy only in degree
\(2r\geq2j\). By descending induction on \(r\), the cofiber sequences
\(\mathrm{Fil}_{\mathbb A}^{r+1}KH/N\to \mathrm{Fil}_{\mathbb A}^rKH/N\to(\Z(r)^{\mathbb A}/N)[2r]\) show that
\(\mathrm{Fil}_{\mathbb A}^jKH(\mathcal O^{\mathrm{sh}})/N\) is \(2j\)-connective. The cofiber
of \(\mathrm{Fil}_{\mathbb A}^jKH(\mathcal O^{\mathrm{sh}})/N\to \mathrm{Fil}_{\mathbb A}^0KH(\mathcal O^{\mathrm{sh}})/N\)
has a finite filtration with grades
\((\Z(r)^{\mathbb A}(\mathcal O^{\mathrm{sh}})/N)[2r]\), \(0\leq r<j\), so its
homotopy groups vanish in degrees greater than \(2j-2\). The long
exact sequence of homotopy groups therefore gives isomorphisms
\(\pi_i(\mathrm{Fil}_{\mathbb A}^jKH(\mathcal O^{\mathrm{sh}})/N)\xrightarrow{\sim}
\pi_i(\mathrm{Fil}_{\mathbb A}^0KH(\mathcal O^{\mathrm{sh}})/N)\) for \(i\geq2j\). With
\eqref{eq:classical-stage-stalk}, and with
\eqref{eq:sheafified-stage-zero} and \eqref{eq:big-small} for the
identification of the target, this proves \textup{(ii)}.
\end{proof}
\begin{lemma}[Normalization of the classical grade map]
\label{lem:classical-grade-normalization}
Let \(P\) be a smooth finitely presented \(\Z_{(\Sigma)}\)-algebra,
let \(j\in\mathbf N\), and let \(a_{\mathrm{\acute et}}\) denote
sheafification of presheaves of abelian groups on
\(P_{\mathrm{\acute et}}\). Write
\[
 q_j:\mathrm{Fil}_{\mathbb A}^jKH/N\longrightarrow (\Z(j)^{\mathbb A}/N)[2j],\qquad
 u_j:\mathrm{Fil}_{\mathbb A}^jKH/N\longrightarrow \mathrm{Fil}_{\mathbb A}^0KH/N
\]
for the grade projection and the stage map on \(\mathcal C_\Sigma\),
restricted to \(P_{\mathrm{\acute et}}\). On the affine basis of
\(P_{\mathrm{\acute et}}\), we identify \(\mathrm{Fil}_{\mathbb A}^0KH/N\) with
\(K^{\mathrm{cn}}/N\) by \cref{sub:classical-filtration-inputs}\textup{(iv)};
after sheafification we use \eqref{eq:sheafified-stage-zero}.
\begin{enumerate}
\item The maps \(\pi_{2j}(q_j)\) and
\(\pi_{2j}(u_j)\) are isomorphisms, and the
sheaf \(L^{\wedge}_{P_{\mathrm{\acute et}}}(\Z(j)^{\mathbb A}/N)\) is the
Eilenberg--Mac Lane sheaf of \(a_{\mathrm{\acute et}}H^0(\Z(j)^{\mathbb A}/N)\).
With the isomorphism \(\lambda_j\) of
\cref{sub:finite-k-etale-coefficients}\textup{(ii)}, the composite
\begin{equation}\label{eq:grade-normalization}
\begin{aligned}
 a_{\mathrm{\acute et}}H^0(\Z(j)^{\mathbb A}/N)
 &\xrightarrow{\,(\pi_{2j}(q_j))^{-1}\,}
       \pi_{2j}(\mathrm{Fil}_{\mathbb A}^jKH/N)\\
 &\xrightarrow{\,\pi_{2j}(u_j)\,}
       \pi_{2j}\mathcal K^{\mathrm{\acute et}}_{P,N}
   \xrightarrow{\,\lambda_j^{-1}\,}\Z/N(j)
\end{aligned}
\end{equation}
is therefore an equivalence
\(L^{\wedge}_{P_{\mathrm{\acute et}}}(\Z(j)^{\mathbb A}/N)\simeq\Z/N(j)\).
\item The unit
\(\Z(j)^{\mathbb A}/N\to L^{\wedge}_{P_{\mathrm{\acute et}}}(\Z(j)^{\mathbb A}/N)\), followed by
\eqref{eq:grade-normalization} and evaluated at \(P\), is a map
\begin{equation}
\label{eq:classical-graded-unit}
 \Z(j)^{\mathbb A}(P)/N\longrightarrow
             R\Gamma_{\mathrm{\acute et}}(P,\Z/N(j)),
\end{equation}
natural in \(P\in\mathrm{Sm}_\Sigma\), which induces isomorphisms on
\(H^q\) for all \(q\leq j\).
\end{enumerate}
\end{lemma}

\begin{proof}
\textup{(i)} The fibre of \(q_j\) is \(\mathrm{Fil}_{\mathbb A}^{j+1}KH/N\). By
\cref{lem:classical-filtration-stalks}\textup{(ii)} and
\eqref{eq:sheafified-stage-small}, the stalks of
\(\pi_i(\mathrm{Fil}_{\mathbb A}^{j+1}KH/N)\) vanish for \(i<2j+2\),
and \(\pi_{2j}(u_j)\) is an isomorphism on
stalks. Hence both maps are isomorphisms
\cite[\stackstag{03PU}]{Stacks}. By
\cref{lem:classical-filtration-stalks}\textup{(i)}, the homotopy
sheaves of \(L^{\wedge}_{P_{\mathrm{\acute et}}}(\Z(j)^{\mathbb A}/N)\) vanish
outside degree zero, where the homotopy sheaf is
\(a_{\mathrm{\acute et}}H^0(\Z(j)^{\mathbb A}/N)\). Since
\(\pi_{2j}((\Z(j)^{\mathbb A}/N)[2j])=a_{\mathrm{\acute et}}H^0(\Z(j)^{\mathbb A}/N)\), the composite
\eqref{eq:grade-normalization} is defined, and it is an isomorphism
of sheaves of abelian groups.

\textup{(ii)} We use the sheafification unit in presheaves of complexes,
viewed as $H\Z$-modules. Forgetting the module structure preserves
limits and colimits and detects equivalences; sheafification and its
unit are therefore compatible with the underlying spectra
\cite[\S4.2.3]{HA}. The target in \textup{(i)} is discrete, so the
isomorphism of abelian sheaves \eqref{eq:grade-normalization} defines
an equivalence of Eilenberg--Mac Lane complexes. The maps in
\eqref{eq:grade-normalization} are the restrictions to
\(P_{\mathrm{\acute et}}\) of maps of abelian sheaves on
\(\mathcal C_\Sigma\), by \cref{sub:big-site-restriction}\textup{(ii)}
and the naturality of \(q_j,u_j\) and \(\lambda_j\)
(\cref{sub:finite-k-etale-coefficients}). These maps are isomorphisms
on \(\mathcal C_\Sigma\), since this can be checked on its affine
objects and their small \'etale sites by \textup{(i)}. Inverting them
on \(\mathcal C_\Sigma\) gives the naturality of
\eqref{eq:classical-graded-unit} for all maps of smooth algebras.
The normalization \eqref{eq:grade-normalization} is the
one constructed in
\cite[\arxivhtml{2508.09915v1}{S6.E8}{Proposition~6.8}, Step~1,
\arxivhtml{2508.09915v1}{S6.E10}{formula (6.10)}]{BEMMotivic25}.
Indeed, on smooth $P$ the maps from the effective cover to its slice
and to $\mathrm{KH}$ in \cite[\arxivhtml{2508.09915v1}{S6.E10}{(6.10)}]{BEMMotivic25} identify with $q_j$ and $u_j$ via
\cref{con:classical-source-filtration,sub:classical-filtration-inputs}.
Regularity identifies $\mathrm{KH}(P)$ with $K^{\mathrm{cn}}(P)$, and
$\lambda_j$ is the same isomorphism from the Tate twist to the
homotopy sheaf by
\cref{sub:finite-k-etale-coefficients}. Thus inverting the left arrow
and the Tate isomorphism after sheafification gives exactly
\eqref{eq:grade-normalization}.
Conditions \((\mathrm B)\) and \((\mathrm G)\) of that reference
imply that the ordinary \'etale sheafification of \(\Z(j)^{\mathbb A}/N\) is
an Eilenberg--Mac Lane sheaf, hence agrees with its
hypercompletion in \textup{(i)}. Condition \((\mathrm C)\),
together with \textup{(i)}, therefore identifies the map
\eqref{eq:classical-graded-unit} with the unit in the equivalence
\[
 \Z(j)^{\mathbb A}/N\xrightarrow{\ \sim\ }
 L_{\mathrm{Nis}}\tau^{\leq j}
 R\Gamma_{\mathrm{\acute et}}(-,\Z/N(j))
\]
on the small Nisnevich site of \(P\); here
\(L_{\mathrm{Nis}}\) denotes Nisnevich sheafification of
presheaves of complexes. Since \'etale cohomology satisfies
Nisnevich descent, the cofiber of the map to
\(R\Gamma_{\mathrm{\acute et}}(-,\Z/N(j))\) is
\[
 L_{\mathrm{Nis}}\tau^{>j}
 R\Gamma_{\mathrm{\acute et}}(-,\Z/N(j)).
\]
Its Nisnevich cohomology sheaves vanish in degrees at most
\(j\), and the same holds after derived global sections over
\(P\). The long exact cohomology sequence proves the assertion.
\end{proof}

\begin{remark}[The external comparison inputs]
\label{rem:classical-comparison-hypotheses}
For every prime \(\ell\), including \(2\), with \(\ell\) invertible on
\(\Z_{(\Sigma)}\), smooth quasi-compact quasi-separated schemes over this mixed-characteristic
Dedekind domain belong to the class \(\mathfrak{BL}^\ell\) of
\cite[\S6.1, \arxivhtml{2508.09915v1}{S6.E21}{Theorem~6.21(2)}]{BEMMotivic25};
this passes to strict localizations by
\cite[\arxivhtml{2508.09915v1}{S6.E2}{Lemma~6.2(3)}]{BEMMotivic25}.
In this invertible-coefficient setting conditions \((\mathrm D)\)--\((\mathrm F)\)
are vacuous. The proof of \cite[Theorem~6.21(2)]{BEMMotivic25} obtains \((\mathrm A)\)--\((\mathrm C)\)
and \((\mathrm G)\) from the Beilinson--Lichtenbaum identification of the
zero slice \cite[\arxivhtml{2508.09915v1}{S6.E21}{Theorem~6.21(1)}]{BEMMotivic25};
the completeness condition \((\mathrm H)\) also follows directly from
\cite[\arxivhtml{2508.09915v1}{S9.E5}{Corollary~9.5}]{BEMMotivic25}, since the henselian local
rings of smooth schemes are noetherian of finite valuative dimension.
The normalization of the comparison map is the one of
\cite[\arxivhtml{2508.09915v1}{S6.E8}{Proposition~6.8}, formula (6.10)]{BEMMotivic25},
as checked in \cref{lem:classical-grade-normalization}. In particular,
\cref{lem:classical-grade-normalization} identifies the comparison
morphism itself, independently of \cref{thm:semilocal-motivic-vanishing}.
\end{remark}

\begin{lemma}[Comparison on smooth algebras]
\label{lem:classical-etale-stage-comparison}
Use the notation of \cref{con:classical-etale-filtration}.
\begin{enumerate}
\item For every \(j\in\mathbf N\), the map from the hypercomplete
sheafification of \(\mathrm{Fil}_{\mathbb A}^jKH/N\) to
\(\mathcal K^{\mathrm{big}}_N|_{\mathcal C_\Sigma}\), given by the stage-zero
map and \eqref{eq:sheafified-stage-zero}, admits a contractible space
of factorizations through the \((2j)\)-connective cover, each giving
an equivalence
\[
 L^{\wedge}_{\mathcal C_\Sigma}(\mathrm{Fil}_{\mathbb A}^jKH/N)
 \xrightarrow{\ \sim\ }
 \tau_{\geq2j}(\mathcal K^{\mathrm{big}}_N|_{\mathcal C_\Sigma}).
\]
These equivalences are compatible with the transition maps of the
sheafified filtration and of the connective covers.
\item For every smooth finitely presented \(\Z_{(\Sigma)}\)-algebra
\(P\), the sheafification units, the equivalences of
\textup{(i)}, and \eqref{eq:big-small-truncated} define a map of
decreasing filtrations
\begin{equation}\label{eq:smooth-filtered-comparison}
 \mathrm{Fil}_{\mathbb A}^{\bullet}KH(P)/N\longrightarrow
 R\Gamma_{\mathrm{\acute et}}\bigl(P,\tau_{\geq 2\bullet}
      \mathcal K^{\mathrm{\acute et}}_{P,N}\bigr),
\end{equation}
natural in \(P\in\mathrm{Sm}_\Sigma\). Under the identification
\(\operatorname{gr}_j\bigl(\tau_{\geq2\bullet}
\mathcal K^{\mathrm{\acute et}}_{P,N}\bigr)\simeq\Z/N(j)[2j]\) given
by \(\lambda_j\), its map on grade \(j\) is the shift by \([2j]\) of
\eqref{eq:classical-graded-unit}.
\end{enumerate}
\end{lemma}

\begin{proof}
\textup{(i)} Let \(i<2j\). By
\cref{lem:classical-filtration-stalks}\textup{(ii)} and
\cite[\stackstag{03PU}]{Stacks}, the sheaf
\(\pi_i((\mathrm{Fil}_{\mathbb A}^jKH/N)|_{P_{\mathrm{\acute et}}})\)
vanishes for every smooth finitely presented \(\Z_{(\Sigma)}\)-algebra
\(P\). By \cref{sub:big-site-restriction}\textup{(ii)}, the homotopy
sheaf \(\pi_i(\mathrm{Fil}_{\mathbb A}^jKH/N)\) on \(\mathcal C_\Sigma\)
has no nonzero section over an affine object; since every object has
a Zariski covering by affine objects, this homotopy sheaf is zero.
The sheafification unit induces isomorphisms on homotopy sheaves,
so \(L^{\wedge}_{\mathcal C_\Sigma}(\mathrm{Fil}_{\mathbb A}^jKH/N)\)
is \((2j)\)-connective. Its map to
\(\mathcal K^{\mathrm{big}}_N|_{\mathcal C_\Sigma}\) therefore factors
through \(\tau_{\geq2j}(\mathcal K^{\mathrm{big}}_N|_{\mathcal C_\Sigma})\)
by the truncation adjunction \cite[\S1.2.1]{HA}, with a contractible
space of factorizations. The same reduction to stalks, with
\cref{lem:classical-filtration-stalks}\textup{(ii)} for \(i\geq2j\),
shows that the factorization is \(\infty\)-connective. Its source is
hypercomplete by construction. Its target is the fibre of
\[
 \mathcal K^{\mathrm{big}}_N|_{\mathcal C_\Sigma}\longrightarrow
 \tau_{\leq2j-1}(\mathcal K^{\mathrm{big}}_N|_{\mathcal C_\Sigma}).
\]
The source of this last map is hypercomplete by
\eqref{eq:sheafified-stage-zero}, and its target is bounded above
for the Postnikov \(t\)-structure, hence hypercomplete by
\cite[Lemma~6.5.2.9]{HTT}, applied after suspension and
\(\Omega^\infty\). Hypercomplete spectral sheaves are closed under
fibres. Hence the factorization is an equivalence
(\cref{not:etale-sites-spectra}). Coherence follows by applying the
functorial connective covers to the augmented diagram
\[
 L^{\wedge}_{\mathcal C_\Sigma}(\mathrm{Fil}_{\mathbb A}^{\bullet}KH/N)
 \longrightarrow \mathcal K^{\mathrm{big}}_N|_{\mathcal C_\Sigma}:
\]
for \(j\geq k\), apply \(\tau_{\geq2j}\) to the transition from stage
\(j\) to stage \(k\), then use \(\tau_{\geq2j}\to\tau_{\geq2k}\).
The counit on stage \(j\) is an equivalence by its \((2j)\)-connectivity.
Thus the factorizations form a map of augmented filtrations.
By \cref{sub:big-site-restriction}\textup{(ii)}, connective covers
commute with restriction from the big site to \(\mathcal C_\Sigma\).

\textup{(ii)} Evaluating at \(P\) the maps of filtered presheaves
\[
 \mathrm{Fil}_{\mathbb A}^{\bullet}KH/N\longrightarrow
 L^{\wedge}_{\mathcal C_\Sigma}(\mathrm{Fil}_{\mathbb A}^{\bullet}KH/N)
 \xrightarrow{\ \sim\ }
 (\tau_{\geq2\bullet}\mathcal K^{\mathrm{big}}_N)|_{\mathcal C_\Sigma}
\]
and applying
\eqref{eq:big-small-truncated} gives
\eqref{eq:smooth-filtered-comparison}, naturally in \(P\). On grade
\(j\), this is the evaluation at \(P\) of
\[
 (\Z(j)^{\mathbb A}/N)[2j]\longrightarrow
 L^{\wedge}_{\mathcal C_\Sigma}\bigl((\Z(j)^{\mathbb A}/N)[2j]\bigr)
 \xrightarrow{\ \sim\ }
 \bigl((\pi_{2j}\mathcal K^{\mathrm{big}}_N)|_{\mathcal C_\Sigma}\bigr)[2j],
\]
where the equivalence is induced by \textup{(i)} via
\eqref{eq:sheafified-grade}. On \(\pi_{2j}\), it is
\(\pi_{2j}(u_j)\circ
(\pi_{2j}(q_j))^{-1}\) on every
\(P_{\mathrm{\acute et}}\). By the definition
\eqref{eq:grade-normalization} of the normalization, the map on
grade \(j\) is the shift of \eqref{eq:classical-graded-unit}.
\end{proof}
\subsection{The filtered comparison on semilocal rings}

\begin{sublemma}[Smooth Kan extension of connective \(K\)-theory]
\label{sub:connective-k-smooth-kan}
Let \(A\) be a \(\Z_{(\Sigma)}\)-algebra. The pullback maps
\(\varphi^*:K^{\mathrm{cn}}(P)\to K^{\mathrm{cn}}(A)\) induce
equivalences
\[
 \operatorname*{colim}_{(P,\varphi)\in\mathrm{Sm}_{\Sigma/A}}
     K^{\mathrm{cn}}(P)\xrightarrow{\ \sim\ }K^{\mathrm{cn}}(A),
 \qquad
 \operatorname*{colim}_{(P,\varphi)\in\mathrm{Sm}_{\Sigma/A}}
     K^{\mathrm{cn}}(P)/N\xrightarrow{\ \sim\ }K^{\mathrm{cn}}(A)/N.
\]
\end{sublemma}

\begin{proof}
The space \(\Omega^\infty K^{\mathrm{cn}}(A)\) is the group completion
of the groupoid of finitely generated projective \(A\)-modules.
The smooth left Kan extension result of Bhatt--Lurie, recorded in
\cite[\arxivhtml{1908.02162v2}{A1.SS0.SSS6}{Example~A.0.6(1)}]{EHK},
applies with the discrete base \(\Z_{(\Sigma)}\): smooth derived
algebras are flat over the base, hence discrete in this case, and
maps between discrete algebras are ordinary ring homomorphisms. Its
indexing category is thus \(\mathrm{Sm}_{\Sigma/A}\). It follows that the map
\[
 \operatorname*{colim}_{(P,\varphi)\in\mathrm{Sm}_{\Sigma/A}}
     \Omega^\infty K^{\mathrm{cn}}(P)\longrightarrow
     \Omega^\infty K^{\mathrm{cn}}(A)
\]
is an equivalence in \(\mathcal S\). The category
\(\mathrm{Sm}_{\Sigma/A}\) is sifted
(\cref{lem:finite-ordinary-truncation-smooth-kan}\textup{(i)}), and
\(\Omega^\infty\) preserves sifted colimits of connective spectra and
is conservative on connective spectra \cite[\S1.4.3]{HA}; compare
\cite[\arxivhtml{1908.02162v2}{A1.SS0.SSS5}{Lemma~A.0.5}]{EHK}. Hence the first map of the
statement is an equivalence. The cofiber of multiplication by \(N\)
commutes with colimits of spectra, which gives the second.
\end{proof}

\begin{definition}[The lisse filtration with finite coefficients]
\label{def:lisse-mod-N}
For a \(\Z_{(\Sigma)}\)-algebra \(A\) and \(j\in\mathbf N\), put
\[
 \mathrm{Fil}_{\mathrm{lis}}^jK^{\mathrm{cn}}(A)
 =\operatorname*{colim}_{(P,\varphi)\in\mathrm{Sm}_{\Sigma/A}}\mathrm{Fil}_{\mathbb A}^jKH(P)
 \in\mathrm{Sp}.
\]
The maps \(\mathrm{Fil}_{\mathbb A}^{j+1}KH\to \mathrm{Fil}_{\mathbb A}^jKH\) make \(\mathrm{Fil}_{\mathrm{lis}}^{\bullet}K^{\mathrm{cn}}(A)\) a
decreasing filtration, functorial in \(A\): it is the left Kan
extension of \(\mathrm{Fil}_{\mathbb A}^{\bullet}KH\) along
\(\mathrm{Sm}_\Sigma\hookrightarrow\mathrm{CAlg}_{\Z_{(\Sigma)}}\)
\cite[Lemma~4.3.2.13]{HTT}, that is, the lisse filtration \emph{over
\(\Z_{(\Sigma)}\)} (the fixed base is suppressed from the subscript
\(\mathrm{lis}\)) of
\cite[\arxivhtml{2508.09915v1}{S7.E1}{Definition~7.1} and \arxivhtml{2508.09915v1}{S7.E3}{Lemma~7.3(2)}]{BEMMotivic25}. Cofibers
commute with colimits, so
\begin{equation}\label{eq:lisse-grades}
\begin{aligned}
 \mathrm{Fil}_{\mathrm{lis}}^jK^{\mathrm{cn}}(A)/N&\simeq
 \operatorname*{colim}_{(P,\varphi)\in\mathrm{Sm}_{\Sigma/A}}\mathrm{Fil}_{\mathbb A}^jKH(P)/N,\\
 \operatorname{gr}_j\bigl(\mathrm{Fil}_{\mathrm{lis}}^\bullet K^{\mathrm{cn}}(A)/N\bigr)&\simeq
 \operatorname*{colim}_{(P,\varphi)\in\mathrm{Sm}_{\Sigma/A}}
     \bigl(\Z(j)^{\mathbb A}(P)/N\bigr)[2j].
\end{aligned}
\end{equation}
and \cref{sub:classical-filtration-inputs}\textup{(iv)} with
\cref{sub:connective-k-smooth-kan} gives
\begin{equation}\label{eq:lisse-stage-zero}
 \mathrm{Fil}_{\mathrm{lis}}^0K^{\mathrm{cn}}(A)/N\simeq K^{\mathrm{cn}}(A)/N .
\end{equation}
\end{definition}

\begin{proposition}[The filtered comparison on semilocal rings]
\label{prop:filtered-etale-comparison}
Let \(A\) be a nonzero semilocal \(\Z_{(\Sigma)}\)-algebra.
Assume that for every \((D,\varphi)\) in
\(\mathrm{EssSm}^{\mathrm{sl}}_{\Sigma/A}\),
\(H^q(\Z(j)^{\mathrm{cla}}(\Spec D))=0\) for \(q>j\) and \(j\in\mathbf N\).
\begin{enumerate}
\item For every \(j\in\mathbf N\), the spectrum
\(\mathrm{Fil}_{\mathrm{lis}}^jK^{\mathrm{cn}}(A)/N\) of \cref{def:lisse-mod-N} is
\(j\)-connective.
\item There is a map of decreasing filtrations
\begin{equation}
\label{eq:filtered-etale-comparison}
 \mathrm{Fil}_{\mathrm{lis}}^{\bullet}K^{\mathrm{cn}}(A)/N\longrightarrow
 R\Gamma_{\mathrm{\acute et}}\bigl(A,\tau_{\geq 2\bullet}
      \mathcal K^{\mathrm{\acute et}}_{A,N}\bigr),
\end{equation}
natural in \(A\), whose stage zero is, under
\eqref{eq:lisse-stage-zero}, the map \eqref{eq:truncated-unit} for
\(a=0\) and \(b=+\infty\).
\item For every \(j\in\mathbf N\), the grades \(j\) of the source and
of the target of \eqref{eq:filtered-etale-comparison} are canonically
identified with the source and the target of the truncation map
\begin{equation}
\label{eq:truncation-on-graded}
 \bigl(\tau^{\leq j}
 R\Gamma_{\mathrm{\acute et}}(A,\Z/N(j))\bigr)[2j]
 \longrightarrow
 R\Gamma_{\mathrm{\acute et}}(A,\Z/N(j))[2j],
\end{equation}
and under these identifications the map on grade \(j\) induced by
\eqref{eq:filtered-etale-comparison} is \eqref{eq:truncation-on-graded}.
\end{enumerate}
\end{proposition}

\begin{proof}
\emph{Construction of \eqref{eq:filtered-etale-comparison}.}
For an object \((P,\varphi)\) of \(\mathrm{Sm}_{\Sigma/A}\), compose
\eqref{eq:smooth-filtered-comparison} with the restriction map of the
presheaf \(\tau_{\geq2\bullet}\mathcal K^{\mathrm{big}}_N\) along
\(\Spec\varphi\), using \eqref{eq:big-small-truncated} for \(P\) and
for \(A\):
\begin{equation}\label{eq:cocone-component}
 \mathrm{Fil}_{\mathbb A}^{\bullet}KH(P)/N\longrightarrow
 \bigl(\tau_{\geq2\bullet}\mathcal K^{\mathrm{big}}_N\bigr)(\Spec P)
 \xrightarrow{\ \varphi^*\ }
 \bigl(\tau_{\geq2\bullet}\mathcal K^{\mathrm{big}}_N\bigr)(\Spec A)
 \simeq
 R\Gamma_{\mathrm{\acute et}}\bigl(A,\tau_{\geq 2\bullet}
      \mathcal K^{\mathrm{\acute et}}_{A,N}\bigr).
\end{equation}
The first arrow is a natural transformation of functors of
\(P\in\mathrm{Sm}_\Sigma\)
(\cref{lem:classical-etale-stage-comparison}\textup{(ii)}), and the
second is given by the functor
\(B\mapsto(\tau_{\geq2\bullet}\mathcal K^{\mathrm{big}}_N)(\Spec B)\)
on \(\mathrm{CAlg}_{\Z_{(\Sigma)}}\). By the universal property of
the left Kan extension \(\mathrm{Fil}_{\mathrm{lis}}^{\bullet}K^{\mathrm{cn}}/N\), the maps
\eqref{eq:cocone-component} induce the map
\eqref{eq:filtered-etale-comparison}, naturally in \(A\). At stage zero, the
first arrow of \eqref{eq:cocone-component} is the unit
\(K^{\mathrm{cn}}(P)/N\to\mathcal K^{\mathrm{big}}_N(\Spec P)\), by
\eqref{eq:sheafified-stage-zero}. The unit is a map of presheaves,
so \eqref{eq:cocone-component} at stage zero is
\(\varphi^*:K^{\mathrm{cn}}(P)/N\to K^{\mathrm{cn}}(A)/N\) followed by
the unit for \(A\). This proves \textup{(ii)}.

\smallskip\noindent
\emph{Proof of \textup{(i)}.}
By \cref{sub:classical-filtration-inputs}\textup{(ii)}, the functor
\(\mathrm{Fil}_{\mathbb A}^jKH/N\) on \(\mathrm{EssSm}_\Sigma\) is the left Kan extension of
its restriction to \(\mathrm{Sm}_\Sigma\). Hence
\cref{lem:semilocal-smooth-source-finality}\textup{(ii)}, applied
with \(\mathcal V=\mathrm{Sp}\) and \(R=A\), and
\eqref{eq:lisse-grades} give
\begin{equation}\label{eq:lisse-semilocal-index}
 \mathrm{Fil}_{\mathrm{lis}}^jK^{\mathrm{cn}}(A)/N\simeq
 \operatorname*{colim}_{(D,\varphi)\in
       \mathrm{EssSm}^{\mathrm{sl}}_{\Sigma/A}}\mathrm{Fil}_{\mathbb A}^jKH(D)/N .
\end{equation}
For every \((D,\varphi)\) in
\(\mathrm{EssSm}^{\mathrm{sl}}_{\Sigma/A}\), the ring \(D\) is a
nonzero semilocal essentially smooth \(\Z_{(\Sigma)}\)-algebra
(\cref{lem:semilocal-smooth-source-finality}\textup{(iii)}), so
\(\mathrm{Fil}_{\mathbb A}^jKH(D)\), and therefore \(\mathrm{Fil}_{\mathbb A}^jKH(D)/N\), is \(j\)-connective by
\cref{lem:classical-slice-source-bounds}. The \(j\)-connective
spectra are closed under arbitrary colimits. In particular, no
filteredness of the category of smooth algebras over the possibly
singular ring \(A\), and no regularity of \(A\), is being asserted here.

\smallskip\noindent
\emph{The grades of the source.}
Let \((D,\varphi)\) be an object of
\(\mathrm{EssSm}^{\mathrm{sl}}_{\Sigma/A}\), and write \(D=W^{-1}C\)
as in \cref{def:essentially-smooth-sigma}. The maps
\eqref{eq:classical-graded-unit} for the algebras \(C[1/f]\),
\(f\in W\), are natural, and their colimit is a map
\begin{equation}\label{eq:graded-unit-essentially-smooth}
 \Z(j)^{\mathbb A}(D)/N\longrightarrow
 R\Gamma_{\mathrm{\acute et}}(D,\Z/N(j)),
\end{equation}
by \cref{sub:classical-filtration-inputs}\textup{(ii)} for the
source and by continuity of \'etale cohomology
\cite[\stackstag{03Q6}]{Stacks} for the target. Filtered colimits
are exact, so \eqref{eq:graded-unit-essentially-smooth} induces
isomorphisms on \(H^q\) for \(q\leq j\)
(\cref{lem:classical-grade-normalization}\textup{(ii)}). By
\cref{lem:classical-slice-source-bounds} and the exact sequence
\(H^q(\Z(j)^{\mathbb A}(D))\to H^q(\Z(j)^{\mathbb A}(D)/N)\to H^{q+1}(\Z(j)^{\mathbb A}(D))\), the source of
\eqref{eq:graded-unit-essentially-smooth} has no cohomology in
degrees \(q>j\). The map \eqref{eq:graded-unit-essentially-smooth}
therefore factors, by the
truncation adjunction \cite[\S1.2.1]{HA}, as
\[
 \Z(j)^{\mathbb A}(D)/N\xrightarrow{\ \sim\ }
 \tau^{\leq j}R\Gamma_{\mathrm{\acute et}}(D,\Z/N(j))\longrightarrow
 R\Gamma_{\mathrm{\acute et}}(D,\Z/N(j)),
\]
where the first map is an equivalence, natural in \((D,\varphi)\),
and the second is the truncation map. The functor
\eqref{eq:truncated-etale-functor}, for \(B=\Z_{(\Sigma)}\), on
\(\mathrm{EssSm}_\Sigma\) is the left Kan extension of its restriction
to \(\mathrm{Sm}_\Sigma\), by
\cref{sub:truncated-etale-finitary-rigid}\textup{(i)} and
\cref{lem:semilocal-smooth-source-finality}\textup{(i)}. Applying
\cref{lem:semilocal-smooth-source-finality}\textup{(ii)} to
\(\Z(j)^{\mathbb A}/N\) and to \eqref{eq:truncated-etale-functor},
with \(\mathcal V=\mathcal D(\Z)\), and
then \cref{lem:finite-ordinary-truncation-smooth-kan}\textup{(ii)},
we obtain from \eqref{eq:lisse-grades}
\begin{equation}\label{eq:lisse-grade-identification}
 \begin{aligned}
 &\operatorname{gr}_j\bigl(\mathrm{Fil}_{\mathrm{lis}}^\bullet K^{\mathrm{cn}}(A)/N\bigr)\\
 &\quad\simeq\operatorname*{colim}_{(D,\varphi)\in
       \mathrm{EssSm}^{\mathrm{sl}}_{\Sigma/A}}
       \bigl(\tau^{\leq j}R\Gamma_{\mathrm{\acute et}}(D,\Z/N(j))\bigr)[2j]\\
 &\quad\simeq\operatorname*{colim}_{(P,\varphi)\in
       \mathrm{Sm}_{\Sigma/A}}
       \bigl(\tau^{\leq j}R\Gamma_{\mathrm{\acute et}}(P,\Z/N(j))\bigr)[2j]\\
 &\quad\simeq\bigl(\tau^{\leq j}R\Gamma_{\mathrm{\acute et}}(A,\Z/N(j))\bigr)[2j].
 \end{aligned}
\end{equation}

\smallskip\noindent
\emph{The grades of the target and the map between them.}
The functor \(R\Gamma_{\mathrm{\acute et}}(A,-)\) preserves cofiber
sequences, and
\(\cofib\bigl(\tau_{\geq2j+2}\mathcal K^{\mathrm{\acute et}}_{A,N}
\to\tau_{\geq2j}\mathcal K^{\mathrm{\acute et}}_{A,N}\bigr)
\simeq\Z/N(j)[2j]\) by
\cref{sub:finite-k-etale-coefficients}\textup{(ii)}. Hence the grade
\(j\) of the target of \eqref{eq:filtered-etale-comparison} is
\(R\Gamma_{\mathrm{\acute et}}(A,\Z/N(j))[2j]\). By
\cref{lem:classical-etale-stage-comparison}\textup{(ii)} and
\cref{sub:finite-k-etale-coefficients}\textup{(iii)}, the map
induced by \eqref{eq:cocone-component} on grade \(j\) is the shift
by \([2j]\) of
\[
 \Z(j)^{\mathbb A}(P)/N\xrightarrow{\ \eqref{eq:classical-graded-unit}\ }
 R\Gamma_{\mathrm{\acute et}}(P,\Z/N(j))
 \xrightarrow{\ \varphi^*\ }
 R\Gamma_{\mathrm{\acute et}}(A,\Z/N(j)).
\]
For \((D,\varphi)\) in \(\mathrm{EssSm}^{\mathrm{sl}}_{\Sigma/A}\),
the corresponding map is, by naturality of the truncation maps, the
shift of
\[
 \begin{aligned}
 \Z(j)^{\mathbb A}(D)/N&\xrightarrow{\ \sim\ }
 \tau^{\leq j}R\Gamma_{\mathrm{\acute et}}(D,\Z/N(j))\\
 &\xrightarrow{\ \varphi^*\ }\tau^{\leq j}R\Gamma_{\mathrm{\acute et}}(A,\Z/N(j))\\
 &\longrightarrow R\Gamma_{\mathrm{\acute et}}(A,\Z/N(j)).
 \end{aligned}
\]
Under \eqref{eq:lisse-grade-identification}, the map on grade \(j\)
induced by \eqref{eq:filtered-etale-comparison} is therefore
\eqref{eq:truncation-on-graded}. This proves \textup{(iii)}.
\end{proof}
\subsection{Passage to the finite comparison}

\begin{sublemma}[A comparison of finite filtrations]
\label{sub:finite-filtration-comparison}
Let \(n\in\mathbf N\), and let
\(\theta_\bullet:D_\bullet\to F_\bullet\) be a natural transformation
of functors \(\{0<1<\cdots<n+1\}^{\mathrm{op}}\to\mathrm{Sp}\) with
\(D_{n+1}=0=F_{n+1}\). For \(0\leq j\leq n\), let \(Q_j\) be the
cofiber of the induced map
\[
 \cofib(D_{j+1}\to D_j)\longrightarrow\cofib(F_{j+1}\to F_j),
\]
and suppose that \(\pi_iQ_j=0\) for all \(i\geq j\). Then
\(\pi_i\cofib(\theta_0)=0\) for all \(i\geq n\). In particular
\(\pi_n(\theta_0):\pi_nD_0\to\pi_nF_0\) is an isomorphism.
\end{sublemma}

\begin{proof}
Put \(C_j=\cofib(\theta_j)\). For \(0\leq j\leq n\), the commutative diagram
\[
\begin{tikzcd}[column sep=large]
 D_{j+1}\arrow[r]\arrow[d,"\theta_{j+1}"'] & D_j\arrow[r]\arrow[d,"\theta_j"'] &
 \cofib(D_{j+1}\to D_j)\arrow[d] \\
 F_{j+1}\arrow[r] & F_j\arrow[r] & \cofib(F_{j+1}\to F_j)
\end{tikzcd}
\]
has cofiber sequences as rows. Taking vertical cofibers gives the
cofiber sequence
\(C_{j+1}\to C_j\to Q_j\), hence an exact sequence
\[
 \pi_iC_{j+1}\longrightarrow\pi_iC_j\longrightarrow\pi_iQ_j
 \qquad(i\in\Z).
\]
Here \(C_{n+1}=0\), and \(\pi_iQ_j=0\) for \(i\geq n\geq j\).
Descending induction on \(j\) gives \(\pi_iC_j=0\) for all
\(i\geq n\) and \(0\leq j\leq n+1\). The exact sequence
\[
 \pi_{n+1}C_0\longrightarrow\pi_nD_0
 \xrightarrow{\ \pi_n(\theta_0)\ }\pi_nF_0\longrightarrow\pi_nC_0
\]
gives the last assertion.
\end{proof}

\begin{lemma}[Criterion for the finite connective \'etale comparison]
\label{lem:finite-comparison-from-source-bound}
Let \(A\) be a nonzero semilocal ring and let \(N\geq 2\)
be an integer invertible in \(A\). Choose \(\Sigma\) to be the set
of positive characteristics of the closed residue fields of \(A\),
or a singleton \(\{p\}\) with \(p\nmid N\) if there are none.
Assume that \(H^q(\Z(j)^{\mathrm{cla}}(\Spec D))=0\) for \(q>j\) and
\(j\in\mathbf N\) for every \((D,\varphi)\in\mathrm{EssSm}^{\mathrm{sl}}_{\Sigma/A}\). For every \(n\in\mathbf N\), the
finite comparison map \eqref{eq:finite-comparison-map} induces an
isomorphism
\[
 K_n^{\mathrm{cn}}(A;\Z/N)\xrightarrow{\ \sim\ }
 \pi_n R\Gamma_{\mathrm{\acute et}}\bigl(A,
        \tau_{\leq 2n}\mathcal K^{\mathrm{\acute et}}_{A,N}\bigr),
\]
natural in \(A\). If \(A\) is regular and noetherian, then
\(K_n^{\mathrm{cn}}(A;\Z/N)=K_n(A;\Z/N)\).
\end{lemma}

\begin{proof}
\emph{The base ring.}
The set \(\Sigma\) is finite, and no element of \(\Sigma\)
divides \(N\). An integer divisible by no element of \(\Sigma\) lies
in no maximal ideal of \(A\), hence is a unit of \(A\). Thus \(A\) is
a \(\Z_{(\Sigma)}\)-algebra, and \(N\) is a unit of
\(\Z_{(\Sigma)}\).

\smallskip\noindent
\emph{Prime powers.}
Let \(N=\ell^\nu\) with \(\ell\) prime and \(\nu\geq1\). Compose
\eqref{eq:filtered-etale-comparison} with the map induced by
\(\mathcal K^{\mathrm{\acute et}}_{A,N}\to
\tau_{\leq2n}\mathcal K^{\mathrm{\acute et}}_{A,N}\). Since
\(\tau_{\geq2n+2}\tau_{\leq2n}=0\) on sheaves, the target filtration
has a zero object at stage $n+1$. Its diagram of maps from that
stage is canonically zero. Taking cofibers of the maps from stage
$n+1$ in both diagrams therefore gives a natural transformation of
functors
\(\{0<1<\cdots<n+1\}^{\mathrm{op}}\to\mathrm{Sp}\)
\begin{equation}\label{eq:finite-filtration-map}
 \theta_j:
 \cofib\bigl(\mathrm{Fil}_{\mathrm{lis}}^{n+1}K^{\mathrm{cn}}(A)/N\to
             \mathrm{Fil}_{\mathrm{lis}}^jK^{\mathrm{cn}}(A)/N\bigr)
 \longrightarrow
 R\Gamma_{\mathrm{\acute et}}\bigl(A,\tau_{\geq 2j}\tau_{\leq2n}
      \mathcal K^{\mathrm{\acute et}}_{A,N}\bigr),
\end{equation}
whose source and target vanish for \(j=n+1\). Let \(0\leq j\leq n\).
The cofiber of the source of \(\theta_{j+1}\) in the source of
\(\theta_j\) is \(\operatorname{gr}_j(\mathrm{Fil}_{\mathrm{lis}}^\bullet K^{\mathrm{cn}}(A)/N)\), and
the cofiber of the target of \(\theta_{j+1}\) in the target of
\(\theta_j\) is
\(R\Gamma_{\mathrm{\acute et}}(A,\Z/N(j))[2j]\), by
\cref{sub:finite-k-etale-coefficients}\textup{(ii)} and
\(2j\leq2n\). By
\cref{prop:filtered-etale-comparison}\textup{(iii)}, the
map induced by \(\theta\) between these cofibers is
\eqref{eq:truncation-on-graded}, whose cofiber is
\[
 Q_j=\bigl(\tau^{>j}
       R\Gamma_{\mathrm{\acute et}}(A,\Z/N(j))\bigr)[2j],
 \qquad
 \pi_iQ_j=H^{2j-i}\bigl(\tau^{>j}
       R\Gamma_{\mathrm{\acute et}}(A,\Z/N(j))\bigr)=0
 \quad(i\geq j).
\]
By \cref{sub:finite-filtration-comparison}, the map
\(\pi_n(\theta_0)\) is an isomorphism.

By \eqref{eq:lisse-stage-zero}, the source of \(\theta_0\) sits in a
cofiber sequence
\begin{equation}\label{eq:stage-zero-cofiber}
 \mathrm{Fil}_{\mathrm{lis}}^{n+1}K^{\mathrm{cn}}(A)/N\longrightarrow K^{\mathrm{cn}}(A)/N
 \xrightarrow{\ \rho\ }
 \cofib\bigl(\mathrm{Fil}_{\mathrm{lis}}^{n+1}K^{\mathrm{cn}}(A)/N\to
             \mathrm{Fil}_{\mathrm{lis}}^0K^{\mathrm{cn}}(A)/N\bigr).
\end{equation}
The spectrum \(\mathrm{Fil}_{\mathrm{lis}}^{n+1}K^{\mathrm{cn}}(A)/N\) is \((n+1)\)-connective
by \cref{prop:filtered-etale-comparison}\textup{(i)}, so
\(\pi_n(\rho)\) is an isomorphism. The sheaf
\(\mathcal K^{\mathrm{\acute et}}_{A,N}\) is connective
(\cref{sub:finite-k-etale-coefficients}\textup{(i)}), so the target
of \(\theta_0\) is
\(R\Gamma_{\mathrm{\acute et}}(A,\tau_{\leq2n}
\mathcal K^{\mathrm{\acute et}}_{A,N})\), and
\(\theta_0\circ\rho\) is the finite comparison map
\eqref{eq:finite-comparison-map} by
\cref{prop:filtered-etale-comparison}\textup{(ii)}. Hence
\eqref{eq:finite-comparison-map} is an isomorphism on \(\pi_n\).

\smallskip\noindent
\emph{General \(N\).}
Write \(N=\prod_{i=1}^rN_i\), where the \(N_i\) are powers of
pairwise distinct prime numbers. By
\cref{sub:moore-primary-decomposition}, the presheaf
\(K^{\mathrm{cn}}/N\) is the direct sum of the presheaves
\(K^{\mathrm{cn}}/N_i\). The functors
\(L^{\wedge}_{A_{\mathrm{\acute et}}}\), \(\tau_{\leq2n}\), and
\(R\Gamma_{\mathrm{\acute et}}(A,-)\) preserve finite direct sums.
Hence \eqref{eq:finite-comparison-map} for \(N\) is the direct sum
of the maps \eqref{eq:finite-comparison-map} for the \(N_i\), and it
is an isomorphism on \(\pi_n\) by the prime-power case.

\smallskip\noindent
\emph{Naturality and regular rings.}
The map \eqref{eq:finite-comparison-map} is natural in \(A\) by
\cref{con:etale-k-sheaves}; it does not depend on \(\Sigma\). If
\(A\) is regular and noetherian, then
\(K^{\mathrm{cn}}(A)\to K(A)\) is an equivalence
\cite[Proposition~6.8(b)]{TT}.
\end{proof}
\begin{theorem}[Finite connective \'etale comparison on local rings]
\label{thm:local-finite-connective-comparison}
Let \(A\) be a nonzero local ring and let \(N\geq 2\)
be an integer invertible in \(A\). For every \(n\in\mathbf N\), the
finite comparison map \eqref{eq:finite-comparison-map} induces an
isomorphism
\[
 K_n^{\mathrm{cn}}(A;\Z/N)\xrightarrow{\ \sim\ }
 \pi_n R\Gamma_{\mathrm{\acute et}}\bigl(A,
        \tau_{\leq 2n}\mathcal K^{\mathrm{\acute et}}_{A,N}\bigr),
\]
natural in \(A\). If \(A\) is regular and noetherian, then
\(K_n^{\mathrm{cn}}(A;\Z/N)=K_n(A;\Z/N)\).
\end{theorem}

\begin{proof}
Choose \(\Sigma\) as in \cref{lem:finite-comparison-from-source-bound}.
By \cref{lem:semilocal-smooth-source-finality}\textup{(iii)}, each ring
\(D\) in the index category for \(A\) is local, say $D=C_{\mathfrak p}$
with $C$ smooth. By cycle continuity (\cref{sub:classical-cycle-spread})
its weight-$j$ cohomology is the stalk at $\mathfrak p$ of the Zariski
cohomology sheaf $\mathcal H^q(\Z(j)^{\mathrm{cla}}_{\Spec C})$:
principal opens containing $\mathfrak p$ form a basis and filtered
colimits commute with cohomology. Geisser's vanishing
\cite[Corollary~4.4]{Geisser04} for $q>j$ therefore gives the required
local bound. Apply
\cref{lem:finite-comparison-from-source-bound}.
\end{proof}

\begin{remark}[Relation to other motivic filtrations]
\label{rem:comparison-literature}
Elmanto--Morrow in equal characteristic
\cite[Theorem~1.1]{ElmantoMorrow}, and Bouis for arbitrary
quasi-compact quasi-separated schemes
\cite[Theorems~C(1) and D(3)]{Bouis}, construct motivic
complexes that occur as the shifted grades of a filtration on
nonconnective \(K\)-theory and, with invertible finite coefficients,
agree with \'etale cohomology through the weight.

For the lisse filtration over \(\Z\) (rather than over
\(\Z_{(\Sigma)}\) as in \cref{def:lisse-mod-N}), the comparison of
its graded pieces on local rings in positive weights is known.
Let \(A\) be a local ring, let \(j\geq1\), and let \(N\geq2\) be a natural
number invertible in \(A\).
In the notation of \cite{BouisWeibel25},
\(\Z(j)^{\mathrm{lisse}}\) is the left Kan extension of
classical motivic cohomology from smooth \(\Z\)-algebras,
and \(\Z(j)^{\mathrm{mot}}\) is Bouis' motivic complex.
The comparison with truncated motivic cohomology is proved
over a field in \cite[Theorem~7.7]{ElmantoMorrow} and for
every local ring in \cite[Corollary~2.12]{BouisWeibel25}.
Together with \cite[Lemma~2.9]{BouisWeibel25} and
\cite[Corollary~5.6]{Bouis}, it gives
\[
 \begin{aligned}
 \Z(j)^{\mathrm{lisse}}(A)/N
 &\simeq \tau^{\leq j}\bigl(\Z(j)^{\mathrm{mot}}(A)/N\bigr)\\
 &\simeq \tau^{\leq j}
       R\Gamma_{\mathrm{\acute et}}(A,\Z/N(j)).
 \end{aligned}
\]
For prime-power \(N\), the surjectivity in
\cite[Lemma~2.9]{BouisWeibel25} and the coefficient exact
sequence give
\(H^{j+1}(\Z(j)^{\mathrm{mot}}(A))[N]=0\).
This proves that reduction modulo \(N\) commutes with
\(\tau^{\leq j}\). The general modulus follows from
\cref{sub:moore-primary-decomposition}. On the other hand,
\cref{prop:filtered-etale-comparison}\textup{(iii)}
identifies the prime-power coefficient grades of the filtration over
\(\Z_{(\Sigma)}\) directly with the shifted truncated \'etale
cohomology complexes in every weight \(j\geq0\), under its upper-bound
hypothesis (which holds for local rings by the proof of
\cref{thm:local-finite-connective-comparison}).
No integral base-change comparison of lisse filtrations is used.
For local rings essentially smooth over a Dedekind domain,
the full motivic complex agrees with the classical one by
\cite[Theorem~6.1]{BouisKundu25}; its vanishing above the
weight follows from \cite[Corollary~4.4]{Geisser04}.

\end{remark}

\section{Finite \'etale effacement and bounded spectra}
\label{sec:finite-etale-effacement}

This section converts two sections and a path between them into
finite Galois data. The unit symbols of \cref{prop:semilocal-unit-symbols}
efface positive-degree cohomology on finite \'etale covers
(\cref{prop:finite-etale-effacement}); bounded
spectral sheaves with finite locally constant \(\ell\)-primary
homotopy sheaves then descend to the finite \'etale site
(\cref{prop:bounded-finite-etale-descent}).
\Cref{prop:finite-galois-data-two-sections} represents the two
sections by a finite-group classifying map and a specified path
by a homotopy lift at a finite quotient.

\subsection{Milnor symbols and \'etale residues}

\begin{notation}[Milnor symbols and tame residues]
\label{not:milnor-tame-residues}
For a commutative local ring $A$, let $K_*^M(A)$ denote the
quotient of the tensor algebra of the abelian group $A^\times$
over $\Z$ by the homogeneous two-sided ideal generated by
$u\otimes(1-u)$, where $u,1-u\in A^\times$
\cite[\arxivhtml{1509.01087v3}{S2.Thmdefin1}{Definition~2.1}]{Dahlhausen18}. Its degree-$j$ part is
$K_j^M(A)$, for $j\in\mathbf N$, and the image of
$u_1\otimes\cdots\otimes u_j$ is denoted by
$\{u_1,\ldots,u_j\}$. In particular, $K_0^M(A)=\Z$, with empty
symbol $1$. For fields this is the usual Milnor $K$-theory
\cite[Chapter~III, Definition~7.1]{Weibel}.
The notation $\widehat K_j^M(A)$ denotes improved Milnor
$K$-theory \cite{Kerz10}. Its construction and the natural map
$K_j^M(A)\to\widehat K_j^M(A)$ are recalled in
\cite[\arxivhtml{1509.01087v3}{S2.Thmdefin2}{Definition~2.2}]{Dahlhausen18}. For a field $E$, this natural
map is an isomorphism \cite[\arxivhtml{1509.01087v3}{S2.Thmdefin3}{Proposition~2.3(v)}]{Dahlhausen18}.
For a discrete valuation ring $A$ with fraction field $E$ and
residue field $\kappa(A)$, we write
\[
 \partial^{\mathrm M}_A:K_j^M(E)\longrightarrow
      K_{j-1}^M(\kappa(A))\qquad(j\geq1)
\]
for the tame residue with the convention of
\cite[proof of Corollary~4.4]{Geisser04}. Thus, for a uniformizer
$\pi\in A$ and units $u_1,\ldots,u_j\in A^\times$,
\[
\begin{aligned}
 \partial^{\mathrm M}_A\{u_1,\ldots,u_{j-1},\pi\}
       &=\{\overline u_1,\ldots,\overline u_{j-1}\},\\
 \partial^{\mathrm M}_A\{u_1,\ldots,u_j\}&=0.
\end{aligned}
\]
\end{notation}

\begin{definition}[Kummer and Galois symbols]\label{def:kummer-symbol}
Let \(m\geq2\) be a natural number, let \(X\) be a scheme on which
\(m\) is invertible, and let \(q\in\mathbf N\).
The finite locally constant sheaf of \(m\)-th roots of unity on
\(X_{\mathrm{\acute et}}\) is denoted by \(\mu_m\), and
\(\mu_m^{\otimes q}\) denotes its \(q\)-fold tensor product over
\(\Z/m\), with \(\mu_m^{\otimes0}=\Z/m\). For \(q\geq1\), a \emph{Kummer symbol of degree \(q\)}
is a cup product of \(q\) images of units under the Kummer boundary
\[
 \operatorname{Kum}_m:\Gamma(X,\mathcal O_X^\times)
       \longrightarrow H^1_{\mathrm{\acute et}}(X,\mu_m).
\]
For a field \(F\) in which \(m\) is
invertible, the \emph{Galois symbol} is the homomorphism
\[
\begin{aligned}
 K_q^M(F)/m&\longrightarrow
 H^q_{\mathrm{\acute et}}(F,\mu_m^{\otimes q}),\\
 \{u_1,\ldots,u_q\}&\longmapsto
 \operatorname{Kum}_m(u_1)\cup\cdots\cup\operatorname{Kum}_m(u_q)
\end{aligned}
\]
for \(q\geq1\); in degree zero it is the identity of \(\Z/m\)
\cite[Chapter~VI, \S4]{Weibel}.
\end{definition}

\begin{notation}[Inverse twist of a finite coefficient sheaf]
\label{not:finite-coefficient-inverse-twist}
Let $m\geq2$ be a natural number invertible on a scheme $X$,
and let $M$ be a finite locally constant abelian sheaf on
$X$ with $mM=0$. Its inverse twist is the \'etale sheaf
\[
 M(-1)=\underline{\operatorname{Hom}}_{\Z/m}(\mu_m,M).
\]
Evaluation defines a map $\mu_m\otimes_{\Z/m}M(-1)\to M$.
\end{notation}

\begin{sublemma}[Residues with locally constant coefficients]
\label{sub:trait-locally-constant-residue}
Let $A$ be a discrete valuation ring with maximal ideal
$\mathfrak m_A$, fraction field $L$, residue field $\kappa(A)$ and
normalized valuation $v_A$. Let $m\geq2$ be a natural number
invertible in $A$, and let $M$ be a finite locally constant abelian
sheaf on $\Spec A$ with $mM=0$.
Cohomology is taken on the small \'etale site.
\begin{enumerate}
\item For every $r\in\Z$ there is an isomorphism
\begin{equation}\label{eq:trait-support-residue}
 H^r_{\mathrm{\acute et},\{\mathfrak m_A\}}(\Spec A,M)
   \xrightarrow{\ \sim\ }H^{r-2}_{\mathrm{\acute et}}(\kappa(A),M(-1)),
\end{equation}
natural in $M$. For $r=2$ and $M=\mu_m$, the localization
boundary of the Kummer class of a uniformizer maps to $1$.
In particular, $H^1_{\mathrm{\acute et},\{\mathfrak m_A\}}(\Spec A,M)=0$.
\item For each integer $r\geq1$, let
\[
 \partial_A:H^r_{\mathrm{\acute et}}(L,M)\longrightarrow H^{r-1}_{\mathrm{\acute et}}(\kappa(A),M(-1))
\]
be the localization boundary followed by
\eqref{eq:trait-support-residue} in degree $r+1$.
For $b\in L^\times$ with Kummer class
$[b]\in H^1_{\mathrm{\acute et}}(L,\mu_m)$ and for
$a\in H^{r-1}_{\mathrm{\acute et}}(A,M(-1))$, evaluation gives
\[
 \partial_A([b]\cup a|_L)=v_A(b)\,a|_{\kappa(A)}.
\]
\item If $d\geq2$ divides $m$ and $dM=0$, the isomorphism
\[
 \underline{\operatorname{Hom}}_{\Z/d}(\mu_d,M)
 \xrightarrow{\ \sim\ }
 \underline{\operatorname{Hom}}_{\Z/m}(\mu_m,M),
 \qquad a\longmapsto a\circ(\,\cdot\,)^{m/d},
\]
is compatible with \eqref{eq:trait-support-residue} and with
$\partial_A$ defined using $d$ and using $m$, in every degree.
\item Let $B$ be a finite \'etale $A$-algebra that is a domain.
For every closed point $y$ of $\Spec B$, the ring $B_y$ is a
discrete valuation ring. For each integer $r\geq1$ and
$\gamma\in H^r_{\mathrm{\acute et}}(\Frac B,M)$,
\begin{equation}\label{eq:residue-corestriction}
 \partial_A\operatorname{Tr}_{\Frac B/L}(\gamma)
 =\sum_{y\in\Spec(B/\mathfrak m_A B)}
   \operatorname{Tr}_{\kappa(y)/\kappa(A)}
      \bigl(\partial_{B_y}\gamma\bigr),
\end{equation}
where $\operatorname{Tr}$ is finite \'etale trace.
\item For each natural number $q\geq1$, the square
\[
\begin{tikzcd}[column sep=large]
 K_q^M(L)/m \arrow[r,"(-1)^{q-1}\partial^{\mathrm M}_A"] \arrow[d] &
 K_{q-1}^M(\kappa(A))/m \arrow[d]\\
 H^q_{\mathrm{\acute et}}(L,\mu_m^{\otimes q}) \arrow[r,"\partial_A"] &
 H^{q-1}_{\mathrm{\acute et}}(\kappa(A),\mu_m^{\otimes(q-1)})
\end{tikzcd}
\]
commutes. The vertical maps are Galois symbols. The Milnor residue
$\partial^{\mathrm M}_A$ uses the last-slot convention of
\cref{not:milnor-tame-residues}; the factor $(-1)^{q-1}$ moves the
uniformizer to the first slot. For $q=1$, the right vertical map
is the identity $\Z/m\to H^0_{\mathrm{\acute et}}(\kappa(A),\Z/m)$.
\end{enumerate}
\end{sublemma}

\begin{proof}
\emph{Purity and normalization.}
Let $i:\Spec\kappa(A)\to\Spec A$ be the closed immersion.
Choose a geometric closed point $\bar z$ and a strict
henselization $A^{\mathrm{sh}}$ at $\bar z$. A uniformizer
$\pi$ of $A$ is a uniformizer of the discrete valuation ring
$A^{\mathrm{sh}}$
\cite[\stackstag{0AP3}, \stackstag{07QM}]{Stacks}.
Every unit of $A^{\mathrm{sh}}$ has a $d$-th root for every
divisor $d\geq2$ of $m$: the polynomial defining such a root
is separable over the separably closed residue field, and
Hensel's lemma lifts a root
\cite[\stackstag{04GF}]{Stacks}.

For $r\geq2$, multilinearity shows that
$K_r^M(\Frac A^{\mathrm{sh}})/d$ is generated by
$\{\pi,\ldots,\pi\}$. This symbol is zero because
$\{\pi,\pi\}=\{\pi,-1\}$ and $-1$ is a $d$-th power.
The field norm-residue theorem
\cite[Chapter~VI, Corollary~4.1.1]{Weibel} gives
\[
 H^r_{\mathrm{\acute et}}(\Frac A^{\mathrm{sh}},\mu_d^{\otimes r})=0
 \qquad(r\geq2).
\]
The pullback of $M$ is constant. Each of its cyclic summands
of order $d$ is isomorphic to $\mu_d^{\otimes r}$ after choosing
a primitive $d$-th root of unity, so
$H^r_{\mathrm{\acute et}}(\Frac A^{\mathrm{sh}},M)=0$ for $r\geq2$.
Moreover, positive-degree cohomology of $\Spec A^{\mathrm{sh}}$
vanishes \cite[\stackstag{03QO}]{Stacks}, and restriction in
degree zero to its fraction field is an isomorphism.
The stalk formula for derived direct image
\cite[\stackstag{03Q9}]{Stacks}, applied to the localization
triangle, therefore shows that $Ri^!M$ is concentrated in
degree $2$ and that
\[
 (R^2i^!M)_{\bar z}=H^1_{\mathrm{\acute et}}(\Frac A^{\mathrm{sh}},M).
\]

Kummer theory \cite[\stackstag{03PL}]{Stacks} and the
divisibility of units give the valuation isomorphism
\[
 (\Frac A^{\mathrm{sh}})^\times/
       ((\Frac A^{\mathrm{sh}})^\times)^d
       \xrightarrow[\ v_{A^{\mathrm{sh}}}\ ]{\ \sim\ }\Z/d.
\]
Consequently
\[
 \operatorname{Hom}_{\Z/m}
       (\mu_m(A^{\mathrm{sh}}),M_{\bar z})
 \xrightarrow{\ \sim\ } H^1_{\mathrm{\acute et}}(\Frac A^{\mathrm{sh}},M),
 \qquad a\longmapsto a_*[\pi].
\]
Indeed, this follows on each cyclic summand from Kummer theory;
the map $\mu_m\to\mu_d$, $\zeta\mapsto\zeta^{m/d}$, takes
$[\pi]_m$ to $[\pi]_d$.
The isomorphism is independent of the uniformizer, natural
in $M$, and equivariant for the residue-field Galois action,
since $[\pi]$ comes from $L$. Thus
$Ri^!M\simeq M|_{\kappa(A)}(-1)[-2]$ with the normalization in
\textup{(i)}. Taking cohomology over $\kappa(A)$ proves
\eqref{eq:trait-support-residue}.

\smallskip\noindent
\emph{Cup products.}
Cup product and evaluation define a morphism
\[
 Ri^!\mu_m\otimes^{\mathbf L}_{\Z/m}i^*M(-1)
       \longrightarrow Ri^!M.
\]
Under the degree-two identification just computed, it sends
$a$ to $a_*[\pi]$. It therefore identifies
$M|_{\kappa(A)}(-1)[-2]$ with $Ri^!M$.
Tensoring the localization triangle on the right and applying
evaluation shows that its boundary $\delta$ satisfies
\[
 \delta([b]\cup a|_L)=\delta([b])\cup a.
\]
Under \eqref{eq:trait-support-residue}, the right side is
$v_A(b)\,a|_{\kappa(A)}$, since
$[b]=v_A(b)[\pi]$ over $\Frac A^{\mathrm{sh}}$.
This proves \textup{(ii)} with the stated order of the factors.

\smallskip\noindent
\emph{Coefficient maps.}
For a morphism $u:M\to M'$ of sheaves killed by $m$,
$u_*(a_*[\pi])=(u\circ a)_*[\pi]$.
If $d\mid m$ and $dM=0$, the kernel of
$\mu_m\xrightarrow{(\,\cdot\,)^{m/d}}\mu_d$ is $d\mu_m$,
which every homomorphism $\mu_m\to M$ kills. Precomposition
therefore gives the isomorphism in \textup{(iii)}.
Because this power map takes $[\pi]_m$ to $[\pi]_d$, the
degree-two support identifications agree. The concentration
of $Ri^!M$ in degree $2$ and naturality of localization give
the compatibility in every degree.

\smallskip\noindent
\emph{Corestriction and symbols.}
Finite \'etale trace is compatible with \'etale localization
and is the sum map on a disjoint union of copies of the base
\cite[\stackstag{03SE} and the characterization following it]{Stacks}.
The finite direct images on abelian sheaves are exact
\cite[\stackstag{03QP}]{Stacks}, so trace induces a morphism of
localization triangles. After strict
henselization at the closed point, the cover splits into such
copies, on each of which $\pi$ remains a uniformizer.
Hence trace on the degree-two support sheaves is the trace on
$M(-1)$. Taking cohomology over $\kappa(A)$ gives
\eqref{eq:residue-corestriction} in every degree. The rings
$B_y$ are discrete valuation rings because $B$ is \'etale over $A$.

For \textup{(v)}, multilinearity, graded commutativity and
$\{\pi,\pi\}=\{\pi,-1\}$ show that $K_q^M(L)$ is generated
by symbols $\{u_1,\ldots,u_q\}$ and
$\{\pi,u_2,\ldots,u_q\}$ with $u_i\in A^\times$.
Both composites in the square vanish on the first kind:
its Galois symbol extends over $A$.
For the second kind, the last-slot convention gives
\[
 \partial^{\mathrm M}_A\{\pi,u_2,\ldots,u_q\}
   =(-1)^{q-1}\{\overline u_2,\ldots,\overline u_q\}.
\]
The factor $(-1)^{q-1}$ in the top arrow cancels this sign.
By \textup{(ii)}, the other composite is
$[\overline u_2]\cup\cdots\cup[\overline u_q]$ as well. For $q=1$,
both composites send $\{\pi\}$ to $1$.
\end{proof}

\subsection{Unit symbols and finite \'etale effacement}

\begin{sublemma}[Unramified Milnor classes on a smooth local ring]
\label{sub:unramified-milnor-smooth-local}
Let $V$ be an excellent discrete valuation ring, let $B$ be a
nonzero local ring essentially smooth over $V$, and let
$m\geq2$ be a natural number invertible in $V$. For every
natural number $q\geq1$, the sequence
\begin{equation}\label{eq:local-symbols-milnor-gersten}
 K_q^M(B)/m\longrightarrow K_q^M(\Frac B)/m
 \xrightarrow{(\partial^{\mathrm M}_{B_z})_z}
 \bigoplus_{z\in(\operatorname{Spec}B)^{(1)}}
                    K_{q-1}^M(\kappa(z))/m
\end{equation}
is exact at the middle term. Here
$\partial^{\mathrm M}_{B_z}$ is the tame residue of the discrete
valuation ring $B_z$, with the last-slot convention of
\cref{sub:trait-locally-constant-residue}\textup{(v)}.
\end{sublemma}

\begin{proof}
The ring $B$ is regular by \cite[\stackstag{07NF}]{Stacks}.
For $q=1$, factoriality of regular local rings
\cite[\stackstag{0AG0}]{Stacks} shows that an element of
$(\Frac B)^\times$ whose valuation at every height-one prime of $B$
is divisible by $m$ is the product of a unit of $B$ and an $m$-th power.

Let $q\geq2$, and write $B=\mathcal O_{X,y}$ for a smooth $V$-scheme
$X$ and a point $y\in X$. The sequence
\eqref{eq:local-symbols-milnor-gersten} with the improved Milnor
group $\widehat K_q^M(B)/m$ in its first term is the stalk at $y$ of
the finite-coefficient Gersten sequence of Zariski sheaves of
\cite[\arxivhtml{2010.02622v4}{S4.Thmtheorem3}{Corollary~4.3} and
\arxivhtml{2010.02622v4}{S4.Thmtheorem4}{Remark~4.4}]{LudersMilnor24}
on $X$; in particular it is exact at the middle term. Here we use
the finite-coefficient assertion in
\cite[\arxivhtml{2010.02622v4}{S4.Thmtheorem4}{Remark~4.4}]{LudersMilnor24},
not merely reduction of the integral sequence modulo $m$.
The natural homomorphism $K_q^M(B)\to\widehat K_q^M(B)$ is surjective
\cite[Theorem~B]{Kerz10}, hence so is its reduction modulo $m$.
\end{proof}

\begin{sublemma}[Unit symbols; {\cite[\arxivhtml{2201.11553v4}{S3.Thmtheorem8}{Corollary~3.8}]{Luders22}}]
\label{sub:local-symbols-via-smooth-sources}
Let $V$ be an excellent discrete valuation ring, let $m\geq 2$ be a
natural number invertible in $V$, and let $B$ be a nonzero local $V$-algebra.
For every natural number $q\geq 1$, the group
$H^q_{\mathrm{\acute et}}(B,\mu_m^{\otimes q})$ is generated by cup products
of $q$ Kummer classes of units of $B$.
\end{sublemma}

\begin{proof}
We follow \cite[proof of \arxivhtml{2201.11553v4}{S3.Thmtheorem8}{Corollary~3.8}]{Luders22}, using
\cite[Theorem~1.2]{Geisser04} for the \'etale comparison in the
smooth case.

\smallskip\noindent
\emph{Injectivity on smooth local rings.}
Suppose first that $B=\mathcal O_{X,y}$, where $X$ is smooth over
$V$ and $y\in X$. The ring $B$ is regular because $X$ is smooth over
the regular scheme $\Spec V$ \cite[\stackstag{07NF}]{Stacks}, and is a domain by
\cite[\stackstag{00NP}]{Stacks}; let $F=\Frac B$. The
norm-residue theorem \cite[Chapter~VI, Theorem~4.1]{Weibel} supplies
the Bloch--Kato hypothesis of \cite[Theorem~1.2]{Geisser04}, and
\cite[Theorem~1.2(2), (4), (5) and p.~775]{Geisser04} give a Gersten
resolution of the Zariski sheaf on $X$ associated with
$U\mapsto H^q_{\mathrm{\acute et}}(U,\mu_m^{\otimes q})$, whose first
map goes to
$\bigoplus_{x\in X^{(0)}}i_{x*}H^q_{\mathrm{\acute et}}(\kappa(x),\mu_m^{\otimes q})$.
\'Etale cohomology commutes with the limit
$\Spec B=\lim_U U$ over the affine open neighbourhoods $U$ of $y$
\cite[\stackstag{09YQ}]{Stacks}, so the stalk at $y$ of this first
map is the restriction map
\begin{equation}\label{eq:local-symbols-injection}
 H^q_{\mathrm{\acute et}}(B,\mu_m^{\otimes q})
     \longrightarrow H^q_{\mathrm{\acute et}}(F,\mu_m^{\otimes q}),
\end{equation}
which is therefore injective.

\smallskip\noindent
\emph{Symbols on smooth local rings.}
Let $\alpha\in H^q_{\mathrm{\acute et}}(B,\mu_m^{\otimes q})$. By the
norm-residue isomorphism for fields
\cite[Chapter~VI, Corollary~4.1.1]{Weibel}, there is
$\beta\in K_q^M(F)/m$ whose Galois symbol is $\alpha|_F$. Let $z$ be
a height-one prime ideal of $B$. The class $\alpha|_F$ comes from
$H^q_{\mathrm{\acute et}}(B_z,\mu_m^{\otimes q})$, so its residue at
$z$, defined by the localization sequence of the discrete valuation
ring $B_z$ with purity, is zero. By \cref{sub:trait-locally-constant-residue}\textup{(v)},
the Galois symbol of the tame residue of $\beta$ at $z$ is zero.
For $q\geq2$, this symbol is injective by the field norm-residue
theorem \cite[Chapter~VI, Corollary~4.1.1]{Weibel}; for $q=1$, it
is the identity $K_0^M(\kappa(z))/m=\Z/m\to H^0_{\mathrm{\acute et}}(\kappa(z),\Z/m)$.
Thus the tame residue is zero, and $\beta$ lies in the kernel of
the second map of \eqref{eq:local-symbols-milnor-gersten}.
By \cref{sub:unramified-milnor-smooth-local}, there is an element
$\sum_kn_k\{u_{k,1},\ldots,u_{k,q}\}\in K^M_q(B)/m$ mapping to
$\beta$, with $n_k\in\Z$ and $u_{k,i}\in B^\times$. Put
\[
 \alpha'=\sum_kn_k\operatorname{Kum}_m(u_{k,1})\cup\cdots\cup
          \operatorname{Kum}_m(u_{k,q})
 \in H^q_{\mathrm{\acute et}}(B,\mu_m^{\otimes q}).
\]
Restriction to $F$ commutes with Kummer boundaries and cup products,
so $\alpha'|_F$ is the Galois symbol of $\beta$, that is $\alpha|_F$.
By the injectivity of \eqref{eq:local-symbols-injection},
$\alpha=\alpha'$.

\smallskip\noindent
\emph{Henselian approximation.}
For general $B$, let $\mathfrak m_B$ be its maximal ideal. Let $P$ be the
polynomial $V$-algebra on the set $B$, with the surjection $P\to B$,
$X_b\mapsto b$, let $P_0$ be the localization of $P$ at the inverse
image of $\mathfrak m_B$, and put $I=\ker(P_0\to B)$. Let
$(P_0^h,IP_0^h)$ be the henselization of the pair $(P_0,I)$. By
\cite[\stackstag{0A02}]{Stacks} and its proof, $P_0^h$ is the filtered
colimit of the \'etale $P_0$-algebras $E$ with $P_0/I\simeq E/IE$,
and its map to $B$ induces an isomorphism $P_0^h/IP_0^h\xrightarrow{\sim}B$.
The ideal $IP_0^h$ is contained in the Jacobson
radical of $P_0^h$ \cite[\stackstag{09XE}]{Stacks}, so $P_0^h$ is local.
By Gabber's theorem \cite[\stackstag{09ZI}]{Stacks}, restriction
\[
 H^q_{\mathrm{\acute et}}(P_0^h,\mu_m^{\otimes q})
 \longrightarrow H^q_{\mathrm{\acute et}}(B,\mu_m^{\otimes q})
\]
is an isomorphism, and it maps Kummer symbols of units to Kummer
symbols of units. It therefore suffices to write every
$\alpha\in H^q_{\mathrm{\acute et}}(P_0^h,\mu_m^{\otimes q})$ as a sum
of Kummer symbols of units of $P_0^h$.

The ring $P_0$ is the filtered colimit of the smooth $V$-algebras
$Q_\nu=V[X_b\mid b\in\lambda][1/s]$, indexed by the pairs
$\nu=(\lambda,s)$ of a finite subset $\lambda\subseteq B$ and an
element $s$ of $V[X_b\mid b\in\lambda]$ outside the inverse image of
$\mathfrak m_B$, ordered by $(\lambda,s)\leq(\lambda',s')$ if
$\lambda\subseteq\lambda'$ and $s$ divides $s'$ in
$V[X_b\mid b\in\lambda']$ \cite[\stackstag{00CR}]{Stacks}. By continuity of
\'etale cohomology \cite[\stackstag{03Q6}]{Stacks}, the class $\alpha$
comes from a class on an \'etale $P_0$-algebra $E$ with
$P_0/I\simeq E/IE$ in the filtered colimit defining $P_0^h$.
This algebra is of finite presentation over $P_0$, so there are an
index $\nu$ and an \'etale $Q_\nu$-algebra $E_\nu$ with
$E\simeq E_\nu\otimes_{Q_\nu}P_0$
\cite[\stackstag{05N9}, \stackstag{07RI}]{Stacks}. By continuity
again, $\alpha$ comes from a class on $E_\nu\otimes_{Q_\nu}Q_{\nu'}$
for some $\nu'\geq\nu$; this ring is smooth of finite type over $V$
\cite[\stackstag{00U9}, \stackstag{00TD}]{Stacks}. Let $C$ be its
localization at the inverse image of the maximal ideal of $P_0^h$. Then
$C$ is the local ring of a scheme smooth over $V$, the map $C\to P_0^h$
is local, and $\alpha$ is the image of a class on $C$. By the second
step, that class is a sum of Kummer symbols of units of $C$, and its
image $\alpha$ is the sum of the Kummer symbols of their images in
$P_0^h$.
\end{proof}
\begin{lemma}[Milnor excision with Tate coefficients]
\label{lem:tate-milnor-excision}
Let \(m\geq2\) be a natural number, let \(j\in\mathbf N\), and let
\[
 \begin{tikzcd}[column sep=large]
 A' \arrow[r] \arrow[d,two heads] & A \arrow[d,two heads]\\
 B' \arrow[r] & B
 \end{tikzcd}
\]
be a cartesian square of \(\Z[1/m]\)-algebras whose vertical
maps are surjective. The induced square
\[
 \begin{tikzcd}[column sep=large]
 R\Gamma_{\mathrm{\acute et}}(A',\mu_m^{\otimes j})
       \arrow[r] \arrow[d] &
 R\Gamma_{\mathrm{\acute et}}(A,\mu_m^{\otimes j}) \arrow[d]\\
 R\Gamma_{\mathrm{\acute et}}(B',\mu_m^{\otimes j}) \arrow[r] &
 R\Gamma_{\mathrm{\acute et}}(B,\mu_m^{\otimes j})
 \end{tikzcd}
\]
is homotopy cartesian.
\end{lemma}

\begin{proof}
The cartesian square identifies
\(\ker(A'\to B')\) with \(\ker(A\to B)\); hence it is a
Milnor square in the sense of
\cite[Definition~1.11]{BhattMathew21}.
The assertion follows from the excision statement of
\cite[\arxivhtml{1807.04725v3}{S5.Thmtheorem4}{Theorem~5.4}]{BhattMathew21},
applied to the torsion sheaf \(\mu_m^{\otimes j}\) on
\(\operatorname{Spec}\Z[1/m]\), whose pullbacks are the indicated coefficients.

Alternatively, one can use excision for constant coefficients and
finite \'etale descent. Let \(D=\Z[1/m][t]/(\Phi_m(t))\), where \(\Phi_m\) is the
\(m\)-th cyclotomic polynomial. The morphism
\[
 \operatorname{Spec}D\longrightarrow
                    \operatorname{Spec}\Z[1/m]
\]
is the finite \'etale torsor of primitive \(m\)-th roots under
\((\Z/m)^\times\). Indeed, \(D\) is finite free over
\(\Z[1/m]\), the primitive roots in every geometric fibre are
simple roots of \(t^m-1\), and the power action of
\((\Z/m)^\times\) is simply transitive on those roots.
Tensoring the square with \(D\) over \(\Z[1/m]\) preserves
the pullback and the identified kernels, because \(D\) is flat.
The constant-coefficient excision theorem
\cite[Corollary~1.17]{BhattMathew21} therefore gives a cartesian
square of derived \'etale sections with coefficient \(\Z/m\)
on the four base changes.

For each of the four rings \(T\in\{A',A,B',B\}\), the universal
root defines an isomorphism of sheaves on
\(\operatorname{Spec}(T\otimes_{\Z[1/m]}D)\):
\[
 \Z/m\xrightarrow{\ \sim\ }\mu_m^{\otimes j},
 \qquad 1\longmapsto t^{\otimes j}.
\]
For \(j=0\), this is the identity of the constant sheaf
\(\Z/m\). For \(a\in(\Z/m)^\times\), the automorphism
\(t\mapsto t^a\) sends \(t^{\otimes j}\) to
\(a^j t^{\otimes j}\). Thus the transported descent action
on constant-coefficient sections is the pullback action
multiplied by \(a^j\). All four ring maps preserve the same
universal root, so the cartesian square is equivariant for
these transported actions.

The \v{C}ech nerve of the finite \'etale torsor over
\(\operatorname{Spec}T\) has degree-\(r\) term
\[
 \operatorname{Spec}(T\otimes_{\Z[1/m]}D)
                       \times((\Z/m)^\times)^r.
\]
Descent for derived \'etale sections therefore identifies
\[
 R\Gamma_{\mathrm{\acute et}}(T,\mu_m^{\otimes j})
 \simeq
 R\Gamma_{\mathrm{\acute et}}
   (T\otimes_{\Z[1/m]}D,\Z/m)^{h((\Z/m)^\times)},
\]
where the right side uses the transported action just described.
Taking homotopy fixed points preserves cartesian squares,
because both operations are limits. Applied to the equivariant
constant-coefficient square, this gives the asserted square.
\end{proof}

\begin{proposition}[Symbols on a semilocal ring]
\label{prop:semilocal-unit-symbols}
Let \(B\) be a nonzero semilocal ring whose closed residue
fields have at most one positive characteristic. Let \(m\geq2\)
be a natural number invertible in \(B\). For every integer \(q\geq1\), the group
\(H^q_{\mathrm{\acute et}}(B,\mu_m^{\otimes q})\) is generated by
Kummer symbols of degree \(q\) in the sense of
\cref{def:kummer-symbol}.
\end{proposition}

\begin{proof}
If a closed residue field has positive characteristic, denote
this characteristic by \(p\). If all the closed residue
fields have characteristic zero, choose a prime \(p\) not
dividing \(m\). Every integer prime to \(p\) is a unit in
\(B\), since its image in every closed residue field is nonzero.
Thus the structure homomorphism factors through
\(\Z_{(p)}\longrightarrow B\), and \(m\) is invertible in
\(\Z_{(p)}\). Let \(\mathfrak m_1,\ldots,\mathfrak m_s\)
be the maximal ideals of \(B\), and let
\(J=\bigcap_i\mathfrak m_i\) be its Jacobson radical.
The quotient maps induce an isomorphism
\[
 B/J\xrightarrow{\ \sim\ }\prod_{i=1}^s B/\mathfrak m_i
\]
by the Chinese remainder theorem \cite[\stackstag{00DT}]{Stacks}.

First, Kummer symbols suffice in degree one. The Kummer
sequence \cite[\stackstag{03PL}]{Stacks} gives
\[
 0\longrightarrow B^\times/(B^\times)^m
   \longrightarrow H^1_{\mathrm{\acute et}}(B,\mu_m)
   \longrightarrow \operatorname{Pic}(B)[m]
   \longrightarrow0.
\]
Every invertible \(B\)-module is free of rank one by
\cref{sub:semilocal-projective-free}, so
\(\operatorname{Pic}(B)=0\). The map
\(B^\times/(B^\times)^m\to H^1_{\mathrm{\acute et}}(B,\mu_m)\) is
therefore surjective, proving the degree-one assertion.

For any \(q\geq2\), let \(V_0\) be the image of \(\Z_{(p)}\) in
\(B/J\). This image is \(\mathbf F_p\) if all the closed
residue fields have characteristic \(p\). If some have
characteristic zero, the map to any such factor is injective
and \(V_0\simeq \Z_{(p)}\). Put \(B'=B\times_{B/J} V_0\). We may view
\(B'\) as the subring of \(B\) of elements whose residues lie
in \(V_0\). The ideal \(J\) is contained in the Jacobson
radical of \(B'\): for \(x\in J\) and \(y\in B'\), the
inverse of \(1-xy\) in \(B\) has residue \(1\) in every
factor, hence belongs to \(B'\). Since \(B'/J=V_0\) is local,
\(B'\) is local. It is a \(\Z_{(p)}\)-algebra and may be singular.
The DVR \(\Z_{(p)}\) is excellent, so
\cref{sub:local-symbols-via-smooth-sources} applies to \(B'\).

The cartesian square of rings
\[
 \begin{tikzcd}[column sep=large]
 B' \arrow[r] \arrow[d,two heads] & B \arrow[d,two heads]\\
 V_0 \arrow[r] & B/J
 \end{tikzcd}
\]
has surjective vertical maps. By \cref{lem:tate-milnor-excision},
it induces a cartesian square of derived \'etale sections
with coefficient \(\mu_m^{\otimes q}\).
In degree \(q\), the resulting fibre sequence gives the exact sequence
\[
 \begin{aligned}
 H^q_{\mathrm{\acute et}}(B',\mu_m^{\otimes q})&\longrightarrow
       H^q_{\mathrm{\acute et}}(B,\mu_m^{\otimes q})
        \oplus H^q_{\mathrm{\acute et}}(V_0,\mu_m^{\otimes q})\\
      &\longrightarrow H^q_{\mathrm{\acute et}}(B/J,\mu_m^{\otimes q}).
 \end{aligned}
\]
Let \(\alpha\in H^q_{\mathrm{\acute et}}(B,\mu_m^{\otimes q})\).
By the field norm-residue theorem
\cite[Chapter~VI, Corollary~4.1.1]{Weibel}, each component of
its restriction to \(B/J\) is a sum of unit symbols. Lift every
entry of each symbol over \(B/\mathfrak m_i\) to a unit of \(B\) with the
prescribed residue in \(B/\mathfrak m_i\) and residue \(1\) in all
other factors. The Chinese remainder theorem \cite[\stackstag{00DT}]{Stacks} gives such a
lift, which is a unit because all its residues are nonzero.
The resulting sum \(\beta\) of symbols on \(B\) has
restriction \(\alpha|_{B/J}\). The pair \((\alpha-\beta,0)\)
therefore lifts to cohomology on \(B'\). By
\cref{sub:local-symbols-via-smooth-sources}, this lift is a
sum of unit symbols. Restriction carries these symbols to
unit symbols on \(B\). Adding \(\beta\) gives the asserted
expression for \(\alpha\).
\end{proof}

\begin{proposition}[Effacement by finite \'etale covers]
\label{prop:finite-etale-effacement}
Let \(B\) be a nonzero semilocal ring whose closed residue fields
have at most one positive characteristic, as in
\cref{prop:semilocal-unit-symbols}, and let \(\ell\) be a prime
number invertible in \(B\). Let \(T\to\Spec B\) be a finite \'etale
morphism, let \(M\) be a finite locally constant sheaf of abelian
groups of \(\ell\)-power order on \(T_{\mathrm{\acute et}}\), let
\(q>0\), and let \(\alpha\in H^q_{\mathrm{\acute et}}(T,M)\). Then
there is a finite \'etale surjective morphism \(g:T'\to T\) such
that
\[
 g^*\alpha=0\quad\text{in}\quad
 H^q_{\mathrm{\acute et}}(T',g^*M).
\]
\end{proposition}

\begin{proof}
If \(T=\varnothing\), take \(T'=\varnothing\); assume
\(T\neq\varnothing\), and choose \(\nu\geq1\) such that \(\ell^\nu\) annihilates \(M\).

\smallskip\noindent
\emph{Connected finite \'etale \(T\)-schemes.}
Let \(P\) be a connected finite \'etale \(T\)-scheme. Then
\(P=\Spec B_P\) for a finite \'etale \(B\)-algebra \(B_P\), which is
semilocal because it is finite over \(B\). A maximal ideal
\(\mathfrak n\) of \(B_P\) lies over a maximal ideal \(\mathfrak m\)
of \(B\) \cite[\stackstag{00GR}]{Stacks}, and \(B_P/\mathfrak n\) is
a quotient of the finite \(B/\mathfrak m\)-algebra
\(B_P/\mathfrak mB_P\). Hence the closed residue fields of \(B_P\)
have at most one positive characteristic, and \(B_P\) satisfies the
hypotheses of \cref{prop:semilocal-unit-symbols} for every divisor
\(m\geq2\) of \(\ell^\nu\). The connected components of the semilocal
scheme \(T\) are closed \cite[\stackstag{004T}]{Stacks} and each
contains a closed point (a maximal ideal of \(\mathcal O(T)\)
containing the prime ideal of any point in the component).
Thus they are finitely many, and hence open as well.

\smallskip\noindent
\emph{Trivializing the coefficients.}
Let \(T_c\) be a connected component of \(T\) and let \(\bar t\) be a
geometric point of \(T_c\). By \cite[\stackstag{03RV},
\stackstag{0BND}]{Stacks}, the finite locally constant sheaf
\(M\oplus\mu_{\ell^\nu}\) on \(T_c\) corresponds to a finite abelian
group with a continuous action of \(\pi_1^{\mathrm{\acute et}}(T_c,\bar t)\) by group
automorphisms. Its kernel \(J\) is an open subgroup; let
\(P_c\to T_c\) be the connected finite \'etale cover corresponding to
the \(\pi_1^{\mathrm{\acute et}}(T_c,\bar t)\)-set \(\pi_1^{\mathrm{\acute et}}(T_c,\bar t)/J\)
\cite[\stackstag{0BND}]{Stacks}. For a geometric point \(\bar p\) of
\(P_c\) over \(\bar t\), the image of \(\pi_1^{\mathrm{\acute et}}(P_c,\bar p)\) in
\(\pi_1^{\mathrm{\acute et}}(T_c,\bar t)\) is \(J\), so \(M|_{P_c}\) and
\(\mu_{\ell^\nu}|_{P_c}\) are constant sheaves of abelian groups
\cite[\stackstag{0BND}]{Stacks}. Put \(T_1=\coprod_cP_c\); the
morphism \(g_1:T_1\to T\) is finite \'etale and surjective. On a
connected component \(P\) of \(T_1\), choose a generator \(\zeta\) of
the cyclic group \(\mu_{\ell^\nu}(P)\) of order \(\ell^\nu\) and an
isomorphism \(M|_P\simeq\bigoplus_k\Z/\ell^{b_k}\) with
\(1\leq b_k\leq \nu\). For \(m=\ell^b\) with \(1\leq b\leq \nu\), the
root \(\zeta^{\ell^{\nu-b}}\) identifies the constant sheaf \(\Z/m\)
with \(\mu_m\), hence with \(\mu_m^{\otimes q}\), and therefore
\begin{equation}\label{eq:effacement-components}
 H^q_{\mathrm{\acute et}}(P,M|_P)\simeq
 \bigoplus_kH^q_{\mathrm{\acute et}}\bigl(P,\mu_{\ell^{b_k}}^{\otimes q}\bigr).
\end{equation}

\smallskip\noindent
\emph{Symbols and their roots.}
By \cref{prop:semilocal-unit-symbols} for the ring \(\mathcal O(P)\),
the component of \((g_1^*\alpha)|_P\) in the summand of index \(k\) of
\eqref{eq:effacement-components} is a finite sum of Kummer symbols
\(\operatorname{Kum}_m(u_1)\cup\cdots\cup\operatorname{Kum}_m(u_q)\),
where \(m=\ell^{b_k}\) and \(u_1,\ldots,u_q\in\mathcal O(P)^\times\).
For such a unit \(u\), the \(\mathcal O(P)\)-algebra
\(\mathcal O(P)[z]/(z^m-u)\) is finite free of rank \(m\), hence
faithfully flat, and it is standard \'etale, because \(z\) and \(m\)
are units in it and so is the derivative \(mz^{m-1}\)
\cite[\stackstag{00UB}, \stackstag{00UC}]{Stacks}. In it
\(u=z^m\), so \(\operatorname{Kum}_m(u)\) vanishes there by the
Kummer sequence \cite[\stackstag{03PL}]{Stacks}. Let \(P'\) be the
spectrum of the tensor product over \(\mathcal O(P)\) of these
algebras, for all units occurring in all components, and put
\(T'=\coprod_PP'\). The composite \(g:T'\to T_1\to T\) is finite
\'etale and surjective. Every factor of each Kummer symbol in the
chosen decompositions of the components of
\eqref{eq:effacement-components} vanishes on \(T'\), so
\(g^*\alpha=0\).
\end{proof}
\subsection{The bounded spectral sheaf on the finite \'etale site}
\label{subsec:bounded-descent}

\begin{notation}[Finite \'etale sites and sheaves of spaces]
\label{not:finite-etale-topoi}
A \emph{space} is an object of the \(\infty\)-category \(\mathcal S\)
of \cref{sec:conventions}. For a site \(\mathcal C\), let
\(\mathrm{Sh}(\mathcal C)\) be the \(\infty\)-topos of sheaves of
spaces on \(\mathcal C\) \cite[Definition~6.2.2.6]{HTT}; a
\emph{sheaf of spaces on \(\mathcal C\)} is an object of
\(\mathrm{Sh}(\mathcal C)\). The constant sheaf
\(\underline{\mathsf A}\in\mathrm{Sh}(\mathcal C)\) with value a space \(\mathsf A\) is
the sheaf associated with the constant presheaf, and
\(\Gamma(\mathcal C,Z)=\operatorname{Map}(\underline{\ast},Z)\) is
the space of global sections of \(Z\in\mathrm{Sh}(\mathcal C)\).

For a scheme \(S\), let \(S_{\mathrm{\acute et}}\) be its small
\'etale site (\cref{not:etale-sites-spectra}), without imposing the
quasi-compactness and quasi-separatedness restrictions of
\cref{sec:finite-connective-comparison}, and let
\(S_{\mathrm{fet}}\subseteq S_{\mathrm{\acute et}}\) be the
\emph{finite \'etale site}: the full subcategory of finite \'etale
\(S\)-schemes, whose coverings are the jointly surjective families.
Restriction of sheaves along
\(S_{\mathrm{fet}}\subseteq S_{\mathrm{\acute et}}\) is the direct
image \(\rho_*\) of a geometric morphism
\[
 \rho:\mathrm{Sh}(S_{\mathrm{\acute et}})\longrightarrow
      \mathrm{Sh}(S_{\mathrm{fet}}).
\]
On both sites we use hypercomplete sheaves of spectra
(\cref{not:etale-sites-spectra}); \(\rho_*\) and \(\rho^*\) also
denote the induced direct image and hypercompleted inverse image on
them. We denote by \(HM\) the Eilenberg--Mac Lane sheaf of spectra
of a sheaf of abelian groups \(M\), placed in degree zero. For a
profinite group \(\Pi\), let \(\Pi\text{-}\mathrm{FSet}\) be the site
of finite continuous \(\Pi\)-sets, whose coverings are the jointly
surjective families.
\end{notation}

\begin{sublemma}[The finite \'etale site]\label{sub:finite-site-basics}
Let \(S\) be a connected scheme and let \(\bar s\) be a geometric
point of \(S\) (\cref{sec:conventions}).
\begin{enumerate}
\item For a sheaf of spaces \(Z\) on \(S_{\mathrm{fet}}\), put
\begin{equation}\label{eq:finite-site-stalk}
 Z_{\bar s}=\varinjlim_{(P,\bar p)} Z(P),
\end{equation}
where the indexing category is the opposite of the category of
connected finite \'etale \(S\)-schemes \(P\) with a geometric point
\(\bar p\) over \(\bar s\), whose morphisms are the pointed
\(S\)-morphisms. This indexing category is filtered, and its full
subcategory of pointed connected Galois covers is cofinal. A morphism \(\alpha\) of sheaves of sets, of
truncated sheaves of spaces, or of hypercomplete sheaves of spectra
on \(S_{\mathrm{fet}}\) is an equivalence if and only if
\(\alpha_{\bar s}\) is an equivalence.
\item Let \(f\) be a geometric morphism between \(\infty\)-topoi of
the form \(\mathrm{Sh}(\mathcal C)\) in
\cref{not:finite-etale-topoi}. On hypercomplete sheaves of spectra,
the inverse image \(f^*\) is exact and satisfies
\(\pi_i(f^*\mathcal F)=f^*\pi_i(\mathcal F)\) for all \(i\in\Z\), and the direct image
\(f_*\) is exact and preserves hypercomplete sheaves.
\item For a hypercomplete sheaf of spectra \(\mathcal F\) on
\(S_{\mathrm{\acute et}}\) and a finite \'etale \(S\)-scheme \(P\),
\begin{equation}\label{eq:finite-site-direct-image}
 (\rho_*\mathcal F)(P)=R\Gamma_{\mathrm{\acute et}}(P,\mathcal F|_P).
\end{equation}
\item If \(S=\Spec L\) for a field \(L\), then \(\rho\) is an
equivalence.
\end{enumerate}
\end{sublemma}
\begin{proof}
\textup{(i)} The Galois category theorem identifies the site
\(S_{\mathrm{fet}}\) with the site
\(\pi_1^{\mathrm{\acute et}}(S,\bar s)\text{-}\mathrm{FSet}\) of
\cref{not:finite-etale-topoi} \cite[\stackstag{0BND}]{Stacks}.
The underlying-set functor preserves finite limits and sends covers
to surjections. It defines a point with stalk \eqref{eq:finite-site-stalk}.
Connected Galois covers are cofinal. Indeed, the pointed connected covers
correspond to the coset sets \(\pi_1^{\mathrm{\acute et}}(S,\bar s)/J\),
pointed at \(J\), where \(J\) is an open subgroup of
\(\pi_1^{\mathrm{\acute et}}(S,\bar s)\). Each \(J\) contains an open
normal subgroup, and finite intersections of open normal subgroups
refine any finite collection of these covers. The category of sheaves
of sets is therefore the category of discrete continuous
\(\pi_1^{\mathrm{\acute et}}(S,\bar s)\)-sets. The stalk
functor \eqref{eq:finite-site-stalk} therefore detects isomorphisms
of sheaves of sets. It detects equivalences of truncated space
sheaves by applying this fact to the sheaves of components and
of homotopy groups, with all local base points. It also detects
equivalences of hypercomplete spectral sheaves by testing the
homotopy sheaves in all spectrum degrees.

\textup{(ii)} and \textup{(iii)}. An inverse image functor of
\(\infty\)-topoi preserves finite limits and colimits.
It preserves \(\infty\)-connective morphisms, since it preserves
the successive connectivity tests by iterated diagonals and
effective epimorphisms. Its right adjoint thus takes hypercomplete
objects to hypercomplete objects. On the hypercomplete topoi the
inverse image is the ordinary inverse image followed by
hypercompletion. The induced inverse image on spectra is exact and satisfies
\(\pi_i(f^*\mathcal F)=f^*\pi_i(\mathcal F)\) for any inverse-image morphism \(f^*\). The spectral direct
image is exact, since it is a right adjoint between stable
categories. A finite \'etale object represents the same sheaf under inverse
image. The sections adjunction tested on that object proves
\eqref{eq:finite-site-direct-image}.

\textup{(iv)} For a field \(L\), an affine \'etale \(L\)-scheme of finite type is
a finite disjoint union of spectra of finite separable extensions
of \(L\). Such affine schemes form a basis of the small \'etale
site. Its restriction to this basis is the finite \'etale site,
which proves the last assertion.
\end{proof}

\begin{proposition}[Descent of a bounded spectrum]\label{prop:bounded-finite-etale-descent}
Let \(S\) be a connected scheme and let \(\ell\) be a prime number.
Suppose that the following effacement property holds: for every
connected finite \'etale \(S\)-scheme \(P\), every finite locally
constant sheaf \(M\) of abelian groups of \(\ell\)-power order on
\(P_{\mathrm{\acute et}}\), every \(q>0\), and every
\(\alpha\in H^q_{\mathrm{\acute et}}(P,M)\), there is a finite
\'etale surjective morphism \(g:P'\to P\) with \(g^*\alpha=0\).
Let \(\mathcal F\) be a hypercomplete sheaf of spectra on
\(S_{\mathrm{\acute et}}\) such that \(\pi_i\mathcal F=0\) for all but
finitely many \(i\in\Z\) and every \(\pi_i\mathcal F\) is a finite locally
constant sheaf of abelian groups of \(\ell\)-power order.
\begin{enumerate}
\item The sheaf \(\rho_*\mathcal F\) on \(S_{\mathrm{fet}}\) has only
finitely many nonzero homotopy sheaves, and the counit
\(\rho^*\rho_*\mathcal F\to \mathcal F\) is an equivalence.
\item For every \(i\in\Z\), the sheaf \(\pi_i(\rho_*\mathcal F)\) on
\(S_{\mathrm{fet}}\) is the restriction of the sheaf \(\pi_i\mathcal F\) to
\(S_{\mathrm{fet}}\); it corresponds to the finite continuous
\(\pi_1^{\mathrm{\acute et}}(S,\bar s)\)-module
\((\pi_i\mathcal F)_{\bar s}\), for any geometric point \(\bar s\) of
\(S\).
\end{enumerate}
\end{proposition}

\begin{proof}
Fix a geometric point \(\bar s\) of \(S\), and let \(M\) be a
finite locally constant sheaf of abelian groups of \(\ell\)-power
order on \(S_{\mathrm{\acute et}}\). For
\(i\leq0\), the stalk of \(\pi_i(\rho_*HM)\) at \(\bar s\) is
\begin{equation}\label{eq:finite-site-em-stalk}
  \varinjlim_{(P,\bar p)} H^{-i}_{\mathrm{\acute et}}(P,M|_P)
\end{equation}
by \eqref{eq:finite-site-stalk} and
\eqref{eq:finite-site-direct-image}. Homotopy groups of spectra
commute with filtered colimits. There is no homotopy in positive
degrees in the sections of the Eilenberg--Mac Lane sheaf, so the
stalks vanish for \(i>0\).

For \(i<0\) the colimit \eqref{eq:finite-site-em-stalk} is zero. A representative is
a positive-degree cohomology class on one pointed connected
\(P\). By the effacement property, there is a finite \'etale surjective
cover of \(P\) on which it vanishes. Choose a lift of \(\bar p\), and then take a connected
pointed Galois cover of \(S\) refining the connected component
through that lift. The resulting map to \(P\) factors through the
effacing cover. Thus the original representative is already zero
at one later stage. For \(i=0\), evaluating sections at \(\bar p\)
identifies the colimit with \(M_{\bar s}\). Covers trivializing
the finite local system are cofinal, and on a connected
trivializing cover every fibre element gives a unique constant
section. The action on these stalk sections is the given
monodromy action. Consequently \(\pi_0(\rho_*HM)\) is the finite
local system with that action.

The direct image is hypercomplete. Its calculated homotopy
sheaves identify it with the Eilenberg--Mac Lane object for this
finite local system on \(S_{\mathrm{fet}}\). Its inverse image
has the corresponding local system \(M\) in degree zero,
because inverse image takes a finite \'etale representative to
the same representative on \(S_{\mathrm{\acute et}}\). The same pointed-cover calculation at any geometric point of
\(S\) identifies the counit in degree zero with the evaluation
isomorphism at that point. Hence the counit is an equivalence for
\(HM\).

If the homotopy sheaves of \(\mathcal F\) occur between degrees \(a\) and
\(b\), its Postnikov filtration has the successive fibre sequences
\begin{equation}\label{eq:finite-site-postnikov-sequences}
 H(\pi_i\mathcal F)[i]\longrightarrow \tau_{\leq i}\mathcal F
            \longrightarrow \tau_{\leq i-1}\mathcal F
 \qquad (a\leq i\leq b).
\end{equation}
Its initial stage at \(a-1\) is zero and its last stage is \(\mathcal F\).
By exactness of the functors, the counit on the middle term of
\eqref{eq:finite-site-postnikov-sequences} is an equivalence once it
is an equivalence on the first and last terms. This proves the
counit assertion by induction on \(i\). The sequences
\eqref{eq:finite-site-postnikov-sequences} also show that
\(\rho_*\mathcal F\) is bounded. The counit induces an isomorphism of abelian
sheaves \(\rho^*\pi_i(\rho_*\mathcal F)\xrightarrow{\sim}\pi_i\mathcal F\) by
\cref{sub:finite-site-basics}\textup{(ii)}. Its adjoint is the canonical map
\(\pi_i(\rho_*\mathcal F)\to\rho_*\pi_i\mathcal F\) on the finite \'etale site.
The sheaf \(\pi_i(\rho_*\mathcal F)\) is finite locally constant by
induction along the Postnikov fibre sequences in
\eqref{eq:finite-site-postnikov-sequences}, using the computation of
\eqref{eq:finite-site-em-stalk} and the stability of finite continuous
modules under extensions, subobjects and quotients. The unit on finite locally constant sheaves is an
isomorphism by the Galois category equivalence of
\cref{sub:finite-site-basics}\textup{(i)}. Thus the canonical map is
an isomorphism, proving \textup{(ii)}.
\end{proof}

\subsection{Finite locally constant spaces}

\begin{definition}\label{def:pi-finite-space}
Let \(d\in\mathbf N\). A space is \emph{\(d\)-truncated} if its
homotopy groups in degrees greater than \(d\) vanish at every base
point. A space is \emph{\(\pi\)-finite} if it is truncated, has
finitely many connected components, and has finite homotopy groups.
For a space \(\mathsf A\) and an \(\infty\)-topos \(\mathcal X\),
\(\underline{\mathsf A}\in\mathcal X\) denotes the constant object with value
\(\mathsf A\), that is, the image of \(\mathsf A\) under the inverse image functor
of the geometric morphism \(\mathcal X\to\mathcal S\)
\cite[Proposition~6.3.4.1]{HTT}; for
\(\mathcal X=\mathrm{Sh}(\mathcal C)\) it is the constant sheaf of
\cref{not:finite-etale-topoi}.
\end{definition}

\begin{sublemma}[Compactness among truncated spaces]
\label{sub:compactness-of-finite-truncated-spaces}
A \(\pi\)-finite \(d\)-truncated space \(\mathsf A\), for \(d\geq0\),
has a simplicial model whose \((d+1)\)-skeleton is finite.
For every filtered diagram of \(d\)-truncated spaces
\((Z_i)\), the canonical map
\[
 \varinjlim_i\operatorname{Map}(\mathsf A,Z_i)
       \longrightarrow \operatorname{Map}(\mathsf A,\varinjlim_i Z_i)
\]
is an equivalence. The space \(\operatorname{Aut}(\mathsf A)\) of
self-equivalences of \(\mathsf A\), with composition as group operation,
is a \(\pi\)-finite \(d\)-truncated space. For constant copies
of \(\mathsf A\) the internal group of self-equivalences is
\(\underline{\operatorname{Aut}(\mathsf A)}\).
\end{sublemma}

\begin{proof}
The assertion about the simplicial model follows from
\cite[Proposition~3.2.3]{AnelPiFinite}.

Choose a minimal Kan model of \(\mathsf A\)
\cite[\href{https://kerodon.net/tag/03TF}{Tag~03TF}]{Kerodon}. In a minimal Kan model
there is one vertex for each connected component. Suppose
simplices in degrees less than \(r\) are known to form finite
sets. There are finitely many maps from the boundary
\(\partial\Delta^r\) into the model, because the boundary is a
finite simplicial set in those degrees. For a fixed boundary
\(b\), the simplices extending \(b\) are the vertices of the
fibre of the Kan fibration
\[
 \operatorname{Map}(\Delta^r,\mathsf A)\longrightarrow
 \operatorname{Map}(\partial\Delta^r,\mathsf A)
\]
over \(b\). If this fibre is nonempty, one filler gives a
homotopy from \(b\) to a constant map. Transport along the
homotopy identifies the fibre with the one over that constant
map, namely the space of pointed maps from
\(\Delta^r/\partial\Delta^r\simeq \mathbb S^r\) to \(\mathsf A\).
It has finitely many components, as they form the group
\(\pi_r(\mathsf A,a)\) at the chosen constant value \(a\). In a minimal
Kan model each component of the fibre contains at most one
vertex. This proves finiteness in degree \(r\) and completes
the induction. This is the argument for cardinal bounds in
\cite[\href{https://kerodon.net/tag/03TK}{Tag~03TK}]{Kerodon}, applied degree by degree to
finite sets.

Let \(I\) be the \((d+1)\)-skeleton of the resulting model.
The inclusion \(I\to \mathsf A\) is an equivalence after \(d\)-truncation,
since the remaining simplicial cells have dimension at least
\(d+2\). Thus maps into a \(d\)-truncated space can be computed
from \(I\). Attaching a simplex from its boundary is a finite
homotopy-colimit operation on spaces. It follows by induction
on the finite cells of \(I\) that \(\operatorname{Map}(I,-)\)
requires finite homotopy limits of copies of the target.
Filtered colimits of spaces commute with finite limits and
preserve \(d\)-truncation. This proves the compactness assertion.

The same computation shows that \(\operatorname{Map}(\mathsf A,\mathsf A)\)
is \(d\)-truncated and \(\pi\)-finite; the latter assertion
is also a special case of
\cite[Proposition~3.4.4]{AnelPiFinite}.
To see the finiteness
in the limit steps, a homotopy pullback of \(\pi\)-finite
spaces lies over the product of two \(\pi\)-finite spaces.
Over a pair of points its fibre is the space of paths between
their images in the third space. The fibre is empty or
equivalent to a loop space of the third space, and in the latter
case has finitely many components and finite homotopy groups.
The homotopy exact sequences show the required finiteness of
the total space. The self-equivalences are a union of components
of \(\operatorname{Map}(\mathsf A,\mathsf A)\), giving the assertion about
\(\operatorname{Aut}(\mathsf A)\).

The constant-object functor preserves finite limits and
colimits. In the topos, maps from \(\underline{\mathsf A}\) into
\(\underline{\mathsf A}\) can again be computed from \(\underline I\).
Indeed \(\underline I\to\underline{\mathsf A}\) becomes an equivalence
after \(d\)-truncation, including after pullback to any slice.
The finite-cell computation of the internal mapping object
therefore gives
\(\underline{\operatorname{Map}}(\underline{\mathsf A},\underline{\mathsf A})
  \simeq\underline{\operatorname{Map}(\mathsf A,\mathsf A)}\).
A section of this object is an equivalence exactly when
its value belongs locally to the union of components of
self-equivalences. Taking these components identifies the
internal group of equivalences with
\(\underline{\operatorname{Aut}(\mathsf A)}\). The identification
respects composition because it comes from the evaluation map.
\end{proof}

\begin{corollary}\label{cor:finite-site-local-trivialization}
Under the hypotheses of \cref{prop:bounded-finite-etale-descent},
let \(n\in\Z\) and put
\(Z=\Omega^\infty\Sigma^{-n}\rho_*\mathcal F\in\mathrm{Sh}(S_{\mathrm{fet}})\).
\begin{enumerate}
\item The counit of \cref{prop:bounded-finite-etale-descent}
induces equivalences
\[
 \rho^*Z\simeq\Omega^\infty\Sigma^{-n}\mathcal F,\qquad
 \Gamma(S_{\mathrm{fet}},Z)\simeq
 \Gamma(S_{\mathrm{\acute et}},\Omega^\infty\Sigma^{-n}\mathcal F).
\]
\item The sheaf \(Z\) is truncated, and there are a connected finite
\'etale cover \(P\to S\), a \(\pi\)-finite space \(\mathsf A\), and an
equivalence \(\underline{\mathsf A}\simeq Z|_{P_{\mathrm{fet}}}\) of sheaves
of spaces on \(P_{\mathrm{fet}}\).
\end{enumerate}
\end{corollary}

\begin{proof}
In the model of spectra as loop spectra,
\(\Omega^\infty\) is evaluation at level zero. The inverse
image and direct image commute with it, as they preserve
finite limits. For the hypercompleted inverse image, this follows
from the same calculation before hypercompletion: all the space
levels of the shifted bounded spectrum are truncated, and inverse
image preserves truncation, so hypercompletion does not change them
\cite[Lemma~6.5.2.9]{HTT}. Exactness of the induced spectral functors
(\cref{sub:finite-site-basics}\textup{(ii)}) also gives compatibility
with the shift \(\Sigma^{-n}\). The counit of
\cref{prop:bounded-finite-etale-descent} identifies
\(\rho^*Z\) with \(\Omega^\infty\Sigma^{-n}\mathcal F\).
The adjunction between \(\rho^*\) and \(\rho_*\), evaluated on
the terminal object, identifies
\(\Gamma(S_{\mathrm{fet}},Z)\) with
\(\Gamma(S_{\mathrm{\acute et}},\Omega^\infty\Sigma^{-n}\mathcal F)\); this
identification is induced by inverse image followed by the counit.
This proves \textup{(i)}.

The space sheaf \(Z\) is \(d\)-truncated for a sufficiently
large \(d\geq0\), since \(\rho_*\mathcal F\) is bounded. Its stalk \(\mathsf A\) is \(\pi\)-finite: the nonnegative
homotopy groups after taking \(\Omega^\infty\) are the
corresponding groups of the shifted stalk spectrum, including
the finite group of components. In \eqref{eq:finite-site-stalk},
write \(\mathsf A=\varinjlim_{(P,\bar p)} Z(P)\). The identity of \(\mathsf A\)
represents a map from \(\mathsf A\) to this colimit. By
\cref{sub:compactness-of-finite-truncated-spaces}
it is represented up to homotopy by a map
\(\mathsf A\to Z(P)\) for one pointed connected finite \'etale \(P\).
The stalk of \(Z|_{P_{\mathrm{fet}}}\) at \(\bar p\) is canonically
\(\mathsf A\): in the indexing category of \eqref{eq:finite-site-stalk}, the
pointed connected covers factoring through \((P,\bar p)\) form a
cofinal subcategory. Indeed, under the Galois category equivalence
of \cref{sub:finite-site-basics}\textup{(i)}, if \(P\) corresponds to
an open subgroup \(J\), this subcategory consists of the open
subgroups contained in \(J\); intersecting any open subgroup with
\(J\) gives a refinement, and the refinement categories are filtered.
Equivalently one takes the pointed connected component of a fibre
product with \(P\). Under this stalk identification the adjoint map
\(\underline{\mathsf A}\to Z|_{P_{\mathrm{fet}}}\) induces the identity of
\(\mathsf A\) up to homotopy, by the choice of \(\mathsf A\to Z(P)\). Since \(P\)
is connected, it is an equivalence by
\cref{sub:finite-site-basics}\textup{(i)}. This proves
\textup{(ii)}.
\end{proof}

\begin{sublemma}[Truncated totalizations]\label{sub:truncated-totalizations}
Let \(d\in\mathbf N\), and let \(Z^\bullet:\Delta\to\mathcal S\)
be a cosimplicial space all of whose terms \(Z^j\) are
\(d\)-truncated. For \(m\in\mathbf N\), let
\(\Delta_{\leq m}\subseteq\Delta\) be the full subcategory on the
objects \([0],\ldots,[m]\), and put
\[
 \operatorname{Tot}(Z^\bullet)=\varprojlim_{[j]\in\Delta}Z^j,
 \qquad
 \operatorname{Tot}_m(Z^\bullet)
   =\varprojlim_{[j]\in\Delta_{\leq m}}Z^j .
\]
\begin{enumerate}
\item The restriction map
\(\operatorname{Tot}(Z^\bullet)\to
\operatorname{Tot}_{d+1}(Z^\bullet)\) is an equivalence.
\item For every filtered diagram
\((Z^\bullet_\lambda)_{\lambda\in\Lambda}\) of cosimplicial spaces
with \(d\)-truncated terms, the canonical map
\(\varinjlim_\lambda\operatorname{Tot}(Z^\bullet_\lambda)\to
\operatorname{Tot}(\varinjlim_\lambda Z^\bullet_\lambda)\) is an
equivalence.
\end{enumerate}
\end{sublemma}

\begin{proof}
Assertion \textup{(i)} is
\cite[Lemma~3.4.3]{AnelPiFinite},
with \(n=d+1\).

Use a Reedy fibrant model, so that the canonical map
\(Z^j\to M^jZ\) to the \(j\)-th matching space is a Kan
fibration \cite[Example~A.2.9.5, Notation~A.2.9.7 and
Proposition~A.2.9.19(F)]{HTT}.
The matching space \(M^jZ\) is the limit over the proper
surjections \([j]\twoheadrightarrow[i]\) in the simplex
category, using \(Z^i\). This indexing category
is finite and nonidentity arrows strictly decrease the target
degree, so the derived matching construction is a finite
homotopy limit. In particular \(M^jZ\) and the fibres of
\(Z^j\to M^jZ\) are \(d\)-truncated.

A point of \(\operatorname{Tot}_{j-1}(Z^\bullet)\)
specifies a map \(\partial\Delta^j\to Z^j\) by the coface
compatibilities and a map \(\Delta^j\to M^jZ\) by the
codegeneracy compatibilities. The two maps agree over
\(\partial\Delta^j\) after the matching projection.
The space of extensions to level \(j\) is the space of
fillers \(\Delta^j\to Z^j\) over the matching space, with
this prescribed boundary. Thus there is a homotopy pullback square
\begin{equation}\label{eq:partial-totalization-square}
\begin{tikzcd}
 \operatorname{Tot}_j(Z^\bullet)\arrow[r]\arrow[d] &
 \operatorname{Map}(\Delta^j,Z^j)\arrow[d] \\
 \operatorname{Tot}_{j-1}(Z^\bullet)\arrow[r] &
 \operatorname{Map}(\partial\Delta^j,Z^j)
   \mathbin{\times_{\operatorname{Map}(\partial\Delta^j,M^jZ)}}
 \operatorname{Map}(\Delta^j,M^jZ).
\end{tikzcd}
\end{equation}
For \(j>d+1\), pull back the matching fibration to
the contractible simplex \(\Delta^j\). Up to equivalence
over \(\Delta^j\) it is a product with its fibre \(F\).
The supplied boundary section gives a map
\(\mathbb S^{j-1}\to F\). This map extends, since
\(j-1>d\) and \(F\) is \(d\)-truncated. Its space of
extensions is equivalent, after choosing an extension,
to \(\Omega^j F\) at the chosen fibre point, and is
contractible. The projection
\(\operatorname{Tot}_j(Z^\bullet)\to
\operatorname{Tot}_{j-1}(Z^\bullet)\) is therefore
an equivalence for \(j>d+1\). Since
\(\operatorname{Tot}(Z^\bullet)=
\varprojlim_j\operatorname{Tot}_j(Z^\bullet)\), this proves
\textup{(i)}.

For \(0\leq j\leq d+1\), the matching space \(M^jZ\) and the
square \eqref{eq:partial-totalization-square} involve only finite
homotopy limits, since the simplices and their boundaries are
finite, and \(\operatorname{Tot}_{-1}(Z^\bullet)\) is a point.
Filtered colimits of spaces commute with finite limits
\cite[Proposition~5.3.3.3]{HTT}, so
\(\operatorname{Tot}_{d+1}\) commutes with filtered colimits.
A filtered colimit of \(d\)-truncated spaces is \(d\)-truncated,
so \textup{(i)} applies to \(\varinjlim_\lambda Z^\bullet_\lambda\)
and gives \textup{(ii)}.
\end{proof}

\begin{sublemma}[Constant spaces on a finite \'etale topos]
\label{sub:finite-site-sections-are-finite-stages}
Let \(\Pi\) be a profinite group and let \(Y\) be a truncated
space with constant object
\(\underline Y\in\mathrm{Sh}(\Pi\text{-}\mathrm{FSet})\)
(\cref{not:finite-etale-topoi,def:pi-finite-space}).
There is a canonical equivalence
\begin{equation}\label{eq:finite-site-constant-sections}
 \varinjlim_{V\triangleleft_o\Pi}\operatorname{Map}(B(\Pi/V),Y)
       \xrightarrow{\sim}
 \Gamma(\Pi\text{-}\mathrm{FSet},\underline Y).
\end{equation}
The index is the filtered poset of open normal subgroups
\(V\triangleleft_o\Pi\), with a morphism \(V\to V'\) when
\(V'\subseteq V\). The space \(B(\Pi/V)\) is the classifying space
of the finite group \(\Pi/V\), and the transition map is inflation
along \(\Pi/V'\to\Pi/V\). This equivalence is natural in \(Y\) and
for homomorphisms of profinite groups acting by restriction.
\end{sublemma}

\begin{proof}
\smallskip\noindent
\emph{A presheaf from finite quotients.}
For a finite continuous \(\Pi\)-set \(I\), consider the
presheaf
\begin{equation}\label{eq:finite-site-constant-presheaf}
 I\longmapsto
 \varinjlim_{\substack{V\triangleleft_o\Pi\\ V\ \mathrm{fixes}\ I}}
   \operatorname{Map}\bigl(I\mathbin{/\mkern-6mu/}(\Pi/V),Y\bigr).
\end{equation}
For a finite group \(G\) acting on \(I\), the notation
\(I\mathbin{/\mkern-6mu/}G\) denotes the classifying space of
the action groupoid. The morphisms in that groupoid from
\(i\) to \(i'\) are the \(g\in G\) with \(gi=i'\).
For a morphism of finite \(\Pi\)-sets choose a subgroup fixing
both sets. Restricting to such subgroups is cofinal and gives
the contravariant maps in \eqref{eq:finite-site-constant-presheaf}. A smaller subgroup
gives a map of the action groupoids to the one at the
original quotient, defining the transition maps of \eqref{eq:finite-site-constant-presheaf}.

\smallskip\noindent
\emph{Descent.}
The presheaf takes finite coproducts to products, since
action groupoids of invariant disjoint unions are disjoint
unions, and filtered colimits commute with finite products.
It suffices to check descent for one finite surjection
\(I'\to I\) of finite \(\Pi\)-sets. Indeed a jointly
surjective family of finite sets has a finite subfamily
which is jointly surjective, and its finite disjoint
union gives such a surjection. Choose an open normal
subgroup fixing \(I'\) and \(I\). At every refinement
\(V\) the augmented simplicial diagram of finite sets
\[
 (I')^{\times_I(\bullet+1)}\longrightarrow I
\]
has geometric realization \(I\): over \(i\in I\) it is
the nerve of the groupoid with object set the nonempty
fibre of \(I'\to I\) at \(i\) and exactly one arrow for
each ordered pair of objects. This groupoid is contractible.
This is an equivalence of spaces with \(\Pi/V\)-action,
as equivalences there are tested on the underlying spaces.
Taking homotopy quotients by \(\Pi/V\) commutes with
geometric realization, since homotopy quotient is the
colimit of the action diagram. Mapping to \(Y\) gives the
descent equivalence at that refinement. Its cosimplicial
mapping spaces are truncated by the same bound as \(Y\).
By \cref{sub:truncated-totalizations}, their totalization
commutes with the filtered colimit over the refinements,
which proves descent for \eqref{eq:finite-site-constant-presheaf}.

\smallskip\noindent
\emph{Comparison with the constant sheaf.}
The constant maps on action groupoids induce a map from
\(\underline Y\) to the resulting sheaf. It is an equivalence
on the stalk. To check this, a stalk representative at the
pointed transitive set \((\Pi/J,J)\) is computed at some
\(V\triangleleft_o\Pi\) with \(V\subseteq J\). Its
action groupoid at that stage is equivalent, using \(J\),
to \(B(J/V)\). Upon restriction to the pointed orbit
\(\Pi/V\), the same quotient stage has the free transitive
action groupoid, whose classifying space is contractible.
Thus the restriction acts by evaluation at the base
point of the groupoid, followed by inclusion of constants.
For a map from any finite sphere into the stalk, it suffices
to compute at one such finite stage, since finite spheres
are compact. If a base point is prescribed, this statement
follows by taking the finite homotopy fibre of evaluation at
the base point of the sphere and commuting that fibre with
the filtered colimit. The same restriction therefore shows
that evaluation at the point and inclusion of constants
are inverse on homotopy groups and on components. Since the
sheaves are truncated, stalk detection makes the map from
\(\underline Y\) an equivalence. Evaluating \eqref{eq:finite-site-constant-presheaf} on
the singleton proves \eqref{eq:finite-site-constant-sections}.

\smallskip\noindent
\emph{Naturality.}
For a homomorphism \(u:\Pi'\to\Pi\) of profinite groups,
the stage \(V\) on \(\Pi\) pulls back to the stage
\(u^{-1}(V)\) on \(\Pi'\). The homomorphism
\(\Pi'/u^{-1}(V)\to\Pi/V\) induces a map from the action
groupoid of the restricted finite set to its action groupoid
for \(\Pi/V\). Precomposition gives the restriction of
mapping spaces in \eqref{eq:finite-site-constant-presheaf}.
These maps commute with refinement and with the inclusions
of constants, proving naturality of the section comparison.
\end{proof}

\subsection{Finite quotients and two sections}
\label{subsec:finite-etale-marked-quotients}

\begin{definition}[Finite-group torsors and marked lifts]
\label{def:marked-lift}
Let \(H\) be a finite group. We use the same letter for the constant
group scheme over the base, whose underlying scheme is a disjoint
union of \(|H|\) copies of the base. A right \(H\)-torsor over a scheme
\(T\) is a scheme \(Q\to T\) with a right \(H\)-action that is locally
isomorphic in the fppf topology to \(H\times T\) with its right
translation action. In particular, \(Q\to T\) is finite \'etale of
degree \(|H|\). If \(\phi:H\to G\) is a homomorphism of finite
groups, the contracted
product \(Q\times^H G\) is the fppf quotient of \(Q\times_T G_T\) by
\[
 (q,g)\cdot h=(q\cdot h,\phi(h)^{-1}g).
\]
Right multiplication on \(g\) makes the quotient a right
\(G\)-torsor. A \emph{marked lift along \(\phi\)} of a right \(G\)-torsor \(P/T\)
is a right \(H\)-torsor \(Q/T\) together with an isomorphism
\(Q\times^H G\to P\) of \(G\)-torsors.

\end{definition}

\begin{construction}[Moduli of one and two sections]
\label{con:marked-moduli}
Let \(d\in\mathbf N\) and let \(\mathsf A\) be a \(\pi\)-finite
\(d\)-truncated space as in \cref{def:pi-finite-space}.
For a group space \(D\)
acting on a space \(V\), the homotopy quotient
\(V\mathbin{/\mkern-6mu/}D\) is the geometric realization of the
action groupoid, whose space of \(n\)-simplices is \(V\times D^n\).
For \(r\in\{1,2\}\), define
\[
 \mathcal M_r(\mathsf A)=\mathsf A^r\mathbin{/\mkern-6mu/}\operatorname{Aut}(\mathsf A),
\]
using the diagonal action induced by evaluation.
The diagonal \(\mathsf A\to \mathsf A\times \mathsf A\) induces
\(\mathcal M_1(\mathsf A)\to\mathcal M_2(\mathsf A)\). Both spaces are \(\pi\)-finite
\((d+1)\)-truncated spaces.
\end{construction}

\Needspace{8\baselineskip}
\begin{proof}
By \cref{sub:compactness-of-finite-truncated-spaces},
\(\operatorname{Aut}(\mathsf A)\) is \(\pi\)-finite and \(d\)-truncated.
Its classifying space is connected, with
\(\pi_1(B\operatorname{Aut}(\mathsf A))=\pi_0\operatorname{Aut}(\mathsf A)\) and
\(\pi_{i+1}(B\operatorname{Aut}(\mathsf A))=\pi_i\operatorname{Aut}(\mathsf A)\)
for \(i\geq1\). The fibre sequences
\[
 \begin{gathered}
 \mathsf A\longrightarrow\mathcal M_1(\mathsf A)\longrightarrow B\operatorname{Aut}(\mathsf A),\\
 \mathsf A\times \mathsf A\longrightarrow\mathcal M_2(\mathsf A)
                         \longrightarrow B\operatorname{Aut}(\mathsf A)
 \end{gathered}
\]
therefore give the asserted finiteness and truncation.
\end{proof}

\begin{definition}[Parameters and finite stages]\label{def:finite-stage}
Let \(\mathsf A\) be a \(\pi\)-finite space, let
\(\mathcal M_1(\mathsf A)\to\mathcal M_2(\mathsf A)\) be as in
\cref{con:marked-moduli}, and let \(S\) be a scheme.
A \emph{parameter} is a section
\(c\in\Gamma(S_{\mathrm{fet}},\underline{\mathcal M_2(\mathsf A)})\).
A \emph{homotopy lift} of \(c\) along
\(\mathcal M_1(\mathsf A)\to\mathcal M_2(\mathsf A)\) is a point of
\[
 \Gamma(S_{\mathrm{fet}},\underline{\mathcal M_1(\mathsf A)})
 \mathbin{\times^h_{\Gamma(S_{\mathrm{fet}},\underline{\mathcal M_2(\mathsf A)})}}
 \{c\}.
\]
Thus a lift includes the homotopy from its image to \(c\).
If \(S\) is connected and locally noetherian with a geometric point
\(\bar s\), a \emph{finite \(G\)-stage} of \(c\) consists of a finite quotient
\(\pi_1^{\mathrm{\acute et}}(S,\bar s)\twoheadrightarrow G\), a map
\(\gamma:BG\to\mathcal M_2(\mathsf A)\), and a path from its image under
\eqref{eq:finite-site-constant-sections} to \(c\).
The same term \emph{homotopy lift} is used for any map to
\(\mathcal M_2(\mathsf A)\), with the analogous homotopy-fibre definition.
\end{definition}

\begin{sublemma}[Forms and homotopy quotients]
\label{sub:forms-and-homotopy-quotients}
Let \(S\) be a scheme, and let \(\mathsf A\) be a \(\pi\)-finite space.
\begin{enumerate}
\item If a space sheaf \(Z\) on \(S_{\mathrm{fet}}\) becomes
      \(\underline{\mathsf A}\) on a finite \'etale cover, then the internal
      space \(P=\operatorname{Isom}(\underline{\mathsf A},Z)\) is a right
      \(\underline{\operatorname{Aut}(\mathsf A)}\)-torsor, and evaluation
      induces an equivalence
      \[
       P\times^{\underline{\operatorname{Aut}(\mathsf A)}}\underline{\mathsf A}
            \xrightarrow{\ \sim\ } Z.
      \]
\item For a group space \(D\) acting on a space \(V\), sections of
      \(\underline{V\mathbin{/\mkern-6mu/}D}\) classify right
      \(\underline D\)-torsors together with a section of their
      associated \(\underline V\)-bundle. This is an equivalence
      of spaces, natural under inverse image of topoi.
\end{enumerate}
Here a right torsor in the topos is a sheaf of spaces \(P\) with
an effective epimorphism \(P\to1\) and an action for which
\(P\times\underline D\to P\times P\),
\((p,g)\mapsto(p,pg)\), is an equivalence.
\end{sublemma}

\begin{proof}
The description of \(Z\) by an automorphism torsor in
\textup{(i)} is a case of
\cite[\arxivhtml{1207.0248v3}{S4.Thmtheorem10}{Proposition~4.10}]{NikolausSchreiberStevenson}.
The classification in \textup{(ii)} follows from
\cite[\arxivhtml{1207.0248v3}{S3.Thmtheorem17}{Theorem~3.17}
and \arxivhtml{1207.0248v3}{S4.Thmtheorem17}{Proposition~4.17}]{NikolausSchreiberStevenson}.

Let \(P\) be the internal space of equivalences from
\(\underline{\mathsf A}\) to \(Z\).
By \cref{sub:compactness-of-finite-truncated-spaces}, the
internal group of automorphisms of \(\underline{\mathsf A}\) is the
constant group space \(\underline{\operatorname{Aut}(\mathsf A)}\).
Act on \(P\) on the right by composition. The internal map
between two objects restricts to the internal map between
their restrictions in the slice over a covering object:
for a test object in the slice, the exponential adjunction
uses its product with the first object over that slice.
Equivalences of the maps are tested in the same slice.
The trivializing cover for \(Z\) thus makes \(P\) a trivial
torsor for \(\underline{\operatorname{Aut}(\mathsf A)}\).
In particular \(P\) covers the terminal object and the map
\[
 P\times\underline{\operatorname{Aut}(\mathsf A)}
    \longrightarrow P\times P,\qquad (p,g)\longmapsto(p,p\circ g)
\]
is an equivalence. These assertions may be tested on the
trivializing cover, where they are assertions about the
regular right action of the group on itself. They are the
torsor conditions in a topos. Evaluation
\(P\times\underline{\mathsf A}\to Z\) is invariant under the diagonal
right action \((p,a)g=(p\circ g,g^{-1}a)\). It descends to an
equivalence from the associated space bundle
\(P\times^{\underline{\operatorname{Aut}(\mathsf A)}}\underline{\mathsf A}\) to \(Z\). This can be checked under the same trivialization,
where evaluation with the identity frame is an inverse
after taking the quotient.

For the second assertion, effectiveness of action groupoids and
stability of their colimits under pullback follow from
\cite[Theorem~6.1.0.6(3)(ii), (iv)]{HTT}.
For a group space \(D\) acting on a space \(V\),
a map to \(V\mathbin{/\mkern-6mu/}D\) classifies a right
\(D\)-torsor together with a section of its associated
\(V\)-bundle. The projection \(V\mathbin{/\mkern-6mu/}D\to BD\) is the
associated \(V\)-bundle of the universal \(D\)-torsor.
For a map to the total space, pull back the torsor along
the composite to \(BD\). The map to the total space then
gives a section of its associated bundle. In the other direction, the
section defines a map from the torsor to \(V\), by expressing
the section in each frame. Its value at \(pg\) is \(g^{-1}\)
applied to its value at \(p\). The functorial
induced map on homotopy quotients is the classifying map.
These constructions are inverse, as pulling back to a
torsor trivialization recovers the frame transitions and
the original points of \(V\). This description applies to
constant \(D,V\) in the topos, since the constant-object
functor commutes with colimits and finite products, hence
with the bar construction of their homotopy quotient.

The constructions use pullbacks and action-groupoid colimits, so
they commute with inverse image. Their inverse identifications
also hold for families over every test object; consequently they
give an equivalence of spaces, including all homotopies.
\end{proof}

\begin{proposition}[Finite group stages of a family with two sections]
\label{prop:finite-group-stages}
Let \(S\) be a scheme, and let \(Z\) be a space sheaf on
\(S_{\mathrm{fet}}\) which becomes constant with value a \(\pi\)-finite
\(d\)-truncated space \(\mathsf A\) on a finite \'etale surjection. Let
\(z_0,z_1\) be sections of \(Z\).
Then \((Z,z_0,z_1)\) has a parameter \(c\) in the sense of
\cref{def:finite-stage}. The space of paths from \(z_0\) to \(z_1\)
is equivalent to the space of homotopy lifts of \(c\) along
\(\mathcal M_1(\mathsf A)\to\mathcal M_2(\mathsf A)\), naturally under inverse image on
finite \'etale topoi.
If moreover \(S\) is connected and locally noetherian with a geometric
point, then \(c\) admits a finite \(G\)-stage.
\end{proposition}

\begin{proof}
Apply \cref{sub:forms-and-homotopy-quotients} with
\(P=\operatorname{Isom}(\underline{\mathsf A},Z)\) and \(V=\mathsf A\times \mathsf A\).
In a frame \(p\), the two sections give
\((p^{-1}z_0,p^{-1}z_1)\). The resulting section of
\(\underline{\mathcal M_2(\mathsf A)}\) is the parameter \(c\), with its
identification with the original family and both sections.

Apply \cref{sub:forms-and-homotopy-quotients} with
\(V=\operatorname{Map}(\Delta^1,\mathsf A)\), with its map to
\(\mathsf A\times\mathsf A\) given by evaluation at the two endpoints.
The associated bundle parametrizes paths in \(Z\). The homotopy
quotient of \(V\) by \(\operatorname{Aut}(\mathsf A)\) is equivalent to
\(\mathcal M_1(\mathsf A)\):
the equivariant map from \(\mathsf A\) to its path space sending a
point to the constant path is an equivalence, by evaluation
at the initial point, and its composite with the endpoints
is the diagonal. Thus the equivalence commutes with both
maps to \(\mathcal M_2(\mathsf A)\). A lift in the homotopy fibre retains the
identification of its family and endpoints with the given
parameter. This proves the path assertion and its naturality.

Suppose now that \(S\) is connected and locally noetherian with a
geometric point. Formula \eqref{eq:finite-site-constant-sections}
applied to \(\mathcal M_2(\mathsf A)\), which is truncated by
\cref{con:marked-moduli}, gives the required finite quotient \(G\)
and representative \(\gamma\). We use the connected right
\(G\)-torsor associated to this quotient of the fundamental
group. The map from the \(G\)-stage in
\eqref{eq:finite-site-constant-sections} sends \(\gamma\)
to the composite of the classifying section of that torsor
with the constant map induced by \(\gamma\). Indeed, for the finite group \(G\), the final quotient in
\eqref{eq:finite-site-constant-sections} is \(G/1\). The
\v{C}ech nerve of the cover \(G\to *\) of finite left \(G\)-sets
has orbit sets \(G^r\) in degree \(r\); descent for a constant space
identifies its sections with maps from the bar construction \(BG\).
The classifying section of the torsor \(G\) with its commuting right
action corresponds to the identity of \(BG\) under this identification.
Its pullback along the finite quotient of the fundamental group
is the chosen right \(G\)-torsor on \(S\).
This describes the map on the original site with the
pointed Galois convention.

\end{proof}

\begin{lemma}[Finite representation of a specified homotopy lift]
\label{lem:finite-stage-homotopy-lift}
Let \(\mathsf A\) be a \(\pi\)-finite space and let
\(\mathcal M_1(\mathsf A)\to\mathcal M_2(\mathsf A)\) be as in
\cref{con:marked-moduli}. Let \(G\) be a finite group and let
\(\gamma:BG\to\mathcal M_2(\mathsf A)\). Let \(\lambda:\Pi'\to G\) be a
continuous homomorphism from a profinite group. Every homotopy lift
of the restricted parameter on \(\mathrm{Sh}(\Pi'\text{-}\mathrm{FSet})\) is represented,
up to homotopy, on \(BH\) for a finite quotient \(\Pi'\twoheadrightarrow H\)
through which \(\lambda\) factors. The representative includes a
map \(\delta:BH\to\mathcal M_1(\mathsf A)\) and a homotopy from its image to
\(\gamma\circ B\lambda_H\), where \(\lambda_H:H\to G\) is induced
by \(\lambda\).
\end{lemma}

\begin{proof}
Take the cofinal system of finite
quotients \(H=\Pi'/V\) with \(V\subseteq\ker(\lambda)\) and
write \(\lambda_H:H\to G\) for the induced map. At each stage
form the space
\[
 \operatorname{Map}(BH,\mathcal M_1(\mathsf A))
 \mathbin{\times^h_{\operatorname{Map}(BH,\mathcal M_2(\mathsf A))}}
 \{\gamma\circ B\lambda_H\}.
\]
The superscript \(h\) indicates the homotopy pullback.
The colimit of these spaces identifies with the space of
homotopy lifts after restriction to \(\mathrm{Sh}(\Pi'\text{-}\mathrm{FSet})\),
by \eqref{eq:finite-site-constant-sections} and commutation
of filtered colimits with finite homotopy limits. On
components a point of this colimit is represented at a
finite stage. Thus \(H\) and \(\delta\) can be chosen with a
specified homotopy making the diagram
\[
\begin{tikzcd}
 BH\arrow[r,"\delta",dashed]\arrow[d,"B\lambda_H"'] &
 \mathcal M_1(\mathsf A)\arrow[d] \\
 BG\arrow[r,"\gamma"'] & \mathcal M_2(\mathsf A)
\end{tikzcd}
\]
commute, representing the specified homotopy lift up to homotopy.
\end{proof}

\begin{corollary}[Transport of a path by a lifting torsor]
\label{cor:transport-marked-path-by-torsor}
Keep the space \(\mathsf A\) and a finite \(G\)-stage
\(\gamma:BG\to\mathcal M_2(\mathsf A)\) of
\cref{prop:finite-group-stages} and let \(P_0\to S\)
be its right \(G\)-torsor. Let \(T\) be an \(S\)-scheme,
\(H\) a finite group, and \(\phi:H\to G\) a homomorphism. Suppose there is a finite \'etale right
\(H\)-torsor \(Q\to T\) and an isomorphism
\[
                Q\times^H G \ \simeq\  T\times_S P_0
\]
of right \(G\)-torsors. The contracted product uses left
multiplication by \(\phi(h)\) for the left \(H\)-action
on \(G\). If \(\gamma\circ B\phi\) has a specified
homotopy lift to \(\mathcal M_1(\mathsf A)\), the two sections pulled back to
\(T_{\mathrm{fet}}\) admit a path. Their images on
\(T_{\mathrm{\acute et}}\) under inverse image also admit a path.
\end{corollary}

\begin{proof}
The classifying section of the \(G\)-torsor
\(T\times_SP_0\simeq Q\times^HG\) over
\(T\) is the composite of the classifying section of \(Q\)
with the map induced by \(B\phi\). The assumed
isomorphism identifies this composite with the pullback
of the original \(G\)-torsor parameter. Hence the
composite of the classifying section of \(Q\) with the
constant map induced by the lift \(BH\to \mathcal M_1(\mathsf A)\) has image
homotopic to the parameter for the pulled-back family with
its two sections. The trivializing finite \'etale surjection for the
family pulls back to one over \(T\). The path assertion of
\cref{prop:finite-group-stages} therefore gives the
path on \(T_{\mathrm{fet}}\); no connectedness or noetherianity
assumption on \(T\) is needed.

Base change of finite \'etale schemes is compatible with
their inclusion in the small \'etale site. Thus the two
inverse image composites from \(S_{\mathrm{fet}}\) to
\(T_{\mathrm{\acute et}}\), through \(S_{\mathrm{\acute et}}\) and through
\(T_{\mathrm{fet}}\), agree: they agree on every representable
finite \'etale object and preserve colimits. Inverse image
takes the constructed path to a path with the asserted
endpoints on the small \'etale site.
\end{proof}

\subsection{Base change and finite Galois data}
\label{subsec:postnikov-generic-base-change}

\begin{proposition}[Finite Galois data for two sections]
\label{prop:finite-galois-data-two-sections}
Let \(S=\operatorname{Spec}B\) be connected, where \(B\) is a
nonzero noetherian semilocal ring whose closed residue fields have at most
one positive characteristic. Let \(\ell\) be a prime invertible in
\(B\), and let \(\mathcal F\) be a bounded hypercomplete spectral sheaf on
\(S_{\mathrm{\acute et}}\) with finite locally constant
\(\ell\)-primary homotopy sheaves. Given two points
\(z_0,z_1\in\Omega^\infty R\Gamma_{\mathrm{\acute et}}(S,\mathcal F)\), there
are a \(\pi\)-finite fibre \(\mathsf A\), a connected right \(G\)-torsor
\(P_0\to S\) for a finite group \(G\), and a finite stage
\(\gamma:BG\to\mathcal M_2(\mathsf A)\) representing the family
\((\Omega^\infty\rho_*\mathcal F,z_0,z_1)\).
For its parameter \(c\), the space of paths is identified by
\[
 \begin{aligned}
 &\{z_0\}\mathbin{\times^h_{\Omega^\infty
              R\Gamma_{\mathrm{\acute et}}(S,\mathcal F)}}\{z_1\}\\
 &\qquad\simeq
 \Gamma(S_{\mathrm{fet}},\underline{\mathcal M_1(\mathsf A)})
 \mathbin{\times^h_{\Gamma(S_{\mathrm{fet}},\underline{\mathcal M_2(\mathsf A)})}}
 \{c\}.
 \end{aligned}
\]
For every field \(L\) over \(S\), the same identification holds for
the pulled-back sections. Every specified path over \(L\) is
represented by a homotopy lift on a finite quotient of the absolute
Galois group of \(L\) refining its homomorphism to \(G\).
\end{proposition}

\begin{proof}
By \cref{prop:finite-etale-effacement}, positive-degree
cohomology with the required coefficients is effaced on finite
\'etale covers of \(S\). By
\cref{prop:bounded-finite-etale-descent}, the counit
\(\rho^*\rho_*\mathcal F\to \mathcal F\) is an equivalence and \(\rho_*\mathcal F\) has the
corresponding finite locally constant homotopy sheaves.
The sections adjunction gives
\begin{equation}\label{eq:finite-site-sections-adjunction}
 \Omega^\infty R\Gamma(S_{\mathrm{fet}},\rho_*\mathcal F)
       \simeq\Omega^\infty R\Gamma_{\mathrm{\acute et}}(S,\mathcal F).
\end{equation}
By \cref{cor:finite-site-local-trivialization},
\(\Omega^\infty\rho_*\mathcal F\) becomes a constant \(\pi\)-finite space
on one connected finite \'etale cover. Apply
\cref{prop:finite-group-stages} to this family and
the two sections under the equivalence
\eqref{eq:finite-site-sections-adjunction}. It supplies
\(\mathsf A\), \(P_0\), \(\gamma\), and the asserted path identification.

For \(\operatorname{Spec}L\to S\), base change of finite \'etale
schemes commutes with their inclusion in the small \'etale site.
The two inverse-image composites from \(S_{\mathrm{fet}}\) to
\((\operatorname{Spec}L)_{\mathrm{\acute et}}\) agree, by the
representable-object argument in the proof of
\cref{cor:transport-marked-path-by-torsor}. The small and finite
\'etale topoi of \(L\) are equivalent by
\cref{sub:finite-site-basics}. Thus the path identification remains
an equivalence over \(L\). Choose a point of the pulled-back \(G\)-torsor over a separable
closure of \(L\), giving its classifying homomorphism from the
absolute Galois group of \(L\) to \(G\)
\cite[\stackstag{0BND}]{Stacks}. Naturality in
\cref{sub:finite-site-sections-are-finite-stages} identifies the
parameter with the restriction of \(\gamma\) along this homomorphism.
Now \cref{lem:finite-stage-homotopy-lift} represents the specified
homotopy lift at the claimed finite quotient.
\end{proof}

\section{A smooth projective model for lifts of a finite torsor}
\label{sec:tame-lifting-model}

\Cref{con:lifting-model,prop:marked-lifting-model-field-versality}
construct a family of marked lifts of a finite torsor that realizes
every such lift over a field. When the kernel of the finite-group
extension has invertible order, the base of this family is a
fibrewise dense open subscheme of a smooth projective scheme
with integral geometric fibres
(\cref{prop:projective-marked-lifting-model}); projectivity uses
\cite[\arxivhtml{1905.00872v1}{S1.SSx2.3}{Remark~1.3}]{BerghRydh}
and \cite[\S3.3]{EdidinRydh}.
\Cref{lem:destackification-fixes-scheme-open} shows that the
destackification sequence of \cite[\arxivhtml{1905.00872v1}{theoremx2}{Theorem~B}]{BerghRydh} is an
isomorphism over the open subscheme carrying the family: on a
scheme, the divisorialification of \cite[\arxivhtml{1905.00872v1}{S6.Thmtheorem1}{Theorem~6.6}]{BerghRydh} is
trivial and \cite[\arxivhtml{1409.5713v1}{S8}{\S8}, Algorithm~E]{Bergh} terminates without performing
a blow-up.

\subsection{The marked lifting model}

\begin{construction}[A model for marked lifts]\label{con:lifting-model}
Let \(S=\operatorname{Spec}R\) for a noetherian ring \(R\), let
\(\phi:H\twoheadrightarrow G\) be a surjection of finite groups,
and let \(P\to S\) be a right \(G\)-torsor.

Let \(E_0\) be the free \(\mathcal O_S\)-module with basis
\((e_g)_{g\in H}\), and let
\[
 \rho_H:H\longrightarrow\operatorname{GL}(E_0),\qquad
 \rho_H(a)e_g=e_{ag},
\]
be the regular representation. Let
\(V=\operatorname{Spec}R[(X_{g,g'})_{g,g'\in H}]\) be the affine space
of endomorphisms \(v=(X_{g,g'})\) of \(E_0\), and let
\(V^\circ=\operatorname{GL}(E_0)\subseteq V\) be the open subscheme on
which \(\det(X_{g,g'})\) is invertible. The group \(H\) acts on \(V\)
on the right by \(v\cdot a=v\circ\rho_H(a)\), and \(V^\circ\) is
stable under this action. Put
\begin{equation}\label{eq:lifting-model}
 Q=P\times_SV^\circ,\qquad
 (p,v)\cdot h=(p\cdot\phi(h),v\cdot h),\qquad
 U=Q/H,
\end{equation}
where \(U\) is the quotient fppf sheaf, and define
\begin{equation}\label{eq:lifting-model-marking}
 \theta:Q\times^H G\longrightarrow P_U=U\times_SP,\qquad
 [(p,v),g]\longmapsto\bigl([(p,v)],p\cdot g\bigr).
\end{equation}
We will verify that \(U\) is a scheme and that \((Q,\theta)\) is a
marked lift of \(P_U\) along \(\phi\) in the sense of \cref{def:marked-lift}.
\end{construction}

\begin{sublemma}[The model carries a marked lift]
\label{sub:lifting-model-marked-lift}
In the situation of \cref{con:lifting-model}, the following
assertions hold.
\begin{enumerate}
\item The group \(H\) acts freely on \(V^\circ\) and on \(Q\), and
for every geometric point \(\bar s\) of \(S\) the open subscheme
\(V^\circ_{\bar s}\subseteq V_{\bar s}\) is nonempty and dense.
\item The quotient \(U=Q/H\) is a scheme, and \(Q\to U\) is a right
\(H\)-torsor.
\item The map \(\theta\) of \eqref{eq:lifting-model-marking} is well
defined and is an isomorphism of right \(G\)-torsors over \(U\).
\end{enumerate}
\end{sublemma}

\begin{proof}
\textup{(i)} Let \(T\) be a nonempty \(S\)-scheme, let
\(M\in V^\circ(T)\), and let \(a\in H\) satisfy
\(M\circ\rho_H(a)=M\). Since \(M\) is invertible,
\(\rho_H(a)=1\), so \(e_a=\rho_H(a)e_1=e_1\) in the free
\(\mathcal O(T)\)-module \(E_0\otimes\mathcal O(T)\), and \(a=1\).
A section of the constant group scheme \(H_S\) over an \(S\)-scheme
is locally constant, so the action of \(H_S\) on \(V^\circ\) is free.
If \((p,v)\cdot h=(p,v)\), then \(v\cdot h=v\); hence the action on
\(Q\) is free. The scheme \(V_{\bar s}\) is an affine space, and
\(V^\circ_{\bar s}\) is the complement of the zero locus of the
nonzero polynomial \(\det(X_{g,g'})\).

\textup{(ii)} The scheme \(Q\) is affine, the group scheme \(H_S\)
is finite locally free over \(S\), and the action is free by
\textup{(i)}. Hence \(U\) is a scheme and \(Q\to U\) is an fppf
\(H\)-torsor \cite[\stackstag{07S7}]{Stacks}.

\textup{(iii)} In \(Q\times^HG\) we have
\([(p,v)\cdot h,\phi(h)^{-1}g]=[(p,v),g]\), and
\[
 \bigl([(p,v)\cdot h],\,p\cdot\phi(h)\cdot\phi(h)^{-1}g\bigr)
 =\bigl([(p,v)],\,p\cdot g\bigr).
\]
Thus \(\theta\) is well defined. It commutes with right
multiplication by \(G\) on \(g\). A morphism of right \(G\)-torsors
is an isomorphism: fppf locally both torsors are trivial, and a
\(G\)-equivariant map \(G_T\to G_T\) is a left translation.
\end{proof}

By \cref{sub:lifting-model-marked-lift}, the quotient in
\cref{con:lifting-model} carries the prescribed marked lift.

\begin{proposition}[Every marked lift over a field occurs in the model]
\label{prop:marked-lifting-model-field-versality}
Keep the data of \cref{con:lifting-model}. For every field \(L\)
over \(S\) and every marked lift \((Q_0,\theta_0)\) of \(P_L\)
along \(\phi\), there is a point \(u\in U(L)\) such that
\(u^*(Q,\theta)\simeq(Q_0,\theta_0)\) as marked lifts.
\end{proposition}

\begin{proof}
We adapt the lifting-variety construction of
\cite[\arfiv{1711.06585}{S2.Thmtheorem1}{Proposition~2.1}]{DemarcheFlorence} to finite groups over \(S\).
The use of the whole invertible-matrix locus allows us to work
also over finite fields. By descent along the finite \'etale cover \(Q_0\to\Spec L\), the
\(H\)-torsor \(Q_0\) and the representation \(\rho_H\) define the
\(L\)-vector space
\[
 E_{Q_0}=Q_0\times^H(E_0\otimes_RL),\qquad
 [q\cdot a,v]=[q,\rho_H(a)v],
\]
of dimension \(|H|\). Choose an \(L\)-linear isomorphism
\(\beta:E_{Q_0}\to E_0\otimes_RL\). For an \(L\)-scheme \(T\) and a
point \(q\in Q_0(T)\), put
\[
 M_q=\beta_T\circ\bigl(v\mapsto[q,v]\bigr)\in
 \operatorname{GL}(E_0)(T)=V^\circ(T).
\]
The relation defining \(E_{Q_0}\) gives
\(M_{q\cdot a}=M_q\circ\rho_H(a)\), and
\(\theta_0([q\cdot a,1])=\theta_0([q,\phi(a)])
=\theta_0([q,1])\cdot\phi(a)\). Hence
\begin{equation}\label{eq:versality-equivariant-map}
 Q_0\longrightarrow P_L\times_LV^\circ_L=Q_L,\qquad
 q\longmapsto\bigl(\theta_0([q,1]),M_q\bigr),
\end{equation}
is \(H\)-equivariant for the action \eqref{eq:lifting-model}. On
quotients by \(H\) it induces a point
\(u:\Spec L=Q_0/H\to Q_L/H=U_L\), and
\eqref{eq:versality-equivariant-map} is a morphism of \(H\)-torsors
\(Q_0\to u^*Q\), hence an isomorphism. Under this isomorphism the
markings \(u^*\theta\) and \(\theta_0\) are two \(G\)-equivariant
maps \(Q_0\times^HG\to P_L\). They agree on the classes \([q,1]\),
which both send to \(\theta_0([q,1])\) by
\eqref{eq:lifting-model-marking} and
\eqref{eq:versality-equivariant-map}. Their pullbacks to the
fppf cover \(Q_0\times G\to Q_0\times^HG\) therefore agree by
\(G\)-equivariance. Hence the isomorphism
\(Q_0\simeq u^*Q\) is an isomorphism of marked lifts.
\end{proof}

\subsection{The projective quotient and its modification}

\begin{sublemma}[Finite generation of graded invariants]
\label{sub:graded-invariant-projectivity}
Let \(R\) be a noetherian ring and let \(B\) be a finitely
generated nonnegatively graded \(R\)-algebra. Let a finite group
\(H\) act on \(B\) by graded \(R\)-algebra automorphisms, and
suppose that \(B_0^H=R\). Then \(B\) is finite over \(B^H\),
the invariant ring \(B^H\) is a finitely generated graded
\(R\)-algebra, and \(\operatorname{Proj}(B^H)\) is projective
over \(R\).
\end{sublemma}

\begin{proof}
Every \(b\in B\) satisfies the monic polynomial
\[
                  \prod_{h\in H}(T-h(b))\in B^H[T].
\]
A finite set of \(R\)-algebra generators of \(B\) also
generates \(B\) as a \(B^H\)-algebra. Since these generators
are integral over \(B^H\), the \(B^H\)-module \(B\) is finite.

We follow the coefficient argument of the Artin--Tate lemma
\cite[\stackstag{00IS}]{Stacks}. Choose \(B^H\)-module
generators \(b_1=1,b_2,\ldots,b_u\) of \(B\), and choose
\(R\)-algebra generators \(x_1,\ldots,x_v\). There exist
coefficients in \(B^H\) such that
\[
 b_i b_j=\sum_{k=1}^u a_{ijk}b_k,
 \qquad
 x_l=\sum_{k=1}^u c_{lk}b_k.
\]
Let \(R'\subseteq B^H\) be the \(R\)-subalgebra generated
by the finitely many coefficients \(a_{ijk}\) and \(c_{lk}\).
The \(R'\)-submodule \(\sum_{k=1}^u R'b_k\) contains \(1\),
is closed under multiplication by the first set of identities,
and contains every \(x_l\) by the second. It is therefore a
subring containing the chosen \(R\)-algebra generators of
\(B\), so it equals \(B\). The ring \(R'\) is noetherian;
hence the \(R'\)-submodule \(B^H\subseteq B\) is finite.
It follows that \(B^H\) is of finite type over \(R\).
Taking homogeneous components of a finite set of algebra
generators gives homogeneous generators of \(B^H\), since
the \(H\)-action preserves degrees.

The degree-zero part of \(B^H\) is \(B_0^H=R\). We follow the
Veronese argument of \cite[\stackstag{0EGH}]{Stacks}.
Let \(f_1,\ldots,f_r\) be positive-degree homogeneous
generators of \(B^H\) over \(R\), with
\(\deg(f_i)=d_i\). If \(r=0\), then \(B^H=R\) is
concentrated in degree zero, and \(\operatorname{Proj}(B^H)\)
is empty. Suppose \(r\geq1\), and choose a positive common
multiple \(e\) of \(d_1,\ldots,d_r\).
For a monomial \(\prod_i f_i^{u_i}\) of degree \(nre\),
where \(n\geq1\), write
\[
 u_i=q_i(e/d_i)+v_i,
 \qquad 0\leq v_i<e/d_i.
\]
The residual monomial \(\prod_i f_i^{v_i}\) has degree
\(se\) for an integer \(0\leq s<r\), because
\[
 \sum_i v_i d_i
       =nre-e\sum_i q_i,
 \qquad 0\leq\sum_i v_i d_i<re.
\]
There are \(\sum_i q_i=nr-s\) pure-power factors
\(f_i^{e/d_i}\), each of degree \(e\). Since
\(nr-s\geq r-s\), multiply the residual monomial by
\(r-s\) of these factors to obtain a monomial of degree
\(re\). The remaining \((n-1)r\) factors form \(n-1\)
groups of \(r\), each of degree \(re\). Thus every monomial
of degree \(nre\) is a product of \(n\) monomials of
degree \(re\).

Consequently the Veronese algebra
\(\bigoplus_{n\geq0}(B^H)_{nre}\), with \((B^H)_{nre}\)
in degree \(n\), is generated in degree one. Its degree-one
part is a finite \(R\)-module: it is spanned by the finitely
many monomials in the \(f_i\) of total weighted degree
\(re\). Its Proj therefore admits a closed immersion into
projective space over \(R\). The canonical isomorphism
\[
 \operatorname{Proj}(B^H)\xrightarrow{\ \sim\ }
 \operatorname{Proj}\!\left(\bigoplus_{n\geq0}(B^H)_{nre}\right)
\]
of \cite[\stackstag{0B5J}]{Stacks} proves projectivity.
\end{proof}

\begin{sublemma}[A coarse quotient in homogeneous coordinates]
\label{sub:projective-graded-coarse-quotient}
Let \(R\) be a noetherian ring, let \(B\) be a nonnegatively graded
\(R\)-algebra with \(B_0\) finite over \(R\), and suppose that \(B\)
is generated over \(B_0\) by finitely many elements of degree one.
Let a finite group \(H\) act on \(B\) by graded \(R\)-algebra
automorphisms, and suppose that \(B_0^H=R\).
\begin{enumerate}
\item Every point of \(\operatorname{Proj}B\) lies in \(D_+(f)\) for
some homogeneous element \(f\in B^H\) of positive degree. Hence the
inclusion \(B^H\subseteq B\) induces a morphism
\[
 \pi:\operatorname{Proj}B\longrightarrow\operatorname{Proj}(B^H),
 \qquad \pi^{-1}\bigl(D_+(f)\bigr)=D_+(f)
\]
for every homogeneous \(f\in B^H\) of positive degree.
\item For every such \(f\), we have \((B_{(f)})^H=(B^H)_{(f)}\), and
\(B_{(f)}\) is a finite \((B^H)_{(f)}\)-algebra. In particular
\(\pi\) is finite.
\item The morphism \(\pi\) is a geometric quotient of
\(\operatorname{Proj}B\) by \(H\), the scheme
\(\operatorname{Proj}(B^H)\) is the coarse space of the quotient
stack \([\operatorname{Proj}B/H]\), and it is projective over \(R\).
\end{enumerate}
\end{sublemma}

\begin{proof}
The algebra \(B\) is finitely generated over \(R\). By
\cref{sub:graded-invariant-projectivity}, the ring \(B^H\) is a
finitely generated graded \(R\)-algebra, \(B\) is a finite
\(B^H\)-module, and \(\operatorname{Proj}(B^H)\) is projective over
\(R\).

\textup{(i)} Let \(\mathfrak p\) be a point of
\(\operatorname{Proj}B\), that is, a homogeneous prime ideal of
\(B\) not containing \(B_+\). For \(h\in H\), the ideal
\(h(\mathfrak p)\) is a homogeneous prime ideal not containing
\(h(B_+)=B_+\). By graded prime avoidance
\cite[\stackstag{00JS}]{Stacks}, there is a homogeneous element
\(b\in B_+\) of positive degree with \(b\notin h(\mathfrak p)\) for
all \(h\in H\). The element
\[
 f=\prod_{h\in H}h(b)\in B^H
\]
is homogeneous of positive degree, and \(f\notin\mathfrak p\),
because \(h(b)\notin\mathfrak p\) for every \(h\) and \(\mathfrak p\)
is prime. Thus \(\mathfrak p\in D_+(f)\). The existence of \(\pi\)
and the formula for \(\pi^{-1}(D_+(f))\) now follow from
\cite[\stackstag{01MY}]{Stacks}; on \(D_+(f)\), the morphism \(\pi\)
corresponds to the ring map \((B^H)_{(f)}\to B_{(f)}\).

\textup{(ii)} Localization is exact, so
\((B^H)_{(f)}\subseteq(B_{(f)})^H\). Let \(b/f^i\) be an invariant
element of \(B_{(f)}\). For every \(h\in H\) there is an exponent
\(a\geq0\) with \(f^a(h(b)-b)=0\); choose \(a\) valid for all \(h\).
Then \(f^ab\in B^H\) and \(b/f^i=(f^ab)/f^{a+i}\in(B^H)_{(f)}\).
Every \(x\in B_{(f)}\) is a root of the monic polynomial
\(\prod_{h\in H}(T-h(x))\) with coefficients in
\((B_{(f)})^H=(B^H)_{(f)}\). Moreover \(B_{(f)}\) is a finitely
generated \(R\)-algebra: it is the coordinate ring of an affine open
subscheme of \(\operatorname{Proj}B\), which is a closed subscheme of
a projective space over \(B_0\). Hence \(B_{(f)}\) is finite over
\((B^H)_{(f)}\). The charts \(D_+(f)\) cover
\(\operatorname{Proj}(B^H)\), and by \textup{(i)} their inverse
images are the affine schemes \(\Spec B_{(f)}\); hence \(\pi\) is
finite.

\textup{(iii)} By \textup{(ii)} and
\cite[Theorem~4.1]{RydhGeometricQuotients}, the morphism
\(\Spec B_{(f)}\to\Spec(B_{(f)})^H=D_+(f)\) is a geometric quotient
by \(H\). The scheme \(D_+(f)\) is the coarse space of
\([\Spec B_{(f)}/H]\) by
\cite[Definition~6.8 and Remark~6.9]{RydhGeometricQuotients}. Both
properties are local on the target: a morphism from
\([\operatorname{Proj}B/H]\) to an algebraic space is an
\(H\)-invariant morphism from \(\operatorname{Proj}B\), and the
morphisms obtained by descent on the charts \(D_+(f)\) agree on the
intersections \(D_+(f)\cap D_+(f')=D_+(ff')\).
\end{proof}
\begin{definition}[Finite tame inertia]\label{def:finite-tame-inertia}
An algebraic stack \(\mathcal W\) has \emph{finite tame inertia}
if its inertia morphism \(I_{\mathcal W}\to\mathcal W\) is finite
and every geometric stabilizer is linearly reductive
(cf.\ \cite[Definition~3.1 and Theorem~3.2]{AOV08}). Here a finite
group scheme over a field is \emph{linearly reductive} if taking
invariants is exact on its finite-dimensional representations.
\end{definition}

\begin{lemma}[Destackification fixes an open subscheme]
\label{lem:destackification-fixes-scheme-open}
Let \(S\) be an affine noetherian scheme, let \(\mathcal W\) be a
smooth finitely presented algebraic stack over \(S\) with finite
tame inertia \textup{(\cref{def:finite-tame-inertia})}, and let
\(U\subseteq\mathcal W\) be an open substack that is a scheme. Let
\begin{equation}\label{eq:destackification-sequence}
 \Pi:(\mathcal Y,E)=(X_n,D_n)\xrightarrow{\ \pi_n\ }\cdots
 \xrightarrow{\ \pi_1\ }(X_0,D_0)=(\mathcal W,\varnothing)
\end{equation}
be the smooth stacky blow-up sequence that
\cite[\arxivhtml{1905.00872v1}{theoremx2}{Theorem~B}]{BerghRydh} associates with the standard pair
\((\mathcal W,\varnothing)\) over \(S\)
\cite[\arxivhtml{1905.00872v1}{S3.6}{Definition~3.5}]{BerghRydh}. For \(0\leq i<n\), the centre of
\(\pi_{i+1}\) is disjoint from the inverse image of \(U\) in
\(X_i\). Consequently \(\Pi\) restricts to an isomorphism
\(\Pi^{-1}(U)\to U\), and \(E\cap\Pi^{-1}(U)=\varnothing\).
\end{lemma}

\begin{proof}
If \(U=\varnothing\), the assertions are immediate. Suppose that
\(U\ne\varnothing\).

\emph{Restriction to \(U\).}
The stack \(\mathcal W\) is smooth and of finite presentation over
\(S\), and the empty family is an ordered simple normal crossings
divisor; hence \((\mathcal W,\varnothing)\) is a standard pair. The
stack \(\mathcal W\) is locally noetherian
\cite[\stackstag{06R6}]{Stacks}, so the open immersion
\(j:U\to\mathcal W\) is quasi-compact
\cite[\stackstag{0CPM}]{Stacks} and thus of finite presentation.
It is smooth, and it is stabilizer
preserving in the sense of \cite[\arxivhtml{1905.00872v1}{S2.SS1}{\S2.1}]{BerghRydh}: the map
\(I_U\to j^*I_{\mathcal W}\) is an isomorphism because \(j\) is a
monomorphism. The base change of \((\mathcal W,\varnothing)\) along
\(j\) is the standard pair \((U,\varnothing)\). By the functoriality
assertion of \cite[\arxivhtml{1905.00872v1}{theoremx2}{Theorem~B}]{BerghRydh}, in the sense of
\cite[\arxivhtml{1905.00872v1}{S3.SS6}{\S3.6}]{BerghRydh}, the base change of \(\Pi\) along \(j\) is
the sequence \(\Pi_U\) that the construction associates with
\((U,\varnothing)\); here two sequences are identified if they differ
by the insertion or deletion of stacky blow-ups with empty centres.
The centres of the base change of \(\Pi\) along \(j\) are the
intersections of the centres of \(\Pi\) with the inverse images of
\(U\). Hence it suffices to prove that all centres of \(\Pi_U\) are
empty.

\smallskip\noindent
\emph{The sequence of a scheme.}
By the proof of \cite[\arxivhtml{1905.00872v1}{S7.Thmtheorem1}{Theorem~7.1}]{BerghRydh} in the case of a base
scheme, \(\Pi_U\) is the divisorialification sequence of
\cite[\arxivhtml{1905.00872v1}{S6.Thmtheorem1}{Theorem~6.6}]{BerghRydh} for \((U,\varnothing)\), followed by
the sequence that \cite[\arxivhtml{1409.5713v1}{S8}{\S8}, Algorithm~E]{Bergh} produces for the
rigidification of its target. The inertia of \(U\) over \(S\) is
trivial, hence flat. By \cite[\arxivhtml{1905.00872v1}{S6.7}{Remark~6.7}]{BerghRydh}, the
divisorial index of \((U,\varnothing)\) is zero at every point of
\(U\). Thus the construction in the proof of
\cite[\arxivhtml{1905.00872v1}{S6.Thmtheorem1}{Theorem~6.6}]{BerghRydh}
has already reached its stopping condition and contributes no
blow-up centres. Its target is \((U,\varnothing)\), which is its
own rigidification.

Let \(\xi\) be a point of \(U\). Its stabilizer is trivial, so every
one-dimensional component of the conormal representation at \(\xi\)
is the trivial representation, and its divisorial part is zero
because no component of the divisor passes through \(\xi\)
\cite[\arxivhtml{1409.5713v1}{S6}{\S6}]{Bergh}. Hence the divisorial index and the toroidal
index at \(\xi\) vanish
\cite[\arxivhtml{1409.5713v1}{S7.Thmtheorem6}{Definition~7.6} and \arxivhtml{1409.5713v1}{S7.Thmtheorem2}{Definition~7.2}]{Bergh}, every
component of the conormal representation is independent, so that
the independency index at \(\xi\) vanishes
\cite[\arxivhtml{1409.5713v1}{S7.Thmtheorem1}{Definition~7.1}]{Bergh}, and the divisorial type at \(\xi\)
\cite[\arxivhtml{1409.5713v1}{S7.Thmtheorem4}{Definition~7.4}]{Bergh} is the one attached to the zero group and
the empty divisor. In particular \((U,\varnothing)\) is a divisorial
stack, as Algorithm~E requires of its input, and the invariant of
step~E1, formed by the independency index, the toroidal index and the
divisorial type, is constant on \(U\). The locus where it attains its
maximum is \(U\), so Algorithm~E terminates in step~E1 without a
blow-up. The proof of \cite[\arxivhtml{1409.5713v1}{S1.Thmtheorem2}{Theorem~1.2}]{Bergh} ends, after
Algorithm~E, with the splitting of components of the divisor by
trivial blow-ups \cite[\arxivhtml{1409.5713v1}{S8}{end of \S8}]{Bergh}; the divisor of
\((U,\varnothing)\) has no components. Thus \(\Pi_U\) has no
nonempty centres, and the transforms of the empty boundary remain
empty on the inverse image of \(U\).
\end{proof}
\begin{sublemma}[Density under stacky blow-ups]
\label{sub:stacky-modifications-preserve-density}
Let \(\mathcal W\) be a noetherian algebraic stack, let
\(U\subseteq\mathcal W\) be an open substack that is dense in the
underlying topological space \(|\mathcal W|\), and let
\(b:\mathcal W'\to\mathcal W\) be one of the following morphisms:
\begin{enumerate}
\item the blow-up of a closed substack \(Z\subseteq\mathcal W\) with
\(Z\cap U=\varnothing\);
\item a root stack along an effective Cartier divisor
\(D\subseteq\mathcal W\) with \(D\cap U=\varnothing\).
\end{enumerate}
Then \(b\) restricts to an isomorphism \(b^{-1}(U)\to U\), and
\(b^{-1}(U)\) is dense in \(|\mathcal W'|\).
\end{sublemma}

\begin{proof}
\textup{(i)} Blowing up commutes with flat base change
\cite[\stackstag{0805}]{Stacks}, so the assertions of
\cite[\stackstag{02OS}]{Stacks} hold for \(b\): it is an isomorphism
over \(\mathcal W\setminus Z\), and \(b^{-1}(Z)\) is an effective
Cartier divisor on \(\mathcal W'\). Let \(V'\to\mathcal W'\) be a
smooth surjective morphism from a scheme. The inverse image of
\(b^{-1}(Z)\) in \(V'\) is an effective Cartier divisor, so its
complement is dense in \(V'\) \cite[\stackstag{07ZU}]{Stacks}. The
map \(|V'|\to|\mathcal W'|\) is open and surjective; hence
\(|\mathcal W'|\setminus|b^{-1}(Z)|\) is dense in \(|\mathcal W'|\).
This open substack is isomorphic to \(\mathcal W\setminus Z\), in
which \(U\) is dense.

\smallskip\noindent
\textup{(ii)} The morphism \(b\) is a coarsening
\cite[\arxivhtml{1905.00872v1}{S3.SS5}{\S3.5}]{BerghRydh}, hence a universal homeomorphism
\cite[\arxivhtml{1905.00872v1}{S2.SS5}{\S2.5}]{BerghRydh}. Let \(r\) be the order of the root stack.
Over \(\mathcal W\setminus D\) the canonical section \(s_D\) of
\(\mathcal O_{\mathcal W}(D)\) is nowhere zero. Hence, for an
\(r\)-th root \((L,s)\) of \((\mathcal O_{\mathcal W}(D),s_D)\) over
a scheme mapping to \(\mathcal W\setminus D\), the section \(s\) is
nowhere zero and identifies \((L,s)\) with \((\mathcal O,1)\) by a
unique isomorphism. Hence \(b\) is an isomorphism over
\(\mathcal W\setminus D\), and density is preserved because
\(|b|\) is a homeomorphism.
\end{proof}

\begin{lemma}[Destackification along a dense open subscheme]
\label{lem:destackification-preserves-scheme-open}
Let \(S\) be an affine noetherian scheme, let \(\mathcal W\) be a
smooth finitely presented algebraic stack over \(S\) with finite
tame inertia, and suppose that its coarse space \(W_0\) is a
projective \(S\)-scheme. Let \(U\subseteq\mathcal W\) be an open
substack that is a scheme, and let \(r\in\mathbf N\). Suppose that,
for every geometric point \(\bar s\) of \(S\), the scheme
\(U_{\bar s}\) is nonempty, integral, of dimension \(r\), and dense
in \(|\mathcal W_{\bar s}|\). Then there are a smooth projective
\(S\)-scheme \(X\) of relative dimension \(r\) and an open immersion
\(U\to X\) of \(S\)-schemes such that, for every geometric point
\(\bar s\) of \(S\), the fibre \(X_{\bar s}\) is integral and
\(U_{\bar s}\) is dense in \(X_{\bar s}\).
\end{lemma}

\begin{proof}
If \(S=\varnothing\), take \(X=\varnothing\).

\smallskip\noindent
\emph{The scheme \(X\).}
Suppose \(S\ne\varnothing\), and let \(\Pi\) be the sequence \eqref{eq:destackification-sequence},
let \(\pi:\mathcal Y\to X\) be the coarse space of \(\mathcal Y\), and
let \(E_{\mathrm{cs}}\subseteq X\) be the image of \(E\). By
\cite[\arxivhtml{1905.00872v1}{theoremx2}{Theorem~B}]{BerghRydh}, \((X,E_{\mathrm{cs}})\) is a standard
pair over \(S\); in particular \(X\) is an algebraic space, smooth
and of finite presentation over \(S\).

\emph{Projectivity.}
We verify the assertion about coarse-space blow-ups in
\cite[\arxivhtml{1905.00872v1}{S1.SSx2.3}{Remark~1.3}]{BerghRydh}
using \cite[\S3.3 and Propositions~3.4, 3.7]{EdidinRydh}.
The diagonal of \(\mathcal W\) over \(S\) is finite: the coarse
space morphism has finite diagonal by
\cite[\arxivhtml{1905.00872v1}{S2.SS5}{\S2.5}]{BerghRydh}, and
\(W_0\) is separated over \(S\). Each blow-up is representable and
separated, and each root stack has finite diagonal by the
coarsening property in
\cite[\arxivhtml{1905.00872v1}{S3.SS5}{\S3.5}]{BerghRydh}.
Finiteness of the diagonal is preserved by composition (factor the
diagonal of the composite through the pullback of the diagonal of
the second morphism). Thus every \(X_i\) in
\eqref{eq:destackification-sequence} is noetherian with finite
(hence affine) diagonal over \(S\), and therefore over \(\Spec\Z\)
as \(S\) is affine; in particular it has finite inertia.
The geometric stabilizers of a blow-up are subgroup schemes of
those of its target. For a root stack of order \(e\), the
homomorphism of geometric stabilizers is an isomorphism away from
the root divisor, and over the divisor it is surjective as an fppf
sheaf map with kernel \(\mu_e\): an automorphism of the original
line with section lifts by taking an \(e\)-th root of its scalar on
the line. Since \(\mu_e\) is diagonalizable, it is linearly
reductive, even in residue characteristics dividing \(e\).
Linearly reductive finite group schemes over a field are closed
under subgroups and extensions \cite[Proposition~2.5(a), (c)]{AOV08}.
Inductively all the stacks \(X_i\) are therefore tame.

For \(0\leq i\leq n\), write \(c_i:X_i\to W_i\) for the coarse
space morphism, so \(W_n=X\). For \(0\leq i<n\), the universal
property of \(c_{i+1}\) gives a commutative square
\[
\begin{tikzcd}
 X_{i+1}\arrow[r,"\pi_{i+1}"]\arrow[d,"c_{i+1}"'] &
 X_i\arrow[d,"c_i"]\\
 W_{i+1}\arrow[r] & W_i .
\end{tikzcd}
\]
If \(\pi_{i+1}\) is a root stack, it is a coarsening and the map
\(W_{i+1}\to W_i\) is an isomorphism by
\cite[\arxivhtml{1905.00872v1}{S2.SS5}{\S2.5}]{BerghRydh}.
Otherwise let \(\mathcal I\subseteq\mathcal O_{X_i}\) be the ideal
of its centre. Since \(c_i\) is a tame coarse space morphism,
\(c_{i,*}\) is exact on quasi-coherent sheaves
\cite[Theorem~3.2]{AOV08} and \(c_{i,*}\mathcal O_{X_i}=\mathcal O_{W_i}\).
In particular it is a good moduli space morphism in the sense of
\cite[Definition~2.1 and Remark~2.4]{EdidinRydh}. The ordinary and saturated Proj
coincide for a tame coarse space morphism by
\cite[\S3.3 and proof of Proposition~3.12]{EdidinRydh}.
By \cite[Propositions~3.4 and 3.7(1)--(2)]{EdidinRydh}, the coarse
space of the blow-up is consequently
\[
 W_{i+1}\simeq
 \operatorname{Proj}_{W_i}\!\left(\bigoplus_{m\geq0}
 c_{i,*}\mathcal I^m\right)
 \simeq \operatorname{Bl}_{c_{i,*}\mathcal I^d}(W_i)
\]
for a sufficiently divisible positive integer \(d\).
Here \(c_{i,*}\mathcal I^d\) is an ideal in
\(\mathcal O_{W_i}\); the second isomorphism uses the Veronese
algebra generated in degree one. Starting from the projective
\(S\)-scheme \(W_0\), each \(W_{i+1}\) is thus a scheme projective
over \(W_i\) (the ideal is of finite type since \(W_i\) is
noetherian, cf.~\cite[\stackstag{02NS}]{Stacks}). By induction and
\cite[\stackstag{0C4P}]{Stacks}, both \(X=W_n\to W_0\) and
\(X\to S\) are projective.

\smallskip\noindent
\emph{The open immersion.}
By \cref{lem:destackification-fixes-scheme-open}, \(\Pi\) restricts
to an isomorphism from \(U'=\Pi^{-1}(U)\) to \(U\). The coarse space
morphism \(\pi\) is a universal homeomorphism, and the formation of
coarse spaces commutes with flat base change
\cite[\arxivhtml{1905.00872v1}{S2.SS2}{\S2.2} and \arxivhtml{1905.00872v1}{S2.SS5}{\S2.5}]{BerghRydh}. Hence \(\pi(U')\) is an open subscheme of \(X\), and its inverse
image under \(\pi\) is precisely \(U'\). By flat base change along
this open immersion, \(U'\to\pi(U')\) is the coarse space of the
scheme \(U'\), that is, an isomorphism. The composite
\(U\simeq U'\simeq\pi(U')\subseteq X\) is the open immersion of the
statement.

\smallskip\noindent
\emph{The geometric fibres.}
Let \(\bar s\) be a geometric point of \(S\). The base change
\(\Pi_{\bar s}\) of \(\Pi\) along \(\bar s\to S\) is a smooth stacky
blow-up sequence starting at \((\mathcal W_{\bar s},\varnothing)\),
whose centres are the base changes of the centres of \(\Pi\)
\cite[\arxivhtml{1905.00872v1}{S3.SS6}{\S3.6}]{BerghRydh}. By
\cref{lem:destackification-fixes-scheme-open}, these centres are
disjoint from the inverse images of \(U_{\bar s}\). Applying
\cref{sub:stacky-modifications-preserve-density} to the successive
steps of \(\Pi_{\bar s}\), starting from the dense open substack
\(U_{\bar s}\subseteq\mathcal W_{\bar s}\), we find that
\(U'_{\bar s}\) is dense in \(|\mathcal Y_{\bar s}|\). Since
\(\pi\) is a universal homeomorphism, its geometric base change
\(\pi_{\bar s}:\mathcal Y_{\bar s}\to X_{\bar s}\) is a
homeomorphism, so \(U_{\bar s}\) is dense in \(X_{\bar s}\).
The argument uses the base-change property of
\eqref{eq:destackification-sequence}, not arbitrary base change for
coarse spaces or for blow-ups of nonflat centres. The scheme
\(X_{\bar s}\) is smooth over the field \(\kappa(\bar s)\), hence
reduced, and it is irreducible because its dense open subscheme
\(U_{\bar s}\) is; thus \(X_{\bar s}\) is integral, of dimension
\(r\). Hence \(X\) has relative dimension \(r\) over \(S\).
\end{proof}
\subsection{A smooth projective compactification}

\begin{construction}[The homogeneous compactification]
\label{con:lifting-compactification}
Keep the data of \cref{con:lifting-model}. The ring
\(\Gamma(P,\mathcal O_P)\) is a finite \'etale \(R\)-algebra. Consider the
graded polynomial algebra
\[
 B=\Gamma(P,\mathcal O_P)\bigl[X_0,(X_{g,g'})_{g,g'\in H}\bigr],
\]
where the coefficients have degree zero and the variables have degree
one. Put
\[
 W=\operatorname{Proj}B
  =P\times_S\mathbb P_S\bigl(\mathcal O_S\oplus
      \operatorname{End}(E_0)\bigr).
\]
The group \(H\) acts on \(B\) by graded \(R\)-algebra automorphisms:
on \(\Gamma(P,\mathcal O_P)\) through \(\phi\) and the action of
\(G\) on \(P\), trivially on \(X_0\), and on the matrix
\((X_{g,g'})\) by right multiplication by \(\rho_H(a)\). The induced
right action of \(H\) on \(W\) restricts to the action
\eqref{eq:lifting-model} on the open subscheme
\begin{equation}\label{eq:lifting-model-chart}
 Q=D_+(f)\subseteq W,\qquad
 f=X_0\det\bigl((X_{g,g'})_{g,g'\in H}\bigr)^2 .
\end{equation}
\end{construction}

\begin{lemma}[Tameness and geometric fibres of the compactification]
\label{lem:lifting-compactification-tame-fibres}
Keep the data of \cref{con:lifting-model,con:lifting-compactification}.
\begin{enumerate}
\item We have \(B_0^H=R\) and \(f\in B^H\). Consequently the
quotient stack \([W/H]\) has a coarse space that is projective over
\(S\), and \(U=[Q/H]\) is an open substack of \([W/H]\) that is a
scheme.
\item Assume that \(|\ker\phi|\) is invertible in \(R\). Then
\([W/H]\) is smooth and finitely presented over \(S\) and has finite
tame inertia (\cref{def:finite-tame-inertia}).
\item For every geometric point \(\bar s\) of \(S\), the scheme
\(U_{\bar s}\) is nonempty, integral, of dimension \(|H|^2\), and
dense in the underlying topological space of \([W/H]_{\bar s}\).
\end{enumerate}
\end{lemma}

\begin{proof}
\textup{(i)} Since \(\phi\) is surjective,
\(\Gamma(P,\mathcal O_P)^H=\Gamma(P,\mathcal O_P)^G\). A function on
\(P\) is \(G\)-invariant if and only if its two pullbacks to
\(P\times_SP\simeq P\times G\) agree, that is, if and only if it
descends along the fppf covering \(P\to S\)
\cite[\stackstag{023Q}]{Stacks}. Hence \(B_0^H=R\). For \(a\in H\),
the matrix \(\rho_H(a)\) is a permutation matrix, so
\[
 \det\bigl((X_{g,g'})\rho_H(a)\bigr)^2
 =\det\bigl((X_{g,g'})\bigr)^2\det\bigl(\rho_H(a)\bigr)^2
 =\det\bigl((X_{g,g'})\bigr)^2,
\]
and \(f\in B^H\). By \cref{sub:projective-graded-coarse-quotient},
the scheme \(\operatorname{Proj}(B^H)\) is the coarse space of
\([W/H]\) and is projective over \(S\). The open subscheme
\(Q=D_+(f)\) of \eqref{eq:lifting-model-chart} is \(H\)-stable, and
\([Q/H]=U\) is a scheme by
\cref{sub:lifting-model-marked-lift}\textup{(ii)}.

\textup{(ii)} The atlas \(W\to[W/H]\) is finite \'etale and
surjective, and \(W\) is smooth and finitely presented over \(S\);
hence so is \([W/H]\). The pullback of the inertia stack to this
atlas is
\[
 I_{[W/H]}\times_{[W/H]}W=\coprod_{h\in H}W^h,
\]
where \(W^h\subseteq W\) is the equalizer of the action of \(h\) and
the identity. Each \(W^h\) is a closed subscheme, because \(W\) is
separated over \(S\). Hence the inertia is finite. Let \((p,x)\) be
a geometric point of \(W\) fixed by \(h\in H\). Then
\(p\cdot\phi(h)=p\), and \(\phi(h)=1\) because \(G\) acts freely on
\(P\). Thus every geometric stabilizer is a subgroup of
\(\ker\phi\). Its order is invertible in the residue field, so
averaging over the group makes the functor of invariants exact.

\textup{(iii)} A point of \(P_{\bar s}\) identifies \(P_{\bar s}\)
with \(G\times\bar s\), and then
\[
 W_{\bar s}=\coprod_{g\in G}
   \mathbb P_{\bar s}\bigl(\mathcal O\oplus\operatorname{End}(E_0)\bigr),
 \qquad
 Q_{\bar s}=\coprod_{g\in G}V^\circ_{\bar s}.
\]
The group \(H\) permutes these components transitively through
\(\phi\), and \(Q_{\bar s}\) is dense in \(W_{\bar s}\) by
\cref{sub:lifting-model-marked-lift}\textup{(i)}. The atlas
\(W_{\bar s}\to[W/H]_{\bar s}\) is surjective, so the image
\(U_{\bar s}\) of \(Q_{\bar s}\) is dense in \([W/H]_{\bar s}\). By
transitivity, the component \(V^\circ_{\bar s}\) indexed by
\(1\in G\) maps onto \(U_{\bar s}\), which is therefore irreducible
and nonempty. The scheme \(U_{\bar s}\) is smooth of dimension
\(|H|^2\) over \(\bar s\), because \(Q_{\bar s}\to U_{\bar s}\) is
finite \'etale and surjective \cite[\stackstag{05B5}]{Stacks}. Hence
it is integral.
\end{proof}

\begin{proposition}[A smooth projective model for lifts of a finite torsor]
\label{prop:projective-marked-lifting-model}
Let \(S=\operatorname{Spec}R\) for a noetherian ring \(R\), let
\(\phi:H\twoheadrightarrow G\) be a surjection of finite groups such
that \(|\ker\phi|\) is invertible in \(R\), and let \(P\to S\)
be a right \(G\)-torsor. There are
\begin{enumerate}
\item a smooth projective \(S\)-scheme \(X\) of relative dimension
\(|H|^2\) with integral geometric fibres,
\item an open subscheme \(U\subseteq X\) such that \(U_{\bar s}\) is
dense in \(X_{\bar s}\) for every geometric point \(\bar s\) of
\(S\),
\item a right \(H\)-torsor \(Q\to U\) and an isomorphism of right
\(G\)-torsors
\(\theta:Q\times^H G\xrightarrow{\sim}P_U\),
\end{enumerate}
with the following property: for every field \(L\) over \(S\) and
every marked lift \((Q_0,\theta_0)\) of \(P_L\) along \(\phi\), there
is a point \(u\in U(L)\) such that
\(u^*(Q,\theta)\simeq(Q_0,\theta_0)\) as marked lifts.
\end{proposition}

\begin{proof}
Let \((U,Q,\theta)\) be as in \cref{con:lifting-model}; the property
for fields is \cref{prop:marked-lifting-model-field-versality}. By
\cref{lem:lifting-compactification-tame-fibres}, the stack
\(\mathcal W=[W/H]\) and its open substack \(U\) satisfy the
hypotheses of \cref{lem:destackification-preserves-scheme-open} with
\(r=|H|^2\). \Cref{lem:destackification-preserves-scheme-open} gives the smooth projective \(S\)-scheme
\(X\) containing \(U\).
\end{proof}
\subsection{The two Sylow reductions}

\begin{sublemma}[A prime-to-\(\ell\) change of base for a finite torsor]
\label{sub:finite-sylow-base-change}
Let \(S=\operatorname{Spec}R\) be a nonempty connected normal
noetherian affine scheme with fraction field \(F\). Let \(P\to S\)
be a connected right torsor under a finite group \(G\), let
\(\ell\) be a prime number, and let \(G_\ell\subseteq G\) be a Sylow
\(\ell\)-subgroup. Put \(S'=P/G_\ell\).
\begin{enumerate}
\item The scheme \(S'\) is affine, connected and normal, the
morphism \(S'\to S\) is finite \'etale of degree \([G:G_\ell]\), and
\(P\to S'\) is a right \(G_\ell\)-torsor.
\item The schemes \(S'\) and \(P\) are integral. If \(F'\) is the
fraction field of \(S'\), then
\(S'\times_S\operatorname{Spec}F=\operatorname{Spec}F'\),
\([F':F]=[G:G_\ell]\), and
\(P\times_{S'}\operatorname{Spec}F'\) is the spectrum of a field.
\end{enumerate}
\end{sublemma}

\begin{proof}
\textup{(i)} The group \(G_\ell\) acts freely on the affine scheme
\(P\), so \(S'\) is a scheme and \(P\to S'\) is a
\(G_\ell\)-torsor \cite[\stackstag{07S7}]{Stacks}. \'Etale locally on
\(S\), the torsor \(P\) is \(S\times G\) and \(S'\) is
\(S\times(G/G_\ell)\); hence \(S'\to S\) is finite \'etale of degree
\([G:G_\ell]\), and \(S'\) is affine. It is connected, as the image
of \(P\). The schemes \(S'\) and \(P\) are normal, because normality
is local in the \'etale topology \cite[\stackstag{034F}]{Stacks}.

\textup{(ii)} A nonempty connected noetherian normal scheme is
integral \cite[\stackstag{033M}]{Stacks}. Let \(Y'\to Y\) be a finite
dominant morphism of integral schemes, and let \(\eta\) be the
generic point of \(Y\). The coordinate ring of
\(Y'\times_Y\eta\) is a localization of a domain and is finite over
the field \(\kappa(\eta)\), so it is a field
\cite[\stackstag{00GS}]{Stacks}; it is the function field of \(Y'\).
Apply this to \(S'\to S\) and to \(P\to S'\). The degree
\([F':F]\) is the degree of the finite \'etale morphism \(S'\to S\).
\end{proof}

\begin{sublemma}[A pro-\(\ell\) Sylow fixed field]
\label{sub:pro-ell-sylow-fixed-field}
Let \(F'\) be a field with separable closure
\((F')^{\mathrm{sep}}\), put
\(\Gamma=\operatorname{Gal}((F')^{\mathrm{sep}}/F')\), and let
\(\ell\) be a prime number.
\begin{enumerate}
\item There is a closed subgroup \(\Gamma_\ell\subseteq\Gamma\)
whose image in \(\Gamma/J\) is a Sylow \(\ell\)-subgroup for every
open normal subgroup \(J\subseteq\Gamma\). It is a pro-\(\ell\)
group.
\item Put \(L=((F')^{\mathrm{sep}})^{\Gamma_\ell}\). Then
\(\operatorname{Gal}((F')^{\mathrm{sep}}/L)=\Gamma_\ell\), and every
finite extension \(E\) of \(F'\) contained in \(L\) has degree prime
to \(\ell\).
\item Let \(G_\ell\) be a finite \(\ell\)-group, and let \(P'\) be a
connected right \(G_\ell\)-torsor over \(F'\), classified by a
homomorphism \(\chi:\Gamma\to G_\ell\). Then
\(\chi(\Gamma_\ell)=G_\ell\). For every open normal subgroup \(J\)
of \(\Gamma_\ell\) contained in \(\ker(\chi|_{\Gamma_\ell})\), with
quotient \(H=\Gamma_\ell/J\) and induced surjection
\(\phi:H\twoheadrightarrow G_\ell\), the spectrum of
\(((F')^{\mathrm{sep}})^{J}\) is a right \(H\)-torsor over \(L\)
that is a marked lift of \(P'_L\) along \(\phi\).
\end{enumerate}
\end{sublemma}

\begin{proof}
\textup{(i)} This is the existence of Sylow pro-\(\ell\)-subgroups
of profinite groups \cite[Chapter~I, \S1.4]{SerreGC}; the argument
is as follows. For an open normal subgroup \(J\subseteq\Gamma\), let
\(q_J:\Gamma\to\Gamma/J\) be the projection and let
\(\operatorname{Syl}_\ell(\Gamma/J)\) be the finite nonempty set of
Sylow \(\ell\)-subgroups of \(\Gamma/J\). For \(J'\subseteq J\), the
surjection \(\Gamma/J'\to\Gamma/J\) carries a Sylow
\(\ell\)-subgroup onto an \(\ell\)-subgroup of index prime to
\(\ell\), that is, onto a Sylow \(\ell\)-subgroup. The sets
\(\operatorname{Syl}_\ell(\Gamma/J)\) therefore form an inverse
system of finite nonempty sets, and its limit is nonempty
\cite[\stackstag{086J}]{Stacks}. Choose an element \((Y_J)_J\) of
the limit, and put
\[
 \Gamma_\ell=\bigcap_{J}q_J^{-1}(Y_J)=\varprojlim_JY_J .
\]
This is a closed subgroup and a pro-\(\ell\) group. For \(y\in Y_J\),
the fibres over \(y\) of the surjections \(Y_{J'}\to Y_J\),
\(J'\subseteq J\), form an inverse system of finite nonempty sets;
a point of its limit \cite[\stackstag{086J}]{Stacks} is an element
of \(\Gamma_\ell\) with image \(y\). Hence \(q_J(\Gamma_\ell)=Y_J\).

\textup{(ii)} Since \(L\) is an intermediate field of the separable
algebraic extension \((F')^{\mathrm{sep}}/F'\), the field
\((F')^{\mathrm{sep}}\) is also a separable closure of \(L\).
The equality
\(\operatorname{Gal}((F')^{\mathrm{sep}}/L)=\Gamma_\ell\) is the
Galois correspondence for the closed subgroup \(\Gamma_\ell\)
\cite[\stackstag{0BML}]{Stacks}. Let \(E\subseteq L\) be finite over
\(F'\), let \(D\subseteq(F')^{\mathrm{sep}}\) be a finite Galois
extension of \(F'\) containing \(E\), and let
\(J=\operatorname{Gal}((F')^{\mathrm{sep}}/D)\). The subgroup
\(\operatorname{Gal}(D/E)\) of \(\Gamma/J=\operatorname{Gal}(D/F')\)
contains \(q_J(\Gamma_\ell)=Y_J\). Hence
\[
 [E:F']=\bigl[\operatorname{Gal}(D/F'):\operatorname{Gal}(D/E)\bigr]
 \quad\text{divides}\quad[\Gamma/J:Y_J],
\]
which is prime to \(\ell\).

\textup{(iii)} The torsor \(P'\) is connected, so \(\Gamma\) acts
transitively on \(P'((F')^{\mathrm{sep}})\) and \(\chi\) is
surjective. By \textup{(i)}, applied to \(J=\ker\chi\), the image
\(\chi(\Gamma_\ell)\) is a Sylow \(\ell\)-subgroup of the
\(\ell\)-group \(G_\ell\), so \(\chi(\Gamma_\ell)=G_\ell\). By
\textup{(ii)} and \cite[\stackstag{0BML}]{Stacks}, the field
\(((F')^{\mathrm{sep}})^{J}\) is a finite Galois extension of \(L\)
with group \(H\); its spectrum is the \(H\)-torsor over \(L\)
classified by \(\Gamma_\ell\twoheadrightarrow H\). The associated
\(G_\ell\)-torsor is classified by the composite
\(\Gamma_\ell\twoheadrightarrow H\xrightarrow{\phi}G_\ell\), which is
\(\chi|_{\Gamma_\ell}\), the homomorphism classifying \(P'_L\).
The equivalence between finite \'etale schemes and finite sets with
continuous Galois action \cite[\stackstag{0BND}]{Stacks} therefore
gives an isomorphism of right \(G_\ell\)-torsors
\(Q_0\times^H G_\ell\simeq P'_L\), where
\(Q_0=\Spec\bigl(((F')^{\mathrm{sep}})^J\bigr)\). Choose such an isomorphism
as the marking in the sense of \cref{def:marked-lift}.
\end{proof}

\begin{lemma}[A closed point of prime-to-\(\ell\) degree]
\label{lem:closed-point-from-sylow-field}
Let \(S=\operatorname{Spec}R\), \(F\), \(S'\) and \(F'\) be as in
\cref{sub:finite-sylow-base-change}, and let \(L\) be the field of
\cref{sub:pro-ell-sylow-fixed-field}\textup{(ii)} for \(F'\). Let
\(X\to S'\) be a smooth projective morphism of relative dimension
\(r\) with integral geometric fibres, and let \(U\subseteq X\) be an
open subscheme whose geometric fibres over \(S'\) are dense in those
of \(X\). Regard \(X\) and \(U\) as \(S\)-schemes through
\(S'\to S\).
\begin{enumerate}
\item The scheme \(X\) is smooth and projective over \(S\) of
relative dimension \(r\), and for every geometric point \(\bar s\) of
\(S\) the open subscheme \(U_{\bar s}\) is dense in \(X_{\bar s}\).
In particular \(\dim(X_{\bar s}\setminus U_{\bar s})\leq r-1\),
with the convention that the empty scheme has dimension \(-1\).
\item If \(U_{F'}(L)\) is nonempty, then \(U_F\) contains a closed point
whose residue field has degree prime to \(\ell\) over \(F\).
\end{enumerate}
\end{lemma}

\begin{proof}
\textup{(i)} The morphism \(S'\to S\) is finite \'etale, so the
composite \(X\to S'\to S\) is smooth of relative dimension \(r\). A
finite morphism is projective \cite[\stackstag{0B3I}]{Stacks}, and a
composite of projective morphisms over an affine base is projective
\cite[\stackstag{0C4P}]{Stacks}. For a geometric point \(\bar s\) of
\(S\), the scheme \(S'\times_S\bar s\) is a finite disjoint union of
copies \(\bar s'\) of \(\bar s\), and
\[
 X_{\bar s}=\coprod_{\bar s'}X_{\bar s'},\qquad
 U_{\bar s}=\coprod_{\bar s'}U_{\bar s'} .
\]
Each \(X_{\bar s'}\) is integral of dimension \(r\) and contains
\(U_{\bar s'}\) as a dense open subscheme. Its complement is a
proper closed subset and has dimension at most \(r-1\) (it is
empty if \(r=0\)). The assertions follow by taking the disjoint
union over \(\bar s'\).

\textup{(ii)} By \cref{sub:finite-sylow-base-change}\textup{(ii)},
\(X_F=X\times_{S'}\Spec F'\) and \(U_F=U\times_{S'}\Spec F'\). Let
\(u\in U_{F'}(L)\), and let \(\Spec D\subseteq U_F\) be an affine open
subscheme through which \(u\) factors. The image \(E\) of the
\(F'\)-algebra homomorphism \(D\to L\) is a finitely generated
\(F'\)-algebra and a domain algebraic over \(F'\); hence it is a
field, finite over \(F'\) \cite[\stackstag{00GS}]{Stacks}. The
kernel of \(D\to E\) is a maximal ideal, that is, a closed point
\(x\) of \(U_F\) with residue field \(E\); it is closed in
\(X_F\) because its residue field is finite over \(F\)
\cite[\stackstag{01TF}]{Stacks}. Its degree
over \(F\) is
\[
 [E:F]=[E:F']\,[F':F],
\]
and both factors are prime to \(\ell\), by
\cref{sub:pro-ell-sylow-fixed-field}\textup{(ii)} and
\cref{sub:finite-sylow-base-change}.
\end{proof}
\section{Regular finite flat multisections}
\label{sec:regular-finite-flat-multisections}

Let \(R\) be a connected noetherian regular semilocal ring with
fraction field \(F\), let \(\ell\) be a prime number, let
\(X\) be a regular flat projective \(R\)-scheme whose geometric
closed and generic fibres are pure of the same dimension, and let
\(U\subseteq X\) be an open subscheme meeting every irreducible
component of every geometric closed fibre.
\Cref{thm:regular-multisections} constructs, from a
regular projective \(R\)-scheme \(C\to X\) whose generic fibre is
finite of degree prime to \(\ell\), a regular scheme \(T\),
finite locally free over \(R\) of rank prime to \(\ell\), with an
\(R\)-morphism \(T\to U\). Each step embeds \(C\) in
\(\mathbb P^a_X\), blows up \(\mathbb P^a_X\) along \(C\) (with
blowup \(\pi:W\to\mathbb P^a_X\) and exceptional divisor \(E\)), and replaces
\(C\) by the common zero scheme of sections of the very ample sheaf
\(\mathcal L_b(\mathcal A)=\pi^*\mathcal A^{\otimes b}\otimes\mathcal O_W(-E)\), where
\(\mathcal A\) is very ample on \(\mathbb P^a_X\) and \(b\) is a sufficiently large
multiple of \(\ell\) (\cref{con:multisection-replacement}). After the
possible finite \'etale base change, the generic degree of the new
zero scheme over each generic residue field is congruent to the
negative of the previous degree modulo \(\ell\)
(\cref{prop:multisection-generic-chow-degree}). Simultaneous choices of
sections on the closed fibres ensure regularity of the zero scheme
(\cref{sub:multisection-simultaneous-parameter-choice,sub:multisection-absolute-cotangent}).
Each step decreases the maximum dimension of the closed fibres until
it is at most zero, without subsequently increasing it. The maximum
dimension of the closed fibres of the inverse image of \(X\setminus U\)
decreases until these fibres are empty
(\cref{prop:multisection-decrease-closed-fibres}). In the iteration
of \cref{thm:regular-multisections}, the generic
degrees over all components of each finite \'etale base extension
have a common nonzero residue modulo \(\ell\).
\Cref{cor:arithmetic-multisections} obtains the
initial \(C\) from an alteration
(\cref{lem:multisection-arithmetic-alteration}) when \(R\) is a
regular semilocal domain essentially of finite type over a field or
over \(\Z\).

\begin{notation}[Dimensions and degrees]
\label{not:multisection-dimensions-degrees}
Let \(\ell\) be a prime number, let \(R\) be a noetherian
regular semilocal ring, and write \(S=\operatorname{Spec}R\).
For a closed point \(s\in S\), put \(d_s=\dim\mathcal O_{S,s}\)
and let \(\bar s\) be a geometric point with residue field an
algebraic closure of \(\kappa(s)\). In this section the dimension of
an empty scheme is \(-1\). If \(W\) is a finite scheme over
a field \(F\), we write
\[
 \deg_F(W)=\dim_F \Gamma(W,\mathcal O_W).
\]
Equivalently, this is the degree of its fundamental
zero-cycle, with both local lengths and residue field degrees
included.
\end{notation}

Gabber--Liu--Lorenzini construct finite flat quasi-sections by
successive hypersurface sections under Cohen--Macaulay hypotheses
\cite[Theorem~6.3]{GLLHypersurfaces15}. Their
\cite[Theorem~3.4]{GLLHypersurfaces15} gives conditions under which
successive sections preserve a prescribed regular immersion;
\cite[Theorem~5.1 and Remark~6.6]{GLLHypersurfaces15} give avoidance
variants.
Related moving results for horizontal cycles are proved in
\cite[\arfiv{1209.2828}{S2.Thmthm3}{Theorem~2.3}]{GLLIndex13} and
\cite[Theorem~7.2]{GLLHypersurfaces15}.
The finite-support moving statement for a regular immersion
over a semilocal noetherian base is
\cite[Proposition~2.4]{GLLIndex13}.

We first give the parameter choice and the two calculations
used in \cref{con:multisection-replacement}.

\subsection{Choosing sections on all closed fibres}

\begin{sublemma}[Lifting choices from the closed residue fields]
\label{sub:multisection-simultaneous-parameter-choice}
Let \(R\) be a semilocal ring, let \(\ell\) be a prime
number, and let \(a\in\mathbf N\). For each maximal ideal
\(\mathfrak m\) of \(R\), let
\(V_{\mathfrak m}\subseteq\mathbb A^a_{\kappa(\mathfrak m)}\)
be a nonempty open subscheme. There are a finite \'etale
\(R\)-algebra \(R'\), free of rank prime to \(\ell\) as an
\(R\)-module, and an element \(x\in(R')^a\) such that, for every
maximal ideal \(\mathfrak m'\) of \(R'\) with contraction
\(\mathfrak m\) in \(R\), the image of \(x\) in
\(\kappa(\mathfrak m')^a\) belongs to
\(V_{\mathfrak m}(\kappa(\mathfrak m'))\). If all the fields
\(\kappa(\mathfrak m)\) are infinite, one can take \(R'=R\).
\end{sublemma}

\begin{proof}
If \(R=0\), take \(R'=R\) with free rank one and \(x=0\); all
conditions on residue fields are vacuous. Assume \(R\ne0\).
The polynomial avoidance step and its combination with the Chinese
remainder theorem also appear in
\cite[proof of
\arxivhtml{1404.5366v1}{S2.Thmthm1}{Theorem~2.1}]{GLLHypersurfaces15}.

\emph{Values of a polynomial.}
Let \(k\) be a field and let \(g\in k[t_1,\ldots,t_a]\) be a nonzero
polynomial of degree at most \(\delta\) in each variable. If
\(|k|>\delta\), there is \(y\in k^a\) with \(g(y)\neq0\). This is
clear for \(a=0\). For \(a\geq1\), write
\(g=\sum_jg_j(t_1,\ldots,t_{a-1})t_a^j\), choose by induction
\((y_1,\ldots,y_{a-1})\) at which some \(g_j\) does not vanish, and
choose \(y_a\in k\) outside the at most \(\delta\) roots of
\(g(y_1,\ldots,y_{a-1},t_a)\).

\smallskip\noindent
\emph{The algebra \(R'\).}
For each \(\mathfrak m\), choose a nonzero polynomial
\(g_{\mathfrak m}\in\kappa(\mathfrak m)[t_1,\ldots,t_a]\) with
\(D(g_{\mathfrak m})\subseteq V_{\mathfrak m}\), and let \(\delta\)
bound the degrees of all the \(g_{\mathfrak m}\) in each variable.
Choose an integer \(e\geq1\) prime to \(\ell\) such that
\(|\kappa(\mathfrak m)|^e>\delta\) for every finite field
\(\kappa(\mathfrak m)\); if there is none, take \(e=1\). For each
\(\mathfrak m\), let \(f_{\mathfrak m}\in\kappa(\mathfrak m)[t]\) be
the following monic separable polynomial of degree \(e\): the
minimal polynomial of a generator of the extension of degree \(e\)
of \(\kappa(\mathfrak m)\) if this field is finite, and a product of
\(e\) distinct monic linear polynomials otherwise. By the Chinese
remainder theorem \cite[\stackstag{00DT}]{Stacks}, there is a monic
polynomial \(f\in R[t]\) of degree \(e\) with image
\(f_{\mathfrak m}\) in \(\kappa(\mathfrak m)[t]\) for every
\(\mathfrak m\). Put
\[
 R'=R[t]/(f),\qquad
 R'/\mathfrak mR'=\kappa(\mathfrak m)[t]/(f_{\mathfrak m}).
\]
The \(R\)-module \(R'\) is free of rank \(e\). Every maximal ideal
of \(R'\) lies over a maximal ideal \(\mathfrak m\) of \(R\), and
the derivative \(f'\) is a unit of
\(\kappa(\mathfrak m)[t]/(f_{\mathfrak m})\), because
\(f_{\mathfrak m}\) is separable. Hence \(f'\) is a unit of \(R'\),
and \(R\to R'\) is standard \'etale
\cite[\stackstag{00UB}, \stackstag{00UC}]{Stacks}.

\smallskip\noindent
\emph{The element \(x\).}
Let \(\mathfrak m'\) be a maximal ideal of \(R'\) over
\(\mathfrak m\). The field \(\kappa(\mathfrak m')\) is a factor of
\(\kappa(\mathfrak m)[t]/(f_{\mathfrak m})\); it is infinite or has
\(|\kappa(\mathfrak m)|^e>\delta\) elements. By the first step
there is \(y_{\mathfrak m'}\in\kappa(\mathfrak m')^a\) with
\(g_{\mathfrak m}(y_{\mathfrak m'})\neq0\). The Chinese remainder
theorem \cite[\stackstag{00DT}]{Stacks} for the finitely many
maximal ideals of \(R'\) gives \(x\in(R')^a\) with image
\(y_{\mathfrak m'}\) in \(\kappa(\mathfrak m')^a\) for every
\(\mathfrak m'\).
\end{proof}
\begin{sublemma}[Dimensions of general hyperplane sections]
\label{sub:multisection-hyperplane-dimensions}
Let \(k\) be a field, let \(b,c\) be natural numbers with
\(c\geq1\), and let \(Z_1,\ldots,Z_v\subseteq\mathbb P_k^b\)
be closed subschemes. Identify \(\mathbb A_k^{c(b+1)}\) with the
affine space of ordered \(c\)-tuples \(h=(h_1,\ldots,h_c)\) of
linear forms in the homogeneous coordinates \(x_0,\ldots,x_b\) of
\(\mathbb P_k^b\). For a field extension \(k'/k\), a point
\(h\in\mathbb A_k^{c(b+1)}(k')\), and a closed subscheme
\(Z\subseteq\mathbb P_k^b\), let \(Z_{k'}\cap V_+(h)\) be the closed
subscheme of \(Z_{k'}\) defined by \(h_1=\cdots=h_c=0\). There is a
nonempty open subset \(U\subseteq\mathbb A_k^{c(b+1)}\) such that,
for every field extension \(k'/k\), every \(h\in U(k')\), and every
\(1\leq j\leq v\),
\[
 \dim\bigl((Z_j)_{k'}\cap V_+(h)\bigr)\leq\dim Z_j-c
 \quad\text{if }\dim Z_j\geq c,
 \qquad
 (Z_j)_{k'}\cap V_+(h)=\varnothing
 \quad\text{if }\dim Z_j<c .
\]
\end{sublemma}

\begin{proof}
If every \(Z_j\) is empty (including the case \(v=0\)), take
\(U=\mathbb A_k^{c(b+1)}\). Otherwise discard the empty subschemes,
which impose no condition. The scheme \(\mathbb A_k^{c(b+1)}\) is integral,
so a finite intersection of nonempty open subsets is nonempty. Hence we
may assume that \(v=1\) and that \(Z=Z_1\) is nonempty. Let
\(I\subseteq Z\times_k\mathbb A_k^{c(b+1)}\) be the common zero
scheme of the \(c\) universal linear forms, with projections
\(\operatorname{pr}_1:I\to Z\) and
\(\operatorname{pr}_2:I\to\mathbb A_k^{c(b+1)}\). For a point
\(h\in\mathbb A_k^{c(b+1)}(k')\), the fibre of
\(\operatorname{pr}_2\) over \(h\) is \(Z_{k'}\cap V_+(h)\). The
morphism \(\operatorname{pr}_1\) is the vector bundle associated
with the kernel of
\[
 \mathcal O_{Z}^{\oplus c(b+1)}
       \longrightarrow \mathcal O_{Z}(1)^{\oplus c},
 \qquad
 (a_{im})_{i,m}\longmapsto
 \Bigl(\sum_{m=0}^{b}a_{im}x_m\Bigr)_{1\leq i\leq c}.
\]
This map is surjective, because at every point of \(Z\) one
homogeneous coordinate is nonzero. Hence
\begin{equation}\label{eq:incidence-dimension}
 \dim I=\dim Z+c(b+1)-c .
\end{equation}
The morphism \(\operatorname{pr}_2\) is projective. If
\(\dim Z<c\), then \(\dim I<c(b+1)\) by
\eqref{eq:incidence-dimension}, so the closed subset
\(\operatorname{pr}_2(I)\) is proper, and we take for \(U\) its
complement. Suppose that \(\dim Z\geq c\). By the dimension formula
\cite[\stackstag{02JU}]{Stacks} and
\eqref{eq:incidence-dimension}, the generic fibre of every
irreducible component of \(I\) dominating \(\mathbb A_k^{c(b+1)}\)
has dimension at most \(\dim Z-c\). The fibre dimension of the
proper morphism \(\operatorname{pr}_2\) is upper semicontinuous
\cite[\stackstag{0D4I}]{Stacks}, so the set \(U\) of points of
\(\mathbb A_k^{c(b+1)}\) over which the fibre of
\(\operatorname{pr}_2\) has dimension at most \(\dim Z-c\) is open,
and it contains the generic point. In both cases, the dimension of
the fibre over a point \(h\in U(k')\) is the dimension of the fibre
over the image of \(h\) in \(U\)
\cite[\stackstag{05F8}]{Stacks}, which proves the assertion for
every extension \(k'/k\).
\end{proof}
\begin{sublemma}[Cotangent classes of linear forms at a point]
\label{sub:multisection-pointwise-incidence}
Let \(\Omega\) be an algebraically closed field, let
\(Z\subseteq\mathbb P^b_\Omega\) be a closed subscheme, let
\(z\in Z(\Omega)\), let \(\mathfrak m_z\) be the maximal ideal of
\(\mathcal O_{Z,z}\), and put
\(t=\dim_\Omega\mathfrak m_z/\mathfrak m_z^2\). Fix an index \(m\)
with \(x_m(z)\neq0\). For a linear form \(h\) in
\(x_0,\ldots,x_b\) with \(h(z)=0\), let
\([h]_z\in\mathfrak m_z/\mathfrak m_z^2\) be the class of the
function \(h/x_m\).
\begin{enumerate}
\item The linear forms vanishing at \(z\) form a hyperplane in the
space \(\Omega^{b+1}\) of linear forms, and the \(\Omega\)-linear map
\(h\mapsto[h]_z\) from this hyperplane to
\(\mathfrak m_z/\mathfrak m_z^2\) is surjective. Another choice of
\(m\) multiplies all the classes \([h]_z\) by the same nonzero
scalar.
\item Let \(1\leq i\leq c\) be integers, and let
\(B^{(i)}_z\subseteq(\Omega^{b+1})^c=\mathbb A^{c(b+1)}(\Omega)\) be
the set of \(c\)-tuples \((h_1,\ldots,h_c)\) such that
\(h_1(z)=\cdots=h_i(z)=0\) and the classes
\([h_1]_z,\ldots,[h_i]_z\) are linearly dependent. Then
\(B^{(i)}_z\) is a closed subset of the affine space, and
\[
 \dim B^{(i)}_z=c(b+1)-(t+1)\ \text{ if }t\geq i,
 \qquad
 \dim B^{(i)}_z=c(b+1)-i\ \text{ if }t<i .
\]
\end{enumerate}
\end{sublemma}

\begin{proof}
Assertion \textup{(i)} and the case \(i=1\) of assertion
\textup{(ii)} are the computation at a closed point in the proof of
the Bertini theorem \cite[\stackstag{0FD6}]{Stacks}.

\textup{(i)} Put \(y_j=x_j/x_m\). The scheme \(Z\cap D_+(x_m)\) is a
closed subscheme of \(\Spec\Omega[y_j:j\neq m]\), so the ideal
\(\mathfrak m_z\) is generated by the functions \(y_j-y_j(z)\),
\(j\neq m\). The linear form \(x_j-y_j(z)x_m\) vanishes at \(z\),
and its class is the class of \(y_j-y_j(z)\). This proves the
surjectivity. The vanishing at \(z\) is the single nontrivial linear
condition \((h/x_m)(z)=0\). For another index \(m'\) with
\(x_{m'}(z)\neq0\), we have
\(h/x_{m'}=(x_m/x_{m'})\cdot h/x_m\), and the class of this product
of an element of \(\mathfrak m_z\) by a unit is
\((x_m/x_{m'})(z)\,[h]_z\).

\textup{(ii)} By \textup{(i)}, the map
\((h_1,\ldots,h_i)\mapsto([h_1]_z,\ldots,[h_i]_z)\) from the
\(i\)-tuples of linear forms vanishing at \(z\) to the space
\((\mathfrak m_z/\mathfrak m_z^2)^i\) of \(t\times i\) matrices is
linear and surjective, and \(B^{(i)}_z\) is the product of
\((\Omega^{b+1})^{c-i}\) with the inverse image of the closed set of
matrices of rank less than \(i\). If \(t<i\), every matrix has rank
less than \(i\), and \(B^{(i)}_z\) is defined by the \(i\)
independent conditions \(h_1(z)=\cdots=h_i(z)=0\). If \(t\geq i\),
the matrices of rank \(u\) are parametrized by a \(u\)-dimensional
subspace of \(\Omega^t\) and a surjection of \(\Omega^i\) onto it,
so they form a set of dimension \(u(t-u)+ui\); for \(u=i-1\) this
is \(ti-(t-i+1)\), the largest value for \(u<i\). Hence the
matrices of rank less than \(i\) have codimension \(t-i+1\), and
\(B^{(i)}_z\) has codimension \(i+(t-i+1)=t+1\).
\end{proof}

\Needspace{7\baselineskip}
\begin{lemma}[Hyperplane values and cotangent vectors]
\label{lem:multisection-geometric-incidence}
Let \(k\) be a field, let \(Z\subseteq\mathbb P_k^b\) be a
closed subscheme, and let \(Z_1,\ldots,Z_v\) be closed
subschemes of \(Z\). Let \(c\geq1\) be a natural number, and use the
notation \(\mathbb A_k^{c(b+1)}\), \(h=(h_1,\ldots,h_c)\), and
\(V_+(h)\) of \cref{sub:multisection-hyperplane-dimensions}. There is
a nonempty open subset \(U\subseteq\mathbb A_k^{c(b+1)}\) such that,
for every field extension \(k'/k\) with algebraic closure
\(\overline{k'}\) and every \(h\in U(k')\), the following hold.
\begin{enumerate}
\item For each \(1\leq i\leq c\) and every closed point \(z\) of
\(Z_{\overline{k'}}\) with \(h_1(z)=\cdots=h_i(z)=0\), the classes
\([h_1]_z,\ldots,[h_i]_z\) of
\cref{sub:multisection-pointwise-incidence} are linearly independent
in the cotangent space \(\mathfrak m_z/\mathfrak m_z^2\) of
\(Z_{\overline{k'}}\) at \(z\).
\item For each \(1\leq j\leq v\),
\(\dim\bigl((Z_j)_{k'}\cap V_+(h)\bigr)\leq\dim Z_j-c\) if
\(\dim Z_j\geq c\), and \((Z_j)_{k'}\cap V_+(h)=\varnothing\) if
\(\dim Z_j<c\).
\end{enumerate}
\end{lemma}

\begin{proof}
The proof adapts the cotangent-space incidence argument in
\cite[proof of \stackstag{0FD6}]{Stacks} to strata of constant
cotangent dimension and to ordered tuples of hyperplanes.

Fix \(1\leq i\leq c\), and let
\(I^{(i)}\subseteq Z\times_k\mathbb A_k^{c(b+1)}\) be the common zero
scheme of the first \(i\) universal linear forms, with projections
\(\operatorname{pr}_Z\) and \(\operatorname{pr}_2\).

\smallskip\noindent
\emph{A constructible set of failures.}
The universal forms \(h_1,\ldots,h_i\) are sections of
\(\operatorname{pr}_Z^*\mathcal O_Z(1)\) vanishing on \(I^{(i)}\).
Their differentials in the direction of \(Z\) define a homomorphism
of coherent sheaves on \(I^{(i)}\)
\begin{equation}\label{eq:incidence-differential}
 \delta^{(i)}:\mathcal O_{I^{(i)}}^{\oplus i}\longrightarrow
 \mathcal F^{(i)}=\left.\operatorname{pr}_{Z}^{*}
  \bigl(\Omega^1_{Z/k}\otimes\mathcal O_Z(1)\bigr)\right|_{I^{(i)}};
\end{equation}
it does not depend on local trivializations of \(\mathcal O_Z(1)\),
because a change of trivialization changes the differential of a
section by a multiple of the section, which vanishes on
\(I^{(i)}\). Let \(\mathcal B^{(i)}\subseteq I^{(i)}\) be the set of
points \(x\) at which \(\delta^{(i)}\otimes\kappa(x)\) is not
injective, that is, at which
\[
 \dim_{\kappa(x)}\bigl(\mathcal F^{(i)}\otimes\kappa(x)\bigr)
 -\dim_{\kappa(x)}\bigl(\operatorname{coker}(\delta^{(i)})
      \otimes\kappa(x)\bigr)<i .
\]
The two dimensions are upper semicontinuous functions of \(x\)
\cite[\stackstag{0C3F}]{Stacks}, so \(\mathcal B^{(i)}\) is
constructible. Let \(\Omega\) be an algebraically closed field
containing \(k\), and let \((z,h)\) be an \(\Omega\)-point of
\(I^{(i)}\). By \cite[\stackstag{0B2E}]{Stacks}, the fibre of
\(\Omega^1_{Z/k}\) at \(z\) is the cotangent space
\(\mathfrak m_z/\mathfrak m_z^2\) of \(Z_\Omega\) at \(z\), and
\(\delta^{(i)}\) at \((z,h)\) sends the \(j\)-th basis vector to
\([h_j]_z\). Hence, in the notation of
\cref{sub:multisection-pointwise-incidence}\textup{(ii)} for
\(Z_\Omega\),
\begin{equation}\label{eq:incidence-failures}
 (z,h)\in\mathcal B^{(i)}\iff h\in B^{(i)}_z .
\end{equation}

\smallskip\noindent
\emph{The dimension of the set of failures.}
Let \(\bar k\) be an algebraic closure of \(k\). For
\(t\in\mathbf N\), let \(Z_t\subseteq Z_{\bar k}\) be the locally
closed subset of points \(z\) with
\(\dim\Omega^1_{Z/k}\otimes\kappa(z)=t\)
\cite[\stackstag{0C3F}]{Stacks}. For each irreducible component of
the locally closed stratum \(Z_t\), choose a closed point of that component, which is also
closed in \(Z_{\bar k}\). At such a point \(z\), the dimension of
that component is at most \(\dim\mathcal O_{Z_{\bar k},z}\)
\cite[\stackstag{00P1}]{Stacks}, and the latter is at most the
number \(t\) of generators of \(\mathfrak m_z\)
\cite[\stackstag{00KQ}]{Stacks}. Hence \(\dim Z_t\leq t\). By
\eqref{eq:incidence-failures} and
\cref{sub:multisection-pointwise-incidence}\textup{(ii)}, the fibre
of \(\mathcal B^{(i)}_{\bar k}\) over such a point \(z\) has
dimension \(c(b+1)-(t+1)\) if \(t\geq i\), and \(c(b+1)-i\) if
\(t<i\). In both cases
\[
 \dim\bigl(\mathcal B^{(i)}_{\bar k}\cap
      \operatorname{pr}_Z^{-1}(Z_t)\bigr)
 \leq t+c(b+1)-\max(t+1,i)\leq c(b+1)-1 ,
\]
by the dimension formula \cite[\stackstag{02JU}]{Stacks}. Indeed,
stratify this constructible incidence set into finitely many integral
locally closed pieces with their reduced structures. For each piece,
over a nonempty open of its image closure the fibre dimension is the
generic fibre dimension \cite[\stackstag{05F7}]{Stacks}; a closed
point of that open in the image bounds it by the computed fibre
dimension over \(Z_t\). This bounds the dimension of every piece,
hence also of the closure of the constructible incidence set. There
are finitely many strata \(Z_t\). The dimension of a constructible
subset of a finite-type scheme over a field is unchanged by field
extension: decompose it into finitely many locally closed pieces
and apply \cite[\stackstag{00P3}]{Stacks} on affine charts. Thus
\(\dim\mathcal B^{(i)}\leq c(b+1)-1\).

\smallskip\noindent
\emph{The open set.}
By Chevalley's theorem \cite[\stackstag{054K}]{Stacks}, the set
\(\operatorname{pr}_2(\mathcal B^{(i)})\) is constructible. Applying
the dimension formula to the integral locally closed pieces and
the closures of their images \cite[\stackstag{02JU}]{Stacks}, its
closure \(C^{(i)}\) has dimension at most \(c(b+1)-1\); it is a
proper closed subset of \(\mathbb A_k^{c(b+1)}\). Let \(U'\) be the
open subset given by
\cref{sub:multisection-hyperplane-dimensions} for
\(Z_1,\ldots,Z_v\), and put
\[
 U=U'\setminus\bigl(C^{(1)}\cup\cdots\cup C^{(c)}\bigr).
\]
It is nonempty because \(\mathbb A_k^{c(b+1)}\) is irreducible. Let
\(h\in U(k')\). Assertion \textup{(ii)} holds by
\cref{sub:multisection-hyperplane-dimensions}. For a closed point
\(z\) of \(Z_{\overline{k'}}\) with
\(h_1(z)=\cdots=h_i(z)=0\), the \(\overline{k'}\)-point \((z,h)\) of
\(I^{(i)}\) does not lie in \(\mathcal B^{(i)}\), since its image
\(h\) is not in \(C^{(i)}\). By \eqref{eq:incidence-failures}, the
classes \([h_1]_z,\ldots,[h_i]_z\) are linearly independent. This is
\textup{(i)}.
\end{proof}

\begin{sublemma}[Regular cuts from geometric cotangent independence]
\label{sub:multisection-absolute-cotangent}
Let \(W\) be a regular projective \(S\)-scheme, let \(\mathcal L\)
be an invertible sheaf on \(W\), and let
\(s_1,\ldots,s_c\in\Gamma(W,\mathcal L)\). For \(1\leq i\leq c\), let
\(Y_i\subseteq W\) be the common zero scheme of \(s_1,\ldots,s_i\).
Suppose that, for every closed point \(s\) of \(S\), every
\(1\leq i\leq c\), and every closed point \(\bar w\) of
\((Y_i)_{\bar s}\), local equations of \(s_1,\ldots,s_i\) at
\(\bar w\) have linearly independent classes in the cotangent space
\(\mathfrak m_{\bar w}/\mathfrak m_{\bar w}^2\) of \(W_{\bar s}\) at
\(\bar w\). Then, for every \(1\leq i\leq c\):
\begin{enumerate}
\item at every point \(w\) of \(Y_i\), local equations of
\(s_1,\ldots,s_i\) form a regular sequence in \(\mathcal O_{W,w}\);
\item the scheme \(Y_i\) is regular, and
\(\dim\mathcal O_{Y_i,w}=\dim\mathcal O_{W,w}-i\) for every closed
point \(w\) of \(Y_i\).
\end{enumerate}
\end{sublemma}

\begin{proof}
The scheme \(Y_i\) is projective over \(S\). Hence its closed points
lie over closed points of \(S\), and every point of \(Y_i\)
specializes to a closed point of \(Y_i\).

Let \(w\) be a closed point of \(Y_i\) with image \(s\). The
extension \(\kappa(w)/\kappa(s)\) is finite, so there is a closed
point \(\bar w\) of \((Y_i)_{\bar s}\) over \(w\). Let
\(a_1,\ldots,a_i\in\mathfrak m_w\) be local equations of
\(s_1,\ldots,s_i\) in a trivialization of \(\mathcal L\) near
\(w\). The local homomorphism
\(\mathcal O_{W,w}\to\mathcal O_{W_{\bar s},\bar w}\) induces a
\(\kappa(w)\)-linear map
\[
 \mathfrak m_w/\mathfrak m_w^2
 \longrightarrow \mathfrak m_{\bar w}/\mathfrak m_{\bar w}^2 ,
\]
which sends the classes of \(a_1,\ldots,a_i\) to the classes of
local equations of \(s_1,\ldots,s_i\) at \(\bar w\). The latter are
linearly independent by hypothesis, so the classes of
\(a_1,\ldots,a_i\) in \(\mathfrak m_w/\mathfrak m_w^2\) are linearly
independent over \(\kappa(w)\). Put \(D=\dim\mathcal O_{W,w}\), and
complete them to a basis by the classes of
\(a_{i+1},\ldots,a_D\in\mathfrak m_w\). By Nakayama's lemma
\cite[\stackstag{00DV}]{Stacks}, the elements \(a_1,\ldots,a_D\)
generate \(\mathfrak m_w\); they form a regular system of parameters
of the regular local ring \(\mathcal O_{W,w}\). Hence
\(a_1,\ldots,a_i\) is a regular sequence, and
\(\mathcal O_{Y_i,w}=\mathcal O_{W,w}/(a_1,\ldots,a_i)\) is a
regular local ring of dimension \(D-i\)
\cite[\stackstag{00NQ}]{Stacks}.

Every point of \(Y_i\) specializes to a closed point of \(Y_i\).
At such a generalization, localization preserves the injectivity of
multiplication by each successive equation; the final quotient is
nonzero because the point belongs to \(Y_i\). Thus the localized
sequence is regular, and localizations of regular local rings are regular
\cite[\stackstag{0AFS}]{Stacks}. This proves \textup{(i)} and
\textup{(ii)} at all points.
\end{proof}
\subsection{Blowing up a regular centre and computing the degree}

\begin{notation}[The blowup of a regular centre]
\label{not:multisection-blowup-setup}
With \(S\) and \(d_s\) as in
\cref{not:multisection-dimensions-degrees}, let \(c\geq1\) be a
natural number, and let \(P\) be a regular flat projective
\(S\)-scheme whose fibres at all closed points of \(S\) are pure of
dimension \(c\). Let \(C\subseteq P\) be a regular closed subscheme,
with ideal sheaf \(\mathcal I_C\), such that
\(\dim\mathcal O_{C,y}=d_s\) for every closed point \(s\) of \(S\)
and every closed point \(y\) of \(C\) over \(s\). Let
\(\pi:W\to P\) be the blowup of \(P\) along \(C\), and let
\(E=\pi^{-1}(C)\). For an invertible sheaf \(\mathcal A\) on \(P\) and
\(b\in\mathbf N\), put
\[
 \mathcal L_b(\mathcal A)=\pi^*\mathcal A^{\otimes b}\otimes\mathcal O_W(-E).
\]
Projective bundles parametrize invertible quotients.
\end{notation}

\begin{sublemma}[The regular centre]
\label{sub:multisection-regular-centre}
In the situation of \cref{not:multisection-blowup-setup}, let \(y\)
be a closed point of \(P\), and let \(s\) be its image in \(S\).
\begin{enumerate}
\item The point \(s\) is closed, and
\(\dim\mathcal O_{P,y}=d_s+c\).
\item If \(y\in C\), there is a regular system of parameters
\(x_1,\ldots,x_{c+d_s}\) of \(\mathcal O_{P,y}\) such that
\(\mathcal I_{C,y}=(x_1,\ldots,x_c)\).
\item The immersion \(C\hookrightarrow P\) is regular of codimension
\(c\), and \(\mathcal I_C/\mathcal I_C^2\) is a locally free
\(\mathcal O_C\)-module of rank \(c\).
\end{enumerate}
\end{sublemma}

\begin{proof}
\textup{(i)} The morphism \(P\to S\) is proper, so \(s\) is closed,
and \(y\) is a closed point of the fibre \(P_s\), which is pure of
dimension \(c\). The dimension formula for the flat local
homomorphism \(\mathcal O_{S,s}\to\mathcal O_{P,y}\) gives
\(\dim\mathcal O_{P,y}=d_s+c\) \cite[\stackstag{00ON}]{Stacks}.

\textup{(ii)} We follow the argument of
\cite[\stackstag{00NR}]{Stacks}. The surjection
\(\mathcal O_{P,y}\to\mathcal O_{C,y}\) of regular local rings of
dimensions \(c+d_s\) and \(d_s\) induces a surjection of cotangent
spaces over \(\kappa(y)\), with kernel of dimension \(c\). Choose
\(x_1,\ldots,x_c\in\mathcal I_{C,y}\) whose classes form a basis of
this kernel, and lifts \(x_{c+1},\ldots,x_{c+d_s}\) of a basis of the
cotangent space of \(\mathcal O_{C,y}\). These \(c+d_s\) elements form
a regular system of parameters of \(\mathcal O_{P,y}\). The ring
\(\mathcal O_{P,y}/(x_1,\ldots,x_c)\) is a regular local domain of
dimension \(d_s\) \cite[\stackstag{00NQ}]{Stacks}, and it surjects onto
\(\mathcal O_{C,y}\), which also has dimension \(d_s\). The quotient of
a local domain of dimension \(d_s\) by a nonzero ideal has dimension
at most \(d_s-1\), because every chain of prime ideals containing a
nonzero ideal can be extended by the prime ideal \((0)\). Hence the
surjection is an isomorphism, and
\(\mathcal I_{C,y}=(x_1,\ldots,x_c)\).

\textup{(iii)} By \textup{(ii)}, the ideal \(\mathcal I_C\) is
generated by a regular sequence of length \(c\) at every closed point
of \(C\). Every point of \(C\) specializes to a closed point of
\(C\), and a localization of a regular sequence is a regular
sequence. Hence \(C\hookrightarrow P\) is a regular immersion of
codimension \(c\), and its conormal sheaf is locally free of rank
\(c\) \cite[\stackstag{063M}]{Stacks}.
\end{proof}

\begin{sublemma}[The blowup of the regular centre]
\label{sub:multisection-regular-ambient-blowup}
In the situation of \cref{not:multisection-blowup-setup}, the
following assertions hold.
\begin{enumerate}
\item The closed subscheme \(E\subseteq W\) is an effective
Cartier divisor, the projection \(E\to C\) identifies \(E\) with
\(\mathbb P_C(\mathcal I_C/\mathcal I_C^2)\), and under this
identification
\[
 \mathcal O_W(E)|_E\simeq\mathcal O_E(-1).
\]
\item The scheme \(W\) is regular, and
\(\dim\mathcal O_{W,w}=d_s+c\) for every closed point \(w\) of
\(W\) with image \(s\) in \(S\).
\end{enumerate}
\end{sublemma}

\begin{proof}
\textup{(i)} The inverse image \(E\) is an effective Cartier divisor,
and \(\mathcal I_C\mathcal O_W=\mathcal O_W(-E)\) is the
tautological invertible sheaf \(\mathcal O_W(1)\) of the blowup
\cite[\stackstag{02OS}]{Stacks}. The associated graded algebra of
the ideal of a regular immersion is the symmetric algebra of its
conormal sheaf
\cite[\stackstag{063K}, \stackstag{063M}]{Stacks}. With
\cref{sub:multisection-regular-centre}\textup{(iii)}, this gives
\[
 E=\operatorname{Proj}_C
       \Bigl(\bigoplus_{j\geq0}\mathcal I_C^j/\mathcal I_C^{j+1}\Bigr)
 \simeq \mathbb P_C(\mathcal I_C/\mathcal I_C^2),
 \qquad
 \mathcal O_W(-E)|_E\simeq\mathcal O_E(1).
\]

\textup{(ii)} Let \(w\) be a closed point of \(W\). Its image
\(y=\pi(w)\) is closed because \(\pi\) is proper. Let \(s\) be the
image of \(y\) in \(S\). If \(y\notin C\),
then \(\mathcal O_{W,w}=\mathcal O_{P,y}\), and the assertion is
\cref{sub:multisection-regular-centre}\textup{(i)}. Suppose that
\(y\in C\). Put \(B=\mathcal O_{P,y}\) and let
\(x_1,\ldots,x_{c+d_s}\) be as in
\cref{sub:multisection-regular-centre}\textup{(ii)}. Blowing up
commutes with the flat base change \(\Spec B\to P\)
\cite[\stackstag{0805}]{Stacks}, and the chart of the blowup of
\(\Spec B\) along \((x_1,\ldots,x_c)\) corresponding to \(x_1\) is
\begin{equation}\label{eq:blowup-chart}
 \Spec B[t_2,\ldots,t_c]/(x_2-x_1t_2,\ldots,x_c-x_1t_c);
\end{equation}
the other charts have the analogous presentations
\cite[\stackstag{0BIQ}]{Stacks}. Let \(\mathfrak n\) be a maximal
ideal of \(B[t_2,\ldots,t_c]\) over the maximal ideal of \(B\). The
local ring \(B[t_2,\ldots,t_c]_{\mathfrak n}\) is regular of
dimension \(d_s+2c-1\), by smoothness over \(B\) and the dimension
formula \cite[\stackstag{00ON}]{Stacks}. The images of
\(x_1,\ldots,x_{c+d_s}\), together with lifts of \(c-1\) parameters
of the local ring of the fibre
\(\kappa(y)[t_2,\ldots,t_c]\) at \(\mathfrak n\), generate
\(\mathfrak n\). Their number is \(d_s+2c-1\), so their classes form a
basis of the cotangent space. In this basis, the classes of
\(x_i-x_1t_i\), \(2\leq i\leq c\), are linearly independent: their
coefficients on the classes of \(x_2,\ldots,x_c\) form the identity
matrix. Hence the local ring of \eqref{eq:blowup-chart} at
\(\mathfrak n\) is regular of dimension \(d_s+c\)
\cite[\stackstag{00NQ}]{Stacks}. For \(c=1\) there are no variables
\(t_i\) and no equations, and \(W=P\). Every point of \(W\)
specializes to a closed point, and a localization of a regular local
ring is regular \cite[\stackstag{0AFS}]{Stacks}, so \(W\) is regular.
\end{proof}

\begin{sublemma}[A very ample sheaf on the blowup]
\label{sub:multisection-blowup-very-ample}
In the situation of \cref{not:multisection-blowup-setup}, let \(\mathcal A\)
be a very ample invertible sheaf for \(P\to S\). There is
\(b_1\in\mathbf N\) such that, for every integer \(b\geq b_1\), the
sheaf \(\mathcal L_b(\mathcal A)\) is very ample for \(W\to S\). In particular, \(W\) is projective over \(S\).
\end{sublemma}

\begin{proof}
We twist the argument of \cite[\stackstag{02NS}]{Stacks} by
\(\mathcal A^{\otimes(b-1)}\).
The sheaf \(\mathcal A\) is ample and \(P\) is noetherian, so there is
\(b_1\geq1\) such that \(\mathcal I_C\otimes \mathcal A^{\otimes(b-1)}\) is
globally generated for every \(b\geq b_1\)
\cite[\stackstag{01Q3}]{Stacks}. Fix such a \(b\). Since \(P\) is
quasi-compact, there are \(h\in\mathbf N\) and a surjection
\(\mathcal O_P^{\oplus(h+1)}\to
\mathcal I_C\otimes \mathcal A^{\otimes(b-1)}\). It induces a surjection of
graded \(\mathcal O_P\)-algebras
\[
 \operatorname{Sym}\bigl(\mathcal O_P^{\oplus(h+1)}\bigr)
 \longrightarrow
 \bigoplus_{j\geq0}\mathcal I_C^j\otimes \mathcal A^{\otimes(b-1)j}.
\]
Twisting the piece of degree \(j\) of a graded algebra by the
\(j\)-th power of an invertible sheaf does not change its relative
\(\operatorname{Proj}\), and twists the sheaf \(\mathcal O(1)\) by
the pullback of the invertible sheaf. Hence the relative
\(\operatorname{Proj}\) of the target is \(W\), with
\(\mathcal O(1)=\mathcal O_W(-E)\otimes\pi^*(\mathcal A^{\otimes(b-1)})\),
and the surjection defines a closed immersion
\(W\hookrightarrow P\times_S\mathbb P_S^h\) under which the
hyperplane bundle of the second factor pulls back to this sheaf.
The sheaf
\(\operatorname{pr}_1^*\mathcal A\otimes\operatorname{pr}_2^*\mathcal O(1)\)
on \(P\times_S\mathbb P_S^h\) is very ample over \(S\), by the Segre
embedding composed with the embedding of \(P\) defined by \(\mathcal A\). Its
restriction to \(W\) is \(\mathcal L_b(\mathcal A)\).
\end{proof}
\begin{proposition}[The degree of the generic intersection]
\label{prop:multisection-generic-chow-degree}
In the situation of \cref{not:multisection-blowup-setup}, let \(\mathcal A\)
be an invertible sheaf on \(P\), let \(b\in\mathbf N\), and let
\(F\) be the residue field of a generic point of \(S\). Assume that
\(P_F\) is pure of dimension \(c\) and
that \(C_F\) is finite over \(F\). Let
\(s_1,\ldots,s_c\in\Gamma(W,\mathcal L_b(\mathcal A))\) be sections such that, for
\(1\leq i\leq c\), the common zero scheme of \(s_1,\ldots,s_i\) is
regularly immersed in \(W\) in codimension \(i\). Let
\(Y\) be the common zero scheme of \(s_1,\ldots,s_c\). Then \(Y_F\) is finite
over \(F\), and its fundamental cycle satisfies
\begin{equation}\label{eq:multisection-generic-degree}
 \pi_{F*}[Y_F]
   = b^{c} c_1(\mathcal A_F)^c\cap[P_F] - [C_F]
       \quad \text{in } \operatorname{CH}_0(P_F).
\end{equation}
Here the cycles of \(Y_F\) and \(C_F\) are viewed on \(W_F\)
and \(P_F\), respectively, by their closed immersions.
In particular, if \(\ell\mid b\), then
\(\deg_F(Y_F)\equiv-\deg_F(C_F)\pmod{\ell}\).
\end{proposition}

\begin{proof}
\emph{The generic fibre.}
The morphism \(\Spec F\to S\) is flat, so
\(W_F\to P_F\) is the blowup of \(P_F\) along \(C_F\), with
exceptional divisor \(E_F\) \cite[\stackstag{0805}]{Stacks}. The
scheme \(P_F\) is regular and pure of dimension \(c\), and \(C_F\)
is a finite regular \(F\)-scheme, that is, a finite disjoint union
of spectra of fields; in particular \(C_F\) contains no irreducible
component of \(P_F\), because \(c\geq1\). Hence \(\pi_F\) is an
isomorphism over the dense open subscheme \(P_F\setminus C_F\),
whose inverse image is the complement of the effective Cartier divisor
\(E_F\), hence dense in \(W_F\)
\cite[\stackstag{02OS}, \stackstag{07ZU}]{Stacks}; so \(W_F\) is pure of dimension
\(c\) and
\begin{equation}\label{eq:blowup-birational}
 \pi_{F*}[W_F]=[P_F].
\end{equation}

\smallskip\noindent
\emph{The class of \(Y_F\).}
For \(0\leq i\leq c\), let \(Y_{i,F}\subseteq W_F\) be the common
zero scheme of \(s_1,\ldots,s_i\), so that \(Y_{0,F}=W_F\) and
\(Y_{c,F}=Y_F\). By hypothesis and flat base change, \(Y_{i,F}\)
is empty or is regularly immersed in \(W_F\) in codimension \(i\),
so it is pure of dimension \(c-i\) when nonempty. For \(0\leq i<c\)
and a point of \(Y_{i+1,F}\), choose a trivialization of \((\mathcal L_b(\mathcal A))_F\)
there. In the regular local ring of \(W_F\)
(regular by \cref{sub:multisection-regular-ambient-blowup}\textup{(ii)}),
the ideal of the first \(i+1\) equations has quotient of dimension
\(i+1\) less than that of the ring, by the codimension hypothesis
and catenarity. This ring is Cohen--Macaulay
\cite[\stackstag{00NQ}]{Stacks}; by
\cite[\stackstag{02JN}]{Stacks} these equations form a regular
sequence in the given order. Hence \(s_{i+1}\) is a nonzerodivisor
on \(Y_{i,F}\) at the points where it vanishes, and elsewhere it is
locally a unit. Thus it defines the effective Cartier divisor
\(Y_{i+1,F}\) on \(Y_{i,F}\), including the case of an empty zero
scheme. In particular \(Y_F\) is
projective of dimension at most zero, hence finite over \(F\). The
formula for the first Chern class of the invertible sheaf of an
effective Cartier divisor
\cite[\stackstag{02T9}, \stackstag{02SU}]{Stacks}, applied \(c\)
times, gives
\begin{equation}\label{eq:chern-class-of-zero-scheme}
 [Y_F]=c_1((\mathcal L_b(\mathcal A))_F)^c\cap[W_F]
            \quad\text{in }\operatorname{CH}_0(W_F).
\end{equation}

\smallskip\noindent
\emph{Expansion.}
First Chern classes are additive and commute
\cite[\stackstag{02SP}, \stackstag{02TJ}]{Stacks}, so
\[
 c_1((\mathcal L_b(\mathcal A))_F)
 =b\,c_1(\pi_F^*\mathcal A_F)-c_1(\mathcal O_{W_F}(E_F)).
\]
The binomial formula gives
\begin{equation}\label{eq:chern-expansion}
\begin{aligned}
 c_1((\mathcal L_b(\mathcal A))_F)^c\cap[W_F]
 &=\sum_{i=0}^{c}\binom ci(-1)^ib^{c-i}\,\cdot\\
 &\qquad c_1(\pi_F^*\mathcal A_F)^{c-i}
    c_1(\mathcal O_{W_F}(E_F))^{i}\cap[W_F].
\end{aligned}
\end{equation}
By the projection formula and \eqref{eq:blowup-birational}, the
term \(i=0\) has image
\[
 \pi_{F*}\bigl(b^cc_1(\pi_F^*\mathcal A_F)^c\cap[W_F]\bigr)
 =b^c\,c_1(\mathcal A_F)^c\cap[P_F].
\]
Let \(1\leq i\leq c\), and let \(j:E_F\hookrightarrow W_F\) be the
inclusion. The identity for an effective Cartier divisor, the
projection formula, and
\(\mathcal O_{W_F}(E_F)|_{E_F}\simeq\mathcal O_{E_F}(-1)\)
(\cref{sub:multisection-regular-ambient-blowup}\textup{(i)}) give
\[
 c_1(\mathcal O_{W_F}(E_F))^{i}\cap[W_F]
 =j_*\Bigl(c_1\bigl(\mathcal O_{E_F}(-1)\bigr)^{i-1}\cap[E_F]\Bigr)
 \quad\text{in }\operatorname{CH}_{c-i}(W_F).
\]
The morphism \(\pi_F\circ j\) factors through the finite
\(F\)-scheme \(C_F\), and \(\operatorname{CH}_{c-i}(C_F)=0\) for
\(i<c\). With the projection formula for \(c_1(\pi_F^*\mathcal A_F)^{c-i}\), the terms
\(1\leq i<c\) of \eqref{eq:chern-expansion} have image zero in
\(\operatorname{CH}_0(P_F)\). For \(i=c\), the scheme \(E_F\) is the
projective bundle \(\mathbb P_{C_F}((\mathcal I_C/\mathcal I_C^2)|_{C_F})\)
of relative dimension \(c-1\) over \(C_F\), and for every point
\(z\) of \(C_F\) the class
\(c_1(\mathcal O(1))^{c-1}\cap[\mathbb P^{c-1}_{\kappa(z)}]\) has
degree one over \(\kappa(z)\). Hence
\[
 \pi_{F*}\bigl(c_1(\mathcal O_{W_F}(E_F))^{c}\cap[W_F]\bigr)=(-1)^{c-1}[C_F],
\]
and the term \(i=c\) of \eqref{eq:chern-expansion}, whose
coefficient is \((-1)^c\), has image \(-[C_F]\). Together with
\eqref{eq:chern-class-of-zero-scheme}, this proves
\eqref{eq:multisection-generic-degree}.

\smallskip\noindent
\emph{The congruence.}
The degree map \(\operatorname{CH}_0(P_F)\to\Z\) is the proper
pushforward to \(\Spec F\); it sends \(\pi_{F*}[Y_F]\) to
\(\deg_F(Y_F)\) and \([C_F]\) to \(\deg_F(C_F)\). If \(\ell\mid b\),
the first term on the right of
\eqref{eq:multisection-generic-degree} has degree divisible by
\(\ell\), because \(c\geq1\).
\end{proof}
\subsection{Replacing the projective scheme}

\begin{notation}[The two closed-fibre dimensions]
\label{not:fibre-invariants}
Assume \(R\ne0\). Let \(X\) be a projective \(S\)-scheme,
let \(U\subseteq X\) be open, and give \(D=X\setminus U\)
its reduced closed subscheme structure. For a projective
\(S\)-scheme \(C\) equipped with an \(S\)-morphism \(C\to X\),
define
\[
 f(C)=\max_{s\in S\text{ closed}}\dim C_{\bar s},\qquad
 g(C)=\max_{s\in S\text{ closed}}\dim(C\times_XD)_{\bar s}.
\]
The maxima are over a nonempty finite set, and an empty fibre
has dimension \(-1\). After a nonzero finite \'etale base change
\(R'/R\), the same notation uses the closed points of
\(\operatorname{Spec}R'\) and the base changes \(X_{R'}\),
\(U_{R'}\) and \(D_{R'}\).
\end{notation}

\begin{construction}[The blowup and cutting replacement]
\label{con:multisection-replacement}
Use the base \(S=\operatorname{Spec}R\), the prime \(\ell\), and the
integers \(d_s\) of \cref{not:multisection-dimensions-degrees},
with \(R\ne0\). Let \(r\in\mathbf N\), and let \(X\) be a regular,
flat, projective \(S\)-scheme whose geometric fibres over the closed
points and over the generic points of \(S\) are pure of dimension
\(r\). Let \(U\subseteq X\) be an open subscheme that meets every
irreducible component of \(X_{\bar s}\) for every closed point \(s\)
of \(S\). Let \(C\) be a nonempty regular projective \(S\)-scheme
with an \(S\)-morphism \(C\to X\), such that
\(\dim\mathcal O_{C,y}=d_s\) for every closed point \(y\) of \(C\)
over \(s\). Use \(D\), \(f(C)\), and \(g(C)\) as in
\cref{not:fibre-invariants}.

\begin{enumerate}
\item \emph{The ambient scheme.} Choose a closed immersion
\(C\hookrightarrow\mathbb P^a_S\) with \(a\geq1\), and put
\(c=a+r\geq1\). The morphism
\(C\to\mathbb P^a_X=\mathbb P^a_S\times_S X\) with components the
chosen immersion and \(C\to X\) is a closed immersion, because \(X\)
is separated over \(S\). The scheme \(\mathbb P^a_X\) is smooth over
\(X\), hence regular; it is flat and projective over \(S\), and its
geometric fibres over the closed and the generic points of \(S\) are
pure of dimension \(c\). Thus \(\mathbb P^a_X\) and \(C\) satisfy the
hypotheses of \cref{not:multisection-blowup-setup}.

\item \emph{The blowup.} Let \(\pi:W\to\mathbb P^a_X\) be the blowup along
\(C\), with exceptional divisor \(E\), and let
\(W_D=W\times_XD\) for the composite \(W\to\mathbb P^a_X\to X\).
Choose a very ample invertible sheaf \(\mathcal A\) for
\(\mathbb P^a_X\to S\) and a positive multiple \(b\) of \(\ell\) such
that the sheaf \(\mathcal L_b(\mathcal A)\) of
\cref{not:multisection-blowup-setup} is very ample for \(W\to S\)
(\cref{sub:multisection-blowup-very-ample}), and a closed immersion
\(W\hookrightarrow\mathbb P^m_S\) with
\(\mathcal O(1)|_W\simeq \mathcal L_b(\mathcal A)\).

\item \emph{The sections.} For each closed point \(s\) of \(S\),
let \(V_s\subseteq\mathbb A^{c(m+1)}_{\kappa(s)}\) be the nonempty
open subset given by \cref{lem:multisection-geometric-incidence} for
\(Z=W_s\subseteq\mathbb P^m_{\kappa(s)}\), the closed subschemes
\(Z_1=W_s\) and \(Z_2=(W_D)_s\), and \(c\) linear forms. By
\cref{sub:multisection-simultaneous-parameter-choice}, there are a
finite \'etale \(R\)-algebra \(R'\), free of rank prime to \(\ell\)
as an \(R\)-module, and a \(c\)-tuple \(h=(h_1,\ldots,h_c)\) of
linear forms with coefficients in \(R'\) whose reduction at every
closed point \(s'\) of \(\operatorname{Spec}R'\) over \(s\) lies in
\(V_s(\kappa(s'))\). Let \(s_i=h_i|_{W_{R'}}\in
\Gamma(W_{R'},(\mathcal L_b(\mathcal A))_{R'})\).

\item \emph{The replacement.} Let \(Y\subseteq W_{R'}\) be the
common zero scheme of \(s_1,\ldots,s_c\), with the morphism
\[
 Y\hookrightarrow W_{R'}\xrightarrow{\ \pi_{R'}\ }\mathbb P^a_{X_{R'}}
       \longrightarrow X_{R'}.
\]
\end{enumerate}
\end{construction}

\begin{proposition}[The dimensions and degree of the replacement]
\label{prop:multisection-decrease-closed-fibres}
Let \(R'\) and \(Y\to X_{R'}\) be obtained from
\cref{con:multisection-replacement}.
\begin{enumerate}
\item The scheme \(Y\) is regular and projective over \(R'\), and
\(\dim\mathcal O_{Y,y}=\dim\mathcal O_{\operatorname{Spec}R',s'}\)
for every closed point \(y\) of \(Y\) with image \(s'\).
\item The two dimensions of \cref{not:fibre-invariants} satisfy
\begin{equation}\label{eq:multisection-dimension-decrease}
 f(Y)\leq\max\{0,f(C)-1\},\qquad
 g(Y)\leq\max\{-1,g(C)-1\}.
\end{equation}
\item If the fibres of \(C\) over the generic points of \(S\) are
finite, then the fibres of \(Y\) over the generic points of
\(\operatorname{Spec}R'\) are finite, and for every generic point of
\(\operatorname{Spec}R'\), with residue field \(F'\),
\[
 \deg_{F'}(Y_{F'})\equiv-\deg_{F'}(C_{F'})\pmod{\ell}.
\]
\end{enumerate}
\end{proposition}

\begin{proof}
\emph{The closed fibres of \(W\) and \(W_D\).}
Let \(s\) be a closed point of \(S\). We use the decomposition of
\(W_{\bar s}\) into \(E_{\bar s}\) and its complement. By the
definition of \(E\), the square
\[
\begin{tikzcd}
 E_{\bar s}\arrow[r,hook]\arrow[d] & W_{\bar s}\arrow[d,"\pi_{\bar s}"]\\
 C_{\bar s}\arrow[r,hook] & \mathbb P^a_{X_{\bar s}}
\end{tikzcd}
\]
is cartesian; we do not identify \(W_{\bar s}\) with the blowup of
\(\mathbb P^a_{X_{\bar s}}\) along \(C_{\bar s}\). The scheme
\(W_{\bar s}\setminus E_{\bar s}\) is isomorphic to
\(\mathbb P^a_{X_{\bar s}}\setminus C_{\bar s}\), of dimension at most \(c\),
and \(E_{\bar s}\) is a projective bundle of relative dimension
\(c-1\) over \(C_{\bar s}\)
(\cref{sub:multisection-regular-ambient-blowup}\textup{(i)}). Hence
\begin{equation}\label{eq:blowup-fibre-dimension}
 \dim W_{\bar s}\leq\max\{c,\;f(C)+c-1\}.
\end{equation}
The scheme \((W_D)_{\bar s}\setminus E_{\bar s}\) is an open
subscheme of \(\mathbb P^a_{D_{\bar s}}\), and
\(\dim D_{\bar s}\leq r-1\), because \(X_{\bar s}\) is pure of
dimension \(r\) and \(U_{\bar s}\) meets each of its irreducible
components. The scheme \((W_D)_{\bar s}\cap E_{\bar s}\) is the
restriction of the projective bundle \(E_{\bar s}\) to
\((C\times_XD)_{\bar s}\); it is empty if \(g(C)=-1\). Hence
\begin{equation}\label{eq:boundary-fibre-dimension}
 \dim(W_D)_{\bar s}\leq\max\{c-1,\;g(C)+c-1\}.
\end{equation}
The fibres of \(W_{R'}\) and \((W_D)_{R'}\) over a closed point
\(s'\) of \(\operatorname{Spec}R'\) over \(s\) are obtained from
\(W_s\) and \((W_D)_s\) by the extension
\(\kappa(s')/\kappa(s)\), so the same bounds hold for them
\cite[\stackstag{05F8}]{Stacks}.

\smallskip\noindent
\emph{Proof of \textup{(i)}.}
The morphism \(\operatorname{Spec}R'\to S\) is finite \'etale.
Hence \(W_{R'}\) is regular, it is the blowup of \(\mathbb P^a_{X_{R'}}\) along
\(C_{R'}\) \cite[\stackstag{0805}]{Stacks}, and
\(\dim\mathcal O_{W_{R'},w}=d_s+c\) for every closed point \(w\) over
a closed point \(s\) of \(S\)
(\cref{sub:multisection-regular-ambient-blowup}\textup{(ii)}). By
the choice of the sections and
\cref{lem:multisection-geometric-incidence}\textup{(i)}, the
hypothesis of \cref{sub:multisection-absolute-cotangent} holds for
\(W_{R'}\), \((\mathcal L_b(\mathcal A))_{R'}\) and \(s_1,\ldots,s_c\). Hence \(Y\) is
regular, and for a closed point \(y\) of \(Y\) over \(s'\) and \(s\),
\[
 \dim\mathcal O_{Y,y}=d_s+c-c=d_s
 =\dim\mathcal O_{\operatorname{Spec}R',s'},
\]
the last equality because \(\operatorname{Spec}R'\to S\) is
\'etale. The scheme \(Y\) is closed in \(W_{R'}\), hence projective
over \(R'\).

\smallskip\noindent
\emph{Proof of \textup{(ii)}.}
Let \(s'\) be a closed point of \(\operatorname{Spec}R'\) over
\(s\). By construction,
\(Y_{s'}=(W_s)_{\kappa(s')}\cap V_+(h)\) and
\((Y\times_XD)_{s'}=((W_D)_s)_{\kappa(s')}\cap V_+(h)\), where
\(h\in V_s(\kappa(s'))\). By
\cref{lem:multisection-geometric-incidence}\textup{(ii)}, each of
these schemes is empty if the dimension of \(W_s\), respectively
\((W_D)_s\), is less than \(c\), and has dimension at most that
dimension minus \(c\) otherwise. With
\eqref{eq:blowup-fibre-dimension} and
\eqref{eq:boundary-fibre-dimension}, this gives
\[
 \dim Y_{\bar s'}\leq\max\{0,f(C)-1\},\qquad
 \dim(Y\times_XD)_{\bar s'}\leq\max\{-1,g(C)-1\},
\]
that is, \eqref{eq:multisection-dimension-decrease}.

\smallskip\noindent
\emph{Proof of \textup{(iii)}.}
Let \(\eta'\) be a generic point of \(\operatorname{Spec}R'\), with
residue field \(F'\). By \textup{(i)} and
\cref{sub:multisection-absolute-cotangent}\textup{(i)}, the common
zero scheme of \(s_1,\ldots,s_i\) is a regular closed immersion of
codimension \(i\) in \(W_{R'}\), for \(1\leq i\leq c\). The scheme
\(\mathbb P^a_{X_{F'}}\) is pure of dimension \(c\) and \(C_{F'}\) is finite over \(F'\).
Moreover, the hypotheses on the closed local dimensions of \(C\) in
\cref{not:multisection-blowup-setup} are preserved by the finite
\'etale base change (use the dimension formula
\cite[\stackstag{00ON}]{Stacks}). Thus
\cref{prop:multisection-generic-chow-degree} applies over \(R'\) with
\(\ell\mid b\); it gives the finiteness of \(Y_{F'}\) and the
congruence.
\end{proof}

\begin{remark}[The dimension reduction over a surface]
Suppose the base is local of dimension two and $C$ is integral of
dimension two, dominant and generically finite of degree prime to
$\ell$ over the base. If the closed fibre of $C$ is positive-dimensional,
its dimension is at most
one, since it is a proper closed subset of the integral surface.
In this case $f(C)=1$ and $g(C)\leq1$. The first replacement of
\cref{prop:multisection-decrease-closed-fibres} makes the closed fibres
zero-dimensional, with $g(Y)\leq0$. If a point still maps to the
boundary, a second replacement makes $g=-1$ without increasing the
closed-fibre dimension. Throughout, the generic degree changes sign
modulo $\ell$ and stays prime to $\ell$. This example illustrates why
two dimension estimates, rather than just quasi-finiteness, are needed.
\end{remark}

\begin{sublemma}[Closed local dimensions of the initial source]
\label{sub:multisection-closed-local-dimensions}
Assume that \(R\) is nonzero and connected, and let \(F\)
be its fraction field. Let \(C\) be a regular projective
\(S\)-scheme with finite generic fibre \(C_F\).
\begin{enumerate}
\item The irreducible components of \(C\) are open and closed. The
union \(C_0\) of those that dominate \(S\) satisfies
\((C_0)_F=C_F\).
\item For every closed point \(y\) of \(C_0\), with image \(s\) in
\(S\), we have \(\dim\mathcal O_{C,y}=d_s\).
\end{enumerate}
\end{sublemma}

\begin{proof}
A noetherian regular scheme is the disjoint union of its
irreducible components, which are integral
\cite[\stackstag{033M}]{Stacks}. This applies to \(S\), so \(R\) is
a domain, and to \(C\). A component of \(C\) that does not dominate
\(S\) has empty generic fibre. This proves \textup{(i)}.

Let \(C'\) be a component of \(C\) dominating \(S\), with function
field \(\kappa(C')\), and let \(y\) be a closed point of \(C'\) with
image \(s\). The extension \(\kappa(C')/F\) is finite because
\(C_F\) is finite. The point \(s\) is closed because \(C\to S\) is
proper, and \(\kappa(y)/\kappa(s)\) is finite because \(y\) is a
closed point of the fibre \(C_s\). The regular ring \(R\) is
universally catenary \cite[\stackstag{00NM}]{Stacks}, and the
dimension formula \cite[\stackstag{02JU}]{Stacks} gives
\[
 \dim\mathcal O_{C,y}
 =\dim\mathcal O_{S,s}+\operatorname{trdeg}_F\kappa(C')
   -\operatorname{trdeg}_{\kappa(s)}\kappa(y)=d_s .
\]
\end{proof}

\begin{sublemma}[Closed fibres and finite local freeness]
\label{sub:multisection-finite-flatness}
Let \(Y\) be a regular projective \(S\)-scheme such that
\(\dim Y_s\leq0\) for every closed point \(s\) of \(S\) and
\(\dim\mathcal O_{Y,y}=d_s\) for every closed point \(y\) of \(Y\)
over \(s\).
\begin{enumerate}
\item The scheme \(Y\) is finite locally free over \(S\).
\item Let \(X\), \(U\), \(D\) and a morphism \(Y\to X\) be as in
\cref{not:fibre-invariants}. If \(g(Y)=-1\), then \(Y\to X\) factors
through \(U\).
\end{enumerate}
\end{sublemma}

\begin{proof}
\textup{(i)} The set of points \(s\in S\) with \(\dim Y_s\geq1\) is
closed, because \(Y\to S\) is proper
\cite[\stackstag{0D4I}]{Stacks}. It contains no closed point of
\(S\), so it is empty. Hence \(Y\to S\) is proper with finite
fibres, and therefore finite
\cite[\stackstag{02LS}]{Stacks}. Let \(y\) be a closed point of
\(Y\) with image \(s\). The local homomorphism
\(\mathcal O_{S,s}\to\mathcal O_{Y,y}\) has a regular source and a
regular, hence Cohen--Macaulay, target, and
\[
 \dim\mathcal O_{Y,y}=d_s
 =\dim\mathcal O_{S,s}
  +\dim\mathcal O_{Y,y}/\mathfrak m_s\mathcal O_{Y,y},
\]
because the fibre \(Y_s\) is finite. By miracle flatness
\cite[\stackstag{00R4}]{Stacks}, the homomorphism is flat. Put
\(B=\Gamma(Y,\mathcal O_Y)\), so that \(Y=\Spec B\). The local
homomorphism \(R_{\mathfrak m}\to B_{\mathfrak n}\) is flat
for every maximal ideal \(\mathfrak n\) of \(B\), where
\(\mathfrak m=\mathfrak n\cap R\). By
\cite[\stackstag{00HT}\textup{(7)}]{Stacks}, the finite
\(R\)-module \(B\) is flat, hence finite locally free
\cite[\stackstag{00NX}]{Stacks}.

\textup{(ii)} The scheme \(Y\times_XD\) is proper over \(S\), so its
image in \(S\) is closed. If \(g(Y)=-1\), this image contains no
closed point of \(S\), so it is empty. Hence
\(Y\times_XD=\varnothing\), and \(Y\to X\) factors through
\(U=X\setminus D\).
\end{proof}

\begin{theorem}[A regular multisection of degree prime to \(\ell\)]
\label{thm:regular-multisections}
Let \(R\) be a nonzero connected noetherian regular semilocal
ring with fraction field \(F\), and let \(\ell\) be a prime number.
Let \(X\) be a regular, flat, projective \(R\)-scheme whose
geometric fibres over the closed points of \(\operatorname{Spec}R\)
and over its generic point are pure of dimension \(r\). Let
\(U\subseteq X\) be an open subscheme that meets every irreducible
component of every geometric fibre of \(X\) over a closed point.
\begin{enumerate}
\item Suppose that there are a regular projective \(R\)-scheme
\(C\) and an \(R\)-morphism \(C\to X\) such that \(C_F\) is finite
over \(F\) and \(\deg_F(C_F)\) is prime to \(\ell\). Then there are
a regular scheme \(T\), finite locally free over \(R\) of rank prime
to \(\ell\), and an \(R\)-morphism \(T\to U\).
\item Let \(x\) be a closed point of \(X_F\) with
\([\kappa(x):F]\) prime to \(\ell\), let \(Z\subseteq X\) be its
closure with the reduced structure, and let \(C\to Z\) be a
projective, surjective, generically finite morphism from a regular
integral scheme \(C\) with \([\kappa(C):\kappa(x)]\) prime to
\(\ell\). Then \(C\to Z\subseteq X\) satisfies the hypothesis of
\textup{(i)}.
\end{enumerate}
\end{theorem}

\begin{proof}
\textup{(ii)} The scheme \(Z\) is integral and projective over
\(R\), with generic fibre \(Z_F=\Spec\kappa(x)\). Hence \(C\) is
projective over \(R\), and \(C_F\) is an integral scheme of finite
type over \(F\) whose function field \(\kappa(C)\) is finite over
\(F\); thus \(C_F=\Spec\kappa(C)\) and
\[
 \deg_F(C_F)=[\kappa(C):\kappa(x)]\,[\kappa(x):F] .
\]

\textup{(i)} We use \cref{not:multisection-dimensions-degrees}
and \cref{not:fibre-invariants} with \(S=\operatorname{Spec}R\).

\smallskip\noindent
\emph{The initial source.}
Let \(C_0\subseteq C\) be the union of the irreducible components
of \(C\) that dominate \(S\). By
\cref{sub:multisection-closed-local-dimensions}, the scheme \(C_0\)
is regular and projective over \(R\), it satisfies
\((C_0)_F=C_F\), and \(\dim\mathcal O_{C_0,y}=d_s\) for every closed
point \(y\) over \(s\). It is nonempty because \(\deg_F(C_F)\neq0\).

\smallskip\noindent
\emph{Iteration.}
We construct, for every \(j\in\mathbf N\), a finite \'etale
\(R\)-algebra \(R_j\), free of rank prime to \(\ell\) as an
\(R\)-module, and a nonempty regular projective \(R_j\)-scheme
\(C_j\) with an \(R_j\)-morphism \(C_j\to X_{R_j}\), such that:
\begin{enumerate}
\item[\textup{(a)}]
\(\dim\mathcal O_{C_j,y}=\dim\mathcal O_{\operatorname{Spec}R_j,s}\)
for every closed point \(y\) of \(C_j\) with image \(s\);
\item[\textup{(b)}] for every generic point of
\(\operatorname{Spec}R_j\), with residue field \(F_j\), the scheme
\((C_j)_{F_j}\) is finite and
\begin{equation}\label{eq:multisection-iterated-degree}
 \deg_{F_j}\bigl((C_j)_{F_j}\bigr)
       \equiv(-1)^j\deg_F(C_F)\pmod{\ell};
\end{equation}
\item[\textup{(c)}] \(f(C_j)\leq\max\{0,f(C_0)-j\}\) and
\(g(C_j)\leq\max\{-1,g(C_0)-j\}\), where the maxima defining
\(f\) and \(g\) are taken over all closed points of
\(\operatorname{Spec}R_j\).
\end{enumerate}
For \(j=0\), take \(R_0=R\) and the scheme \(C_0\). Suppose that
\(R_j\) and \(C_j\) are constructed. The ring \(R_j\) is noetherian,
regular and semilocal, and \(X_{R_j}\), \(U_{R_j}\) satisfy the
hypotheses of \cref{con:multisection-replacement} over \(R_j\):
the scheme \(X_{R_j}\) is \'etale over \(X\), hence regular, and
its geometric fibres over \(\operatorname{Spec}R_j\) are geometric
fibres of \(X\) over \(\operatorname{Spec}R\). Apply
\cref{con:multisection-replacement} to \(C_j\to X_{R_j}\); it gives
a finite \'etale \(R_j\)-algebra \(R_{j+1}\), free of rank prime to
\(\ell\), and a scheme \(C_{j+1}=Y\) over \(R_{j+1}\). Properties
\textup{(a)} and \textup{(c)} for \(j+1\) are
\cref{prop:multisection-decrease-closed-fibres}\textup{(i)} and
\textup{(ii)}. Let \(F_{j+1}\) be the residue field of a generic
point of \(\operatorname{Spec}R_{j+1}\), lying over a generic point
of \(\operatorname{Spec}R_j\) with residue field \(F_j\). The degree
of a finite scheme over a field is invariant under extension of the
field, so
\cref{prop:multisection-decrease-closed-fibres}\textup{(iii)} gives
\[
 \deg_{F_{j+1}}\bigl((C_{j+1})_{F_{j+1}}\bigr)
 \equiv-\deg_{F_j}\bigl((C_j)_{F_j}\bigr)\pmod{\ell},
\]
which is \textup{(b)} for \(j+1\). The right-hand side of
\eqref{eq:multisection-iterated-degree} is nonzero modulo \(\ell\),
so \(C_{j+1}\) is nonempty.

\smallskip\noindent
\emph{Termination.}
Put \(n=\max\{f(C_0),g(C_0)+1\}\). By \textup{(b)}, the proper
morphism \(C_n\to\operatorname{Spec}R_n\) is surjective, so
\(f(C_n)\geq0\); with \textup{(c)}, \(f(C_n)=0\) and
\(g(C_n)=-1\). By \textup{(a)} and
\cref{sub:multisection-finite-flatness}, the scheme \(C_n\) is
finite locally free over \(R_n\), and \(C_n\to X_{R_n}\) factors
through \(U_{R_n}\).

\smallskip\noindent
\emph{The rank over \(R\).}
Put \(T=C_n\), with the composite morphisms
\(T\to\operatorname{Spec}R_n\to\operatorname{Spec}R\) and
\(T\to U_{R_n}\to U\). It is regular and finite locally free over
\(R\). Write \(R_n=\prod_iR_{n,i}\) with \(R_{n,i}\) connected, let
\(e_i\) be the rank of the finite locally free \(R\)-module
\(R_{n,i}\), and let \(b_i\) be the rank of \(T\) over \(R_{n,i}\),
which is the degree of its fibre over the generic point of
\(\operatorname{Spec}R_{n,i}\). By
\eqref{eq:multisection-iterated-degree},
\[
 \operatorname{rank}_R(T)=\sum_i e_i b_i
   \equiv(-1)^n\deg_F(C_F)\sum_i e_i
   =(-1)^n\deg_F(C_F)\operatorname{rank}_R(R_n)\pmod{\ell},
\]
and both factors are prime to \(\ell\).
\end{proof}
\subsection{The initial alteration over a field or over the integers}

\begin{lemma}[The projective regular source]
\label{lem:multisection-arithmetic-alteration}
Let \(\ell\) be a prime number, and let \(R\) be a regular semilocal
domain essentially of finite type over a field or over \(\Z\)
\textup{(\cref{def:arithmetic-semilocal})}, with fraction field
\(F\), such that \(\ell\in R^\times\). Use
\cref{not:multisection-dimensions-degrees} for
\(S=\operatorname{Spec}R\). For every projective \(R\)-scheme \(X\)
and every closed point \(x\in X_F\), there are a regular integral
projective \(R\)-scheme \(C\) and an \(R\)-morphism \(C\to X\) with
the following properties.
\begin{enumerate}
\item The morphism \(C\to X\) maps \(C\) onto the closure of \(x\)
in \(X\), the generic fibre \(C_F\) is the spectrum of the function
field \(\kappa(C)\) of \(C\), and
\([\kappa(C):\kappa(x)]\) is finite and prime to \(\ell\). In
particular, if \([\kappa(x):F]\) is prime to \(\ell\), then
\(\deg_F(C_F)\) is prime to \(\ell\).
\item For every closed point \(y\) of \(C\), with image \(s\) in
\(S\), we have \(\dim\mathcal O_{C,y}=d_s\).
\end{enumerate}
\end{lemma}

\begin{proof}
\emph{A presentation of \(R\).}
By \cref{sub:arithmetic-finite-etale}\textup{(i)}, there are a ring
\(B\) that is a field or \(\Z\), an integral finite-type
\(B\)-algebra \(A\) with \(B\to A\) injective, and pairwise
incomparable prime ideals \(\mathfrak p_1,\ldots,\mathfrak p_t\) of
\(A\) such that \(R\) is the semilocalization of \(A\) at
\(\mathfrak p_1,\ldots,\mathfrak p_t\). The fraction field of \(A\)
is \(F\). If \(B\) is a field, its characteristic is different from
\(\ell\), because \(\ell\in R^\times\). If \(B=\Z\), then
\(\ell\notin\mathfrak p_i\) for every \(i\), and
\[
 R=\Bigl(A[1/\ell]\setminus
       \bigcup_{i=1}^t\mathfrak p_iA[1/\ell]\Bigr)^{-1}
 A[1/\ell];
\]
the ring \(A[1/\ell]\) is an integral domain of finite type over
\(\Z[1/\ell]\) with fraction field \(F\). Replacing \(A\) by
\(A[1/\ell]\) in this case, we may assume that \(B\) is a field of
characteristic different from \(\ell\) or that \(B=\Z[1/\ell]\).
\smallskip\noindent
\emph{A model of the closure of \(x\).}
Choose a closed immersion \(X\hookrightarrow\mathbb P^n_R\), defined
by finitely many homogeneous polynomials. The ring \(R\) is the
filtered colimit of the rings \(A[1/a]\),
\(a\in A\setminus\bigcup_i\mathfrak p_i\). Hence there are such an
element \(a\) and a closed subscheme
\(X_0\subseteq\mathbb P^n_{A[1/a]}\) with
\(X_0\times_{A[1/a]}R=X\); in particular \((X_0)_F=X_F\). Let
\(Z_0\subseteq X_0\) be the closure of \(x\) with the reduced
structure. It is integral and projective over \(A[1/a]\), its
function field is \(\kappa(x)\), and
\((Z_0)_F=\operatorname{Spec}\kappa(x)\).

\smallskip\noindent
\emph{The alteration.}
If \(B\) is a field, Gabber's theorem
\cite[Expos\'e~X, Theorem~2.1, p.~185]{ILO14}, applied to \(Z_0\) with the
empty boundary, gives a projective \(\ell'\)-alteration
\(\widetilde Z\to Z_0\) whose source is smooth over a finite
extension of \(B\), hence regular. If \(B=\Z[1/\ell]\), the base
\(\operatorname{Spec}B\) is separated, integral, noetherian,
excellent, and regular of dimension one, and \(Z_0\) is separated
and of finite type over \(B\); it is flat over \(B\), because it is
integral and dominates \(\operatorname{Spec}B\). In this case
\cite[Expos\'e~X, Theorem~2.4, p.~190]{ILO14}, applied with the
empty boundary, gives a projective \(\ell'\)-alteration
\(\widetilde Z\to Z_0\) with regular source. By
definition \cite[Expos\'e~X, \S2, p.~185]{ILO14}, an \(\ell'\)-alteration
is proper, surjective, generically finite and maximally dominating
(it sends generic points of irreducible components to generic points
of irreducible components); all residue field degrees over maximal
points of the target are prime to \(\ell\). The irreducible components of the regular
noetherian scheme \(\widetilde Z\) are open and closed. Let \(C_0\)
be one whose generic point lies over the generic point of \(Z_0\)
with \([\kappa(C_0):\kappa(x)]\) prime to \(\ell\). Then
\(C_0\to Z_0\) is projective, surjective and generically finite.

\smallskip\noindent
\emph{The scheme \(C\).}
Put \(C=C_0\times_{A[1/a]}R\), with the composite
\(C\to Z_0\times_{A[1/a]}R\hookrightarrow X\). The scheme \(C\) is
a localization of \(C_0\), so it is regular and integral with
function field \(\kappa(C)=\kappa(C_0)\); it is projective over
\(R\), and it maps onto the closure of \(x\) in \(X\). The scheme
\(C_F\) is integral and of finite type over \(F\), and its function
field \(\kappa(C_0)\) is finite over \(\kappa(x)\), hence over
\(F\). Therefore \(C_F=\operatorname{Spec}\kappa(C_0)\) and
\[
 \deg_F(C_F)=[\kappa(C_0):\kappa(x)]\,[\kappa(x):F].
\]
This proves \textup{(i)}.

\smallskip\noindent
\emph{Local dimensions.}
The regular integral projective \(R\)-scheme \(C\) dominates \(S\)
and has finite generic fibre by \textup{(i)}. Thus
\cref{sub:multisection-closed-local-dimensions}\textup{(ii)} gives
\(\dim\mathcal O_{C,y}=d_s\) at every closed point \(y\) over \(s\).
\end{proof}

\begin{corollary}[Multisections over semilocal domains essentially of finite type]
\label{cor:arithmetic-multisections}
Let \(\ell\) be a prime number, and let \(R\) be a regular semilocal
domain essentially of finite type over a field or over \(\Z\), with
fraction field \(F\), such that \(\ell\in R^\times\). Let \(X\) be a
regular, flat, projective \(R\)-scheme whose geometric fibres over
the closed points of \(\operatorname{Spec}R\) and over its generic
point are pure of dimension \(r\), and let \(U\subseteq X\) be an
open subscheme that meets every irreducible component of every
geometric fibre of \(X\) over a closed point. Suppose that \(X_F\)
has a zero-cycle of degree prime to \(\ell\). Then there are a
regular scheme \(T\), finite locally free over \(R\) of rank prime
to \(\ell\), and an \(R\)-morphism \(T\to U\).
\end{corollary}

\begin{proof}
Write the zero-cycle as \(\sum_xn_x[x]\), with closed points
\(x\in X_F\) and \(n_x\in\Z\). Its degree
\(\sum_xn_x[\kappa(x):F]\) is prime to \(\ell\), so
\([\kappa(x):F]\) is prime to \(\ell\) for some \(x\). For this
point, \cref{lem:multisection-arithmetic-alteration} gives an
\(R\)-morphism \(C\to X\) with \(C\) regular and projective over
\(R\), \(C_F\) finite over \(F\), and \(\deg_F(C_F)\) prime to
\(\ell\). Now
\cref{thm:regular-multisections}\textup{(i)}
gives \(T\).
\end{proof}
\begin{remark}[Comparison with finite quasi-sections]
\label{rem:multisection-quasisections}
The finite flat quasi-sections of Gabber--Liu--Lorenzini
\cite[Theorem~6.3, Theorem~5.1 and Remark~6.6]{GLLHypersurfaces15}
are obtained by successive hypersurface sections. In the situation of
\cref{thm:regular-multisections}, with $a$, $c$ and
$\mathcal A$ as in \cref{con:multisection-replacement}, a complete
intersection of $c$ hypersurfaces in the linear system of
$\mathcal A_F^{\otimes b}$, when zero-dimensional, has degree
\[
 b^c\deg_F\bigl(c_1(\mathcal A_F)^c\cap[\mathbb P^a_{X_F}]\bigr)
\]
by \cite[\stackstag{02T9},\stackstag{02SU}]{Stacks}, applied as in
\cref{prop:multisection-generic-chow-degree}, where $\deg_F$ is the
proper pushforward of zero-cycles to $\Spec F$. If the degree of
$c_1(\mathcal A_F)^c\cap[\mathbb P^a_{X_F}]$ is divisible by $\ell$,
this cannot give a finite scheme
of degree prime to $\ell$. The correction term $-[C_F]$ in
\eqref{eq:multisection-generic-degree} gives the required nonzero
congruence even in this case. Regularity and avoidance of the boundary
are obtained simultaneously by
\cref{sub:multisection-absolute-cotangent,prop:multisection-decrease-closed-fibres}.
\end{remark}

\section{Returning the generic path}\label{sec:generic-return}

We prove \cref{thm:generic-path-return} and apply it to finite
Postnikov bands in \cref{prop:finite-type-band-injection}. The trace of
\cref{lem:postnikov-band-trace}, induced by restriction of
scalars on connective \(K\)-theory, makes the rank prime to \(\ell\)
sufficient for injectivity; no application of
\cref{thm:local-finite-connective-comparison,thm:semilocal-finite-connective-comparison}
to the multisection ring is needed.

\subsection{Returning two sections}

\begin{theorem}[Returning two sections from the generic point]
\label{thm:generic-path-return}
Let $\ell$ be a prime number, and let $R$ be a regular semilocal
domain essentially of finite type over a field or over $\Z$
\textup{(\cref{def:arithmetic-semilocal})}, with fraction field $F$,
such that $\ell\in R^\times$. Suppose that the closed residue fields
of $R$ have at most one positive characteristic. Let $\mathcal F$ be a
hypercomplete sheaf of spectra on the small \'etale site of
$\Spec R$ such that $\pi_i\mathcal F=0$ for all but finitely many $i\in\Z$
and every $\pi_i\mathcal F$ is a finite locally constant sheaf of abelian
groups of $\ell$-power order. Let
$z_0,z_1\in\Omega^\infty R\Gamma_{\mathrm{\acute et}}(R,\mathcal F)$ be two
points whose images in
$\Omega^\infty R\Gamma_{\mathrm{\acute et}}(F,\mathcal F|_F)$ lie in the same
connected component. Then there is a regular scheme $T$, finite
locally free over $R$ of rank prime to $\ell$, such that the images
of $z_0$ and $z_1$ in
$\Omega^\infty R\Gamma_{\mathrm{\acute et}}(T,\mathcal F|_T)$ lie in the same
connected component.
\end{theorem}

\begin{proof}
\emph{Finite Galois data.}
The ring $R$ is a noetherian semilocal domain in which $\ell$ is
invertible. By \cref{prop:finite-galois-data-two-sections}, there
are a $\pi$-finite space $\mathsf A$, a finite group $G$, a connected right
$G$-torsor $P_0\to\Spec R$, and a finite stage
$\gamma:BG\to\mathcal M_2(\mathsf A)$ representing the family
$(\Omega^\infty\rho_*\mathcal F,z_0,z_1)$, with $\mathcal M_1(\mathsf A)\to\mathcal M_2(\mathsf A)$
as in \cref{con:marked-moduli} for $\mathsf A$. By \cref{prop:finite-galois-data-two-sections},
for every field $E$ over $R$ and the parameter
$c_E\in\Gamma((\Spec E)_{\mathrm{fet}},\underline{\mathcal M_2(\mathsf A)})$ of
the pulled-back family,
\begin{equation}\label{eq:return-paths-as-lifts}
 \{z_0|_E\}\mathbin{\times^h_{\Omega^\infty
    R\Gamma_{\mathrm{\acute et}}(E,\mathcal F|_E)}}\{z_1|_E\}
 \simeq
 \Gamma((\Spec E)_{\mathrm{fet}},\underline{\mathcal M_1(\mathsf A)})
 \mathbin{\times^h_{\Gamma((\Spec E)_{\mathrm{fet}},
    \underline{\mathcal M_2(\mathsf A)})}}\{c_E\}.
\end{equation}

\smallskip\noindent
\emph{Reduction to an $\ell$-group.}
Let $G_\ell\subseteq G$ be a Sylow $\ell$-subgroup and put
$\Spec R'=P_0/G_\ell$. The domain $R$ is regular, hence normal. By
\cref{sub:finite-sylow-base-change}\textup{(i)}, the $R$-algebra
$R'$ is finite \'etale of rank $[G:G_\ell]$, and $P_0\to\Spec R'$ is a
right $G_\ell$-torsor; in particular $R'$ is regular, noetherian
and semilocal. By
\cref{sub:finite-sylow-base-change}\textup{(ii)}, the ring $R'$ is a
domain, and for its fraction field $F'$ the scheme
$P_0\times_{\Spec R'}\Spec F'$ is the spectrum of a field. The map
\[
 P_0\times^{G_\ell}G\longrightarrow P_0\times_{\Spec R}\Spec R',
 \qquad [p,g]\longmapsto(p\cdot g,[p]),
\]
is a morphism, hence an isomorphism, of right $G$-torsors over
$\Spec R'$. Thus the inverse image of the family
$(\Omega^\infty\rho_*\mathcal F,z_0,z_1)$ on $(\Spec R')_{\mathrm{fet}}$ is
represented by the finite stage
$\gamma_\ell:BG_\ell\to BG\xrightarrow{\gamma}\mathcal M_2(\mathsf A)$, with
the $G_\ell$-torsor $P_0\to\Spec R'$: the isomorphism of torsors
identifies their classifying sections after extension of structure
group, so $\gamma_\ell$ is a finite stage in the sense of
\cref{def:finite-stage}.

\smallskip\noindent
\emph{The lift over the Sylow fixed field.}
Let $\Gamma$ be the absolute Galois group of $F'$, and let
$\Gamma_\ell\subseteq\Gamma$ and
$L=((F')^{\mathrm{sep}})^{\Gamma_\ell}$ be as in
\cref{sub:pro-ell-sylow-fixed-field}. Let $\chi:\Gamma\to G_\ell$
be the homomorphism classifying the connected torsor
$P_0\times_{\Spec R'}\Spec F'$, defined by a point of $P_0$ with values in
$(F')^{\mathrm{sep}}$, and put
$\lambda=\chi|_{\Gamma_\ell}:\Gamma_\ell\to G_\ell$. By
\cref{sub:pro-ell-sylow-fixed-field}\textup{(iii)}, $\lambda$ is
surjective. By hypothesis there is a path from $z_0|_F$ to
$z_1|_F$; its pullback is a path from $z_0|_L$ to $z_1|_L$, and by
\eqref{eq:return-paths-as-lifts} for $E=L$ it gives a homotopy lift
to $\mathcal M_1(\mathsf A)$ of the parameter $c_L$. By
\cref{sub:pro-ell-sylow-fixed-field}\textup{(ii)}, the absolute
Galois group of $L$ is $\Gamma_\ell$, and its homomorphism to $G$
defined by the chosen point of $P_0$ is the composite of $\lambda$
with the inclusion $G_\ell\subseteq G$. By the last assertion of
\cref{prop:finite-galois-data-two-sections}, in the form of
\cref{lem:finite-stage-homotopy-lift}, there are an open normal
subgroup $J\subseteq\Gamma_\ell$ contained in $\ker\lambda$, with
quotient $H=\Gamma_\ell/J$ and induced surjection
$\phi:H\twoheadrightarrow G_\ell$, a map
$\delta:BH\to\mathcal M_1(\mathsf A)$, and a homotopy from the image of
$\delta$ in $\mathcal M_2(\mathsf A)$ to
$\gamma\circ B(G_\ell\subseteq G)\circ B\phi=\gamma_\ell\circ B\phi$.
The group $H$ is a quotient of the pro-$\ell$ group $\Gamma_\ell$
(\cref{sub:pro-ell-sylow-fixed-field}\textup{(i)}), hence a finite
$\ell$-group, and $|\ker\phi|$ is invertible in $R'$.
By \cref{sub:pro-ell-sylow-fixed-field}\textup{(iii)}, the spectrum
$Q_0$ of $((F')^{\mathrm{sep}})^{J}$ is a right $H$-torsor over $L$
and a marked lift of $P_0\times_{\Spec R'}\Spec L$ along $\phi$.

\smallskip\noindent
\emph{The path over the space of lifts.}
Apply \cref{prop:projective-marked-lifting-model} to
$S=\Spec R'$, the surjection $\phi:H\twoheadrightarrow G_\ell$, and
the torsor $P_0\to\Spec R'$. It gives a smooth projective
$R'$-scheme $X$ of relative dimension $|H|^2$ with integral
geometric fibres, an open subscheme $U\subseteq X$ that is dense in
every geometric fibre, a right $H$-torsor $Q\to U$, and an
isomorphism $\theta:Q\times^HG_\ell\simeq P_0\times_{\Spec R'}U$. By
\cref{cor:transport-marked-path-by-torsor}, applied to the finite
stage $\gamma_\ell$ over $\Spec R'$, the morphism $U\to\Spec R'$,
the torsor $(Q,\theta)$, and the homotopy lift $\delta$ of
$\gamma_\ell\circ B\phi$, there is a path between the pulled-back
sections on $U_{\mathrm{fet}}$. Its inverse image on
$U_{\mathrm{\acute et}}$ is a path from $z_0|_U$ to $z_1|_U$, by the
counit equivalence of \cref{prop:bounded-finite-etale-descent}\textup{(i)}
and the compatibility of inverse images in
\cref{cor:transport-marked-path-by-torsor}.

\smallskip\noindent
\emph{The multisection.}
By the last assertion of
\cref{prop:projective-marked-lifting-model}, the marked lift
$Q_0$ of $P_0\times_{\Spec R'}\Spec L$ is the pullback of $(Q,\theta)$ along a point of
$U_{F'}(L)$. By \cref{lem:closed-point-from-sylow-field}, the
scheme $X$ is smooth and projective over $R$ of relative dimension
$|H|^2$, the open subscheme $U_{\bar s}$ is dense in $X_{\bar s}$
for every geometric point $\bar s$ of $\Spec R$, and $U_F$ contains
a closed point of degree prime to $\ell$. In particular $X$ is
regular and flat over $R$, its geometric fibres are pure of
dimension $|H|^2$, and $X_F$ has a zero-cycle of degree prime to
$\ell$. By \cref{cor:arithmetic-multisections},
there are a regular scheme $T$, finite locally free over $R$ of
rank prime to $\ell$, and an $R$-morphism $T\to U$. The pullback
along $T\to U$ of the path over $U$ is a path from $z_0|_T$ to
$z_1|_T$.
\end{proof}

\subsection{First maps for finite Postnikov bands}
\label{subsec:postnikov-band-injectivity}

\begin{notation}[Finite Postnikov bands]
\label{not:finite-postnikov-band}
Let \(N=\ell^\nu\), where \(\ell\) is a prime number and
\(\nu\geq1\) is a natural number. For a scheme \(X\) on which
\(N\) is invertible, we use the sheaf
\(\mathcal K^{\mathrm{\acute et}}_{X,N}\) of
\cref{con:etale-k-sheaves}. For natural numbers \(a\leq b\), its
\emph{finite Postnikov band} is
\[
 \tau_{[a,b]}\mathcal K^{\mathrm{\acute et}}_{X,N}
   =\tau_{\geq a}\tau_{\leq b}\mathcal K^{\mathrm{\acute et}}_{X,N}.
\]
For \(n\in\Z\) and a \(\Z[1/N]\)-algebra \(B\), put
\begin{equation}\label{eq:band-functor}
 \Phi^{[a,b]}_n(B)=\pi_nR\Gamma_{\mathrm{\acute et}}
   \bigl(B,\tau_{[a,b]}\mathcal K^{\mathrm{\acute et}}_{B,N}\bigr).
\end{equation}
By \eqref{eq:big-small-truncated}, this is a functor from
\(\Z[1/N]\)-algebras to abelian groups, and for \(a<c\leq b\) the
cofibre sequence
\(\tau_{[c,b]}\mathcal K^{\mathrm{big}}_N\to
\tau_{[a,b]}\mathcal K^{\mathrm{big}}_N\to
\tau_{[a,c-1]}\mathcal K^{\mathrm{big}}_N\) gives a long exact
sequence, natural in \(B\),
\begin{equation}\label{eq:band-long-exact-sequence}
 \cdots\longrightarrow\Phi^{[c,b]}_n(B)\longrightarrow
 \Phi^{[a,b]}_n(B)\longrightarrow\Phi^{[a,c-1]}_n(B)
 \longrightarrow\Phi^{[c,b]}_{n-1}(B)\longrightarrow\cdots.
\end{equation}
\end{notation}

\begin{sublemma}[Continuity and torsion of a finite Postnikov band]
\label{sub:finite-band-continuity}
Let \(a\leq b\) be natural numbers and let \(n\in\Z\).
\begin{enumerate}
\item The functor \(\Phi^{[a,b]}_n\) of \eqref{eq:band-functor}
commutes with filtered colimits of \(\Z[1/N]\)-algebras.
\item For every \(\Z[1/N]\)-algebra \(B\), the group
\(\Phi^{[a,b]}_n(B)\) is killed by \(N^2\) (and by \(N\) if \(a=b\)).
\item For every morphism \(g:Y\to X\) of affine schemes on which
\(N\) is invertible, the morphism
\begin{equation}\label{eq:band-base-change}
 g^*\tau_{[a,b]}\mathcal K^{\mathrm{\acute et}}_{X,N}
   \longrightarrow\tau_{[a,b]}\mathcal K^{\mathrm{\acute et}}_{Y,N}
\end{equation}
obtained by applying \(\tau_{[a,b]}\) to
\eqref{eq:k-sheaf-base-change} is an equivalence.
\end{enumerate}
\end{sublemma}

\begin{proof}
\textup{(iii)} Ordinary and hypercompleted inverse images are
\(t\)-exact (\cref{sub:finite-site-basics}\textup{(ii)} and
\cite[Proposition~5.5.6.28]{HTT}); for bounded bands their values agree.
By \cref{sub:finite-k-etale-coefficients}\textup{(ii)} and
\textup{(iii)}, the morphism \eqref{eq:band-base-change} induces an
isomorphism on every homotopy sheaf. Its source and its target are bounded above by \(b\) for the
Postnikov \(t\)-structure, so they are hypercomplete
\cite[Lemma~6.5.2.9]{HTT}, and
\eqref{eq:band-base-change} is an equivalence.

\smallskip\noindent
\textup{(i)} Let \(a=b\). If \(a\) is odd, then
\(\Phi^{[a,a]}_n=0\) by
\cref{sub:finite-k-etale-coefficients}\textup{(ii)}. If \(a=2j\),
then \cref{sub:finite-k-etale-coefficients}\textup{(ii)} and
\textup{(iii)} give an isomorphism, natural in \(B\),
\[
 \Phi^{[2j,2j]}_n(B)\simeq H^{2j-n}_{\mathrm{\acute et}}(B,\Z/N(j)).
\]
The right-hand side commutes with filtered colimits
\cite[\stackstag{09YQ}]{Stacks} and is killed by \(N\). For \(a<b\),
continuity follows by induction on \(b-a\) from
\eqref{eq:band-long-exact-sequence} for \(c=b\), exactness of filtered
colimits of abelian groups and the five lemma.

\smallskip\noindent
\textup{(ii)} In the homotopy category of a stable category, consider
any cofiber triangle
\[
 E\xrightarrow{N}E\xrightarrow{u}C\xrightarrow{v}\Sigma E.
\]
Since $(N\,\mathrm{id}_C)u=u(N\,\mathrm{id}_E)=0$, exactness of
$\operatorname{Hom}(-,C)$ on the triangle
\cite[\stackstag{0149}]{Stacks} gives a morphism $w:\Sigma E\to C$
with $N\,\mathrm{id}_C=wv$. Hence
$N^2\,\mathrm{id}_C=wv(N\,\mathrm{id}_C)=w(N\,\mathrm{id}_{\Sigma E})v=0$,
the last equality being the vanishing of the composite of consecutive
arrows of the triangle \cite[\stackstag{0146}]{Stacks}.
Apply this in the stable category of presheaves of spectra on the
small \'etale site of $\Spec B$ to $E=K^{\mathrm{cn}}$.
Hypercomplete sheafification and derived global sections preserve
finite biproducts as adjoints between stable categories. The truncations,
regarded as endofunctors, also preserve finite biproducts: each is an
adjoint to the inclusion of a truncation subcategory, which is closed
under finite biproducts. Thus these functors are additive on homotopy
categories. It follows that multiplication by $N^2$
on $R\Gamma_{\mathrm{\acute et}}(B,\tau_{[a,b]}\mathcal K^{\mathrm{\acute et}}_{B,N})$
is zero, which proves the uniform bound. The single-grade assertion
follows from the calculation in \textup{(i)}.
\end{proof}

\begin{sublemma}[Restriction followed by restriction of scalars]
\label{sub:finite-flat-multiplication}
Let $R$ be a semilocal ring and $B$ a finite locally free $R$-algebra
of constant rank $e$. For every $R$-algebra $R'$, the composite
of extension of scalars along $R'\to B\otimes_R R'$ and restriction
of scalars back to $R'$ is homotopic to multiplication by $e$ on
$K(R')$ and on $K^{\mathrm{cn}}(R')$, naturally in $R'$. In particular,
for every natural number $N\geq2$ and every $n\in\Z$ the composite
\[
 K_n(R;\Z/N)\xrightarrow{\,-\otimes_R B\,}K_n(B;\Z/N)
       \xrightarrow{\,\operatorname{Res}_{B/R}\,}K_n(R;\Z/N)
\]
is multiplication by $e$. If $\gcd(e,N)=1$, a class whose pullback
to $B$ is zero is itself zero.
\end{sublemma}

\begin{proof}
The identity before taking finite coefficients is the finite
locally free case of the projection formula
\cite[Proposition~3.17]{TT}.

By \cref{sub:semilocal-projective-free}, choose an isomorphism of
$R$-modules $B\simeq R^e$. Its base changes identify $B\otimes_R R'$
with $(R')^e$, compatibly with all homomorphisms of $R$-algebras.
Restriction of scalars preserves perfect complexes because
$B\otimes_R R'$ is finite projective over $R'$. It commutes with
extension of scalars along homomorphisms $R'\to R''$ of
$R$-algebras by associativity of derived tensor products. For every
perfect complex $P$ over $R'$, the chosen isomorphism induces an
equivalence
\[
 \operatorname{Res}_{(B\otimes_R R')/R'}
       (P\otimes_R B)\simeq P^{\oplus e},
\]
natural in $P$ and compatible with extension of scalars in $R'$.
The direct sum functor induces addition on connective and nonconnective
$K$-theory spectra, so this equivalence of functors gives a homotopy
of the corresponding natural transformations of spectra with
multiplication by $e$. Taking the cofiber of multiplication by $N$
and then $\pi_n$ yields the asserted composite on $K_n(R;\Z/N)$. By \eqref{eq:coefficient-sequence}, the group
$K_n(R;\Z/N)$ is an extension of a group killed by $N$ by a group
killed by $N$, so it is killed by $N^2$. If $\gcd(e,N)=1$, choose
integers $u,v$ with $ue+vN^2=1$. For $x\in K_n(R;\Z/N)$ with zero
image in $K_n(B;\Z/N)$ we get $ex=0$ and $x=u\,ex+vN^2x=0$.
\end{proof}

\begin{lemma}[Trace on a finite Postnikov band]
\label{lem:postnikov-band-trace}
Let \(R\) be a semilocal \(\Z[1/N]\)-algebra, let \(B\) be a finite
locally free \(R\)-algebra of constant rank \(e\), and let
\(f:\Spec B\to\Spec R\) be the corresponding morphism. For all
natural numbers \(a\leq b\), there is a morphism of sheaves of
spectra on the small \'etale site of \(\Spec R\)
\[
 \operatorname{tr}_f:
 f_*\tau_{[a,b]}\mathcal K^{\mathrm{\acute et}}_{B,N}
 \longrightarrow\tau_{[a,b]}\mathcal K^{\mathrm{\acute et}}_{R,N}
\]
such that, for every \(n\in\Z\), the composite
\[
 \Phi^{[a,b]}_n(R)\longrightarrow\Phi^{[a,b]}_n(B)
 \xrightarrow{\ \operatorname{tr}_f\ }\Phi^{[a,b]}_n(R)
\]
of the restriction map and the map induced by
\(\operatorname{tr}_f\) is multiplication by \(e\).
\end{lemma}

\begin{proof}
\emph{Presheaves.}
For an \'etale \(R\)-algebra \(R'\) put
\[
 \mathcal P(R')=K^{\mathrm{cn}}(R')/N,\qquad
 \mathcal Q(R')=K^{\mathrm{cn}}(B\otimes_RR')/N.
\]
Base change along \(R'\to B\otimes_RR'\) and restriction of scalars
define morphisms \(\operatorname{res}:\mathcal P\to \mathcal Q\) and \(\operatorname{tr}:\mathcal Q\to \mathcal P\) of presheaves
of spectra on the affine objects of
\((\Spec R)_{\mathrm{\acute et}}\). By \cref{sub:finite-flat-multiplication}, the composite
\(\operatorname{tr}\circ\operatorname{res}\) is homotopic to multiplication by \(e\) as a
morphism of presheaves of spectra.

\smallskip\noindent
\emph{The direct image.}
By \cref{con:etale-k-sheaves},
\(L^{\wedge}_{R_{\mathrm{\acute et}}}\mathcal P=\mathcal K^{\mathrm{\acute et}}_{R,N}\)
and \(f_*\) preserves hypercomplete sheaves. For a hypercomplete
sheaf of spectra \(\mathcal G\) on \((\Spec B)_{\mathrm{\acute et}}\)
with finitely many nonzero homotopy sheaves, one has
\(\pi_i(f_*\mathcal G)\simeq f_*\pi_i(\mathcal G)\) for every \(i\).
Indeed, for an Eilenberg--Mac Lane sheaf this follows from the vanishing
of \(R^qf_*\) for \(q\geq1\) \cite[\stackstag{03QP}]{Stacks}; the
general case follows by applying the exact functor \(f_*\) to the finite
Postnikov filtration of \(\mathcal G\).
In particular $f_*\tau_{\leq b}\mathcal K^{\mathrm{\acute et}}_{B,N}$
is hypercomplete and has homotopy sheaves in degrees $[0,b]$ only,
since $\tau_{\leq b}\mathcal K^{\mathrm{\acute et}}_{B,N}$ is bounded
above \cite[Lemma~6.5.2.9]{HTT} and connective
(\cref{sub:finite-k-etale-coefficients}\textup{(i)}). On presheaves, compose the unit
$\mathcal Q\to f_*\mathcal K^{\mathrm{\acute et}}_{B,N}$ with the map induced
by the truncation $\mathcal K^{\mathrm{\acute et}}_{B,N}\to
\tau_{\leq b}\mathcal K^{\mathrm{\acute et}}_{B,N}$. By the universal
property of hypercomplete sheafification this composite factors through
$L^{\wedge}_{R_{\mathrm{\acute et}}}\mathcal Q$, and by the universal property
of $\tau_{\leq b}$ it factors further through
$\tau_{\leq b}L^{\wedge}_{R_{\mathrm{\acute et}}}\mathcal Q$. We obtain the natural morphism
\begin{equation}\label{eq:band-direct-image}
 \tau_{\leq b}L^{\wedge}_{R_{\mathrm{\acute et}}}\mathcal Q\longrightarrow
 f_*\tau_{\leq b}\mathcal K^{\mathrm{\acute et}}_{B,N}.
\end{equation}
Both sides are hypercomplete, since they are bounded above
\cite[Lemma~6.5.2.9]{HTT}. Let \(\bar x\) be a geometric point of
\(\Spec R\), and index \(\bar y\) by the finitely many points of the geometric fibre
of \(f\) at \(\bar x\), choosing a geometric lift for each. By
\cite[\stackstag{04GH}, \stackstag{05WR}]{Stacks},
\[
 B\otimes_R\mathcal O^{\mathrm{sh}}_{\Spec R,\bar x}
   \simeq\prod_{\bar y}\mathcal O^{\mathrm{sh}}_{\Spec B,\bar y}.
\]
For \(0\leq i\leq b\), the continuity of \(K^{\mathrm{cn}}/N\) used
for \eqref{eq:k-sheaf-stalk}, the formula \eqref{eq:k-sheaf-stalk}
itself, and the stalk formula for the direct image
\cite[\stackstag{03QP}]{Stacks} give
\[
 (\pi_iL^{\wedge}_{R_{\mathrm{\acute et}}}\mathcal Q)_{\bar x}
 =\prod_{\bar y}\pi_i\bigl(K^{\mathrm{cn}}
     (\mathcal O^{\mathrm{sh}}_{\Spec B,\bar y})/N\bigr)
 =\bigl(f_*\pi_i\mathcal K^{\mathrm{\acute et}}_{B,N}\bigr)_{\bar x},
\]
and \eqref{eq:band-direct-image} induces this identification,
because both equalities are induced by the restriction maps to the
rings \(\mathcal O^{\mathrm{sh}}_{\Spec B,\bar y}\). Hence
\eqref{eq:band-direct-image} is an equivalence.

\smallskip\noindent
\emph{The trace.}
Apply \(\tau_{[a,b]}\) to \(L^{\wedge}_{R_{\mathrm{\acute et}}}\operatorname{res}\)
and \(L^{\wedge}_{R_{\mathrm{\acute et}}}\operatorname{tr}\). By the equivalence
\eqref{eq:band-direct-image} and the calculation
of the homotopy sheaves of $f_*$, there are natural
equivalences
\begin{equation}\label{eq:band-direct-image-equivalence}
 \tau_{[a,b]}L^{\wedge}_{R_{\mathrm{\acute et}}}\mathcal Q\simeq
 \tau_{\geq a} f_*\tau_{\leq b}\mathcal K^{\mathrm{\acute et}}_{B,N}
 \simeq f_*\tau_{[a,b]}\mathcal K^{\mathrm{\acute et}}_{B,N}.
\end{equation}
Indeed, apply the exact functor $f_*$ to the fibre sequence defining
$\tau_{\geq a}$ of the bounded sheaf $\tau_{\leq b}\mathcal K^{\mathrm{\acute et}}_{B,N}$;
the finite Postnikov argument identifies the homotopy ranges of the
fibre and the quotient. Hence we obtain morphisms
\[
 \tau_{[a,b]}\mathcal K^{\mathrm{\acute et}}_{R,N}
 \xrightarrow{\ \tau_{[a,b]}L^{\wedge}_{R_{\mathrm{\acute et}}}\operatorname{res}\ }
 f_*\tau_{[a,b]}\mathcal K^{\mathrm{\acute et}}_{B,N}
 \xrightarrow{\ \operatorname{tr}_f=\tau_{[a,b]}L^{\wedge}_{R_{\mathrm{\acute et}}}\operatorname{tr}\ }
 \tau_{[a,b]}\mathcal K^{\mathrm{\acute et}}_{R,N}
\]
whose composite is multiplication by \(e\): hypercomplete sheafification
is exact and the band functor preserves finite biproducts, hence both
preserve the sum of $e$ copies of an identity morphism.
By \eqref{eq:big-small} and naturality of the unit of sheafification
on the big site, the composite
\(\mathcal P\xrightarrow{\operatorname{res}}\mathcal Q\to f_*\mathcal K^{\mathrm{\acute et}}_{B,N}\)
is the composite of the unit \(\mathcal P\to L^{\wedge}_{R_{\mathrm{\acute et}}}\mathcal P\)
and the adjoint of
\eqref{eq:k-sheaf-base-change}. Composing with the map to
\(f_*\tau_{\leq b}\mathcal K^{\mathrm{\acute et}}_{B,N}\), the universal
properties of hypercomplete sheafification and \(\tau_{\leq b}\) identify
\(\tau_{\leq b}L^{\wedge}_{R_{\mathrm{\acute et}}}\operatorname{res}\), followed by
\eqref{eq:band-direct-image},
with the adjoint of the upper truncation of
\eqref{eq:k-sheaf-base-change}. Applying \(\tau_{\geq a}\) and using
\eqref{eq:band-direct-image-equivalence}, the first morphism is the
adjoint of
\eqref{eq:band-base-change} for \(g=f\).
By \eqref{eq:big-small-truncated} it induces on global sections the
restriction map \(\Phi^{[a,b]}_n(R)\to\Phi^{[a,b]}_n(B)\) of
\eqref{eq:band-functor}.
\end{proof}

\begin{remark}[Two features of the band trace]
\label{rem:band-trace-pitfalls}
For an algebraically closed field $k$ of characteristic different
from $\ell$, let $e>1$ be prime to $\ell$ and consider
$R=k[[s]]\to B=k[[t]]$ given by $s=t^e$.
The closed geometric fibre has a single point, but its length is $e$.
The trace of \cref{lem:postnikov-band-trace} at the closed
point is multiplication by $e$ under the identification by rigidity
(\cref{sub:finite-k-etale-coefficients}); it is not an unweighted
sum over the geometric points.

Also, injectivity on the cohomology of each graded term of a
filtration does not suffice for injectivity on the cohomology of the
filtered object. For example, filter the inclusion
$\mathbf F_\ell[-1]\to(\mathbf F_\ell\xrightarrow{1}\mathbf F_\ell)$,
where the latter complex is in cohomological degrees $0,1$, by the
decreasing filtration with stage zero the entire complex, stage one
the subcomplex in degree one, and stage two zero. The induced maps on all graded cohomology
groups are injective, whereas the map on $H^1$ is zero. The trace
identity is proved as a morphism of bands and thus retains their
extension data.
\end{remark}

\begin{proposition}[First maps for finite bands on arithmetic semilocal domains]
\label{prop:finite-type-band-injection}
Let \(R_0\) be a regular semilocal domain essentially of finite type over
a field or over \(\Z\), with fraction field \(F_0\), and suppose that the
closed residue fields of \(R_0\) have at most one positive characteristic.
For \(N=\ell^\nu\) invertible in \(R_0\), \(0\leq a\leq b\), and \(n\in\Z\),
the map \(\Phi^{[a,b]}_n(R_0)\to\Phi^{[a,b]}_n(F_0)\) is injective.
\end{proposition}

\begin{proof}
Let \(x\in\Phi^{[a,b]}_n(R_0)\) have zero image in
\(\Phi^{[a,b]}_n(F_0)\). Put
\(\mathcal F=\Sigma^{-n}\tau_{[a,b]}\mathcal K^{\mathrm{\acute et}}_{R_0,N}\).
By \cref{sub:finite-k-etale-coefficients}\textup{(ii)}, the sheaf
\(\mathcal F\) is hypercomplete with finitely many nonzero homotopy sheaves,
which are finite locally constant of \(\ell\)-power order. By
\cref{sub:finite-band-continuity}\textup{(iii)}, for every affine
\(R_0\)-scheme \(Y=\Spec B\),
\begin{equation}\label{eq:band-sections-components}
 \pi_0\,\Omega^\infty R\Gamma_{\mathrm{\acute et}}(Y,\mathcal F|_Y)
   \simeq\Phi^{[a,b]}_n(B),
\end{equation}
compatibly with restriction. Let
\(z_0\in\Omega^\infty R\Gamma_{\mathrm{\acute et}}(R_0,\mathcal F)\) be a
point of the connected component \(x\), and let \(z_1\) be the zero
section. By \eqref{eq:band-sections-components} for \(B=F_0\), the
points \(z_0|_{F_0}\) and \(z_1|_{F_0}\) lie in the same connected
component. By \cref{thm:generic-path-return}, there is a
regular scheme \(T=\Spec B\), finite locally free over \(R_0\) of
rank \(e\) prime to \(\ell\), such that \(z_0|_T\) and \(z_1|_T\) lie
in the same connected component. By
\eqref{eq:band-sections-components}, the image of \(x\) in
\(\Phi^{[a,b]}_n(B)\) is zero. By
\cref{lem:postnikov-band-trace}, \(e\,x=0\), and by
\cref{sub:finite-band-continuity}\textup{(ii)},
\(N^2x=0\). Since \(\gcd(e,N)=1\), we get \(x=0\).

\end{proof}

\begin{remark}
The proof of \cref{prop:finite-type-band-injection} uses only the
finite local freeness and the rank prime to \(\ell\) of the morphism
in \cref{thm:generic-path-return}, not regularity of its source.
\Cref{thm:regular-multisections} gives more: a regular
finite locally free scheme mapping to the lifting open. It would be
interesting to obtain such a finite locally free morphism directly,
without the regularity requirement.
\end{remark}
\section{Passage to arbitrary noetherian regular local rings}
\label{sec:approximation}

This section proves
\cref{thm:gersten-injectivity}.
\Cref{lem:finite-type-witness} transports an element with zero
generic image, for a functor that commutes with filtered colimits,
from a noetherian regular local ring to a regular local ring
essentially of finite type over a field or over $\Z$, where
\cref{prop:finite-type-band-injection} applies to finite Postnikov
bands. This yields \cref{thm:postnikov-band-injectivity}, and
\cref{thm:local-finite-connective-comparison}
then gives \cref{thm:gersten-injectivity}.
\Cref{lem:finite-type-witness} is also used for
\'etale cohomology in \cref{sec:ordinary-etale-injection}.

\begin{sublemma}[Continuity of $K$-theory with finite coefficients]
\label{sub:k-finite-continuity}
Let $(A_i)_{i\in I}$ be a filtered diagram of rings with colimit
$A$, let $N\geq2$ be a natural number, and let $n\in\Z$. The
canonical map
\begin{equation}\label{eq:k-finite-continuity}
 \varinjlim_{i\in I}K_n(A_i;\Z/N)\longrightarrow K_n(A;\Z/N)
\end{equation}
is an isomorphism.
\end{sublemma}

\begin{proof}
The colimit can be computed over a directed set
\cite[\stackstag{0032}]{Stacks}. The scheme $\Spec A$ is the inverse
limit of the affine schemes $\Spec A_i$, with affine transition
morphisms. These schemes are quasi-compact and quasi-separated,
so \cite[Theorem~7.2, p.~364]{TT} gives the equivalence
$\varinjlim_iK(A_i)\to K(A)$ of spectra. The
cofiber of multiplication by $N$ commutes with colimits of spectra,
and the functor $\pi_n$ from spectra to abelian groups commutes
with filtered colimits \cite[\S1.4.3]{HA}.
\end{proof}

\begin{sublemma}[Descending an element with zero generic image]
\label{sub:principal-open-witness}
Let $B$ be a ring and let $(R_i)_{i\in I}$ be a diagram of
$B$-algebras indexed by a directed set, whose colimit $R$ is a
domain with fraction field $F$.
\begin{enumerate}
\item Let $R_1$ be a $B$-algebra of finite presentation. Every
homomorphism of $B$-algebras $R_1\to R$ is the composite of a
homomorphism of $B$-algebras $R_1\to R_i$, for some $i\in I$, and
the canonical map $R_i\to R$.
\item Let $\Phi$ be a functor from $B$-algebras to abelian groups
that commutes with filtered colimits, and let $x\in\Phi(R)$ be an
element with zero image in $\Phi(F)$. Then there are an index
$i\in I$, an element $f_i\in R_i$ with nonzero image in $R$, and an
element $x_i\in\Phi(R_i)$ such that
\[
 x_i\longmapsto x\in\Phi(R)
 \qquad\text{and}\qquad
 x_i\longmapsto0\in\Phi(R_i[1/f_i]).
\]
If $R_i$ is a domain, then $x_i$ has zero image in
$\Phi(\Frac R_i)$.
\end{enumerate}
\end{sublemma}

\begin{proof}
\textup{(i)} This is \cite[\stackstag{00QO}]{Stacks}.

\smallskip\noindent
\textup{(ii)} By \cite[\stackstag{00CR}]{Stacks},
$F=\varinjlim_{f\in R\setminus\{0\}}R[1/f]$, where
$R\setminus\{0\}$ is a directed preorder under divisibility. Hence
\[
 \Phi(F)=\varinjlim_{f\in R\setminus\{0\}}\Phi(R[1/f]),
\]
and there is $f\in R\setminus\{0\}$ such that $x$ has zero image in
$\Phi(R[1/f])$. Since $R=\varinjlim_iR_i$ and
$\Phi(R)=\varinjlim_i\Phi(R_i)$, there are $i_0\in I$,
$f_{i_0}\in R_{i_0}$ and $x_{i_0}\in\Phi(R_{i_0})$ with images $f$
and $x$. For $i\geq i_0$ let $f_i\in R_i$ and $x_i\in\Phi(R_i)$ be
the images of $f_{i_0}$ and $x_{i_0}$. The identity
$R_i[1/f_i]=R_i[T]/(f_iT-1)$ gives, by passage to the colimit,
\[
 R[1/f]=\varinjlim_{i\geq i_0}R_i[1/f_i],
 \qquad
 \Phi(R[1/f])=\varinjlim_{i\geq i_0}\Phi(R_i[1/f_i]),
\]
and the image of $x_{i_0}$ in $\Phi(R[1/f])$ is zero. Hence
$x_i$ has zero image in $\Phi(R_i[1/f_i])$ for some $i\geq i_0$. If
$R_i$ is a domain, then $f_i\neq0$ gives
$R_i[1/f_i]\subseteq\Frac R_i$, and $x_i$ has zero image in
$\Phi(\Frac R_i)$.
\end{proof}

\begin{sublemma}[Localizing a filtered presentation of a semilocal ring]
\label{sub:localized-filtered-presentation}
Let $(P_i)_{i\in I}$ be a filtered diagram of rings whose colimit
$R$ is a semilocal ring with maximal ideals
$\mathfrak m_1,\ldots,\mathfrak m_t$ (with $t=0$ if $R=0$). For $i\in I$ and
$1\leq j\leq t$, let $\mathfrak p_{i,j}\subset P_i$ be the inverse
image of $\mathfrak m_j$, and put
\[
 S_i=P_i\setminus\bigcup_{j=1}^t\mathfrak p_{i,j},
 \qquad R_i=S_i^{-1}P_i.
\]
\begin{enumerate}
\item The transition maps of $(P_i)_{i\in I}$ induce a filtered
diagram $(R_i)_{i\in I}$, and $\varinjlim_{i\in I}R_i=R$.
\item The ring $R_i$ is the semilocalization of $P_i$ at the
maximal elements of
$\{\mathfrak p_{i,1},\ldots,\mathfrak p_{i,t}\}$
\textup{(\cref{def:arithmetic-semilocal})}. For each such element
$\mathfrak p_{i,j}$, the closed residue field
$R_i/\mathfrak p_{i,j}R_i$ is a subfield of $R/\mathfrak m_j$.
\item If $R$ is local, then $R_i=(P_i)_{\mathfrak p_{i,1}}$ is
local and $R_i\to R$ is a local homomorphism.
\end{enumerate}
\end{sublemma}

\begin{proof}
If $R=0$, then $S_i=P_i$ and $R_i=0$ for every $i$; the assertions
follow from \cref{def:arithmetic-semilocal}, with \textup{(iii)}
vacuous. Assume henceforth that $R\neq0$.
For a local ring $R$, assertion \textup{(i)} is the equality
$R=\varinjlim_i(P_i)_{\mathfrak p_{i,1}}$, used in the proof of
\cite[\stackstag{037D}]{Stacks}; the semilocal case is the same
argument. The units of $R$ are the elements outside
$\bigcup_j\mathfrak m_j$. Hence an element of $P_i$ lies in $S_i$ if
and only if its image in $R$ is a unit.

\smallskip\noindent
\textup{(i)} A transition map $P_i\to P_{i'}$ maps $S_i$ into
$S_{i'}$, and $P_i\to R$ maps $S_i$ into $R^\times$. This gives the
diagram $(R_i)_{i\in I}$ and a map $\varinjlim_iR_i\to R$, which is
surjective because $\varinjlim_iP_i\to R$ is. Let $a\in P_i$ and
$s\in S_i$ be such that $a/s$ has zero image in $R$. Then $a$ has
zero image in $R$, hence in $P_{i'}$ for some $i\to i'$, and $a/s$
has zero image in $R_{i'}$.

\smallskip\noindent
\textup{(ii)} The union $\bigcup_j\mathfrak p_{i,j}$ is the union of
the maximal elements of
$\{\mathfrak p_{i,1},\ldots,\mathfrak p_{i,t}\}$, which are pairwise
incomparable. By \cref{def:arithmetic-semilocal},
$R_i/\mathfrak p_{i,j}R_i=\Frac(P_i/\mathfrak p_{i,j})$ for each of
them, and $P_i/\mathfrak p_{i,j}\to R/\mathfrak m_j$ is injective.

\smallskip\noindent
\textup{(iii)} This is \textup{(ii)} for $t=1$.
\end{proof}

\begin{lemma}[Descent to a regular local ring essentially of finite type]
\label{lem:finite-type-witness}
Let $R$ be a noetherian regular local ring with fraction field $F$.
If $R$ contains a field, let $k$ be its prime field; otherwise let
$k=\Z$. Let $R_1$ be a finite-type $k$-algebra, let $R_1\to R$ be a
ring map, let $\Phi$ be a functor from $R_1$-algebras to abelian
groups that commutes with filtered colimits, and let $x\in\Phi(R)$
be an element with zero image in $\Phi(F)$. Then there are a regular
local ring $R_0$ essentially of finite type over $k$, with fraction
field $F_0$, ring maps $R_1\to R_0\to R$ whose composite is the
given map and whose second map is a local homomorphism, as in the
commutative diagram
\[
\begin{tikzcd}[column sep=large,row sep=large]
 R_1 \arrow[r] \arrow[dr,dashed] & R \\
 & R_0 \arrow[u,dashed, "\text{local}"']
\end{tikzcd}
\]
and an element
\[
 x_0\in\Phi(R_0)
 \qquad\text{with}\qquad
 x_0\longmapsto x\in\Phi(R),\qquad x_0\longmapsto0\in\Phi(F_0).
\]
\end{lemma}

\begin{proof}
\emph{A presentation of $R$.}
We first show that there is a diagram $(P_i)_{i\in I}$ of regular
rings essentially of finite type over $k$, indexed by a directed
set, with colimit $R$. Suppose that $R$ contains a field. Every
finite purely inseparable extension of the perfect field $k$ is
trivial, so the regular noetherian ring $R$ is geometrically regular
over $k$ \cite[\stackstag{0381}, \stackstag{0382}]{Stacks}, and
$k\to R$ is a regular
ring map \cite[\stackstag{07BZ}]{Stacks}. By Popescu's theorem
\cite[\stackstag{07GC}]{Stacks} and
\cite[\stackstag{0032}]{Stacks}, the ring $R$ is the colimit of a
diagram of smooth $k$-algebras indexed by a directed set, and smooth
$k$-algebras are regular \cite[\stackstag{07NF}]{Stacks}. Suppose
that $R$ contains no field. Since $R$ is a domain
\cite[\stackstag{00NP}]{Stacks}, positive characteristic would give
a prime field in $R$. If its residue field had characteristic zero,
every nonzero integer would be invertible in $R$, giving a copy of
$\mathbf Q$. Thus $R$ has mixed characteristic. By the last
sentence of \cite[Theorem~3.8, p.~722]{Popescu19}
(\cite[Theorem~17, p.~7]{Popescu19arXiv}) and
\cite[\stackstag{0032}]{Stacks}, the ring $R$ is the colimit of a
diagram of regular local rings essentially of finite type over $\Z$
indexed by a directed set. In both cases, the transition maps and
structure maps are $k$-algebra homomorphisms, since $k$ is a prime
field or $\Z$.

\smallskip\noindent
\emph{Local stages.}
Let $R_i$ be the localization of $P_i$ at the inverse image of the
maximal ideal of $R$. By
\cref{sub:localized-filtered-presentation}\textup{(i)} and
\textup{(iii)}, $R=\varinjlim_iR_i$ and every $R_i\to R$ is a local
homomorphism. Each $R_i$ is a regular local ring essentially of
finite type over $k$, hence a domain
\cite[\stackstag{00NP}]{Stacks}.

\smallskip\noindent
\emph{Descent.}
The ring $k$ is noetherian, so $R_1$ is a $k$-algebra of finite
presentation. By \cref{sub:principal-open-witness}\textup{(i)},
there are $i_1\in I$ and a factorization $R_1\to R_{i_1}\to R$ of
the given map. On the cofinal tail $i\geq i_1$, the maps
$R_1\to R_{i_1}\to R_i$ make the diagram a diagram of
$R_1$-algebras with colimit $R$ for its given $R_1$-algebra structure. By \cref{sub:principal-open-witness}\textup{(ii)},
there are $i\geq i_1$, $x_i\in\Phi(R_i)$ with image $x$, and
$f_i\in R_i$ with nonzero image $f$ in $R$, such that $x_i$ has zero
image in $\Phi(R_i[1/f_i])$. Put $R_0=R_i$, $x_0=x_i$, and $f_0=f_i$.
The vanishing and the factorization are recorded in the commutative
diagram of rings
\[
\begin{tikzcd}[column sep=large]
 R_0 \arrow[r] \arrow[d] & R_0[1/f_0] \arrow[r,hook] \arrow[d] & F_0 \\
 R \arrow[r] & R[1/f] \arrow[r,hook] & F .
\end{tikzcd}
\]
The image of $x_0$ in $\Phi(R_0[1/f_0])$ is zero, hence so is
its image in $\Phi(F_0)$; its image in $\Phi(R)$ is $x$. In general
$R_0\to R$ does not extend to a homomorphism $F_0\to F$, since it
need not be injective.
\end{proof}

\begin{theorem}[First maps for finite bands]
\label{thm:postnikov-band-injectivity}
Let \(R\) be a noetherian regular local ring with fraction
field \(F\), and let \(N=\ell^\nu\) be invertible in \(R\). For all
natural numbers \(a\leq b\) and all \(n\in\Z\), the restriction map
\[
 \pi_n R\Gamma_{\mathrm{\acute et}}
   (R,\tau_{[a,b]}\mathcal K^{\mathrm{\acute et}}_{R,N})
 \longrightarrow
 \pi_n R\Gamma_{\mathrm{\acute et}}
   (F,\tau_{[a,b]}\mathcal K^{\mathrm{\acute et}}_{F,N})
\]
is injective.
\end{theorem}

\begin{proof}
Let \(x\in\Phi^{[a,b]}_n(R)\) have zero image in
\(\Phi^{[a,b]}_n(F)\). If \(R\) contains a field, let \(k\) be its
prime field and put \(R_1=k\), in which \(N\) is invertible
since \(k\subseteq R\) is a field; otherwise put \(k=\Z\) and
\(R_1=\Z[1/N]\). In both cases \(R_1\) is a finite-type
\(k\)-algebra and a \(\Z[1/N]\)-algebra, and \(R\) is an
\(R_1\)-algebra. By \cref{sub:finite-band-continuity}\textup{(i)},
the functor \(\Phi^{[a,b]}_n\) on \(R_1\)-algebras commutes with
filtered colimits. By \cref{lem:finite-type-witness}, there are a
regular local ring \(R_0\) essentially of finite type over \(k\),
with fraction field \(F_0\), ring maps \(R_1\to R_0\to R\), and
\(x_0\in\Phi^{[a,b]}_n(R_0)\) with image \(x\) in
\(\Phi^{[a,b]}_n(R)\) and zero image in \(\Phi^{[a,b]}_n(F_0)\).
Since \(R_0\) is an \(R_1\)-algebra, \(\ell\in R_0^\times\), and since
it is local it has only one closed residue field. Thus
\cref{prop:finite-type-band-injection} applies to \(R_0\) and gives
\(x_0=0\), hence \(x=0\).
\end{proof}

\begin{proof}[Proof of \cref{thm:gersten-injectivity}]
By \cref{sub:moore-primary-decomposition} it suffices to consider
\(N=\ell^\nu\), with \(\nu\geq1\). Both \(R\) and \(F\) are regular
noetherian local rings, so connective and nonconnective $K$-theory
agree on them by \cite[Proposition~6.8(b), p.~362]{TT}. By
\eqref{eq:coefficient-sequence}, their finite-coefficient groups in
negative degrees vanish. Assume \(n\geq0\). Since
\(\mathcal K^{\mathrm{\acute et}}_{A,N}\) is connective for
\(A=R,F\) by \cref{sub:finite-k-etale-coefficients}\textup{(i)},
its truncations \(\tau_{[0,2n]}\) and \(\tau_{\leq2n}\) agree.

The natural comparison of \cref{thm:local-finite-connective-comparison}
and \eqref{eq:band-functor} therefore give the commutative diagram
\[
\begin{tikzcd}[column sep=large]
 K_n(R;\Z/N) \arrow[r] \arrow[d,"\sim"'] &
 K_n(F;\Z/N) \arrow[d,"\sim"] \\
 \Phi^{[0,2n]}_n(R) \arrow[r,hook] & \Phi^{[0,2n]}_n(F).
\end{tikzcd}
\]
The bottom arrow is injective by \cref{thm:postnikov-band-injectivity}.
\end{proof}

\begin{proposition}[Equicharacteristic rings and arbitrary coefficients]
\label{lem:equicharacteristic-injection}
Let $R$ be a noetherian regular local ring containing a field, with
fraction field $F$. For every natural number $N\geq2$ and every
$n\in\Z$, the restriction map
$K_n(R;\Z/N)\to K_n(F;\Z/N)$ is injective.
\end{proposition}

\begin{proof}
\emph{Rings essentially of finite type over a prime field.}
We first prove the assertion more generally for regular semilocal
domains $R$ essentially of finite type over a prime field $k$, with
fraction field $F$. By \cref{sub:arithmetic-finite-etale}\textup{(i)},
$R$ is the semilocalization of a regular integral finite-type
$k$-algebra, which is smooth over $k$ since $k$ is perfect
\cite[\stackstag{0B8X}]{Stacks}.
The pure-exactness form of Quillen's theorem (see Grayson's proof
in \cite[Chapter~V, Corollary~9.6.4]{Weibel}) shows that the
Gersten sequence of $R$ is pure exact in the sense of
\cite[Chapter~V, Definition~9.6.3]{Weibel}. In particular, for every $m\geq0$, the group
$K_m(R)$ is a pure subgroup of $K_m(F)$; this also holds for $m<0$
by \cite[Proposition~6.8(b)]{TT}, since both groups are zero. Thus
\[
 K_m(R)[N]\hookrightarrow K_m(F)[N],
 \qquad
 K_m(R)/N\hookrightarrow K_m(F)/N.
\]
The restriction maps form a morphism of exact sequences
\eqref{eq:coefficient-sequence}
\[
\begin{tikzcd}[column sep=small]
 0\arrow[r] & K_n(R)/N\arrow[r]\arrow[d,hook]
  & K_n(R;\Z/N)\arrow[r]\arrow[d]
  & K_{n-1}(R)[N]\arrow[r]\arrow[d,hook] & 0\\
 0\arrow[r] & K_n(F)/N\arrow[r]
  & K_n(F;\Z/N)\arrow[r]
  & K_{n-1}(F)[N]\arrow[r] & 0
\end{tikzcd}
\]
whose outer vertical maps are injective. Hence the middle vertical
map is injective.

\smallskip\noindent
\emph{The general case.}
Let $k$ be the prime field of $R$, and let $x\in K_n(R;\Z/N)$ have
zero image in $K_n(F;\Z/N)$. By \cref{sub:k-finite-continuity} and
\cref{lem:finite-type-witness} for $R_1=k$, there are a regular
local ring $R_0$ essentially of finite type over $k$, with fraction
field $F_0$, a ring map $R_0\to R$, and $x_0\in K_n(R_0;\Z/N)$ with
image $x$ in $K_n(R;\Z/N)$ and zero image in $K_n(F_0;\Z/N)$. By the
pure-exactness argument for regular local rings essentially of
finite type over $k$, $x_0=0$. Hence $x=0$.
\end{proof}

\begin{remark}[The equicharacteristic case]
\label{rem:equicharacteristic-case}
The proof of
\cref{thm:gersten-injectivity}
treats regular local rings in equal and in mixed characteristic
by the same argument. For a regular local ring $R$
containing a field, the theorem is classical. Quillen's theorem
\cite[\S7, Theorem~5.11]{Quillen}, in its pure exact form
\cite[Chapter~V, Corollary~9.6.4]{Weibel}, gives the injectivity for
rings essentially of finite type over a field, and Popescu's theorem
extends it to $R$; this is the method of Panin's proof of the
Gersten conjecture for equicharacteristic regular local rings
\cite[Theorem~A]{Panin}. \Cref{lem:equicharacteristic-injection}
records this argument, which applies to every $N\geq2$, invertible
in $R$ or not. For $N$ invertible in $R$, the proof of
\cref{thm:gersten-injectivity} uses
the following inputs of Gersten type: the
vanishing of the cohomology of Bloch's cycle complexes on essentially
smooth local rings above the weight
\cite[Corollary~4.4]{Geisser04}, the Beilinson--Lichtenbaum
comparison for the graded pieces of the slice filtration on smooth
$\Z_{(\Sigma)}$-schemes
(\cref{lem:classical-grade-normalization,rem:classical-comparison-hypotheses}),
and the generation of
$H^q_{\mathrm{\acute et}}(-,\mu_m^{\otimes q})$ by symbols
(\cref{sub:unramified-milnor-smooth-local,sub:local-symbols-via-smooth-sources}); it also uses Gabber's
alteration theorem (\cref{lem:multisection-arithmetic-alteration})
and the destackification theorem of Bergh and Rydh
\cite[\arxivhtml{1905.00872v1}{theoremx2}{Theorem~B}]{BerghRydh}.
\end{remark}

\section{Ordinary \'etale cohomology with finite coefficients}
\label{sec:ordinary-etale-injection}

\Cref{thm:generic-path-return}, applied to Eilenberg--Mac Lane
sheaves,
together with the trace of a finite locally free morphism
(\cref{sub:finite-flat-etale-trace}), gives the injectivity of
restriction to the fraction field for \'etale cohomology with
finite locally constant coefficients of invertible order, first for
regular semilocal domains essentially of finite type over a field or over
\(\Z\) whose closed residue fields have at most one positive
characteristic (\cref{prop:arithmetic-etale-injectivity}) and then for all
noetherian regular local rings
(\cref{thm:etale-injectivity}).
\Cref{subsec:semilocal-etale-injectivity} treats further semilocal rings;
\cref{subsec:height-one-residues} proves injectivity in degree two on
every connected noetherian regular semilocal ring by a residue
calculation independent of \cref{thm:generic-path-return}.

\subsection{Trace and generic injectivity}

\begin{sublemma}[Trace for finite locally free morphisms]
\label{sub:finite-flat-etale-trace}
Let $f:T\to S$ be a finite locally free morphism of schemes of
constant rank $e$, and let $M$ be a finite locally constant sheaf of
abelian groups on $S_{\mathrm{\acute et}}$.
\begin{enumerate}
\item There is a morphism of sheaves
$\operatorname{tr}_f:f_*f^*M\to M$ whose composite with the unit
$M\to f_*f^*M$ is multiplication by $e$. If $f$ is finite \'etale,
this morphism is the usual finite \'etale trace.
\item For every $q\geq0$, the composite
\[
 H^q_{\mathrm{\acute et}}(S,M)\xrightarrow{\ f^*\ }H^q_{\mathrm{\acute et}}(T,f^*M)
 \simeq H^q_{\mathrm{\acute et}}(S,f_*f^*M)
 \xrightarrow{\ \operatorname{tr}_f\ }H^q_{\mathrm{\acute et}}(S,M)
\]
is multiplication by $e$.
\end{enumerate}
\end{sublemma}

\begin{proof}
The trace is the finite locally free case of
\cite[\stackstag{0GKI}]{Stacks}; we give its construction for
finite locally constant coefficients.

\textup{(i)} Let $V\to S$ be an \'etale morphism such that $M|_V$ is
the constant sheaf with value a finite abelian group $A$, and put
$T_V=T\times_SV$. A section $s$ of $f_*f^*M$ over $V$ is a locally
constant function $s:T_V\to A$, and
$T_V=\coprod_{a\in A}s^{-1}(a)$ with $s^{-1}(a)$ open and closed.
The $\mathcal O_V$-algebra of $s^{-1}(a)$ is a direct factor of the
finite locally free $\mathcal O_V$-algebra of $T_V$, so
$s^{-1}(a)\to V$ is finite locally free. Let $e_a(s):V\to\mathbf N$
be its rank, a locally constant function, and put
\begin{equation}\label{eq:finite-flat-trace}
 \operatorname{tr}_V(s)=\sum_{a\in A}e_a(s)\,a\in M(V).
\end{equation}
For two sections $s,s'$, let $e_{a,a'}$ be the rank of
$s^{-1}(a)\cap s'^{-1}(a')$ over $V$. Then
$e_a(s)=\sum_{a'}e_{a,a'}$, $e_{a'}(s')=\sum_ae_{a,a'}$ and
$e_c(s+s')=\sum_{a+a'=c}e_{a,a'}$, so
\[
 \operatorname{tr}_V(s+s')=\sum_{a,a'}e_{a,a'}\,(a+a')
 =\operatorname{tr}_V(s)+\operatorname{tr}_V(s').
\]
The ranks commute with base change along \'etale morphisms $V'\to V$,
and $\operatorname{tr}_V(\sigma\circ s)=\sigma(\operatorname{tr}_V(s))$
for every automorphism $\sigma$ of $A$. A change of trivialization
is locally constant with values in $\operatorname{Aut}(A)$, so the
same identity holds locally on the overlap. Hence the maps
\eqref{eq:finite-flat-trace} do not depend on the trivialization of
$M|_V$ and are compatible with restriction, and they define a
morphism of sheaves $\operatorname{tr}_f$. If $s$ is the constant
function with value $a$, then $e_a(s)=e$ and
$\operatorname{tr}_V(s)=e\,a$. If $f$ is finite \'etale, an
\'etale refinement of $V$ splits $T_V$ into $e$ copies of $V$;
\eqref{eq:finite-flat-trace} is then the sum of the sections on
these copies. This characterizes the finite \'etale trace
\cite[\stackstag{03SE} and the characterization following it]{Stacks}.

\smallskip\noindent
\textup{(ii)} The morphism $f$ is finite, so $R^qf_*(f^*M)=0$ for
$q\geq1$ \cite[\stackstag{03QP}]{Stacks}. This gives the isomorphism
of the statement, under which $f^*$ is induced by the unit
$M\to f_*f^*M$. Apply $H^q_{\mathrm{\acute et}}(S,-)$ to \textup{(i)}.
\end{proof}

\begin{sublemma}[Direct image under a finite \'etale morphism]
\label{sub:finite-etale-coefficient-pushforward}
Let $R$ be a local domain with fraction field $F$, let
$f:\Spec A\to\Spec R$ be a finite \'etale morphism, and let $M$ be a
finite locally constant abelian sheaf on $\Spec A$ whose stalks have
orders invertible in $A$.
\begin{enumerate}
\item The sheaf $f_*M$ on $\Spec R$ is finite locally constant, and
its stalks have orders invertible in $R$.
\item For every $q\geq0$ there is a commutative square
\[
\begin{tikzcd}[column sep=large]
H^q_{\mathrm{\acute et}}(A,M) \ar[r] \ar[d,"\sim"'] &
H^q_{\mathrm{\acute et}}(A\otimes_R F,M) \ar[d,"\sim"] \\
H^q_{\mathrm{\acute et}}(R,f_*M) \ar[r] &
H^q_{\mathrm{\acute et}}(F,(f_*M)|_F)
\end{tikzcd}
\]
whose horizontal maps are the restriction maps and whose vertical
maps are isomorphisms.
\end{enumerate}
\end{sublemma}

\begin{proof}
If $A=0$, then $f_*M=0$ and all four groups in \textup{(ii)} vanish;
suppose that $A\neq0$.

\smallskip\noindent
\textup{(i)} The finite local constancy assertion is
\cite[\stackstag{095B}]{Stacks}. We recall its local description
and check the orders of the stalks.
By \cite[\stackstag{04HN}]{Stacks}, applied at the
closed point of $\Spec R$, there are a surjective \'etale morphism
$V\to\Spec R$ and a decomposition
\[
 \Spec A\times_RV=\coprod_{j=1}^dV_j
\]
in which every $V_j\to V$ is an isomorphism. After refining $V$, the
sheaf $M|_{V_j}$ is constant for every $j$. Then
\begin{equation}\label{eq:finite-etale-pushforward-local}
 (f_*M)|_V\simeq\prod_{j=1}^dM|_{V_j},
\end{equation}
where $M|_{V_j}$ is regarded as a sheaf on $V$. Hence $f_*M$ is
finite locally constant. At a geometric closed point of $\Spec R$,
\eqref{eq:finite-etale-pushforward-local} identifies the stalk of
$f_*M$ with the product of the stalks of $M$ over that point. Each
factor has order invertible in the residue field of $R$, so the
product has order invertible in $R$. Since $\Spec R$ is connected,
the orders of the stalks of $f_*M$ are constant.

\smallskip\noindent
\textup{(ii)} Let $f_F:\Spec(A\otimes_RF)\to\Spec F$ be the base
change of $f$. The base change map
$(f_*M)|_F\to(f_F)_*(M|_{A\otimes_RF})$ is an isomorphism, because
$f$ is finite \cite[\stackstag{0959}]{Stacks}. Since $f$ and $f_F$ are
finite, $R^qf_*M=0$ and $R^q(f_F)_*(M|_{A\otimes_RF})=0$ for
$q\geq1$ \cite[\stackstag{03QP}]{Stacks}. This gives the vertical
isomorphisms, and the square commutes by the naturality of the base
change map.
\end{proof}

\begin{sublemma}[Descending finite locally constant coefficients]
\label{sub:etale-coefficient-witness}
Let $k$ be a ring, let $R$ be a $k$-algebra, let $m\geq2$ be a
natural number invertible in $R$, and let $M$ be a finite locally
constant abelian sheaf on $\Spec R$ with $mM=0$.
\begin{enumerate}
\item There are a finitely generated $k$-subalgebra
$R_1\subseteq R$ with $m\in R_1^\times$, a finite locally constant
abelian sheaf $M_1$ on $\Spec R_1$ with $mM_1=0$, and an
isomorphism $M_1|_{\Spec R}\simeq M$.
\item For $R_1$ and $M_1$ as in \textup{(i)} and every $q\geq0$, the
functor
\begin{equation}\label{eq:etale-coefficient-functor}
 B\longmapsto H^q_{\mathrm{\acute et}}(B,M_1|_{\Spec B})
\end{equation}
from $R_1$-algebras to abelian groups commutes with filtered
colimits.
\end{enumerate}
\end{sublemma}

\begin{proof}
Assertion \textup{(i)} is the case of
\(\Z/m\)-module sheaves in
\cite[\stackstag{0GL2}(3)]{Stacks}. We give the descent in terms of
finite \'etale group schemes.

\textup{(i)} By \cite[\stackstag{03RV}]{Stacks}, the sheaf of sets
underlying $M$ is represented by a finite \'etale $R$-scheme $G$. By
the Yoneda lemma, the addition, the inverse and the zero section of
$M$ make $G$ a commutative group object of the category
$\textit{F\'Et}_R$ of finite \'etale $R$-schemes, and
$m\cdot\mathrm{id}_G$ is the zero endomorphism. The ring $R$ is the
union of the directed set of its finitely generated $k$-subalgebras
$R_\lambda$ that contain $1/m$. By \cite[\stackstag{0BTV}]{Stacks},
\[
 \varinjlim_\lambda\textit{F\'Et}_{R_\lambda}
   \xrightarrow{\ \sim\ }\textit{F\'Et}_R ,
\]
and the base change functors preserve finite products. Hence there
are an index $\lambda$, a finite \'etale $R_\lambda$-scheme
$G_\lambda$ with $G_\lambda\times_{R_\lambda}R\simeq G$, and
morphisms $G_\lambda\times_{R_\lambda}G_\lambda\to G_\lambda$,
$G_\lambda\to G_\lambda$ and $\Spec R_\lambda\to G_\lambda$ inducing
the group structure of $G$. After enlarging $\lambda$, these
morphisms satisfy the axioms of a commutative group object and
$m\cdot\mathrm{id}_{G_\lambda}=0$, because these finitely many
equalities of morphisms hold after base change to $R$. Put
$R_1=R_\lambda$, and let $M_1$ be the sheaf of abelian groups
represented by $G_\lambda$. It is finite locally constant
\cite[\stackstag{03RV}]{Stacks}, $mM_1=0$, and
$M_1|_{\Spec R}$ is represented by $G$.

\smallskip\noindent
\textup{(ii)} This is \cite[\stackstag{09YQ}]{Stacks}, together with
\cite[\stackstag{0032}]{Stacks}.
\end{proof}

\begin{proposition}[Semilocal domains essentially of finite type]
\label{prop:arithmetic-etale-injectivity}
Let $R$ be a regular semilocal domain essentially of finite type
over a field or over $\Z$
\textup{(\cref{def:arithmetic-semilocal})}, with fraction field $F$,
whose closed residue fields have at most one positive
characteristic. Let $M$ be a finite locally constant abelian sheaf
on $\Spec R$ whose stalks have orders invertible in $R$. For every
$q\geq0$, the restriction map
\[
 H^q_{\mathrm{\acute et}}(R,M)\longrightarrow H^q_{\mathrm{\acute et}}(F,M|_F)
\]
is injective.
\end{proposition}

\begin{proof}
\emph{Reduction to primary coefficients.}
If $M=0$, the assertion is immediate. Otherwise the common order of
the stalks of $M$ on the connected scheme $\Spec R$ is an invertible
annihilator $m\geq2$. The
decomposition of $\Z/m$ into its primary factors decomposes $M$ into
a direct sum of subsheaves $M_{(\ell)}$, where $\ell$ runs through
the prime divisors of $m$ and $M_{(\ell)}$ is killed by a power of
$\ell$. Hence we may assume that $\ell^\nu M=0$ for a prime number
$\ell\in R^\times$ and $\nu\geq1$.

\smallskip\noindent
\emph{Two sections.}
Let $HM$ be the Eilenberg--Mac Lane sheaf of spectra of $M$ on the
small \'etale site of $\Spec R$ (\cref{not:finite-etale-topoi}). The
sheaf $\Sigma^qHM$ has homotopy sheaf $M$ in degree $q$ and zero
homotopy sheaves in all other degrees, so it is hypercomplete \cite[Lemma~6.5.2.9]{HTT}. For
every $R$-scheme $Y$,
\begin{equation}\label{eq:em-sections-components}
 \pi_0\,\Omega^\infty R\Gamma_{\mathrm{\acute et}}(Y,(\Sigma^qHM)|_Y)
   =H^q_{\mathrm{\acute et}}(Y,M|_Y).
\end{equation}
Let $x\in H^q_{\mathrm{\acute et}}(R,M)$ have zero image in $H^q_{\mathrm{\acute et}}(F,M|_F)$. Let
$z_0\in\Omega^\infty R\Gamma_{\mathrm{\acute et}}(R,\Sigma^qHM)$ be a point of
the connected component $x$, and let $z_1$ be the zero section. By
\eqref{eq:em-sections-components} for $Y=\Spec F$, the points
$z_0|_F$ and $z_1|_F$ lie in the same connected component.

\smallskip\noindent
\emph{Return and trace.}
By \cref{thm:generic-path-return}, there is a regular scheme
$T$, finite locally free over $R$ of rank $e$ prime to $\ell$, such
that $z_0|_T$ and $z_1|_T$ lie in the same connected component. By
\eqref{eq:em-sections-components} for $Y=T$, the image of $x$ in
$H^q_{\mathrm{\acute et}}(T,M|_T)$ is zero. By
\cref{sub:finite-flat-etale-trace}\textup{(ii)}, $e\,x=0$. Since
$\ell^\nu x=0$ and $\gcd(e,\ell^\nu)=1$, we get $x=0$.
\end{proof}

\begin{theorem}[Finite locally constant coefficients]
\label{thm:etale-injectivity}
Let $R$ be a noetherian regular local ring with fraction field $F$,
and let $A$ be a finite \'etale $R$-algebra. Let $M$ be a finite
locally constant abelian sheaf on $\Spec A$ whose stalks have orders
invertible in $A$. For every $q\geq0$, the restriction map
\[
 H^q_{\mathrm{\acute et}}(A,M)\longrightarrow H^q_{\mathrm{\acute et}}(A\otimes_R F,M|_{A\otimes_RF})
\]
is injective.
\end{theorem}

\begin{proof}
By \cref{sub:finite-etale-coefficient-pushforward}, we may assume
that $A=R$. If $M=0$, there is nothing to prove. Assume $M\neq0$.
As in the proof of
\cref{prop:arithmetic-etale-injectivity}, we may assume that
$\ell^\nu M=0$ for a prime number $\ell\in R^\times$ and $\nu\geq1$. Let
$x\in H^q_{\mathrm{\acute et}}(R,M)$ have zero image in $H^q_{\mathrm{\acute et}}(F,M|_F)$, and let $k$ be as
in \cref{lem:finite-type-witness}. By
\cref{sub:etale-coefficient-witness} for $m=\ell^\nu$, there are a
finitely generated $k$-subalgebra $R_1\subseteq R$ with
$\ell\in R_1^\times$ and a finite locally constant abelian sheaf
$M_1$ on $\Spec R_1$ with $\ell^\nu M_1=0$ and $M_1|_{\Spec R}\simeq M$,
and the functor \eqref{eq:etale-coefficient-functor} commutes with
filtered colimits. By \cref{lem:finite-type-witness} for this
functor, there are a regular local ring $R_0$ essentially of finite
type over $k$, with fraction field $F_0$, ring maps $R_1\to R_0\to R$
and a class
\[
 x_0\in H^q_{\mathrm{\acute et}}(R_0,M_0),\qquad M_0=M_1|_{\Spec R_0},
\]
with image $x$ in $H^q_{\mathrm{\acute et}}(R,M)$ and zero image in
$H^q_{\mathrm{\acute et}}(F_0,M_0|_{F_0})$. The ring $R_0$ and the sheaf $M_0$ satisfy the
hypotheses of \cref{prop:arithmetic-etale-injectivity}, so
$x_0=0$. Hence $x=0$.
\end{proof}

\begin{remark}[Smooth schemes over a Dedekind ring]
\label{rem:etale-dedekind-literature}
Let $B$ be a Dedekind ring, let $X$ be a smooth $B$-scheme, let
$m\geq2$ be a natural number invertible on $B$, and let
$\varepsilon:X_{\mathrm{\acute et}}\to X_{\mathrm{Zar}}$ be the
morphism of sites. We may work Zariski locally on equidimensional
smooth open subschemes, as in the notation of \cite[p.~775]{Geisser04}.
For integers $0\leq q\leq n$, Geisser proves that
the Gersten complex
\begin{equation}\label{eq:geisser-gersten-resolution}
 0\to R^q\varepsilon_*\mu_m^{\otimes n}\to
 \bigoplus_{x\in X^{(0)}}i_{x*}
   H^q_{\mathrm{\acute et}}(\kappa(x),\mu_m^{\otimes n})\to
 \bigoplus_{x\in X^{(1)}}i_{x*}
   H^{q-1}_{\mathrm{\acute et}}(\kappa(x),\mu_m^{\otimes(n-1)})
 \to\cdots
\end{equation}
of Zariski sheaves on $X$ is exact
\cite[Theorem~1.2(2), (4), (5) and p.~775]{Geisser04}, assuming the
Bloch--Kato conjecture, which is the norm-residue theorem of Rost and
Voevodsky \cite{Voevodsky11,HaesemeyerWeibel19}. The restriction
$q\leq n$ can be removed. Let $e$ be the exponent of
$(\Z/m)^\times$. \'Etale locally on $X$, a primitive $m$-th root of
unity $\zeta$ defines a section $\zeta^{\otimes e}$ of
$\mu_m^{\otimes e}$, and $(\zeta^u)^{\otimes e}=u^e\zeta^{\otimes e}
=\zeta^{\otimes e}$ for $u\in(\Z/m)^\times$; these sections glue to
an isomorphism $\Z/m\xrightarrow{\sim}\mu_m^{\otimes e}$ on
$X_{\mathrm{\acute et}}$. Let $b\in\Z$ and $q\geq0$, and choose
$n\geq\max\{q,b\}$ with $n\equiv b\pmod e$. Tensoring with the
$((n-b)/e)$-th power of this isomorphism identifies
$\mu_m^{\otimes(b-c)}$ with $\mu_m^{\otimes(n-c)}$ for every
$c\in\Z$, compatibly with the localization sequences and purity
isomorphisms that define the maps of
\eqref{eq:geisser-gersten-resolution}. Hence the Gersten complex of
$R^q\varepsilon_*\mu_m^{\otimes b}$ is exact for all $q\geq0$ and
$b\in\Z$. In particular the case $A=R=\mathcal O_{X,x}$ of
\cref{thm:etale-injectivity} with coefficients
$\mu_m^{\otimes b}$ follows from this resolution,
since $H^q_{\mathrm{\acute et}}(\mathcal O_{X,x},\mu_m^{\otimes b})$
is the stalk of $R^q\varepsilon_*\mu_m^{\otimes b}$ at $x$
\cite[\stackstag{09YQ}]{Stacks}. Over an excellent discrete valuation
ring, L\"uders treats these coefficients for prime-power moduli within
an axiomatic Bloch--Ogus theory
\cite[\arxivhtml{2201.11553v4}{S2.Thmtheorem9}{Example~2.9} and \arxivhtml{2201.11553v4}{S3.Thmtheorem6}{Theorem~3.6}]{Luders22}.
For finite locally constant abelian sheaves pulled back from a discrete
valuation ring, with order invertible on that ring, see also
\cite[Theorem~1.1]{GaidaiTurlov} for full Gersten exactness over
regular local rings geometrically regular over the base.
\end{remark}
\begin{remark}[Failure for constructible coefficients]
\label{rem:constructible-coefficient-counterexample}
Finite constructible coefficients can fail to have generic injectivity
even in degree one and even when they embed into a finite locally
constant sheaf. Let $R=\mathbf Q[t]_{(t)}$, $X=\Spec R$ and
$Y=\Spec R[z]/(z^2-1-t)$. The morphism $f:Y\to X$ is finite
\'{e}tale because $2z$ is a unit. The polynomial is irreducible
over $\mathbf Q(t)$, since the valuation at $1+t$ takes odd value on
$1+t$. The coordinate ring of $Y$ is free of rank two over $R$,
so it embeds into its generic fibre, which is a field. Thus $Y$
is integral (and regular since it is \'etale over $X$). Its closed
fibre consists of
two $\mathbf Q$-points, given by $z=1$ and $z=-1$.
Let $M=f_*\mathbf F_\ell$, and denote the closed and generic
immersions by
\[
 i:\operatorname{Spec}\mathbf Q\hookrightarrow X,
 \qquad
 j:\operatorname{Spec}\mathbf Q(t)\hookrightarrow X.
\]
The extension by
zero $C=j_!j^*M$ is finite constructible. The sequence
\[
 0\longrightarrow C\longrightarrow M\longrightarrow
 i_*i^*M\longrightarrow0
\]
is exact on geometric stalks. In cohomology, the map
$H^0_{\mathrm{\acute et}}(X,M)\to H^0_{\mathrm{\acute et}}(X,i_*i^*M)$ is the diagonal
$\mathbf F_\ell\to\mathbf F_\ell^2$. Its cokernel therefore
injects into $H^1_{\mathrm{\acute et}}(X,C)$ and restricts to zero at the generic
point. It is the entire kernel. Indeed, generic restriction for
$M$ in degree one is injective: by the finite direct-image calculation
\cite[\stackstag{03QP}]{Stacks}, this is the degree-one restriction
for the constant sheaf on $Y$, and the map of fundamental groups from
the function field of the normal connected scheme $Y$ is surjective
\cite[\stackstag{0BQI}]{Stacks}.
\end{remark}

\subsection{Further semilocal cases}
\label{subsec:semilocal-etale-injectivity}

\begin{lemma}[Approximation of a regular semilocal domain]
\label{lem:semilocal-regular-approximation}
Let \(A\to R\) be a regular homomorphism of noetherian rings, where
\(A\) is a regular ring and \(R\) is a semilocal domain. Then \(R\)
is regular, and there is a diagram \((B_i)_{i\in I}\) of rings,
indexed by a directed set, with colimit \(R\) and with the following
properties for every \(i\in I\).
\begin{enumerate}
\item The ring \(B_i\) is a regular semilocal domain, and it is the
semilocalization of an integral smooth finite-type \(A\)-algebra at
finitely many pairwise incomparable prime ideals. In particular, if
\(A\) is essentially of finite type over a ring \(B_0\), then so is
\(B_i\).
\item Every closed residue field of \(B_i\) is a subfield of a
closed residue field of \(R\).
\end{enumerate}
\end{lemma}

\begin{proof}
The ring \(R\) is regular by \cite[\stackstag{0H7S}]{Stacks}. By
Popescu's theorem \cite[\stackstag{07GC}]{Stacks} and
\cite[\stackstag{0032}]{Stacks}, there is a diagram
\((P_i)_{i\in I}\) of smooth \(A\)-algebras, indexed by a directed
set, with colimit \(R\). Let \(S_i\subseteq P_i\) and
\(B_i=S_i^{-1}P_i\) be as in
\cref{sub:localized-filtered-presentation}. By
\cref{sub:localized-filtered-presentation}\textup{(i)},
\(R=\varinjlim_iB_i\), and
\cref{sub:localized-filtered-presentation}\textup{(ii)} gives
\textup{(ii)}. The ring \(P_i\) is regular
\cite[\stackstag{07NF}]{Stacks} and noetherian, hence a finite
product of regular domains. Since \(R\) is a domain, exactly one of
the corresponding idempotents \(\epsilon\in P_i\) has image \(1\) in
\(R\). Then \(\epsilon\in S_i\), so \(\epsilon=1\) in \(B_i\), and
\(B_i\) is the localization of the regular domain
\(P'_i=\epsilon P_i\) at the image of \(S_i\). The ring \(P'_i\) is
an integral smooth finite-type \(A\)-algebra. By
\cref{sub:localized-filtered-presentation}\textup{(ii)}, the ring
\(B_i\) is the semilocalization of \(P'_i\) at finitely many
pairwise incomparable prime ideals.
\end{proof}

\begin{sublemma}[Unramified regular rings over \(\Z_{(p)}\)]
\label{sub:unramified-regular-over-zp}
Let \(p\) be a prime number and let \(R\) be a noetherian regular
semilocal domain of characteristic zero whose closed residue fields
have characteristic \(0\) or \(p\). Suppose that
\(p\notin\mathfrak m^2R_{\mathfrak m}\) for every maximal ideal
\(\mathfrak m\) of \(R\) containing \(p\). Then every prime number
\(q\neq p\) is a unit of \(R\), and the induced ring map
\(\Z_{(p)}\to R\) is regular.
\end{sublemma}

\begin{proof}
A prime number \(q\neq p\) lies in no maximal ideal of \(R\), because
the closed residue fields have characteristic \(0\) or \(p\); hence
\(q\in R^\times\), and \(\Z\to R\) extends to \(\Z_{(p)}\to R\). The
domain \(R\) has characteristic zero, so it is a torsion-free, hence
flat, module over the discrete valuation ring \(\Z_{(p)}\)
\cite[\stackstag{0539}]{Stacks}. The generic fibre
\(R\otimes_{\Z_{(p)}}\mathbf Q\) is a localization of \(R\), hence
regular \cite[\stackstag{00OD}]{Stacks}. Let \(\mathfrak q\) be a
prime ideal of \(R\) containing \(p\), and let
\(\mathfrak m\supseteq\mathfrak q\) be a maximal ideal. Since
\(p\notin\mathfrak m^2R_{\mathfrak m}\), the element \(p\) belongs to
a minimal system of generators of \(\mathfrak mR_{\mathfrak m}\), so
\(R_{\mathfrak m}/pR_{\mathfrak m}\) is a regular local ring
\cite[\stackstag{00NQ}]{Stacks}, and so is its localization
\((R/pR)_{\mathfrak q}\) \cite[\stackstag{0AFS}]{Stacks}. Hence the
special fibre \(R/pR\) is regular. The fields \(\mathbf Q\) and
\(\mathbf F_p\) are perfect, so both fibres are geometrically regular
\cite[\stackstag{0381}, \stackstag{0382}]{Stacks}, and
\(\Z_{(p)}\to R\) is regular \cite[\stackstag{07BZ}]{Stacks}.
\end{proof}

\begin{corollary}[Semilocal rings with a regular model: \'etale cohomology]
\label{cor:semilocal-regular-model-etale-injectivity}
Let \(R\) be a noetherian semilocal domain with fraction field
\(F\), whose closed residue fields have at most one positive
characteristic. Let \(M\) be a finite locally constant abelian sheaf
on \(\operatorname{Spec}R\) whose stalks have orders invertible in
\(R\).
\begin{enumerate}
\item Suppose that there is a regular homomorphism \(A\to R\) from a
regular ring \(A\) essentially of finite type over a field or over
\(\Z\). Then, for every \(q\geq0\), the restriction map
\(H^q_{\mathrm{\acute et}}(R,M)\to H^q_{\mathrm{\acute et}}(F,M|_F)\) is injective.
\item Let \(p\) be a prime number. Suppose that \(R\) is regular of
characteristic zero, that every closed residue field of \(R\) has
characteristic \(0\) or \(p\), and that
\(p\notin\mathfrak m^2R_{\mathfrak m}\) for every maximal ideal
\(\mathfrak m\) of \(R\) containing \(p\). Then the conclusion of
\textup{(i)} holds.
\end{enumerate}
\end{corollary}

\begin{proof}
\textup{(i)} If $M=0$, the assertion is immediate. Assume $M\neq0$.
As in the proof of
\cref{prop:arithmetic-etale-injectivity}, we may assume that
\(\ell^\nu M=0\) for a prime number \(\ell\in R^\times\) and
\(\nu\geq1\). Let \(x\in H^q_{\mathrm{\acute et}}(R,M)\) have zero image in
\(H^q_{\mathrm{\acute et}}(F,M|_F)\). By \cref{sub:etale-coefficient-witness} for
\(k=\Z\) and \(m=\ell^\nu\), there are a finitely generated subring
\(R_1\subseteq R\) with \(\ell\in R_1^\times\) and a finite locally
constant abelian sheaf \(M_1\) on \(\Spec R_1\) with
\(\ell^\nu M_1=0\) and \(M_1|_{\Spec R}\simeq M\), and the functor
\eqref{eq:etale-coefficient-functor} commutes with filtered
colimits. Let \((B_i)_{i\in I}\) be as in
\cref{lem:semilocal-regular-approximation}. Since $R_1$
is finitely presented over $\Z$, \cref{sub:principal-open-witness}\textup{(i)}
gives a factorization $R_1\to B_{i_0}\to R$ of the inclusion. We
regard $(B_i)_{i\geq i_0}$ as a diagram of $R_1$-algebras with
colimit $R$ via this factorization. By
\cref{sub:principal-open-witness}\textup{(ii)}, there are $i\geq i_0$
and a class
\(x_i\in H^q_{\mathrm{\acute et}}(B_i,M_1|_{\Spec B_i})\) with image \(x\) in
\(H^q_{\mathrm{\acute et}}(R,M)\) and zero image over \(\Frac B_i\). By
\cref{lem:semilocal-regular-approximation}, the ring
\(B_i\) is a regular semilocal domain essentially of finite type
over a field or over \(\Z\), and its closed residue fields have at
most one positive characteristic; moreover \(\ell\in B_i^\times\),
because \(\ell\in R_1^\times\). By
\cref{prop:arithmetic-etale-injectivity}, \(x_i=0\). Hence
\(x=0\).

\smallskip\noindent
\textup{(ii)} By \cref{sub:unramified-regular-over-zp}, the ring map
\(\Z_{(p)}\to R\) is regular, and \(\Z_{(p)}\) is a regular ring
essentially of finite type over \(\Z\). Apply \textup{(i)} to
\(A=\Z_{(p)}\).
\end{proof}

\begin{corollary}[Semilocal rings containing a field]
\label{cor:equicharacteristic-semilocal-etale-injection}
Let \(R\) be a noetherian regular semilocal domain containing a
field, with fraction field \(F\), and let \(M\) be a finite locally
constant abelian sheaf on \(\Spec R\) whose stalks have orders
invertible in \(R\). For every \(q\geq0\), the restriction map
\(H^q_{\mathrm{\acute et}}(R,M)\to H^q_{\mathrm{\acute et}}(F,M|_F)\) is injective.
\end{corollary}

\begin{proof}
Let \(k\) be the prime field of \(R\). As in the proof of
\cref{lem:finite-type-witness}, the map \(k\to R\) is a regular
homomorphism. All closed residue fields of \(R\) have the
characteristic of \(k\). Apply
\cref{cor:semilocal-regular-model-etale-injectivity}\textup{(i)} to
\(A=k\).
\end{proof}
\begin{remark}[Equal characteristic and the Bloch--Ogus theorem]
\label{rem:etale-equicharacteristic-literature}
Let $k$ be a field, let $R$ be the semilocalization of an integral
smooth finite-type $k$-algebra at finitely many prime ideals, and
let $M$ be a finite locally constant abelian sheaf on $\Spec R$ whose
stalks have orders invertible in $k$. Panin and Zainoulline prove
that the Gersten complex of the cohomology of $M$ is exact
\cite[Corollary~1.3]{PaninZainoulline}; for the sheaves
$\mu_n^{\otimes i}$ and local rings this goes back to Bloch and Ogus
\cite{BlochOgus}. By Popescu's theorem, as in the proof of
\cref{lem:finite-type-witness}, their result implies
\cref{thm:etale-injectivity} for regular local rings
containing a field; with
\cref{lem:semilocal-regular-approximation} in place of
\cref{lem:finite-type-witness}, it also implies
\cref{cor:equicharacteristic-semilocal-etale-injection}. The proof of
\cref{thm:etale-injectivity} is the same in
equal and in mixed characteristic. It uses the generation of
$H^q_{\mathrm{\acute et}}(-,\mu_m^{\otimes q})$ by symbols
(\cref{sub:local-symbols-via-smooth-sources}), whose proof rests on
the Gersten resolution of \cite{Geisser04} for the sheaf
$\mu_m^{\otimes q}$ in degree $q$ and on the Milnor input of
\cref{sub:unramified-milnor-smooth-local}.
\end{remark}

\subsection{Height-one residues and degree two}
\label{subsec:height-one-residues}

The construction in
\cref{lem:supported-residue-lifting} prescribes one residue
on a chosen open subset. It gives the exact degree-one
row of \cref{prop:degree-one-residue-sequence} and, through
cohomology with support, the degree-two injectivity of
\cref{prop:semilocal-degree-two-injectivity}.
For a closed subset $Z\subseteq X$, with open complement
$j:U\hookrightarrow X$, let $R\underline{\Gamma}_Z M$ denote the
fibre of $M\to Rj_*j^*M$ in the derived category of abelian sheaves
on $X_{\mathrm{\acute et}}$. Write $\mathcal H^q_Z(M)$ for its
cohomology sheaves and $H^q_{\mathrm{\acute et},Z}(X,M)$ for the
cohomology of $R\Gamma_{\mathrm{\acute et}}(X,R\underline{\Gamma}_Z M)$.

\begin{lemma}[Lifting one residue on a prescribed open]
\label{lem:supported-residue-lifting}
Let $X=\Spec R$, where $R$ is a connected noetherian
regular semilocal ring with fraction field $F$. Let
$m\geq2$ be a natural number invertible in $R$, and let
$M$ be a finite locally constant abelian sheaf on $X$
with $mM=0$. Let $Z\subsetneq X$ be closed and let
$U=X\setminus Z$. For every $x\in Z\cap X^{(1)}$ and
$a\in H^0_{\mathrm{\acute et}}(\kappa(x),M(-1))$, there is a class
$\alpha\in H^1_{\mathrm{\acute et}}(U,M)$ such that
\[
 \partial_{R_{\mathfrak p_{x'}}}(\alpha|_F)=
 \begin{cases}
 a,&x'=x,\\
 0,&x'\in X^{(1)}\setminus\{x\}.
 \end{cases}
\]
Here $\mathfrak p_{x'}\subset R$ is the prime ideal
corresponding to $x'$.
\end{lemma}

\begin{proof}
Choose a geometric point over $x$. Let $G$ be the finite
monodromy image of $M(-1)$ on $X$, and let $H\subseteq G$
be the stabilizer of the value of $a$ in that geometric
fibre. Let $p:S\to X$ be the connected finite \'etale
cover associated with the transitive $G$-set $G/H$
\cite[\stackstag{0DV5}]{Stacks}.
The image of the residue-field Galois group fixes the
coset $H$. Therefore $S\times_X\Spec\kappa(x)$ has a
point $y$ for which the structural map
\[
 \kappa(x)\longrightarrow\kappa(y)
\]
is an isomorphism. With the geometric point of $S$ corresponding
to the coset $H$, its monodromy image in $G$ is $H$. The value of
$a$ is fixed by $H$, and hence defines a section
$a_S\in H^0_{\mathrm{\acute et}}(S,M(-1))$ whose value at $y$ is the image
of $a$ under this isomorphism.

The scheme $S$ is affine, regular and semilocal. It is integral
because it is connected and regular noetherian, and $y$ has
codimension one because $p$ is \'etale. Its local rings are factorial \cite[\stackstag{0AG0}]{Stacks},
so the prime divisor $\overline{\{y\}}\subset S$ is
Cartier. Its ideal is an invertible module, and is free
of rank one by \cref{sub:semilocal-projective-free}. Choose a generator
$w\in\Gamma(S,\mathcal O_S)$. Then
\[
 v_y(w)=1,\qquad v_{y'}(w)=0
       \quad(y'\in S^{(1)}\setminus\{y\}).
\]
Since $p(y)=x\in Z$ and $p^{-1}(Z)$ is closed, it contains
$\overline{\{y\}}$. The generator $w$ of the ideal of
$\overline{\{y\}}$ is therefore a unit on $p^{-1}(U)$. Its Kummer class $[w]\in H^1_{\mathrm{\acute et}}(p^{-1}(U),\mu_m)$
is defined on $p^{-1}(U)$. Using the evaluation of
\cref{not:finite-coefficient-inverse-twist}, define
\[
 \alpha=\operatorname{Tr}_{p^{-1}(U)/U}
       \bigl([w]\cup a_S|_{p^{-1}(U)}\bigr)
       \in H^1_{\mathrm{\acute et}}(U,M).
\]
For every $x'\in X^{(1)}$, the formula
\eqref{eq:residue-corestriction} for the discrete valuation ring
$R_{\mathfrak p_{x'}}$ and
\cref{sub:trait-locally-constant-residue}\textup{(ii)} give
\[
 \partial_{R_{\mathfrak p_{x'}}}(\alpha|_F)
  =\sum_{y'\in S\times_X\Spec\kappa(x')}
       \operatorname{Tr}_{\kappa(y')/\kappa(x')}
          \bigl(v_{y'}(w)\,a_S|_{\kappa(y')}\bigr).
\]
If $x'\neq x$, all summands vanish. For $x'=x$, only the
summand at $y$ can be nonzero; it is $a$ because $v_y(w)=1$ and
$\kappa(x)\to\kappa(y)$ is an isomorphism. This proves
the assertion.
\end{proof}

\begin{sublemma}[Detection of degree-two classes with principal support]
\label{sub:degree-two-support-detection}
Let $X=\Spec R$, where $R$ is a connected noetherian
regular semilocal ring. Let $m\geq2$ be a natural number
invertible in $R$, and let $M$ be a finite locally constant
abelian sheaf on $X$ with $mM=0$. For $0\ne f\in R$,
let $Z=V(f)$. The support sheaves
$\mathcal H^0_Z(M)$ and $\mathcal H^1_Z(M)$ vanish, and
the map
\[
 H^2_{\mathrm{\acute et},Z}(X,M)\longrightarrow
 \bigoplus_{x\in Z\cap X^{(1)}}
 H^2_{\mathrm{\acute et},\{\mathfrak p_x\}}(\Spec R_{\mathfrak p_x},M)
 \simeq
 \bigoplus_{x\in Z\cap X^{(1)}}
 H^0_{\mathrm{\acute et}}(\kappa(x),M(-1))
\]
is injective. The right-hand isomorphism is
\eqref{eq:trait-support-residue} for the discrete valuation rings
$R_{\mathfrak p_x}$.
\end{sublemma}

\begin{proof}
The connected regular noetherian scheme $X$ is integral.
There are finitely many height-one points in $Z$: their
prime ideals are exactly the minimal primes of $fR$.
If $f$ is a unit, $Z$ is empty and the assertions hold
with zero support groups. Otherwise proceed as follows.

\smallskip\noindent
\emph{The stalks of the support sheaves.}
Let \(\bar v\) be a geometric point of \(X\), and let \(A\) be
the strict henselization of the local ring at its image.
The stalk at \(\bar v\) of \(R\underline{\Gamma}_Z M\)
is the cohomology complex with support in \(V(f)\) on
\(\operatorname{Spec}A\). Indeed, the derived direct-image stalk formula
\cite[\stackstag{03Q9}]{Stacks}, applied to $j:X[1/f]\hookrightarrow X$,
identifies $(Rj_*j^*M)_{\bar v}$ with
$R\Gamma_{\mathrm{\acute et}}(A[1/f],M)$. Taking the fibre of
$M_{\bar v}\to (Rj_*j^*M)_{\bar v}$ gives the asserted complex,
since $M_{\bar v}=R\Gamma_{\mathrm{\acute et}}(A,M)$ by
\cite[\stackstag{03QO}]{Stacks}.
The ring \(A\) is noetherian regular local
\cite[\stackstag{06LJ}, \stackstag{06LN}]{Stacks}, and
hence factorial \cite[\stackstag{0AG0}]{Stacks}.
The pulled-back $M$ is constant and has no positive cohomology on
$\Spec A$ by \cite[\stackstag{03QO}]{Stacks}. Since $R\to A$ is flat and $A\neq0$,
the element $f$ has nonzero image in the domain $A$. The open
\(\operatorname{Spec}A[1/f]\) is therefore nonempty and connected,
so restriction on \(H^0_{\mathrm{\acute et}}\) is an isomorphism. It follows from the localization sequence that
the stalks of the support sheaves in degrees zero and one vanish,
and the stalk in degree two is
\[
 H^2_{\mathrm{\acute et},V(f)}(\operatorname{Spec}A,M)
     = H^1_{\mathrm{\acute et}}(A[1/f],M).
\]
In particular the local-to-global sequence gives
\(H^2_{\mathrm{\acute et},Z}(X,M)=\Gamma(X,\mathcal H^2_Z(M))\).

\smallskip\noindent
\emph{Detection by valuations.}
We check that $H^1_{\mathrm{\acute et}}(A[1/f],M)$ is detected by
its residues at the height-one primes of $A$ containing $f$.
Decompose the constant \(M\) into nontrivial cyclic groups of orders
\(b\mid m\), with $b\geq2$ (the empty sum if $M=0$). For each summand choose an identification
with \(\mu_b\) on \(\operatorname{Spec}A\), using a
primitive root in \(A\). Then Kummer theory \cite[\stackstag{03PL}]{Stacks} identifies that
summand of \(H^1_{\mathrm{\acute et}}(A[1/f],M)\) with
\[
     A[1/f]^\times/(A[1/f]^\times)^b.
\]
Here \(\operatorname{Pic}(A[1/f])=0\), since this localization
of a noetherian factorial ring is factorial and normal.
The units in \(A\) have \(b\)-th roots by Hensel's lemma \cite[\stackstag{04GF}]{Stacks}.
An element of \(A[1/f]^\times\) is a unit of \(A\) times a
product of integer powers of the distinct irreducible
factors of \(f\), by unique factorization. Its class is
therefore zero if its valuations at all these primes
are zero modulo \(b\). Via the localization maps to
\(A_{\mathfrak q}\), the trait residue formula of
\cref{sub:trait-locally-constant-residue}\textup{(ii)} identifies those
valuations with the images in
\(H^2_{\mathrm{\acute et},\{\mathfrak q\}}(\Spec A_{\mathfrak q},M)\).

\smallskip\noindent
\emph{Contraction to \(R\).}
Every height-one \(\mathfrak q\subset A\) containing \(f\)
contracts to a height-one prime \(\mathfrak p_x\subset R\).
In fact \(R\to A\) is flat, so the height of the contraction
is at most one by going down, and the contraction contains the
nonzero element \(f\). Restriction to \(A_{\mathfrak q}\)
then factors through \(R_{\mathfrak p_x}\), including
on cohomology with support in \(V(f)\). Thus a class in $H^2_{\mathrm{\acute et},Z}(X,M)$ whose images in
$H^2_{\mathrm{\acute et},\{\mathfrak p_x\}}(\Spec R_{\mathfrak p_x},M)$
all vanish has zero stalk in $\mathcal H^2_Z(M)$ at every
geometric point. It is zero, which proves injectivity.

\end{proof}

\begin{proposition}[The coefficient row of height-one residues]
\label{prop:degree-one-residue-sequence}
Let $R$ be a connected noetherian regular semilocal ring with fraction field
$F$, and let $M$ be a finite locally constant abelian sheaf on
$\operatorname{Spec}R$, annihilated by a natural number $m\geq2$ invertible
in $R$. Use the inverse twist of
\cref{not:finite-coefficient-inverse-twist}.
Then the finite-coefficient residue sequence
\[
\begin{aligned}
0\longrightarrow H^1_{\mathrm{\acute et}}(R,M)\longrightarrow H^1_{\mathrm{\acute et}}(F,M)
\xrightarrow{\;(\partial_{R_{\mathfrak p_x}})_x\;}
\bigoplus_{x\in (\operatorname{Spec}R)^{(1)}}
                  H^0_{\mathrm{\acute et}}(\kappa(x),M(-1))\longrightarrow0
\end{aligned}
\]
is exact.
\end{proposition}

\begin{proof}
Put $X=\Spec R$ and let $G_F$ be the absolute Galois group of $F$.

\smallskip\noindent
\emph{The residue map.}
The maps $\partial_{R_{\mathfrak p_x}}$ are those of
\cref{sub:trait-locally-constant-residue}\textup{(ii)} for the
discrete valuation rings $R_{\mathfrak p_x}$. By \'etale continuity \cite[\stackstag{09YQ}]{Stacks},
a class of $H^1_{\mathrm{\acute et}}(F,M)$ extends over some principal open $U=X[1/f]$,
where $0\ne f\in R$. Its residues vanish at every height-one point
of $U$ by \cref{sub:trait-locally-constant-residue}\textup{(ii)}.
Only the finitely many minimal primes of $fR$ can contribute, so
the residue tuple belongs to the direct sum.

\smallskip\noindent
\emph{Exactness at $H^1_{\mathrm{\acute et}}(R,M)$ and at $H^1_{\mathrm{\acute et}}(F,M)$.}
By the classification of torsors and the finite monodromy description
\cite[\stackstag{03AJ}, \stackstag{0DV5}]{Stacks},
the map $H^1_{\mathrm{\acute et}}(R,M)\to H^1_{\mathrm{\acute et}}(F,M)$
is inflation along the surjection
$G_F\twoheadrightarrow\pi_1^{\mathrm{\acute et}}(\Spec R)$
\cite[\stackstag{0BQM}]{Stacks}, hence injective.
Let $\alpha\in H^1_{\mathrm{\acute et}}(F,M)$ have all its height-one residues zero.
By \'etale continuity \cite[\stackstag{09YQ}]{Stacks}, choose $0\ne f\in R$
and $\widetilde\alpha\in H^1_{\mathrm{\acute et}}(U,M)$ mapping to $\alpha$, where
$U=X[1/f]$ and $Z=V(f)$. Localization gives the exact row and the
commutative square
\begin{equation}\label{eq:degree-one-support-diagram}
\begin{tikzcd}[column sep=large]
 H^1_{\mathrm{\acute et}}(X,M) \arrow[r] & H^1_{\mathrm{\acute et}}(U,M) \arrow[r,"\partial^X_Z"] \arrow[d] &
 H^2_{\mathrm{\acute et},Z}(X,M) \arrow[d,hook] \\
 & H^1_{\mathrm{\acute et}}(F,M) \arrow[r,"(\partial_{R_{\mathfrak p_x}})_x"'] &
 \displaystyle\bigoplus_{x\in Z\cap X^{(1)}}H^0_{\mathrm{\acute et}}(\kappa(x),M(-1)).
\end{tikzcd}
\end{equation}
Here the square commutes by naturality of localization and the right
vertical map is injective by \cref{sub:degree-two-support-detection}.
The bottom arrow sends $\alpha$ to zero, hence
$\partial^X_Z(\widetilde\alpha)=0$ and the top row lifts
$\widetilde\alpha$ to $H^1_{\mathrm{\acute et}}(X,M)$. Conversely a class extending to $X$
has zero residues by the trait localization sequences. This proves
exactness at the middle term.

\smallskip\noindent
\emph{Surjectivity.}
For $x\in(\Spec R)^{(1)}$ and
$a\in H^0_{\mathrm{\acute et}}(\kappa(x),M(-1))$, apply
\cref{lem:supported-residue-lifting} with
$Z=\overline{\{x\}}$. Restriction of the resulting class
to $F$ has residue $a$ at $x$ and zero elsewhere.
An element of this direct sum has finite support.
Adding these classes over that support proves surjectivity.
\end{proof}

\begin{proposition}[The first map in degree two on a regular semilocal ring]
\label{prop:semilocal-degree-two-injectivity}
Let \(R\) be a connected noetherian regular semilocal ring,
let \(F\) be its fraction field, and let \(M\) be a finite
locally constant abelian sheaf on \(X=\operatorname{Spec}R\),
annihilated by a natural number \(m\geq 2\) invertible in \(R\).
Then \(H^2_{\mathrm{\acute et}}(X,M)\to H^2_{\mathrm{\acute et}}(F,M)\)
is injective. Moreover, for every $0\ne f\in R$, putting $Z=V(f)$
and $U=X[1/f]$, the residue map of
\cref{sub:degree-two-support-detection} is an isomorphism, and the
localization boundary
$H^1_{\mathrm{\acute et}}(U,M)\to H^2_{\mathrm{\acute et},Z}(X,M)$
is surjective.
\end{proposition}

\begin{proof}
Fix $0\ne f\in R$ and put $Z=V(f)$ and $U=X[1/f]$.
Localization and restriction to the height-one points of $Z$
give a commutative diagram
\begin{equation}\label{eq:degree-two-support-diagram}
\begin{tikzcd}[column sep=large]
 H^1_{\mathrm{\acute et}}(U,M) \ar[r,"\partial^X_Z"] \ar[d] &
 H^2_{\mathrm{\acute et},Z}(X,M) \ar[r] \ar[d,hook] & H^2_{\mathrm{\acute et}}(X,M)\\
 H^1_{\mathrm{\acute et}}(F,M) \ar[r,"(\partial_{R_{\mathfrak p_x}})_x"'] &
 \displaystyle\bigoplus_{x\in Z\cap X^{(1)}}
 H^0_{\mathrm{\acute et}}(\kappa(x),M(-1)).
\end{tikzcd}
\end{equation}
The first row of \eqref{eq:degree-two-support-diagram} is exact by
localization, and the residue map on $H^2_{\mathrm{\acute et},Z}(X,M)$
is injective by \cref{sub:degree-two-support-detection}.
Given any tuple in
$\bigoplus_{x\in Z\cap X^{(1)}}H^0_{\mathrm{\acute et}}(\kappa(x),M(-1))$,
apply \cref{lem:supported-residue-lifting} to its finitely many
nonzero entries, with the same open $U$ for all of them. The sum
of the resulting classes is a class $\gamma\in H^1_{\mathrm{\acute et}}(U,M)$
whose generic restriction has the prescribed residues. By
\eqref{eq:degree-two-support-diagram}, the composite of $\partial^X_Z$
and the residue map is therefore surjective. Since the residue map
is injective, it is an isomorphism and $\partial^X_Z$ is surjective.

Now let $\alpha\in H^2_{\mathrm{\acute et}}(X,M)$ have zero restriction
to $F$. By \cite[\stackstag{09YQ}]{Stacks}, there is $0\ne f\in R$
such that $\alpha|_{X[1/f]}=0$. By localization, $\alpha$ comes from
$H^2_{\mathrm{\acute et},V(f)}(X,M)$, whose map to
$H^2_{\mathrm{\acute et}}(X,M)$ is zero by surjectivity of
$\partial^X_{V(f)}$ and exactness of
\eqref{eq:degree-two-support-diagram}. Thus $\alpha=0$.
\end{proof}

\begin{remark}[An alternative using purity of the branch locus]
\label{rem:degree-one-purity-proof}
The middle exactness in \cref{prop:degree-one-residue-sequence} also follows
from purity of the branch locus. Put $X=\Spec R$. For a cocycle $c$
representing a class with zero residues, let $\rho$ be the finite
coefficient monodromy. At a height-one point the class extends over the
local discrete valuation ring by the localization sequence and
\cref{sub:trait-locally-constant-residue}\textup{(i),(ii)}; hence its
restriction to the fraction field of the strict henselization is zero.
The monodromy $\rho$ is trivial on inertia, so the restricted cocycle
itself is zero there. The finite homomorphism $(c,\rho)$ to the associated
semidirect product thus kills every height-one inertia group
\cite[\stackstag{0BSW}]{Stacks}.
The normalization $Y$ of $X$ in the corresponding finite Galois extension
of $F$ is finite over $R$ by separable trace duality
\cite[\stackstag{032L}]{Stacks}; in particular no Nagata hypothesis is
needed. It is \'etale above every height-one point of $X$ by the inertia
criterion \cite[\stackstag{0BTF}]{Stacks}. Going down
and incomparability for the integral extension of the normal domain
$R$ give equality of heights of points of $Y$ and of their images.
Purity \cite[\stackstag{0BMB}]{Stacks} now makes $Y\to X$ \'etale
everywhere. Thus the cocycle factors through the \'etale fundamental
group of $X$ (\cite[\stackstag{0BQM}]{Stacks}).
The support argument \eqref{eq:degree-one-support-diagram} avoids both
normalization and branch-locus purity, and uses the same input as
\cref{prop:semilocal-degree-two-injectivity}.
\end{remark}

\Needspace{8\baselineskip}
\section{Continuous cohomology and a discrete valuation ring}
\label{sec:continuous-cohomology}

Let $\ell$ be a prime number. For a noetherian regular local ring on
which $\ell$ is invertible, the injectivity of
\cref{thm:etale-injectivity} at every finite level reduces
the kernel of restriction to the fraction field in continuous
$\ell$-adic cohomology to a kernel between first derived inverse
limits (\cref{cor:continuous-kernel-derived-limit}).
\Cref{subsec:continuous-injectivity-criteria} gives conditions under
which this kernel vanishes, for instance a finite cyclotomic image
(\cref{prop:finite-cyclotomic-image}) or, for a discrete
valuation ring, a local homomorphism from a henselian local ring
whose residue extension is algebraic with finite separable part
(\cref{prop:continuous-henselian-source}), and
\cref{subsec:continuous-valuation-counterexample} constructs an
excellent discrete valuation ring for which it is nonzero in degree
two (\cref{thm:continuous-counterexample}).

\subsection{Continuous cohomology and the Milnor sequence}
\label{subsec:continuous-milnor}

\begin{definition}[Continuous Tate cohomology]
\label{def:continuous-tate-cohomology}
Let $T$ be a scheme on which $\ell$ is invertible, and let $j\in\Z$.
The coefficient-reduction maps define an inverse system of Tate
sheaves $(\Z/\ell^n(j))_{n\geq1}$ on $T_{\mathrm{\acute et}}$, with the
conventions of \cref{sec:conventions}. Its continuous cohomology is
\[
 H^q_{\mathrm{cont}}(T,\Z_\ell(j))
   =H^q\!\left(R\varprojlim_n
       R\Gamma(T_{\mathrm{\acute et}},\Z/\ell^n(j))\right).
\]
For a ring $B$, write $H^q_{\mathrm{cont}}(B,\Z_\ell(j))$
for $H^q_{\mathrm{cont}}(\Spec B,\Z_\ell(j))$.
\end{definition}

This is Jannsen's continuous \'etale cohomology
\cite[\S3, p.~216]{Jannsen}. It agrees
with the cohomology of
$\Z_\ell(j)=\varprojlim_n\nu^*\Z/\ell^n(j)$ on the
pro-\'etale site, where $\nu$ is the morphism from the
pro-\'etale to the small \'etale topos. The coefficient-reduction
maps are surjective as maps of sheaves, so
\cite[Proposition~5.6.2]{BhattScholze} applies.

\begin{corollary}[Tate residues on a trait]
\label{cor:trait-tate-residues}
Let \(A\) be a discrete valuation ring on which \(\ell\)
is invertible, with fraction field \(L\), residue field
\(\kappa(A)\) and normalized valuation \(v_A\). Let \(n\geq1\) be a
natural number and let \(j\in\Z\).
\begin{enumerate}
\item There is an exact sequence
\begin{equation}
\label{eq:trait-tate-residue-sequence}
\begin{aligned}
 0\longrightarrow H^1_{\mathrm{\acute et}}(A,\Z/\ell^{n}(j))
 &\longrightarrow H^1_{\mathrm{\acute et}}(L,\Z/\ell^{n}(j))
 \xrightarrow{\ \partial_{A,n}\ }
 H^0_{\mathrm{\acute et}}(\kappa(A),\Z/\ell^{n}(j-1)) \\
 &\longrightarrow H^2_{\mathrm{\acute et}}(A,\Z/\ell^{n}(j))
 \longrightarrow H^2_{\mathrm{\acute et}}(L,\Z/\ell^{n}(j)),
\end{aligned}
\end{equation}
in which \(\partial_{A,n}\) is the map \(\partial_A\) of
\cref{sub:trait-locally-constant-residue}\textup{(ii)} for
\(M=\Z/\ell^n(j)\) and \(m=\ell^n\).
\item The sequences \eqref{eq:trait-tate-residue-sequence} for \(n+1\) and
for \(n\) are compatible with the reduction maps
\(\Z/\ell^{n+1}(i)\to\Z/\ell^n(i)\), \(i\in\{j,j-1\}\).
\item For \(b\in L^\times\) with Kummer class
\([b]\in H^1_{\mathrm{\acute et}}(L,\mu_{\ell^n})\) and for
\(a\in H^0_{\mathrm{\acute et}}(A,\Z/\ell^n(j-1))\),
\[
 \partial_{A,n}([b]\cup a|_L)=v_A(b)\,a|_{\kappa(A)}.
\]
\item For a finite \'{e}tale \(A\)-algebra \(B\) that is a domain,
the maps \(\partial_{A,n}\) and \(\partial_{B_y,n}\) satisfy
\eqref{eq:residue-corestriction}.
\end{enumerate}
\end{corollary}

\begin{proof}
Apply \cref{sub:trait-locally-constant-residue} to
\(M=\Z/\ell^n(j)\) and \(m=\ell^n\). Tensor evaluation
gives the canonical identification
\[
 \underline{\operatorname{Hom}}_{\Z/\ell^n}
       (\mu_{\ell^n},\Z/\ell^n(j))
 \simeq\Z/\ell^n(j-1).
\]
The localization sequence and
\cref{sub:trait-locally-constant-residue}\textup{(i)} give
\eqref{eq:trait-tate-residue-sequence}, and \textup{(iii)} and
\textup{(iv)} are
\cref{sub:trait-locally-constant-residue}\textup{(ii)} and
\textup{(iv)}.

For \textup{(ii)}, regard both
\(\Z/\ell^{n+1}(j)\) and \(\Z/\ell^n(j)\) as
\(\Z/\ell^{n+1}\)-module sheaves.
\Cref{sub:trait-locally-constant-residue}\textup{(iii)} for
\(d=\ell^n\) and \(m=\ell^{n+1}\) gives
\[
\begin{aligned}
 \underline{\operatorname{Hom}}_{\Z/\ell^{n+1}}
       (\mu_{\ell^{n+1}},\Z/\ell^n(j))
 &\simeq
 \underline{\operatorname{Hom}}_{\Z/\ell^n}
       (\mu_{\ell^n},\Z/\ell^n(j))\\
 &\simeq\Z/\ell^n(j-1).
\end{aligned}
\]
The first isomorphism is inverse to precomposition by
\(\mu_{\ell^{n+1}}\to\mu_{\ell^n}\),
\(\zeta\mapsto\zeta^\ell\). Under these identifications,
postcomposition by coefficient reduction induces the
usual reduction
\(\Z/\ell^{n+1}(j-1)\to\Z/\ell^n(j-1)\).
Indeed, \'etale locally choose a primitive root
\(\zeta_{n+1}\) and put \(\zeta_n=\zeta_{n+1}^{\ell}\).
The tensor powers of these generators, and their duals
for negative twists, identify both maps with reduction
\(\Z/\ell^{n+1}\to\Z/\ell^n\).
Naturality in the coefficient sheaf now proves the
asserted compatibility of the localization sequences.
\end{proof}

\begin{sublemma}[The Milnor sequence {\cite[(3.1), p.~216]{Jannsen}}]
\label{sub:continuous-milnor}
For every $T$ as in \cref{def:continuous-tate-cohomology}
and every pair of integers $q,j$, there is a natural short
exact sequence
\begin{equation}
\label{eq:continuous-milnor-sequence}
\begin{aligned}
0\longrightarrow \varprojlim_n^1 H^{q-1}_{\mathrm{\acute et}}(T,\Z/\ell^n(j))
&\longrightarrow H^q_{\mathrm{cont}}(T,\Z_\ell(j))\\
&\longrightarrow \varprojlim_n H^q_{\mathrm{\acute et}}(T,\Z/\ell^n(j))
\longrightarrow0.
\end{aligned}
\end{equation}
\end{sublemma}

\begin{proof}
Regard the complexes of derived sections as complexes of
\(\Z_\ell\)-modules via coefficient reduction. The spectral
sequence for the derived inverse limit has
\[
 E_2^{a,b}=\varprojlim_n^{\,a}H^b_{\mathrm{\acute et}}(T,\Z/\ell^n(j))
 \quad\Longrightarrow\quad
 H^{a+b}_{\mathrm{cont}}(T,\Z_\ell(j)).
\]
For an inverse system indexed by the positive integers,
$\varprojlim^{\,a}=0$ for $a\geq2$
\cite[\stackstag{07KW}]{Stacks}; the two-term product complex computing
the derived limit is the same for abelian groups and for
\(\Z_\ell\)-modules. Only columns zero and one remain. Their
filtration in total degree $q$ is exactly \eqref{eq:continuous-milnor-sequence}.
The construction is natural in morphisms of towers of complexes of
\(\Z_\ell\)-modules (in particular in restriction in $T$), and its
edge maps are \(\Z_\ell\)-linear.
\end{proof}

\begin{corollary}[The derived-limit obstruction]
\label{cor:continuous-kernel-derived-limit}
Let $R$ be a noetherian regular local ring on which $\ell$ is
invertible, with fraction field $F$, and let $q\in\mathbf N$ and
$j\in\Z$.
\begin{enumerate}
\item The first map of \eqref{eq:continuous-milnor-sequence} induces an
isomorphism
\begin{equation}\label{eq:continuous-kernel}
\begin{aligned}
 &\ker\Bigl(\varprojlim_n^1 H^{q-1}_{\mathrm{\acute et}}(R,\Z/\ell^n(j))
  \longrightarrow\varprojlim_n^1 H^{q-1}_{\mathrm{\acute et}}(F,\Z/\ell^n(j))\Bigr)\\
 &\qquad\xrightarrow{\ \sim\ }
 \ker\Bigl(H^q_{\mathrm{cont}}(R,\Z_\ell(j))\longrightarrow
           H^q_{\mathrm{cont}}(F,\Z_\ell(j))\Bigr).
\end{aligned}
\end{equation}
\item If $q=1$, or if $q=j+1$ and $j\geq1$, the restriction map
$H^q_{\mathrm{cont}}(R,\Z_\ell(j))\to
H^q_{\mathrm{cont}}(F,\Z_\ell(j))$ is injective.
\end{enumerate}
\end{corollary}

\begin{proof}
\textup{(i)} The sequences \eqref{eq:continuous-milnor-sequence} for $R$ and for
$F$ form a commutative diagram with exact rows
\[
\begin{tikzcd}[column sep=small]
 0\arrow[r] & \varprojlim_n^1 H^{q-1}_{\mathrm{\acute et}}(R,\Z/\ell^n(j))\arrow[r]\arrow[d]
  & H^q_{\mathrm{cont}}(R,\Z_\ell(j))\arrow[r]\arrow[d]
  & \varprojlim_n H^q_{\mathrm{\acute et}}(R,\Z/\ell^n(j))\arrow[r]\arrow[d,hook] & 0\\
 0\arrow[r] & \varprojlim_n^1 H^{q-1}_{\mathrm{\acute et}}(F,\Z/\ell^n(j))\arrow[r]
  & H^q_{\mathrm{cont}}(F,\Z_\ell(j))\arrow[r]
  & \varprojlim_n H^q_{\mathrm{\acute et}}(F,\Z/\ell^n(j))\arrow[r] & 0.
\end{tikzcd}
\]
By \cref{thm:etale-injectivity}, the map
$\varprojlim_n H^q_{\mathrm{\acute et}}(R,\Z/\ell^n(j))\to
\varprojlim_n H^q_{\mathrm{\acute et}}(F,\Z/\ell^n(j))$ is injective.
The snake lemma applied to \eqref{eq:continuous-milnor-sequence} for $R$ and $F$
gives \eqref{eq:continuous-kernel}.

\smallskip\noindent
\textup{(ii)} By \textup{(i)}, it suffices to prove that
$\varprojlim_n^1H^{q-1}_{\mathrm{\acute et}}(R,\Z/\ell^n(j))=0$. For $q=1$, the groups
$H^0_{\mathrm{\acute et}}(R,\Z/\ell^n(j))$ are finite, so their inverse system is
Mittag--Leffler and its first derived limit is zero
\cite[\stackstag{07KW}]{Stacks}. Suppose now that $q=j+1$ and
$j\geq1$. By
\cref{prop:semilocal-unit-symbols}, the group $H^j_{\mathrm{\acute et}}(R,\Z/\ell^n(j))$
is generated by Kummer symbols of units of $R$, and the symbol of
the same units in $H^j_{\mathrm{\acute et}}(R,\Z/\ell^{n+1}(j))$ maps to it. Hence the
transition maps are surjective, and the first derived limit is zero
\cite[\stackstag{07KW}]{Stacks}.
\end{proof}

\Needspace{9\baselineskip}
\subsection{Criteria for continuous injectivity}
\label{subsec:continuous-injectivity-criteria}

\begin{proposition}[Finite cyclotomic image]
\label{prop:finite-cyclotomic-image}
Let $R$ be a noetherian regular local ring with fraction
field $F$, and let $\ell$ be a prime invertible in $R$.
Suppose that the $\ell$-adic cyclotomic character
$\chi:\pi_1^{\mathrm{\acute et}}(\Spec R)\to\Z_\ell^\times$ of $R$
has finite image. Then, for all integers $q$ and $j$, the
restriction map
\[
 H^q_{\mathrm{cont}}(R,\Z_\ell(j))
 \longrightarrow H^q_{\mathrm{cont}}(F,\Z_\ell(j))
\]
is injective.
\end{proposition}

\begin{proof}
For a ring $B$ on which $\ell$ is invertible and $i\geq0$, put
\[
 D_B^i=H^i_{\mathrm{\acute et}}(B,\mathbf Q_\ell/\Z_\ell(j))
      =\varinjlim_nH^i_{\mathrm{\acute et}}(B,\Z/\ell^{n}(j)),
 \qquad
 \mathbf Q_\ell/\Z_\ell(j)=\varinjlim_n\Z/\ell^{n}(j);
\]
the second equality is \cite[\stackstag{03Q5}]{Stacks}, where the
colimit of the Tate sheaves uses the inclusions induced by multiplication
by $\ell$ on lifts under coefficient reduction. The group $D_B^i$
is $\ell$-primary torsion.

\smallskip\noindent
\emph{Divisibility over the splitting cover.}
Let $\Spec S\to\Spec R$ be the connected finite \'etale cover
corresponding to $\ker\chi$, and let $d$ be its degree. The sheaves
$\mu_{\ell^n}$ are constant on $\Spec S$, and a compatible system of
generators identifies the inverse systems $(\Z/\ell^n(j))_n$ and
$(\Z/\ell^n(i))_n$ on $\Spec S$ for every $i\in\Z$. The ring $S$ is
semilocal, and its closed residue fields have the characteristic of
the residue field of $R$. For $i\geq1$, the group
$H^i_{\mathrm{\acute et}}(S,\Z/\ell^n(i))$ is generated by Kummer symbols of units
(\cref{prop:semilocal-unit-symbols}), and the symbol of the same
units with coefficients $\Z/\ell^{n+1}(i)$ maps to it; for $i=0$ the
sheaves are constant and $\Spec S$ is connected. Hence the reduction
maps
\[
 H^i_{\mathrm{\acute et}}(S,\Z/\ell^{n+1}(j))\longrightarrow H^i_{\mathrm{\acute et}}(S,\Z/\ell^n(j))
 \qquad(i\geq0,\ n\geq1)
\]
are surjective. The composite of the reduction
$\Z/\ell^{n+1}(j)\to\Z/\ell^n(j)$ with the inclusion
$\Z/\ell^n(j)\to\Z/\ell^{n+1}(j)$ is multiplication by $\ell$.
Therefore every element in the image of $H^i_{\mathrm{\acute et}}(S,\Z/\ell^n(j))$ in
$D_S^i$ is $\ell$ times an element in the image of
$H^i_{\mathrm{\acute et}}(S,\Z/\ell^{n+1}(j))$, and
\begin{equation}\label{eq:splitting-cover-divisible}
 \ell D_S^i=D_S^i\qquad(i\geq0).
\end{equation}

\smallskip\noindent
\emph{The structure of $D_R^i$.}
The traces of \cref{sub:finite-flat-etale-trace} for the sheaves
$\Z/\ell^n(j)$ are compatible with the inclusions
$\Z/\ell^n(j)\to\Z/\ell^{n+1}(j)$, and their colimit is a map
$\operatorname{tr}:D_S^i\to D_R^i$ with
$\operatorname{tr}(x|_S)=d\,x$ for $x\in D_R^i$. By
\eqref{eq:splitting-cover-divisible}, the group
$\operatorname{tr}(D_S^i)$ is divisible, and it contains $dD_R^i$.
Let $D_{\mathrm{div}}\subseteq D_R^i$ be the maximal divisible
subgroup and put $c=v_\ell(d)$. Since $D_R^i$ is $\ell$-primary
torsion, $\ell^cD_R^i=dD_R^i\subseteq D_{\mathrm{div}}$. A divisible
subgroup of an abelian group is a direct summand, so
\begin{equation}\label{eq:divisible-plus-bounded}
 D_R^i=D_{\mathrm{div}}\oplus C,\qquad\ell^cC=0.
\end{equation}

\smallskip\noindent
\emph{Vanishing of the derived limits.}
For each $n\geq1$, the commutative diagram of sheaves with exact rows
\[
\begin{tikzcd}[column sep=small]
0\arrow[r]&\Z/\ell^{n+1}(j)\arrow[r]\arrow[d,"\mathrm{red}"']&
\mathbf Q_\ell/\Z_\ell(j)\arrow[r,"\ell^{n+1}"]\arrow[d,"\ell"]&
\mathbf Q_\ell/\Z_\ell(j)\arrow[r]\arrow[d,"\mathrm{id}"]&0\\
0\arrow[r]&\Z/\ell^n(j)\arrow[r]&
\mathbf Q_\ell/\Z_\ell(j)\arrow[r,"\ell^n"]&
\mathbf Q_\ell/\Z_\ell(j)\arrow[r]&0
\end{tikzcd}
\]
induces the transition maps in the exact sequences
\[
 0\longrightarrow D_R^{i-1}/\ell^n
 \longrightarrow H^i_{\mathrm{\acute et}}(R,\Z/\ell^{n}(j))
 \longrightarrow D_R^i[\ell^n]\longrightarrow0
 \qquad(i\geq0,\ D_R^{-1}=0),
\]
which form an exact sequence of inverse systems: the transition maps
are the projections $D_R^{i-1}/\ell^{n+1}\to D_R^{i-1}/\ell^n$ on the
left, the reductions in the middle, and multiplication by $\ell$ from
$D_R^i[\ell^{n+1}]$ to $D_R^i[\ell^n]$ on the right. The left system
has surjective transition maps. By
\eqref{eq:divisible-plus-bounded}, the right system is the sum of the
system $(D_{\mathrm{div}}[\ell^n])_n$, whose transition maps are
surjective, and of the system $(C[\ell^n])_n$, in which every
composite of $c$ transition maps is zero. Both outer systems are
Mittag--Leffler, so their first derived limits vanish
\cite[\stackstag{07KW}]{Stacks}, and the long exact sequence of
derived limits gives
\[
 \varprojlim_n^1H^i_{\mathrm{\acute et}}(R,\Z/\ell^n(j))=0\qquad(i\geq0).
\]
For $q\geq0$, the assertion now follows from
\cref{cor:continuous-kernel-derived-limit}\textup{(i)}; for $q<0$,
both groups are zero.
\end{proof}
\begin{proposition}[A henselian local source]
\label{prop:continuous-henselian-source}
Let $R$ be a discrete valuation ring with fraction field $F$
and residue field $\kappa$, and let $\ell$ be a prime
invertible in $R$. Suppose that there is a local ring
homomorphism $B\to R$ from a henselian local ring $B$
with residue field $k$. Assume that the induced residue extension
$\kappa/k$ is algebraic and that its maximal separable subextension
$k_s/k$ is finite. Then, for all integers $q$ and $j$, the
restriction map
\[
 H^q_{\mathrm{cont}}(R,\Z_\ell(j))
 \longrightarrow H^q_{\mathrm{cont}}(F,\Z_\ell(j))
\]
is injective.
\end{proposition}

\begin{proof}
\emph{Reduction to compatible sections of the residues.}
Since $B\to R$ is local, $\ell\in B^\times$.
For $i\in\Z$ and $n\geq1$,
\cref{sub:trait-locally-constant-residue}\textup{(i)} and the
localization sequence of the closed point of $\Spec R$ give a
long exact sequence
\begin{equation}\label{eq:trait-purity-sequence}
\begin{aligned}
 \cdots&\to H^i_{\mathrm{\acute et}}(R,\Z/\ell^{n}(j))\to H^i_{\mathrm{\acute et}}(F,\Z/\ell^{n}(j))
 \xrightarrow{\ \partial_n\ }H^{i-1}_{\mathrm{\acute et}}(\kappa,\Z/\ell^{n}(j-1))\\
 &\to H^{i+1}_{\mathrm{\acute et}}(R,\Z/\ell^{n}(j))\to\cdots,
\end{aligned}
\end{equation}
whose residues have the normalization of
\cref{sub:trait-locally-constant-residue}\textup{(ii)}.
Naturality in the coefficient sheaf and
\cref{sub:trait-locally-constant-residue}\textup{(iii)} make these
sequences compatible with reduction in $n$. We construct maps
\[
 s_n:H^{i-1}_{\mathrm{\acute et}}(\kappa,\Z/\ell^{n}(j-1))\longrightarrow
     H^i_{\mathrm{\acute et}}(F,\Z/\ell^{n}(j)),
 \qquad\partial_n\circ s_n=\mathrm{id},
\]
compatible with the reduction maps. Then
\eqref{eq:trait-purity-sequence} breaks into short exact sequences of
inverse systems
\[
 0\to\bigl(H^i_{\mathrm{\acute et}}(R,\Z/\ell^{n}(j))\bigr)_n\to
 \bigl(H^i_{\mathrm{\acute et}}(F,\Z/\ell^{n}(j))\bigr)_n\to
 \bigl(H^{i-1}_{\mathrm{\acute et}}(\kappa,\Z/\ell^{n}(j-1))\bigr)_n\to0
\]
that are split by $(s_n)_n$. Hence
$\varprojlim_n^1H^i_{\mathrm{\acute et}}(R,\Z/\ell^{n}(j))\to
\varprojlim_n^1H^i_{\mathrm{\acute et}}(F,\Z/\ell^{n}(j))$ is injective for every $i$,
and the assertion for $q\geq0$ follows from
\cref{cor:continuous-kernel-derived-limit}\textup{(i)}. For $q<0$,
both continuous cohomology groups vanish by \cref{sub:continuous-milnor}.

For $i\leq0$, the source of $s_n$ is zero, and we take $s_n=0$.
We henceforth assume $i\geq1$.

\smallskip\noindent
\emph{A finite \'etale cover with a rational point.}
Let $B'$ be the finite \'etale local $B$-algebra with residue field
$k_s$ \cite[\stackstag{04GK}]{Stacks}. The $R$-algebra
$R\otimes_BB'$ is finite \'etale, and the multiplication map
$\kappa\otimes_kk_s\to\kappa$ defines a point $y$ of its closed
fibre with $\kappa(y)=\kappa$. Let $\Spec S$ be the connected
component of $\Spec(R\otimes_BB')$ containing $y$, and let
$L=\Frac S$. The ring $S$ is a regular semilocal domain of dimension
one and a $B'$-algebra. The prime ideal of $y$ in $S$ is invertible,
hence free of rank one (\cref{sub:semilocal-projective-free}). Let
$w$ be a generator. Then
\[
 v_y(w)=1,\qquad v_{y'}(w)=0
\]
for every closed point $y'\neq y$ of $\Spec S$.

\smallskip\noindent
\emph{Classes over the residue field.}
The extension $\kappa/k_s$ is purely inseparable, so
$H^{i-1}_{\mathrm{\acute et}}(k_s,-)\to H^{i-1}_{\mathrm{\acute et}}(\kappa,-)$ is an isomorphism
\cite[\stackstag{04DZ}]{Stacks}. The ring $B'$ is henselian local
with residue field $k_s$, so
$H^{i-1}_{\mathrm{\acute et}}(B',-)\to H^{i-1}_{\mathrm{\acute et}}(k_s,-)$ is an isomorphism on torsion
sheaves \cite[\stackstag{09ZI}]{Stacks}. Hence every
$b\in H^{i-1}_{\mathrm{\acute et}}(\kappa,\Z/\ell^{n}(j-1))$ is the image of a unique
class $\widetilde b\in H^{i-1}_{\mathrm{\acute et}}(B',\Z/\ell^{n}(j-1))$, and
$b\mapsto\widetilde b$ commutes with the reduction maps.

\smallskip\noindent
\emph{The sections.}
Let $[w]\in H^1_{\mathrm{\acute et}}(L,\mu_{\ell^n})$ be the Kummer class of $w$, and put
\[
 s_n(b)=\operatorname{Tr}_{L/F}\bigl([w]\cup\widetilde b|_L\bigr)
 \in H^i_{\mathrm{\acute et}}(F,\Z/\ell^{n}(j)).
\]
The ring $S$ and the element $w$ do not depend on $n$.
Reduction takes $[w]_{\ell^{n+1}}$ to $[w]_{\ell^n}$ and
commutes with the lift $b\mapsto\widetilde b$, cup products
and finite \'etale trace; hence the maps $s_n$ commute with
reduction. By
\cref{sub:trait-locally-constant-residue}\textup{(ii), (iv)},
\[
 \partial_ns_n(b)=
 \sum_{y'\in\Spec(S\otimes_R\kappa)}
 \operatorname{Tr}_{\kappa(y')/\kappa}
       \bigl(v_{y'}(w)\,\widetilde b|_{\kappa(y')}\bigr)
 =\widetilde b|_{\kappa(y)}=b.
 \qedhere
\]
\end{proof}
\begin{corollary}[A finitely generated residue extension]
\label{cor:continuous-finitely-generated-residue}
Let $R$ be a discrete valuation ring with fraction field $F$
and residue field $\kappa$, and let $\ell$ be a prime
invertible in $R$. Suppose that $R$ contains a field $k$
and that $\kappa/k$, via $k\hookrightarrow R\twoheadrightarrow\kappa$,
is a finitely generated field extension. Then
\[
 H^q_{\mathrm{cont}}(R,\Z_\ell(j))
 \longrightarrow H^q_{\mathrm{cont}}(F,\Z_\ell(j))
\]
is injective for all integers $q$ and $j$.
In particular, this conclusion holds for discrete valuation
rings essentially of finite type over a field, with $\ell$
invertible in the ring.
\end{corollary}

\begin{proof}
Choose a transcendence basis $t_1,\ldots,t_d$ of
$\kappa/k$, and lift it to
$\widetilde t_1,\ldots,\widetilde t_d\in R$.
Every nonzero polynomial over $k$ in the lifted elements
has nonzero residue, so is a unit of $R$. Consequently the
field $k(\widetilde t_1,\ldots,\widetilde t_d)$ embeds in
$R$ and maps isomorphically to
$k(t_1,\ldots,t_d)\subseteq\kappa$ under reduction.
A field is a henselian local ring, and $\ell$ is invertible
in $k(\widetilde t_1,\ldots,\widetilde t_d)$ because it is
invertible in $R$. The extension
$\kappa/k(t_1,\ldots,t_d)$ is finite because $\kappa/k$
is finitely generated; hence its maximal separable
subextension is finite. Apply
\cref{prop:continuous-henselian-source} to the local map
$k(\widetilde t_1,\ldots,\widetilde t_d)\hookrightarrow R$.
\end{proof}

\Needspace{10\baselineskip}
\begin{corollary}[The degree-two obstruction at height one]
\label{cor:continuous-degree-two-obstruction}
Let \(R\) be a connected noetherian regular semilocal ring on which
\(\ell\) is invertible, with fraction field \(F\). For \(j\in\Z\),
put
\[
 \rho_j:H^2_{\mathrm{cont}}(R,\Z_\ell(j))\longrightarrow
               H^2_{\mathrm{cont}}(F,\Z_\ell(j)).
\]
\begin{enumerate}
\item The kernel of \(\rho_j\) is canonically isomorphic, as a
\(\Z_\ell\)-module, to the cokernel of
\begin{equation}\label{eq:continuous-residue-limit}
 \varprojlim_n H^1_{\mathrm{\acute et}}(F,\Z/\ell^{n}(j))
 \xrightarrow{\;(\partial_{x,n})_x\;}
 \varprojlim_n\bigoplus_{x\in (\operatorname{Spec}R)^{(1)}}
              H^0_{\mathrm{\acute et}}(\kappa(x),\Z/\ell^{n}(j-1)).
\end{equation}
\item Let \(j\neq1\). If the \((j-1)\)-st power of the \(\ell\)-adic
cyclotomic character of \(\kappa(x)\) is nontrivial for every
\(x\in(\Spec R)^{(1)}\), then \(\rho_j\) is injective.
\item If \(R\) is Nagata and the \(\ell\)-adic cyclotomic character
of every closed residue field of \(R\) has infinite image, then
\(\rho_j\) is injective for every \(j\in\Z\).
\item If the \(\ell\)-adic cyclotomic character of \(R\) has finite
image, then \(\rho_j\) is injective for every \(j\in\Z\).
\end{enumerate}
\end{corollary}

\begin{proof}
\textup{(i)}
Apply \cref{prop:degree-one-residue-sequence} to \(M=\Z/\ell^{n}(j)\).
By \cref{cor:trait-tate-residues}\textup{(ii)}, the resulting short
exact sequences
\[
 0\to H^1_{\mathrm{\acute et}}(R,\Z/\ell^{n}(j))\to H^1_{\mathrm{\acute et}}(F,\Z/\ell^{n}(j))\to
 \bigoplus_{x\in(\Spec R)^{(1)}}H^0_{\mathrm{\acute et}}(\kappa(x),\Z/\ell^{n}(j-1))\to0
\]
form an exact sequence of inverse systems. The inverse limit of the
right-hand direct sums in \eqref{eq:continuous-residue-limit} is taken
with finite support at each finite level. It embeds in the product
of the individual limits, but need not equal their direct sum.
The long exact sequence of derived limits of these short exact sequences
identifies the cokernel of
\eqref{eq:continuous-residue-limit} with the kernel of
\(\varprojlim_n^1H^1_{\mathrm{\acute et}}(R,\Z/\ell^{n}(j))\to
\varprojlim_n^1H^1_{\mathrm{\acute et}}(F,\Z/\ell^{n}(j))\). As in the proof of
\cref{cor:continuous-kernel-derived-limit}\textup{(i)}, the Milnor
sequences \eqref{eq:continuous-milnor-sequence} identify this kernel with
\(\ker\rho_j\), because the maps
\(H^2_{\mathrm{\acute et}}(R,\Z/\ell^{n}(j))\to H^2_{\mathrm{\acute et}}(F,\Z/\ell^{n}(j))\) are injective by
\cref{prop:semilocal-degree-two-injectivity}. All these
identifications are \(\Z_\ell\)-linear: the derived inverse limits
are computed by the two-term product complex in \(\Z_\ell\)-modules
as in \cref{sub:continuous-milnor}, and the identifications with
continuous cohomology use its natural edge maps.

\smallskip\noindent
\textup{(ii)}
If the \((j-1)\)-st powers of the cyclotomic characters are
nontrivial, for every \(x\) the limit of
\(H^0_{\mathrm{\acute et}}(\kappa(x),\Z/\ell^{n}(j-1))\) is zero. It is the invariant
submodule of a rank-one \(\Z_\ell\)-module on which some group
element acts by a scalar different from one. The inverse limit
of the direct sums embeds into the product of those limits
and is therefore zero. By \textup{(i)}, \(\rho_j\) is injective.

\smallskip\noindent
\textup{(iii)}
Suppose that \(R\) is Nagata and that the cyclotomic character
of each closed residue field has infinite image. If the
cyclotomic character of a height-one residue field had finite
image, some finite extension of that field would contain all
\(\ell\)-power roots of unity. The normalization of
\(R/\mathfrak p_x\) in this extension is finite by the Nagata
hypothesis. All those roots lie in the normalization by integrality.
Choose a maximal ideal \(\mathfrak m\subset R\) containing
\(\mathfrak p_x\). At any maximal ideal of the normalization
lying above \(\mathfrak m/\mathfrak p_x\), the roots of each
order remain distinct:
the polynomial \(t^{\ell^n}-1\) is separable on the residue field
and is the product of its linear factors in the normalization.
This puts all \(\ell\)-power roots of unity in one finite
extension of the closed residue field, a contradiction. Thus
every height-one cyclotomic image is infinite, and every nonzero
power of its character is nontrivial; for \(j\neq1\) the assertion
follows from \textup{(ii)}. For \(j=1\),
\cref{sub:semilocal-projective-free} gives
\(\operatorname{Pic}(R)=0\). The Kummer sequence \cite[\stackstag{03PL}]{Stacks} gives
\[
 H^1_{\mathrm{\acute et}}(R,\mu_{\ell^n})
 \cong R^\times/(R^\times)^{\ell^n}
\]
for every \(n\geq1\).
Under these isomorphisms, coefficient reduction is the
surjection
\[
 R^\times/(R^\times)^{\ell^{n+1}}
 \longrightarrow R^\times/(R^\times)^{\ell^n}.
\]
Consequently
\(\varprojlim_n^1 H^1_{\mathrm{\acute et}}(R,\mu_{\ell^n})=0\).
The Milnor sequence \eqref{eq:continuous-milnor-sequence} and the
finite-coefficient injectivity of
\cref{prop:semilocal-degree-two-injectivity}
therefore give injectivity in continuous degree two.

\smallskip\noindent
\textup{(iv)}
Suppose that the cyclotomic character of \(R\) has finite
image. For any \(j\), choose the connected finite \'{e}tale
Galois cover corresponding to the kernel of the action on
\(\Z_\ell(j-1)\), and fix a compatible basis over that cover.
At each height-one point \(x\) with nonzero inverse-limit invariant
group, the action on that lattice is trivial. The cover is split
over \(\kappa(x)\), and there is a point \(y\) above \(x\) of
residue degree one. Choose one principal generator of \(y\) on
the cover as in \cref{lem:supported-residue-lifting}. Its Kummer
classes cupped with the compatible basis and transferred to \(F\)
give a compatible sequence whose residue at \(x\) is a basis of
the invariant lattice, and whose residues elsewhere are zero.
An element in the limit of the residue groups is a family of
coefficients in these invariant lattices, with only finitely many
nonzero reductions at each coefficient level. At level \(n\), add the
chosen lifts for the nonzero reductions, with their prescribed
coefficients modulo \(\ell^n\). These sums are finite and compatible
with reduction. They give an element of
\(\varprojlim_nH^1_{\mathrm{\acute et}}(F,\Z/\ell^{n}(j))\) with the prescribed residues.
Hence \eqref{eq:continuous-residue-limit} is surjective, and
\(\rho_j\) is injective by \textup{(i)}.
\end{proof}

\subsection{A discrete valuation ring with a continuous kernel}
\label{subsec:continuous-valuation-counterexample}

\begin{construction}[A cyclotomic power-series valuation]
\label{con:cyclotomic-power-series-valuation}
Let $u$ be an indeterminate. Choose roots of unity
\(\zeta_i\in\mathbf Q(\mu_{\ell^\infty})\) of order \(\ell^i\) for
\(i\geq0\) such that \(\zeta_{i+1}^\ell=\zeta_i\) (thus \(\zeta_0=1\)).
In the power series ring \(\mathbf Q(\mu_{\ell^\infty})[[u]]\), form
\[
 a=\sum_{i\geq 1}\zeta_i u^i.
\]
Let \(v\) be the restriction of the \(u\)-adic valuation of
\(\mathbf Q(\mu_{\ell^\infty})((u))\) to \(\mathbf Q(u,a)\).
\end{construction}

\begin{lemma}[The value group and residue field]
\label{lem:cyclotomic-valuation}
Let $v$ be the valuation of $\mathbf Q(u,a)$ constructed in
\cref{con:cyclotomic-power-series-valuation}.
\begin{enumerate}
\item The valuation $v$ is trivial on $\mathbf Q$, its value group
is $\Z$, and its residue field is $\mathbf Q(\mu_{\ell^\infty})$.
\item The element $a$ is transcendental over $\mathbf Q(u)$. Hence
$\mathbf Q(u,a)$ is purely transcendental of degree two over $\mathbf Q$.
\end{enumerate}
In particular, the valuation ring of $v$ is a discrete valuation
ring containing $\mathbf Q$, with fraction field $\mathbf Q(u,a)$ and
residue field $\mathbf Q(\mu_{\ell^\infty})$.
\end{lemma}

\begin{proof}
\textup{(i)} The valuation \(v\) of
\cref{con:cyclotomic-power-series-valuation} is trivial on
\(\mathbf Q\), and its value group is \(\Z\) because
\(v(u)=1\). Its residue field is the image of
\(\mathbf Q(u,a)\cap\mathbf Q(\mu_{\ell^\infty})[[u]]\) under the
constant-term map, hence a subfield of \(\mathbf Q(\mu_{\ell^\infty})\).

The residue field contains \(\mathbf Q=\mathbf Q(\zeta_0)\).
For \(n\geq1\), suppose inductively that it contains
\(\mathbf Q(\zeta_{n-1})\). Define the truncation
\(s=\sum_{i=1}^{n-1}\zeta_i u^i\in \mathbf Q(\zeta_{n-1})[u]\).
In a second indeterminate $Y$, form the polynomial
\[
 P(u,Y)=\prod_{\sigma\in\operatorname{Gal}(\mathbf Q(\zeta_{n-1})/\mathbf Q)}
            (Y-\sigma(s))\in\mathbf Q[u,Y].
\]
The coefficients are rational because the Galois group permutes the
factors. In the factor corresponding to the identity, the leading
term of \(a-s\) is \(\zeta_nu^n\). For a nonidentity \(\sigma\),
let \(i_\sigma<n\) be the first index with
\(\sigma(\zeta_{i_\sigma})\ne\zeta_{i_\sigma}\). Such an index exists
because the \(\zeta_i\) for \(i<n\) generate \(\mathbf Q(\zeta_{n-1})\). The leading
term of \(a-\sigma(s)\) is
\((\zeta_{i_\sigma}-\sigma(\zeta_{i_\sigma}))u^{i_\sigma}\).
Therefore, with
\[
 d=n+\sum_{\sigma\ne1}i_\sigma,
 \qquad
 b=\prod_{\sigma\ne1}
       (\zeta_{i_\sigma}-\sigma(\zeta_{i_\sigma}))\in \mathbf Q(\zeta_{n-1})^\times,
\]
where the empty sum is zero and the empty product is one,
the element \(u^{-d}P(u,a)\) lies in the valuation ring of \(v\) and its
residue is \(b\zeta_n\). Since \(b\) is a nonzero element of the residue field
by the induction hypothesis, the residue field contains \(\zeta_n\).
The empty truncation gives the same computation for \(n=1\).
Since \(\mathbf Q(\mu_{\ell^\infty})=\bigcup_{n\geq0}\mathbf Q(\zeta_n)\),
this proves the residue field assertion.

\smallskip\noindent
\textup{(ii)} To see that \(\mathbf Q(u,a)/\mathbf Q(u)\) has infinite
degree, take finitely many \(\mathbf Q\)-linearly independent elements
of \(\mathbf Q(\mu_{\ell^\infty})\) and lift them to the valuation ring
of \(v\). A nonzero linear relation between the lifts
with coefficients in \(\mathbf Q(u)\) could be multiplied by a power of
\(u\) so that all coefficients have nonnegative valuation and one has
valuation zero. Reducing the relation would give a nonzero linear
relation over \(\mathbf Q\) between the chosen residues. Since
\([\mathbf Q(\mu_{\ell^\infty}):\mathbf Q]=\infty\), the lifts give arbitrarily large linearly
independent subsets over \(\mathbf Q(u)\). Thus \(a\) is
transcendental over \(\mathbf Q(u)\), and
\(\mathbf Q(u,a)\) is a rational function field in two variables.
\end{proof}

\begin{lemma}[Compatible classes on a rational function field]
\label{lem:rational-function-compatible-classes}
Let $U,Y$ be algebraically independent indeterminates over
$\mathbf Q$, and let
$j\in\Z\setminus\{1\}$. The restriction map
\[
 \varprojlim_n H^1_{\mathrm{\acute et}}(\mathbf Q,\Z/\ell^{n}(j))
 \longrightarrow
 \varprojlim_n H^1_{\mathrm{\acute et}}(\mathbf Q(U,Y),\Z/\ell^{n}(j))
\]
is an isomorphism of \(\Z_\ell\)-modules.
\end{lemma}

\begin{proof}
Let \(\alpha=(\alpha_m)_{m\geq1}\) be a compatible family of classes
\(\alpha_m\in H^1_{\mathrm{\acute et}}(\mathbf Q(U,Y),\Z/\ell^m(j))\).
Fix \(n\geq1\) and put \(E=\mathbf Q(\mu_{\ell^n})\). Restrict the
family to \(E(U,Y)\). For any height-one
prime \(P\) of \(E[U,Y]\), the local ring \(E[U,Y]_P\) is a discrete
valuation ring. Its residue field \(\kappa(P)\) is finitely generated
over \(\mathbf Q\). The algebraic closure of \(\mathbf Q\) in such a
field is finite over \(\mathbf Q\) \cite[\stackstag{037J}]{Stacks}. Hence
\([\kappa(P)(\mu_{\ell^m}):\kappa(P)]\) is unbounded as \(m\) increases.
Indeed, \(\mathbf Q(\mu_{\ell^m})/\mathbf Q\) is Galois of unbounded
degree, and its intersection with \(\kappa(P)\) has bounded degree.
Thus the cyclotomic character of \(\kappa(P)\) has infinite image in
\(\Z_\ell^\times\). Since \(j-1\ne0\), the
\((j-1)\)-st power of this character acts nontrivially on the
rank-one \(\Z_\ell\)-module \(\Z_\ell(j-1)\). Its
invariants are zero because \(\Z_\ell\) is an integral domain.
In particular
\[
 \varprojlim_m H^0_{\mathrm{\acute et}}(\kappa(P),\Z/\ell^{m}(j-1))=0.
\]
By \cref{cor:trait-tate-residues}, the residues of the restricted
family form an element of this inverse limit. Thus the residue of
\(\alpha_n|_{E(U,Y)}\) at every height-one prime of \(E[U,Y]\) is zero.

The field \(E\) contains all \(\ell^n\)-th roots of unity.
Choose a generator of the constant rank-one \(\Z/\ell^n\)-module
sheaf \(\Z/\ell^n(j-1)\) on \(\Spec E\) and use it to identify
\(H^1_{\mathrm{\acute et}}(E(U,Y),\Z/\ell^{n}(j))\) with \(H^1_{\mathrm{\acute et}}(E(U,Y),\mu_{\ell^n})\).
By Kummer theory \cite[\stackstag{03PL}]{Stacks} the latter group is
\(E(U,Y)^\times/(E(U,Y)^\times)^{\ell^n}\). Let \(f\in E(U,Y)^\times\) represent
the class. At a height-one prime \(P\), the residue formula in
\cref{cor:trait-tate-residues} gives \(v_P(f)=0\bmod \ell^n\).
Unique factorization in \(E[U,Y]\) gives \(f=cg^{\ell^n}\) with
\(c\in E^\times\) and \(g\in E(U,Y)^\times\). Hence the class of
\(\alpha_n|_{E(U,Y)}\) is pulled back from \(E\). In particular, it
restricts to zero on \(\overline{\mathbf Q}(U,Y)\).

The field \(\mathbf Q(U,Y)\) is regular over \(\mathbf Q\). For the exact
sequence of absolute Galois groups associated to
\(\mathbf Q\subset\mathbf Q(U,Y)\subset\overline{\mathbf Q}(U,Y)\), the action
on \(\Z/\ell^{n}(j)\) factors through the quotient \(G_{\mathbf Q}\).
The inflation--restriction sequence now shows that \(\alpha_n\)
comes from a unique \(\beta_n\in H^1_{\mathrm{\acute et}}(\mathbf Q,\Z/\ell^{n}(j))\).
Uniqueness and compatibility of \(\alpha\) show that the \(\beta_n\)
form a compatible family. Since inflation is injective at every
coefficient level, restriction on inverse limits is also injective.
All transition and restriction maps are \(\Z_\ell\)-linear.
\end{proof}

\begin{lemma}[Finite-level residue surjectivity]
\label{lem:finite-level-residue-surjectivity}
Let \(R\) be the valuation ring of \cref{lem:cyclotomic-valuation}.
With fraction field \(F\) and residue field \(k\), the residue map
\[
 \partial_{R,n}\colon H^1_{\mathrm{\acute et}}(F,\Z/\ell^{n}(j))
             \longrightarrow H^0_{\mathrm{\acute et}}(k,\Z/\ell^{n}(j-1))
\]
is surjective for every \(n\geq1\) and every \(j\in\Z\).
\end{lemma}

\begin{proof}
Apply \cref{lem:supported-residue-lifting} to
\(X=\Spec R\), with \(Z\) its closed point,
\(M=\Z/\ell^n(j)\), and \(m=\ell^n\). Then
\(X\setminus Z=\Spec F\). Under the canonical
identification
\[
 \underline{\operatorname{Hom}}_{\Z/\ell^n}
       (\mu_{\ell^n},\Z/\ell^n(j))
 \simeq\Z/\ell^n(j-1),
\]
the residue in \cref{lem:supported-residue-lifting}
is exactly \(\partial_{R,n}\), by its definition in
\cref{cor:trait-tate-residues}. Thus every
\(a\in H^0_{\mathrm{\acute et}}(k,\Z/\ell^n(j-1))\) is
\(\partial_{R,n}(\alpha)\) for some
\(\alpha\in H^1_{\mathrm{\acute et}}(F,\Z/\ell^n(j))\), as required.
\end{proof}

\begin{theorem}[A continuous degree-two kernel]
\label{thm:continuous-counterexample}
For every prime \(\ell\), there is an excellent discrete valuation
ring \(R\) containing \(\mathbf Q\), with fraction field \(F\) purely
transcendental of degree two over \(\mathbf Q\), such that, for every
integer \(j\ne 1\), the kernel of
\[
 H^2_{\mathrm{cont}}(R,\Z_\ell(j))
 \longrightarrow H^2_{\mathrm{cont}}(F,\Z_\ell(j))
\]
is a free \(\Z_\ell\)-module of rank one. For $j=1$ the same
continuous restriction map is injective. For every \(j\in\Z\)
and \(n\geq 1\), the corresponding degree-two restriction map with
\(\Z/\ell^{n}(j)\)-coefficients is injective.
\end{theorem}

\begin{proof}
Take the ring \(R\) of \cref{lem:cyclotomic-valuation}, with fraction
field \(F=\mathbf Q(u,a)\) and residue field
\(k=\mathbf Q(\mu_{\ell^\infty})\). It is a discrete valuation ring
and hence a Dedekind domain with fraction field of characteristic zero. Such a ring is excellent
\cite[\stackstag{07QW}\textup{(4)}]{Stacks}. It contains \(\mathbf Q\), so \(\ell\)
is invertible. For any \(n\geq1\) and \(j\in\Z\), finite-coefficient
injectivity in degree two follows from
\cref{cor:trait-tate-residues,lem:finite-level-residue-surjectivity}.

Suppose \(j\ne1\). By \cref{lem:rational-function-compatible-classes},
every compatible family in
\(\varprojlim_n H^1_{\mathrm{\acute et}}(F,\Z/\ell^n(j))\)
comes from \(\mathbf Q\). Since \(R\) is a \(\mathbf Q\)-algebra,
it therefore comes from \(R\). In particular the residue map
\eqref{eq:continuous-residue-limit} is zero. As the only height-one
point of \(\Spec R\) is its closed point,
\cref{cor:continuous-degree-two-obstruction}\textup{(i)} identifies the
continuous degree-two kernel, as a \(\Z_\ell\)-module, with
\[
 \varprojlim_n H^0_{\mathrm{\acute et}}(k,\Z/\ell^n(j-1)).
\]
All \(\ell\)-power roots of unity lie in \(k\), so this group is
free of rank one over \(\Z_\ell\).

For \(j=1\), the continuous degree-two restriction map is injective
by \cref{cor:continuous-kernel-derived-limit}\textup{(ii)}.
\end{proof}

\begin{remark}[The kernel as a cohomology group with support]
\label{rem:continuous-kernel-support}
Let \(R\), \(F\) and \(k\) be as in the proof of
\cref{thm:continuous-counterexample}, let \(j\neq1\), and let
\(z\) be the closed point of \(\operatorname{Spec}R\).
Applying \(R\varprojlim_n\) to the finite
localization triangles gives the continuous localization triangle;
we write \(H^q_{z,\mathrm{cont}}\) for the cohomology of its term
with support in \(z\).
By \cref{sub:trait-locally-constant-residue}\textup{(i)}, the
cohomology with support at every finite level vanishes in degree
one and is identified in degree two with
\(H^0_{\mathrm{\acute et}}(k,\Z/\ell^n(j-1))\), compatibly with
reduction by \cref{cor:trait-tate-residues}\textup{(ii)}. The Milnor sequence
with support then gives
\[
 H^2_{z,\mathrm{cont}}(R,\Z_\ell(j))
   \simeq\varprojlim_n H^0_{\mathrm{\acute et}}(k,\Z/\ell^{n}(j-1)).
\]
The map from \(H^1_{\mathrm{cont}}(F,\Z_\ell(j))\) to
\(H^2_{z,\mathrm{cont}}(R,\Z_\ell(j))\) factors through
\(\varprojlim_n H^1_{\mathrm{\acute et}}(F,\Z/\ell^n(j))\)
by naturality of \eqref{eq:continuous-milnor-sequence} with respect to the
localization boundary. The induced residue map is zero by
\cref{lem:rational-function-compatible-classes} and the factorization
\(\mathbf Q\to R\to F\). Hence the natural
forget-support map takes this rank-one support group isomorphically
onto the kernel. Its composite with every finite coefficient
reduction is zero by \cref{lem:finite-level-residue-surjectivity}.
\end{remark}

\begin{remark}[Rational continuous coefficients]
\label{rem:continuous-rational-counterexample}
Tensoring the degree-two restriction map in
\cref{thm:continuous-counterexample} with
\(\mathbf Q_\ell\) gives a kernel of dimension one for \(j\ne1\).
This is also the kernel for the pro-\'{e}tale sheaf
\(\mathbf Q_\ell(j)\). Bhatt--Scholze's identification
\cite[Lemma~6.8.2]{BhattScholze} in twist zero and a pro-\'etale
local trivialization of the entire Tate system
\cite[Proposition~6.8.4(1)]{BhattScholze} identify it with
\(\Z_\ell(j)[1/\ell]\).
On quasi-compact quasi-separated schemes, derived pro-\'{e}tale
global sections commute with uniformly bounded-below sequential
homotopy colimits \cite[\stackstag{09A3}]{Stacks}. For \(T=\Spec R\)
or \(\Spec F\), apply this to the sequence of multiplications by
\(\ell\) on the pro-\'{e}tale sheaf \(\Z_\ell(j)\) in degree zero.
It gives a natural isomorphism
\[
 H^2_{\mathrm{cont}}(T,\Z_\ell(j))[1/\ell]
 \simeq H^2_{\mathrm{pro\text{-}\acute et}}(T,\mathbf Q_\ell(j)).
\]
For $j=1$ the corresponding rational restriction is injective by
\cref{thm:continuous-counterexample} and flatness of localization.
\end{remark}

\Needspace{11\baselineskip}
\section{Acknowledgment and AI disclosure}\label{sec:ai}

\medskip\noindent\textbf{Acknowledgment.}
I thank Fr\'ed\'eric D\'eglise for introducing me to the Gersten conjecture
and for many subsequent discussions, and Fangzhou Jin, Morten L\"uders,
and Mauro Porta for helpful conversations. This research was supported
in part by the French National Research Agency (ANR) under project
DAG-ARTS (\href{https://anr.fr/Projet-ANR-24-CE40-4098}{ANR-24-CE40-4098}).

\medskip\noindent\textbf{AI disclosure.}
Generative AI tools assisted with mathematical research and proof development, adversarial
readings of the arguments and cited sources, and drafting and revising
this manuscript. I take responsibility for its mathematical claims,
proofs, references, and final text.
I have personally checked every proof in the final version of this manuscript.

\appendix
\crefalias{section}{appendix}
\section{Semilocal comparison and Gersten injectivity}
\label{app:semilocal-moving}

We first prove \cref{thm:semilocal-motivic-vanishing}:
Bloch's cycle complexes of a semilocal essentially smooth
\(\Z_{(\Sigma)}\)-algebra have no cohomology above the weight. For
local rings this is \cite[Corollary~4.4]{Geisser04}, and the proof
follows \cite[Theorem~4.2 and Corollary~4.4]{Geisser04}. Restriction
to the semilocal Dedekind domain obtained by localizing at the generic
points of the closed fibres is injective, because pushforward from a
flat principal divisor vanishes on the homology of cycle complexes
(\cref{prop:flat-principal-divisor-effacement}). Over that Dedekind
domain the bound follows from simultaneous tame symbols
(\cref{sub:semilocal-simultaneous-tame-symbols}). The vanishing of the
pushforward uses the projections of
\cref{lem:relative-curve-projection,lem:simultaneous-curve-projection},
the prism of \cref{con:prism,prop:classical-prism}, and the Quillen graph
of \cref{con:quillen-graph}. The Kan-extension argument of
\cref{lem:semilocal-smooth-source-finality} and the comparison criterion
of \cref{lem:finite-comparison-from-source-bound} then give
\cref{thm:semilocal-finite-connective-comparison}, using the filtered maps of
\cref{lem:classical-etale-stage-comparison,prop:filtered-etale-comparison}
and \cref{sub:finite-filtration-comparison}. We conclude with its
$K$-theoretic applications,
using \cref{prop:finite-type-band-injection,lem:finite-type-witness,lem:semilocal-regular-approximation}
and the degree-two residue calculation of
\cref{prop:semilocal-degree-two-injectivity}.

\subsection{Flat effective principal divisors}

Let \(A\) be a flat \(\Z_{(\Sigma)}\)-algebra. A \emph{flat effective
principal divisor} on \(\Spec A\) is a closed subscheme
\(\Spec(A/fA)\), where \(f\in A\) is a nonzerodivisor and \(A/fA\)
is flat over \(\Z_{(\Sigma)}\).

\begin{sublemma}[Flat effective principal divisors on localizations]
\label{sub:flat-principal-model}
\begin{enumerate}
\item Let \(A\) be a flat \(\Z_{(\Sigma)}\)-algebra and let
\(f\in A\) be a nonzerodivisor. For every \(p\in\Sigma\),
\[
 \ker\bigl(p:A/fA\to A/fA\bigr)\simeq
 \ker\bigl(f:A/pA\to A/pA\bigr).
\]
In particular, \(A/fA\) is flat over \(\Z_{(\Sigma)}\) if and only
if multiplication by \(f\) is injective on \(A/pA\) for every
\(p\in\Sigma\).
\item Let \(C\) be a smooth finitely presented
\(\Z_{(\Sigma)}\)-algebra, let \(S\subseteq C\) be a multiplicative
subset, and let \(\Spec(S^{-1}C/(f))\) be a flat effective principal
divisor on \(\Spec S^{-1}C\). There are \(s\in S\) and \(f_0\in C\)
such that \(f_0S^{-1}C=fS^{-1}C\) and
\(\Spec(C[1/s]/(f_0))\) is a flat effective principal divisor on
\(\Spec C[1/s]\).
\end{enumerate}
\end{sublemma}

\begin{proof}
\emph{Proof of \textup{(i)}.}
The element \(p\) is a nonzerodivisor of the flat
\(\Z_{(\Sigma)}\)-algebra \(A\). The snake lemma
\cite[\stackstag{07JW}]{Stacks} for multiplication by \(p\) on the
exact sequence \(0\to A\xrightarrow{f}A\to A/fA\to0\) gives the
exact sequence
\[
 0\longrightarrow\ker\bigl(p:A/fA\to A/fA\bigr)\longrightarrow
 A/pA\xrightarrow{\ f\ }A/pA .
\]
Every nonzero element of \(\Z_{(\Sigma)}\) is a unit times a product
of elements of \(\Sigma\), and a module over the principal ideal
domain \(\Z_{(\Sigma)}\) is flat if and only if it is torsion free
\cite[\stackstag{0AUW}]{Stacks}. This proves the second assertion.

\smallskip\noindent
\emph{Proof of \textup{(ii)}.}
Write \(f=f_0/s_0\) with \(f_0\in C\) and \(s_0\in S\). The
\(C\)-modules \(\ker(f_0:C\to C)\) and
\(\ker(f_0:C/pC\to C/pC)\), for \(p\in\Sigma\), are finitely
generated. Their localizations at \(S\) are zero: for the first
because \(f\) is a nonzerodivisor of \(S^{-1}C\), and for the
others by \textup{(i)} applied to \(A=S^{-1}C\). Hence some
\(s\in S\) annihilates all of them. Then \(f_0\) is a nonzerodivisor
of \(C[1/s]\) and of \(C[1/s]/pC[1/s]\) for every \(p\in\Sigma\),
and \textup{(i)} applied to \(A=C[1/s]\) shows that
\(C[1/s]/(f_0)\) is flat over \(\Z_{(\Sigma)}\).
\end{proof}
\subsection{Projections for flat principal divisors}

\begin{lemma}[A relative curve projection at one point]
\label{lem:relative-curve-projection}
Let \(X\) be smooth affine over \(\Z_{(\Sigma)}\), of constant relative
dimension \(d\geq1\). Let \(Y\subset X\) be a flat effective
principal divisor. Let \(s\) be a point of \(\operatorname{Spec}\Z_{(\Sigma)}\)
and \(x\) a closed point of the fibre \(X_s\). There are an affine
open neighborhood \(X_0\subseteq X\) of \(x\) and a
\(\Z_{(\Sigma)}\)-morphism
\[
 \pi:X_0\longrightarrow\mathbb A_{\Z_{(\Sigma)}}^{d-1}
\]
such that \(\pi\) is smooth of relative dimension one and
\(\pi|_{Y\cap X_0}\) is quasi-finite.
\end{lemma}

\begin{proof}
The field \(\kappa(s)\) is perfect. Since \(x\) is closed
in the fibre \(X_s\), the extension \(\kappa(x)/\kappa(s)\)
is finite and separable.
The regular local ring \(\mathcal O_{X_s,x}\) has dimension \(d\).
Choose a regular system of parameters \(t_1,\ldots,t_d\) in this
local ring. Their differentials form a basis of
\(\Omega^1_{X_s/\kappa(s)}\otimes \kappa(x)\): the conormal map from
the cotangent space of the fibre local ring onto this space is
surjective because \(\Omega^1_{\kappa(x)/\kappa(s)}=0\), and both spaces have
dimension \(d\). The completion of the fibre local ring is
\(\kappa(x)[[t_1,\ldots,t_d]]\). A coefficient field in the completion
is obtained by lifting the root of a separable polynomial defining
\(\kappa(x)/\kappa(s)\), and the regular parameters give the claimed power
series presentation. Choose an affine open neighborhood
\(X_1\subseteq X\) of \(x\) and elements of
\(\Gamma(X_1,\mathcal O_X)\) lifting \(t_1,\ldots,t_d\); we denote
these lifts again by \(t_1,\ldots,t_d\).

Let \(f\in\Gamma(X,\mathcal O_X)\) be an equation of \(Y\), and let
\(F\) be its image in
\(\widehat{\mathcal O}_{X_s,x}=\kappa(x)[[t_1,\ldots,t_d]]\).
We claim that \(F\neq0\). If \(s=(p)\) is closed, then \(f\) is a
nonzerodivisor on \(\Gamma(X,\mathcal O_X)/p\), by
\cref{sub:flat-principal-model}, hence on its localization
\(\mathcal O_{X_s,x}\). If \(s\) is the generic point, then
\(\mathcal O_{X_s,x}=\mathcal O_{X,x}\), in which \(f\) is a
nonzerodivisor. In both cases the image of \(f\) in
\(\mathcal O_{X_s,x}\) is nonzero, and \(F\neq0\) because the
faithfully flat map
\(\mathcal O_{X_s,x}\to\widehat{\mathcal O}_{X_s,x}\) is injective
\cite[\stackstag{00MC}]{Stacks}.

We adapt the monomial substitution in the proof of Noether
normalization \cite[proof of \stackstag{00OX}]{Stacks} to the power
series \(F\). Among the exponent vectors of nonzero terms of \(F\),
let \(\alpha^*\) be least for the lexicographic order comparing
coordinate \(d\) first and coordinate \(1\) last. This least vector
exists by successively minimizing each nonnegative coordinate.
Put \(w_1=1\). For \(2\leq i\leq d\) choose an integer \(w_i\geq2\)
such that
\[
                   w_i>\sum_{a<i}\alpha_a^*w_a .
\]
For any other occurring exponent vector \(\alpha\), let \(i\)
be the largest index at which \(\alpha\) differs from \(\alpha^*\).
Then \(\alpha_i>\alpha_i^*\). The contribution of coordinate \(i\)
to \(\sum_a (\alpha_a-\alpha_a^*)w_a\) is at least \(w_i\), and
the possible decrease in the lower coordinates is at most
\(\sum_{a<i}\alpha_a^*w_a\). Thus the specialization
\[
                 F(t_1,t_1^{w_2},\ldots,t_1^{w_d})
\]
is a nonzero series. The specialization is defined because all the
weights are positive, and the unique least weighted term has
nonzero coefficient.

Define
\[
 \pi_1=(t_2-t_1^{w_2},\ldots,t_d-t_1^{w_d}):
   X_1\longrightarrow\mathbb A^{d-1}_{\Z_{(\Sigma)}} .
\]
Since \(t_1(x)=0\) and \(w_i\geq2\), the relative differential of
\(t_i-t_1^{w_i}\) at \(x\) is \(dt_i\). The differential criterion
for maps between schemes smooth over \(\Z_{(\Sigma)}\) therefore
shows that \(\pi_1\) is smooth of relative dimension one at \(x\).
The point \(y=\pi_1(x)\) is the origin of the fibre of
\(\mathbb A^{d-1}_{\Z_{(\Sigma)}}\) at \(s\). If \(x\in Y\), the
completed local ring at \(x\) of \(Y\cap\pi_1^{-1}(y)\) is
\[
 \kappa(x)[[t_1]]/
            \bigl(F(t_1,t_1^{w_2},\ldots,t_1^{w_d})\bigr).
\]
This ring has finite length, since the series is nonzero, and
\(\kappa(x)/\kappa(y)=\kappa(x)/\kappa(s)\) is finite. Thus in this
case \(\pi_1|_{Y\cap X_1}\) is quasi-finite at \(x\). If
\(x\notin Y\), we may choose the final neighborhood disjoint from
\(Y\). When \(d=1\) the map \(\pi_1\) is the structure morphism,
and the same argument applies to the nonzero one-variable series
\(F\).
The set of points of \(X_1\) at which \(\pi_1\) is smooth of
relative dimension one is open, and so is the set of points of
\(Y\cap X_1\) at which \(\pi_1|_{Y\cap X_1}\) is quasi-finite
\cite[\stackstag{01TI}]{Stacks}. Let \(X_0\subseteq X_1\) be an
affine open neighborhood of \(x\) contained in the first set and
such that \(Y\cap X_0\) is contained in the second, and put
\(\pi=\pi_1|_{X_0}\).
\end{proof}

\begin{lemma}[A projection on a neighborhood of finitely many points]
\label{lem:simultaneous-curve-projection}
Let \(X\) and \(Y\) be as in
\cref{lem:relative-curve-projection}, and let
\(x_1,\ldots,x_s\in X\). There are a principal affine open
subscheme \(X_0\subseteq X\) containing \(x_1,\ldots,x_s\) and a
\(\Z_{(\Sigma)}\)-morphism
\[
 \pi:X_0\longrightarrow\mathbb A_{\Z_{(\Sigma)}}^{d-1}
\]
such that \(\pi\) is smooth of relative dimension one and
\(\pi|_{Y\cap X_0}\) is quasi-finite.
\end{lemma}

\begin{proof}
Over a discrete valuation ring, the projection at finitely many points
is recalled in \cite[\S2.3.2, Lemma~2.11]{LudersMilnor24}, with
attribution to \cite[Lemma~1]{GilletLevine}. A version over a
local base for finitely many closed points is
\cite[Lemma~2.12, p.~7 of the arXiv version]{LudersMilnor24}.
We give the argument over \(\Z_{(\Sigma)}\). If \(s=0\), take \(X_0=\varnothing\)
and the unique morphism to the affine space. Assume \(s\geq1\).

\smallskip\noindent
\emph{Reduction to closed points.}
Every point of the noetherian scheme \(X\) specializes to a closed
point of \(X\). An open subset of \(X\) containing closed
specializations of \(x_1,\ldots,x_s\) contains
\(x_1,\ldots,x_s\). We may therefore assume that
\(x_1,\ldots,x_s\) are pairwise distinct closed points of \(X\),
with maximal ideals
\(\mathfrak m_{x_1},\ldots,\mathfrak m_{x_s}\subset
\Gamma(X,\mathcal O_X)\). Each
\(x_i\) is a closed point of its fibre over \(\Spec\Z_{(\Sigma)}\).

\smallskip\noindent
\emph{Local data.}
For each \(i\), \cref{lem:relative-curve-projection}
gives an affine open neighborhood of \(x_i\) and a morphism from
it to \(\mathbb A^{d-1}_{\Z_{(\Sigma)}}\), with coordinates
\(f_{i1},\ldots,f_{i,d-1}\), which is smooth of relative dimension
one and quasi-finite on \(Y\). Let \(y_i\) be the image of \(x_i\),
with maximal ideal \(\mathfrak m_{y_i}\) in its local ring. If
\(x_i\in Y\), let \(\mathfrak n_i\) be the maximal ideal of
\(\mathcal O_{Y,x_i}\) and put
\[
 J_i=\mathfrak m_{y_i}\mathcal O_{Y,x_i}.
\]
The ring \(\mathcal O_{Y,x_i}/J_i\) is the local ring at \(x_i\) of
the fibre of \(Y\) over \(y_i\). It has dimension zero by
quasi-finiteness, so \(J_i\) is \(\mathfrak n_i\)-primary. Choose
\(e_i\in\mathbf N\) with \(\mathfrak n_i^{e_i}\subseteq J_i\), and
put \(M_i=\max(2,e_i+1)\). If \(x_i\notin Y\), put \(M_i=2\).

\smallskip\noindent
\emph{Global coordinates.}
The ideals \(\mathfrak m_{x_i}^{M_i}\) of
\(\Gamma(X,\mathcal O_X)\) are pairwise comaximal. The quotient
\(\Gamma(X,\mathcal O_X)/\mathfrak m_{x_i}^{M_i}\) has the unique
maximal ideal \(\mathfrak m_{x_i}/\mathfrak m_{x_i}^{M_i}\), so
localization induces an isomorphism
\[
 \Gamma(X,\mathcal O_X)/\mathfrak m_{x_i}^{M_i}
 \xrightarrow{\sim}
 \mathcal O_{X,x_i}/\mathfrak m_{x_i}^{M_i}\mathcal O_{X,x_i}.
\]
The Chinese remainder theorem \cite[\stackstag{00DT}]{Stacks}
therefore gives \(f_1,\ldots,f_{d-1}\in\Gamma(X,\mathcal O_X)\) with
\[
 f_k\equiv f_{ik}\pmod{\mathfrak m_{x_i}^{M_i}\mathcal O_{X,x_i}}
 \qquad(1\leq i\leq s,\ 1\leq k\leq d-1).
\]
Let \(\pi_1=(f_1,\ldots,f_{d-1}):X\to\mathbb A^{d-1}_{\Z_{(\Sigma)}}\).
Since \(M_i\geq2\), the functions \(f_k\) and \(f_{ik}\) have the
same value in \(\kappa(x_i)\) and the same relative differential at
\(x_i\). Hence \(\pi_1(x_i)=y_i\), and \(\pi_1\) is smooth of
relative dimension one at \(x_i\), by the differential criterion
used in the proof of \cref{lem:relative-curve-projection}.

\smallskip\noindent
\emph{Quasi-finiteness at the points \(x_i\in Y\).}
Let \(u_1,\ldots,u_{d-1}\) be the coordinates of
\(\mathbb A^{d-1}_{\Z_{(\Sigma)}}\), and let
\(\mathfrak p_{y_i}\subset\Z_{(\Sigma)}[u_1,\ldots,u_{d-1}]\) be the
prime ideal of \(y_i\); its elements generate \(\mathfrak m_{y_i}\).
Thus \(J_i\) is generated by the elements
\(g(f_{i1},\ldots,f_{i,d-1})\), for \(g\in\mathfrak p_{y_i}\). Let
\(J'_i\subseteq\mathcal O_{Y,x_i}\) be the ideal generated by the
elements \(g(f_1,\ldots,f_{d-1})\), for \(g\in\mathfrak p_{y_i}\),
so that \(\mathcal O_{Y,x_i}/J'_i\) is the local ring at \(x_i\) of
the fibre of \(\pi_1|_Y\) over \(y_i\). For every
\(g\in\mathfrak p_{y_i}\), the difference
\(g(f_1,\ldots,f_{d-1})-g(f_{i1},\ldots,f_{i,d-1})\) belongs to the
ideal generated by the \(f_k-f_{ik}\), so that
\[
 g(f_1,\ldots,f_{d-1})-g(f_{i1},\ldots,f_{i,d-1})
            \in \mathfrak n_i^{M_i}\subseteq\mathfrak n_iJ_i
\]
in \(\mathcal O_{Y,x_i}\). Consequently \(J'_i\subseteq J_i\) and
\(J_i=J'_i+\mathfrak n_iJ_i\). Nakayama's lemma
\cite[\stackstag{00DV}]{Stacks}, applied to the finite module
\(J_i/J'_i\), gives \(J'_i=J_i\). The local ring
\(\mathcal O_{Y,x_i}/J'_i\) thus has dimension zero. The extension
\(\kappa(x_i)/\kappa(y_i)\) is finite, because \(\kappa(x_i)\) is
finite over the residue field of the image of \(x_i\) in
\(\Spec\Z_{(\Sigma)}\). Hence \(\pi_1|_Y\) is quasi-finite at
\(x_i\).

\smallskip\noindent
\emph{The open subscheme.}
Let \(X_{\mathrm{sm}}\subseteq X\) be the open set of points at
which \(\pi_1\) is smooth of relative dimension one, and let
\(Y_{\mathrm{qf}}\subseteq Y\) be the open set of points at which
\(\pi_1|_Y\) is quasi-finite \cite[\stackstag{01TI}]{Stacks}. The
closed subset
\((X\setminus X_{\mathrm{sm}})\cup(Y\setminus Y_{\mathrm{qf}})\) of
\(X\) contains no \(x_i\). Let
\(J\subseteq\Gamma(X,\mathcal O_X)\) be its radical ideal.
Then \(J\not\subseteq\mathfrak m_{x_i}\) for every \(i\),
and prime avoidance \cite[\stackstag{00DS}]{Stacks} gives
\(a\in J\) outside every \(\mathfrak m_{x_i}\). Put
\(X_0=\Spec\bigl(\Gamma(X,\mathcal O_X)[1/a]\bigr)\) and
\(\pi=\pi_1|_{X_0}\).
\end{proof}

\subsection{A cycle homotopy with prescribed endpoints}
\label{subsec:prism-dimension}
\begin{sublemma}[Finite \'etale points in an open of the affine line]
\label{sub:finite-etale-affine-line-points}
Let $R$ be a semilocal ring whose closed residue fields are
finite, and let $U\subseteq\mathbb A_R^1$ be an open subscheme such that
$U\times_R\kappa(\mathfrak m)$ is nonempty for every maximal
ideal $\mathfrak m$ of $R$.
For every sufficiently large natural number $d$, there are a
finite free \'etale $R$-algebra $R'$ of rank $d$ and an $R$-morphism
$\Spec R'\to U$. In particular, such algebras can be chosen with
two coprime ranks.
\end{sublemma}

\begin{proof}
We adapt the construction by lifting monic polynomials from
the closed residue fields in
\cite[proof of
\arxivhtml{2305.16627v3}{S5.Thmthm2}{Lemma~5.2}]{PaninStavrova25}.
If $R=0$, then $U=\varnothing$ and we may take $R'=R$, which is
free of rank $d$ over $R$ for every $d\geq1$. Assume $R\neq0$.

For each maximal ideal $\mathfrak m\subset R$, the complement of
$U\times_R\kappa(\mathfrak m)$ in
$\mathbb A^1_{\kappa(\mathfrak m)}$ consists of finitely many closed
points. Choose $d\geq1$ larger than the degrees of all these points,
for all the finitely many maximal ideals of $R$.

For a finite field $k$, the roots of $T^{|k|^d}-T$ in an algebraic
closure form a field containing $k$ with exactly $|k|^d$ elements
(the derivative is $-1$). This gives a separable extension of degree
$d$, generated by one element by \cite[\stackstag{030N}]{Stacks}.
Its minimal polynomial is monic, separable and irreducible of degree $d$.
Apply this construction to each $\kappa(\mathfrak m)$ and choose
such a polynomial $g_{\mathfrak m}\in\kappa(\mathfrak m)[T]$.
The corresponding degree-$d$ closed point lies in
$U\times_R\kappa(\mathfrak m)$ by the choice of $d$.

The Chinese remainder theorem
\cite[\stackstag{00DT}]{Stacks}, applied to the coefficients,
gives a monic polynomial $g\in R[T]$ of degree $d$ whose reduction
modulo every $\mathfrak m$ is $g_{\mathfrak m}$.
The algebra $R'=R[T]/(g)$ is free over $R$ with basis given by
the images of $1,T,\ldots,T^{d-1}$.
Every maximal ideal of $R'$ lies over a
maximal ideal $\mathfrak m$ of $R$, and
\[
 R'/\mathfrak mR'\simeq
 \kappa(\mathfrak m)[T]/(g_{\mathfrak m})
\]
is a field. The derivative $g'$ is nonzero in this field because
$g_{\mathfrak m}$ is separable. Thus its image belongs to no
maximal ideal of $R'$ and is a unit. The standard \'etale criterion
\cite[\stackstag{00UB}, \stackstag{00UC}]{Stacks} proves that
$R'$ is \'etale over $R$.

The map $\Spec R'\to\mathbb A_R^1$ sends each closed point to
the point defined by $g_{\mathfrak m}$, hence into $U$.
Its inverse image of $\mathbb A_R^1\setminus U$ is a closed subset
with no closed point. A nonempty closed subset of an affine
scheme has a closed point, so this inverse image is empty.
The map therefore factors through $U$. Choosing two consecutive
integers greater than the degree bound gives coprime ranks.
\end{proof}

\begin{construction}[The prism maps]
\label{con:prism}
Let \(R'/\Z_{(\Sigma)}\) be finite free \'etale, let \(a\in R'\)
with \(a,a-1\in(R')^\times\), and let \(e\in\{0,1\}\).
Let \(m\in\mathbf N\) and \(0\leq r\leq m\). Write
\(w_0,\ldots,w_{m+1}\) for the coordinates of \(\Delta^{m+1}_{R'}\),
and \(t\) and \(v_0,\ldots,v_m\) for those of
\(\mathbb A^1_{R'}\times_{R'}\Delta^m_{R'}\). The prism map
\[
 \theta_{r,e,a}\colon\Delta^{m+1}_{R'}
          \longrightarrow\mathbb A^1_{R'}\times_{R'}\Delta^m_{R'}
\]
is the affine map defined by
\begin{equation}\label{eq:prism-coordinates}
 \theta_{r,e,a}^*(v_i)=
 \begin{cases}
  w_i & (i<r),\\
  w_r+w_{r+1} & (i=r),\\
  w_{i+1} & (i>r),
 \end{cases}
 \qquad
 \theta_{r,e,a}^*(t)=e+(a-e)\sum_{k=r+1}^{m+1}w_k .
\end{equation}
It sends the ordered vertices of \(\Delta^{m+1}_{R'}\) to
\[
 (e,v_0),\ldots,(e,v_r),(a,v_r),(a,v_{r+1}),\ldots,(a,v_m),
\]
where \(v_i\) in a vertex position denotes the vertex of
\(\Delta^m_{R'}\) with \(i\)-th coordinate equal to \(1\). These
are the prism maps used in
\cite[proof of Theorem~4.2]{Geisser04}. Since
\(a-e\in(R')^\times\), the map \(\theta_{r,e,a}\) is an isomorphism:
its inverse is given by \(w_i=v_i\) for \(i<r\), \(w_{i+1}=v_i\)
for \(i>r\), and
\[
 w_{r+1}=\frac{t-e}{a-e}-\sum_{i=r+1}^mv_i,\qquad
 w_r=v_r-w_{r+1}.
\]

Let \(X\) be a smooth affine \(\Z_{(\Sigma)}\)-scheme, let
\(j\in\mathbf N\), and let
\(M\in z^j(X\times_{\Z_{(\Sigma)}}\mathbb A^1_{\Z_{(\Sigma)}},m)\).
We say that \(M\) is \emph{admissible at the endpoints} if, for
every integral summand \(W\) of \(M\), every \(e\in\{0,1\}\) and
every nonempty \(I\subseteq\{0,\ldots,m\}\), every irreducible
component of
\[
 W\cap\bigl(X\times_{\Z_{(\Sigma)}}\{e\}\times_{\Z_{(\Sigma)}}
      \Delta^I_{\Z_{(\Sigma)}}\bigr)
\]
has codimension at least \(j\) in
\(X\times_{\Z_{(\Sigma)}}\{e\}\times_{\Z_{(\Sigma)}}\Delta^I_{\Z_{(\Sigma)}}\).
In that case the Cartier restriction
\[
 i_e^*M=\{t=e\}\cdot M\in z^j(X,m)
\]
of \cref{not:codimension-delta} is defined.
Let \(M_{R'}\) be the flat pullback of \(M\) to \(X_{R'}\). If
\((\operatorname{id}_{X_{R'}}\times\theta_{r,e,a})^*M_{R'}\) is
admissible of codimension \(j\) for every \(0\leq r\leq m\), define
\[
 T_{m,e,a}(M_{R'})=
   \sum_{r=0}^m(-1)^r
        (\operatorname{id}_{X_{R'}}\times\theta_{r,e,a})^*M_{R'}
                   \in z^j(X_{R'},m+1).
\]
\end{construction}

\begin{lemma}[Admissible parameters for both prisms]
\label{lem:prism-admissible-parameters}
Let \(X\) be smooth affine of constant relative dimension \(d\)
over \(\Z_{(\Sigma)}\), let \(j,m\in\mathbf N\), and let
\(M\in z^j(X\times_{\Z_{(\Sigma)}}\mathbb A^1_{\Z_{(\Sigma)}},m)\)
be admissible at the endpoints in the sense of \cref{con:prism}.
There is an open subscheme
\(U\subseteq\mathbb A^1_{\Z_{(\Sigma)}}\setminus\{0,1\}\), with
nonempty fibre over every point of \(\Spec\Z_{(\Sigma)}\), such that
the following holds. Let \(R'/\Z_{(\Sigma)}\) be finite free
\'etale, let \(\Spec R'\to U\) be a \(\Z_{(\Sigma)}\)-morphism, and
let \(a\in R'\) be the image of the coordinate \(t\). Then, for
\(e\in\{0,1\}\),
\[
 (\operatorname{id}_{X_{R'}}\times\theta_{r,e,a})^*M_{R'}
   \in z^j(X_{R'},m+1)\qquad(0\leq r\leq m),
\]
and, if \(m\geq1\),
\[
 (\operatorname{id}_{X_{R'}}\times\theta_{r,e,a})^*(d_i^*M)_{R'}
   \in z^j(X_{R'},m)\qquad(0\leq i\leq m,\ 0\leq r\leq m-1).
\]
\end{lemma}

\begin{proof}
\emph{Faces of the prisms.}
Fix \(e\) and \(r\), and let \(J\subseteq\{0,\ldots,m+1\}\) be
nonempty. By \eqref{eq:prism-coordinates}, the vertex \(k\) of
\(\Delta^{m+1}_{R'}\) lies over the vertex \(\bar k\) of
\(\Delta^m_{R'}\), where \(\bar k=k\) for \(k\leq r\) and
\(\bar k=k-1\) for \(k\geq r+1\). Put \(I=\{\bar k\mid k\in J\}\).
If \(\{r,r+1\}\subseteq J\), then \(\theta_{r,e,a}\) maps the face
\(\Delta^J_{R'}\) isomorphically onto
\(\mathbb A^1_{R'}\times_{R'}\Delta^I_{R'}\), and the components of
the intersection of \(M_{R'}\) with this product have codimension
at least \(j\) because \(M\) is admissible. Otherwise
\(\theta_{r,e,a}\) maps \(\Delta^J_{R'}\) isomorphically onto the
graph
\[
 \{t=e+(a-e)\lambda\}\subseteq
   \mathbb A^1_{R'}\times_{R'}\Delta^I_{R'},
 \qquad
 \lambda=\sum_{k\in J,\ k\geq r+1}v_{\bar k}.
\]
If the sum is empty then \(\lambda=0\), and if it contains all
\(v_i\), \(i\in I\), then \(\lambda=1\). It remains to bound the
intersections with these graphs.

\smallskip\noindent
\emph{The bound to be proved.}
Let \(W\) be an integral summand of \(M\). On
\(W\cap(X\times_{\Z_{(\Sigma)}}\mathbb A^1_{\Z_{(\Sigma)}}\times_{\Z_{(\Sigma)}}\Delta^I_{\Z_{(\Sigma)}})\),
every point has
\[
                         \delta\leq d+|I|+1-j ,
\]
by admissibility of \(M\) and \eqref{eq:codim-delta}.
For admissibility on the graph we need
\(\delta\leq d+|I|-j\), since the graph is isomorphic to
\(\Delta^I\) and
\(X\times_{\Z_{(\Sigma)}}\Delta^I_{\Z_{(\Sigma)}}\) has relative
dimension \(d+|I|-1\).
On the locus \(\lambda=0\), the intersection with the graph is
contained in the intersection with the endpoint \(t=e\),
which satisfies the bound because \(M\) is admissible at the
endpoints.

\smallskip\noindent
\emph{The locus \(\lambda\neq0\).}
Consider the morphism
\[
 \phi=e+\frac{t-e}{\lambda}\colon
           W\cap(X\times_{\Z_{(\Sigma)}}\mathbb A^1_{\Z_{(\Sigma)}}\times_{\Z_{(\Sigma)}}\Delta^I_{\Z_{(\Sigma)}})
               \cap\{\lambda\ne0\}
                          \longrightarrow\mathbb A^1_{\Z_{(\Sigma)}} .
\]
After base change to \(R'\), the intersection with the graph
\(t=e+(a-e)\lambda\) on \(\lambda\ne0\) is the inverse image of the
section \(t=a\) under \(\phi_{R'}\). Let \(s\) be a point of \(\Spec\Z_{(\Sigma)}\).
By the bound on \(\delta\), the fibre over \(s\) of the source of
\(\phi\) is empty or has dimension at most
\(d+|I|-j+\operatorname{ht}(s)\).
On every irreducible component of this fibre which dominates
\(\mathbb A^1_{\kappa(s)}\), the generic fibre of \(\phi\) has
dimension one less than the component, by the
transcendence degree formula over the function field of
\(\mathbb A^1_{\kappa(s)}\). For a component which does not
dominate \(\mathbb A^1_{\kappa(s)}\), we remove the closure of its
image. Applying the generic fibre dimension lemma
\cite[\stackstag{05F7}]{Stacks} to the dominating components (if
there are none, all the remaining fibres are empty), we obtain a
dense open of \(\mathbb A^1_{\kappa(s)}\) over which
the fibres of \(\phi\) are empty or have dimension at most
\(d+|I|-1-j+\operatorname{ht}(s)\). For a closed point \(a_s\) of
this open, its residue field is finite over \(\kappa(s)\).
Thus every point \(z\) of \(\phi^{-1}(a_s)\) satisfies
\[
 \delta(z)=\operatorname{trdeg}_{\kappa(a_s)}\kappa(z)
          +1-\operatorname{ht}(s)\leq d+|I|-j .
\]

\smallskip\noindent
\emph{The open \(U\).}
For each point \(s\) of \(\Spec\Z_{(\Sigma)}\), let
\(U_s\) be the intersection of
\(\mathbb A^1_{\kappa(s)}\setminus\{0,1\}\) with all the dense opens
supplied by the fibre dimension estimate, for the finitely many \(e\),
\(r\), \(J\) and integral summands \(W\) of \(M\).
The complement of \(U_\eta\) in the generic fibre is the zero
set of a nonzero polynomial in \(\mathbf Q[t]\). For every
\(p\in\Sigma\), multiply it by the power of \(p\) for which the
minimum of the \(p\)-adic valuations of its coefficients becomes
zero. The resulting polynomial \(h\) has coefficients in
\(\Z_{(\Sigma)}\) and is nonzero modulo every \(p\in\Sigma\). Put
\[
 U=\mathbb A^1_{\Z_{(\Sigma)}}\setminus
   \Bigl(V(h)\cup\bigcup_{p\in\Sigma}
      \bigl(\mathbb A^1_{\mathbf F_p}\setminus U_{(p)}\bigr)\Bigr).
\]
Each \(\mathbb A^1_{\mathbf F_p}\setminus U_{(p)}\) is a finite set
of closed points, so \(U\) is open. Its fibre over every point
\(s\) is nonempty and contained in \(U_s\).

\smallskip\noindent
\emph{Conclusion.}
Let \(R'\) and \(a\) be as in the statement. For every point
\(s'\) of \(\Spec R'\) over \(s\), the element \(a\) defines a
\(\kappa(s')\)-point of \(U_s\). Since \(R'\) is finite over
\(\Z_{(\Sigma)}\), the extension \(\kappa(s')/\kappa(s)\) is finite.
In particular the image of this point in \(U_s\) is a closed point
\(a_s\). The fibre of \(\phi\otimes_{\kappa(s)}\kappa(s')\) over
the chosen parameter is obtained from \(\phi^{-1}(a_s)\) by the
finite extension \(\kappa(a_s)\subseteq\kappa(s')\). The dimension
bound is preserved by this extension; equivalently, \(\delta\) is
unchanged under the finite base change over \(\Z_{(\Sigma)}\)
(\cref{not:codimension-delta}). Together with the bound on
\(\lambda=0\), this shows that \(\delta\leq d+|I|-j\) holds
on the intersection of \(W_{R'}\) with each graph \(t=e+(a-e)\lambda\), and every
\((\operatorname{id}_{X_{R'}}\times\theta_{r,e,a})^*M_{R'}\) is
admissible of codimension \(j\).
For \(m\geq1\), \(0\leq i\leq m\) and \(0\leq r\leq m-1\), let
\(r'=r+1\) and \(k=i\) if \(i\leq r\), and let \(r'=r\) and
\(k=i+1\) if \(i>r\). By \eqref{eq:prism-coordinates},
\(\theta_{r',e,a}\circ d_k=(\operatorname{id}\times d_i)\circ\theta_{r,e,a}\).
Hence the support of
\((\operatorname{id}_{X_{R'}}\times\theta_{r,e,a})^*(d_i^*M)_{R'}\),
and of its intersection with any face, is contained in the
intersection of the support of
\((\operatorname{id}_{X_{R'}}\times\theta_{r',e,a})^*M_{R'}\) with a
face of \(\Delta^{m+1}_{R'}\). This proves the second assertion.
\end{proof}

\begin{proposition}[A prism on the smooth base]
\label{prop:classical-prism}
Let \(X\) be smooth affine of constant relative dimension over
\(\Z_{(\Sigma)}\), let \(j,m\in\mathbf N\), and let
\(M\in z^j(X\times_{\Z_{(\Sigma)}}\mathbb A^1_{\Z_{(\Sigma)}},m)\)
satisfy \(\partial M=0\). Suppose that \(M\) is admissible at the
endpoints in the sense of \cref{con:prism}. Then
\(\partial(i_e^*M)=0\) for \(e\in\{0,1\}\), and
\[
 [i_0^*M]=[i_1^*M]\quad\text{in }H_m\bigl(z^j(X,\bullet)\bigr).
\]
\end{proposition}

\begin{proof}
\emph{The endpoint restrictions are cycles.}
For \(m=0\) there is no face differential. For \(m\geq1\),
let \(W\) be an integral summand of \(M\), let \(e\in\{0,1\}\) and
\(0\leq i\leq m\), and consider the divisors \(D_1=\{t=e\}\) and
\(D_2=\{v_i=0\}\) of
\(X\times_{\Z_{(\Sigma)}}\mathbb A^1_{\Z_{(\Sigma)}}\times_{\Z_{(\Sigma)}}\Delta^m_{\Z_{(\Sigma)}}\).
The components of \(W\cap D_1\) and of \(W\cap D_2\) have
codimension one in \(W\), and those of \(W\cap D_1\cap D_2\) have
codimension two in \(W\), by admissibility of \(M\), admissibility
at the endpoints, and \eqref{eq:codim-delta}. The hypotheses of
\cref{sub:cartier-restrictions-commute} therefore hold, and
\(d_i^*i_e^*M=i_e^*d_i^*M\). Hence
\(\partial(i_e^*M)=i_e^*(\partial M)=0\).

\smallskip\noindent
\emph{Choice of the parameter.}
Take the open \(U\) from \cref{lem:prism-admissible-parameters}.
By \cref{sub:finite-etale-affine-line-points}, for each of two
sufficiently large consecutive integers \(s\), there is a finite
free \'etale \(\Z_{(\Sigma)}\)-algebra \(R'\) of rank \(s\)
with a morphism \(\Spec R'\to U\). Let \(a\in R'\) be the
image of the coordinate \(t\). Then \(a\) and \(a-1\) are units,
and \cref{con:prism} defines \(T_{m,e,a}(M_{R'})\) for both
\(e=0,1\), as well as \(T_{m-1,e,a}(\partial M_{R'})\) when \(m\geq1\).

\smallskip\noindent
\emph{The chain identity.}
Write \(\theta^m_r=\theta_{r,e,a}\) on \(\Delta^{m+1}_{R'}\).
By \eqref{eq:prism-coordinates}, the faces
\(d_k:\Delta^m_{R'}\to\Delta^{m+1}_{R'}\) satisfy
\begin{align*}
 \theta^m_r\circ d_k&=(\operatorname{id}\times d_k)\circ\theta^{m-1}_{r-1}
      &&(k<r),\\
 \theta^m_r\circ d_k&=(\operatorname{id}\times d_{k-1})\circ\theta^{m-1}_{r}
      &&(k>r+1),\\
 \theta^m_r\circ d_{r+1}&=\theta^m_{r+1}\circ d_{r+1}
      &&(0\leq r<m),\\
 \theta^m_0\circ d_0&=(a,\operatorname{id}_{\Delta^m_{R'}}),
 \qquad
 \theta^m_m\circ d_{m+1}=(e,\operatorname{id}_{\Delta^m_{R'}}).
\end{align*}
In each case the face \(\{w_k=0\}\) of \(\Delta^{m+1}_{R'}\) is the
inverse image under the isomorphism \(\theta^m_r\) of a divisor of
\(\mathbb A^1_{R'}\times_{R'}\Delta^m_{R'}\): the face
\(\{v_k=0\}\), the face \(\{v_{k-1}=0\}\), the graph
\(\{t=e+(a-e)\sum_{i>r}v_i\}\), and the sections \(\{t=a\}\) and
\(\{t=e\}\), respectively. By \eqref{eq:cartier-push-pull}, the
Cartier restrictions of
\((\operatorname{id}_{X_{R'}}\times\theta^m_r)^*M_{R'}\) to the faces
are therefore the pullbacks of the Cartier restrictions of
\(M_{R'}\) to these divisors. Write
\(\theta_r^*\) for \((\operatorname{id}_{X_{R'}}\times\theta^m_r)^*\)
or \((\operatorname{id}_{X_{R'}}\times\theta^{m-1}_r)^*\), and
\(i_a^*M_{R'}=\{t=a\}\cdot M_{R'}\). Then
\begin{align*}
 \partial T_{m,e,a}(M_{R'})
 &=\sum_{r=0}^m\sum_{k=0}^{m+1}(-1)^{r+k}d_k^*\theta_r^*M_{R'}\\
 &=\sum_{r=1}^{m}\sum_{k<r}(-1)^{r+k}\theta_{r-1}^*d_k^*M_{R'}
  +\sum_{r=0}^{m-1}\sum_{k>r+1}(-1)^{r+k}\theta_{r}^*d_{k-1}^*M_{R'}\\
 &\qquad+\sum_{r=0}^m\bigl(d_r^*\theta_r^*M_{R'}
              -d_{r+1}^*\theta_r^*M_{R'}\bigr)\\
 &=-T_{m-1,e,a}(\partial M_{R'})+i_a^*M_{R'}-i_e^*M_{R'} .
\end{align*}
In the last step, the two double sums give
\(-\sum_{r=0}^{m-1}\sum_{i=0}^m(-1)^{r+i}\theta_r^*d_i^*M_{R'}\),
and the third sum telescopes by the third face identity. Thus
\[
  \partial T_{m,e,a}(M_{R'})+T_{m-1,e,a}(\partial M_{R'})
                   =i_a^*M_{R'}-i_e^*M_{R'} ,
\]
where the second term on the left is zero for \(m=0\).
As \(\partial M_{R'}=0\), this gives
\([i_0^*M_{R'}]=[i_a^*M_{R'}]=[i_1^*M_{R'}]\) in
\(H_m(z^j(X_{R'},\bullet))\).

\smallskip\noindent
\emph{Descent from the two coprime ranks.}
Let \(\rho:X_{R'}\to X\) be the projection, which is finite locally
free of rank \(s\). Then \(i_e^*M_{R'}=\rho^*i_e^*M\) by
\eqref{eq:cartier-push-pull}, and \(\rho_*\rho^*\) is multiplication
by \(s\) on \(z^j(X,\bullet)\) \cite[\stackstag{02RH}]{Stacks}.
Applying \(\rho_*\) to \([i_0^*M_{R'}]=[i_1^*M_{R'}]\) gives
\(s([i_0^*M]-[i_1^*M])=0\) in \(H_m(z^j(X,\bullet))\). This holds
for two coprime integers \(s\), so \([i_0^*M]=[i_1^*M]\).
\end{proof}
\subsection{Pushforward from a flat principal divisor}

\begin{sublemma}[A principal Cartier neighborhood in a finite scheme]
\label{sub:finite-cartier-neighborhood}
Let \(X\) be affine noetherian, let \(x_1,\ldots,x_s\in X\), and
let \(b:Z\to X\) be separated, quasi-finite and of finite type.
Let \(h:Y\hookrightarrow Z\) be an effective Cartier immersion
such that \(b\circ h:Y\hookrightarrow X\) is a closed immersion.
There is a principal open subscheme \(X_0\subseteq X\) containing
\(x_1,\ldots,x_s\) such that, putting \(Z_0=Z\times_X X_0\) and
\(Y_0=Y\times_X X_0\), there are a factorization
\[
 Z_0\lhook\joinrel\longrightarrow\overline Z_0
                 \xrightarrow{\,\overline b_0\,}X_0
\]
of the restriction of \(b\), with the first arrow an open immersion
and \(\overline b_0\) finite, and an element
\(g\in\Gamma(\overline Z_0,\mathcal O_{\overline Z_0})\) generating
the ideal of \(Y_0\hookrightarrow Z_0\hookrightarrow\overline Z_0\)
and equal to one on \((\overline Z_0\setminus Z_0)_{\mathrm{red}}\).
\end{sublemma}

\begin{proof}
The construction of \(g\) is a semilocal form of the one in
\cite[proof of Theorem~4.2, pp.~783--784]{Geisser04}.
If \(s=0\), take \(X_0=\varnothing\), with the empty scheme as
finite model and \(g=0\). Assume \(s\geq1\).

Zariski's Main Theorem \cite[\stackstag{05K0}]{Stacks} gives
\(Z\hookrightarrow\overline Z\) with \(\overline Z\) finite
over \(X\). The composite \(\overline h:Y\to\overline Z\)
is closed: it is a section of the separated morphism
\(\overline Z\times_XY\to Y\) followed by the closed immersion
\(\overline Z\times_XY\hookrightarrow\overline Z\).
Its image is contained in \(Z\).

Write \(B=\Gamma(\overline Z,\mathcal O_{\overline Z})\),
let \(I\subseteq B\) be the ideal of \(\overline h(Y)\), and let
\(J\subseteq B\) be the ideal of
\((\overline Z\setminus Z)_{\mathrm{red}}\).
The ideal \(I\) is invertible: on \(Z\) this follows from the
Cartier hypothesis on \(h\), and off \(\overline h(Y)\) it is
the unit ideal. Also \(I+J=B\), since these two closed
subschemes are disjoint.

Discard redundant marked primes contained in another marked prime:
an open containing the remaining points contains their generalizations.
Let \(P\) be the semilocalization of \(\Gamma(X,\mathcal O_X)\)
at the remaining, pairwise incomparable primes, and let
\(B_P=B\otimes_{\Gamma(X,\mathcal O_X)}P\), with ideals
\(I_P\) and \(J_P\). The finite \(P\)-algebra \(B_P\) is
semilocal. Its invertible module \(I_P\) is free of rank one
by \cref{sub:semilocal-projective-free}. Choose a generator
\(g'\in I_P\).

Let \(\mathfrak n_1,\ldots,\mathfrak n_v\) be the maximal
ideals of \(B_P\) not containing \(I_P\). Then
\(I_P+\mathfrak n_i=B_P\) and \(I_P+J_P=B_P\). It follows that
\[
 I_P^2+(J_P\cap\mathfrak n_1\cap\cdots\cap\mathfrak n_v)=B_P.
\]
Indeed each ideal in the intersection maps onto the unit ideal
modulo \(I_P^2\), and their product is contained in their
intersection. The Chinese remainder theorem
\cite[\stackstag{00DT}]{Stacks} supplies \(g_P\in I_P\) with
\[
 g_P\equiv g'\pmod{I_P^2},\qquad
 g_P\equiv1\pmod{J_P},\qquad
 g_P\equiv1\pmod{\mathfrak n_i}\quad(1\leq i\leq v).
\]
At a maximal ideal containing \(I_P\), apply Nakayama's lemma
\cite[\stackstag{00DV}]{Stacks} to \(I_P/(g_P)\) to conclude
that \(g_P\) generates \(I_P\) there. At every other maximal
ideal \(g_P\) is a unit. Thus \(I_P=(g_P)\) on all of
\(\Spec B_P\).

Let \(S\subseteq\Gamma(X,\mathcal O_X)\) be the complement of
the union of the prime ideals corresponding to the \(x_i\), so
that \(P=S^{-1}\Gamma(X,\mathcal O_X)\).
Since \(I_P=S^{-1}I\) and \(J_P=S^{-1}J\), there are \(a_1\in S\)
and \(g\in I[1/a_1]\) with image \(g_P\) in \(I_P\) and with
\(g-1\in J[1/a_1]\). The \(B[1/a_1]\)-module \(I[1/a_1]/(g)\) is
finite and its localization at \(S\) is \(I_P/(g_P)=0\), so it is
annihilated by some \(a_2\in S\). Put
\(X_0=\Spec\bigl(\Gamma(X,\mathcal O_X)[1/(a_1a_2)]\bigr)\).
Then \(I[1/(a_1a_2)]=(g)\) and \(g-1\in J[1/(a_1a_2)]\),
as required over \(X_0\).
\end{proof}

\begin{construction}[The Quillen graph]
\label{con:quillen-graph}
Let \(X'\) be smooth affine of relative dimension \(d\geq1\)
over \(\Z_{(\Sigma)}\), let \(i':Y'\hookrightarrow X'\) be a
flat effective principal divisor, and let
\(\pi:X'\to\mathbb A^{d-1}_{\Z_{(\Sigma)}}\) be smooth of
relative dimension one with \(\pi|_{Y'}\) quasi-finite.
Let \(x_1,\ldots,x_s\in X'\). Form the cartesian square
\[
\begin{tikzcd}
 Z=Y'\times_{\mathbb A_{\Z_{(\Sigma)}}^{d-1}}X'
       \arrow[r,"b"] \arrow[d,"q"'] &
 X' \arrow[d,"\pi"] \\
 Y' \arrow[r,"\pi|_{Y'}"'] & \mathbb A_{\Z_{(\Sigma)}}^{d-1}.
\end{tikzcd}
\]
The maps \(q\) and \(b\) are affine; the first is smooth of
relative dimension one, and the second is quasi-finite. The map
\(h:Y'\to Z\), \(y\mapsto(y,i'(y))\), is a section of the
separated smooth relative curve \(q\), hence an effective Cartier
immersion. Its composite with \(b\) is \(i'\).
Apply \cref{sub:finite-cartier-neighborhood} to \(b\) and \(h\).
It gives a principal open \(X_0\subseteq X'\) containing the marked
points, and, writing \(Z_0=b^{-1}(X_0)\) and \(Y_0=Y'\cap X_0\),
a finite model \(\overline b_0:\overline Z_0=\Spec B\to X_0\) with a function
\(g\) generating the ideal of \(h_0(Y_0)\), where \(h_0=h|_{Y_0}\),
and equal to one on the reduced complement of \(Z_0\).
We retain the original target \(Y'\) of the smooth map: in general
\(q(Z_0)\) is not contained in \(Y_0\).
Put \(W=\Spec B[(g-1)^{-1}]\subseteq Z_0\), and let
\(q_W:W\to Y'\) be the restriction of \(q\), smooth of relative
dimension one. The morphism
\[
 \beta:W\longrightarrow X_0\times_{\Z_{(\Sigma)}}\mathbb A^1_{\Z_{(\Sigma)}},
 \qquad z\longmapsto(\overline b_0(z),g(z)/(g(z)-1))
\]
is finite: its graph in \(\overline Z_0\times\mathbb A^1\) is defined by
\((g-1)t=g\), since this relation makes \((g-1)(t-1)=1\), and
\(\overline Z_0\) is finite over \(X_0\). This is the graph construction
of \cite[proof of Theorem~4.2]{Geisser04}, with the two targets kept distinct.
Writing \(i'_0:Y_0\hookrightarrow X_0\), its scheme-theoretic endpoints are
\begin{equation}\label{eq:graph-endpoints}
 \beta^{-1}(t=0)=V(g)=h_0(Y_0),\qquad
 \beta^{-1}(t=1)=\varnothing,\qquad \beta h_0=(i'_0,0).
\end{equation}
The composite \(q_W\circ h_0\) is the open immersion \(Y_0\to Y'\).
\end{construction}

\begin{proposition}[Vanishing of a flat principal divisor pushforward]
\label{prop:flat-principal-divisor-effacement}
Let \(P\) be a semilocal essentially smooth \(\Z_{(\Sigma)}\)-algebra in the
sense of \cref{def:essentially-smooth-sigma}, let \(X=\operatorname{Spec}P\), and let
\(i:Y\hookrightarrow X\) be a flat effective principal divisor.
For each \(j\geq1\) and \(m\geq0\), the map
\[
       i_*:H_m(z^{j-1}(Y,*))\longrightarrow H_m(z^j(X,*))
\]
is zero.
\end{proposition}

\begin{proof}
The proof adapts that of \cite[Theorem~4.2]{Geisser04} to semilocal
rings. If \(Y=\varnothing\), the source of the pushforward is zero.
Assume \(Y\ne\varnothing\), so that \(P\ne0\).

\smallskip\noindent
\emph{A finite-type model.}
Write \(P=S^{-1}C\), where \(C\) is a smooth finitely presented
\(\Z_{(\Sigma)}\)-algebra, \(\mathfrak p_1,\ldots,\mathfrak p_s\) are
prime ideals of \(C\), and
\(S=C\setminus\bigcup_{i=1}^s\mathfrak p_i\). Let
\(x_1,\ldots,x_s\in\Spec C\) be the corresponding points, and let
\(f\in P\) be an equation of \(Y\). In this proof, a
\emph{shrinking} is the replacement of \(C\) by \(C[1/a]\) for some
\(a\in S\); it does not change \(P\). By
\cref{sub:flat-principal-model}, after a shrinking there is
\(f_0\in C\) with \(f_0P=fP\) such that
\[
 i':Y'=\Spec(C/f_0C)\lhook\joinrel\longrightarrow X'=\Spec C
\]
is a flat effective principal divisor. Then
\(Y=\Spec S^{-1}(C/f_0C)\). Let \(E\in z^{j-1}(Y,m)\) satisfy
\(\partial E=0\). By \cref{sub:classical-cycle-spread}, applied to
\(C/f_0C\), after a shrinking there is \(E'\in z^{j-1}(Y',m)\) with
\(\partial E'=0\) whose restriction to \(Y\) is \(E\). The
restriction of \(i'_*E'\) to \(\Spec P\) is \(i_*E\). It therefore
suffices to prove that, after a shrinking,
\[
 [i'_*E']=0\quad\text{in }H_m\bigl(z^j(X',\bullet)\bigr).
\]
The scheme \(X'\) is the disjoint union of its finitely many
connected components, and \(z^j\) is additive
(\cref{def:bloch-complex}). Components with no marked point may be
discarded by a shrinking, using the idempotent of their complement.
Treating the other components separately, we may assume that \(X'\)
is connected, of constant relative dimension \(d\) over
\(\Z_{(\Sigma)}\). If \(d=0\), every fibre of \(X'\) over
\(\Spec\Z_{(\Sigma)}\) is the spectrum of a finite product of
fields, on which the nonzerodivisor \(f_0\) is a unit
(\cref{sub:flat-principal-model}); thus \(Y'=\varnothing\). Assume
\(d\geq1\).

\smallskip\noindent
\emph{Choosing the relative curve and its graph.}
By \cref{lem:simultaneous-curve-projection}, after a
shrinking there is a smooth map
\(\pi:X'\to\mathbb A^{d-1}_{\Z_{(\Sigma)}}\) of relative
dimension one whose restriction to \(Y'\) is quasi-finite.
Apply \cref{con:quillen-graph} to \(\pi\) and the points
\(x_1,\ldots,x_s\), retaining \(Y'\) and \(E'\) after choosing \(\pi\).
Use the notation of \cref{con:quillen-graph}; in particular, the
endpoints of \(\beta\) are given by \eqref{eq:graph-endpoints}.

\smallskip\noindent
\emph{The cycle.}
Set \(M=\beta_*q_W^*E'\). We verify that
\(M\in z^j(X_0\times\mathbb A^1,m)\) using
\cref{not:codimension-delta,def:bloch-complex}. At every point of
\(Y'\times\Delta^m\) the codimension equals \(d+m-\delta\), because
\(Y'\) is an effective Cartier divisor in \(X'\). Thus each summand
of \(E'\) has generic \(\delta=d+m-j+1\). Smooth pullback of relative
dimension one raises this by one at the generic points of the inverse
image, and finite image preserves it. The dimension constant on
\(X_0\times\mathbb A^1\times\Delta^m\) is \(d+m+2\), so these images
have codimension \(j\).
On a face \(\Delta^I\) the support of \(E'\) has \(\delta\leq d+|I|-j\)
at all points. Smooth pullback and finite image give
\(\delta\leq d+|I|-j+1\); finiteness ensures that intersection of an
image with the face is the image of its intersection with that face.
The ambient dimension constant there is \(d+|I|+1\). This proves
admissibility, and the Cartier identities \eqref{eq:cartier-push-pull}
give \(\partial M=\beta_*q_W^*\partial E'=0\) (also for \(m=0\)).

\smallskip\noindent
\emph{The endpoints.}
For an integral summand \(V\) of \(E'\), put
\(V_0=V\cap(Y_0\times\Delta^m)\) and
\(W_V=(q_W\times\mathrm{id}_{\Delta^m})^{-1}(V)\).
The scheme \(W_V\) is smooth over \(V\); the flat pullback
\(q_W^*[V]\) is the fundamental cycle of this entire inverse image,
not necessarily an integral scheme. Its intersection with \(V(g)\)
is the section over \(V_0\), by \eqref{eq:graph-endpoints}. Each
component of \(W_V\) dominates \(V\), and none is contained in
\(V(g)\) (over the generic point the section is a point on a smooth
curve). If \(V_0\) is nonempty, the generic point of the section over
\(V_0\) lies on a unique component, and its local ring is the discrete
valuation ring of the smooth generic fibre at that rational point.
The ideal of the section is generated by \(g\), so its Cartier
multiplicity is one. If \(V_0\) is empty the intersection is empty.
Consequently \(V(g)\cdot q_W^*[V]=(h_0\times\mathrm{id})_*[V_0]\),
with the empty cycle interpreted as zero. The Cartier pushforward identity \eqref{eq:cartier-push-pull} gives
\[
 i_0^*M=(i'_0)_*(E'|_{Y_0}),\qquad i_1^*M=0.
\]
Moreover, for every face \(\Delta^I\), the intersection of the support
of \(M\) with \(t=0\) on that face is contained in the image of the
support of \(E'|_{Y_0}\) on that face. Indeed the support of a finite
pushforward is contained in the finite image of the original support,
and its inverse image at the endpoint is the section just computed.
The codimension is at least \(j-1\) in \(Y_0\times\Delta^I\), hence
at least \(j\) in \(X_0\times\Delta^I\) by the dimension formula for
the Cartier divisor. The other endpoint is empty. Thus \(M\) is
admissible at the endpoints in the sense of \cref{con:prism}.
The prism of \cref{prop:classical-prism} gives
\([(i'_0)_*(E'|_{Y_0})]=0\) in \(H_m(z^j(X_0,\bullet))\).
Restriction to \(\Spec P\) proves the proposition.
\end{proof}
\Needspace{6\baselineskip}
\subsection{The upper cohomological bound}

\begin{sublemma}[Simultaneous tame symbols]
\label{sub:semilocal-simultaneous-tame-symbols}
Let $V$ be a semilocal Dedekind domain of dimension one, with
maximal ideals $\mathfrak m_1,\ldots,\mathfrak m_t$.
For every natural number $j\geq1$, the homomorphism
\[
 K_j^M(\Frac V)\longrightarrow
 \bigoplus_{i=1}^t K_{j-1}^M(V/\mathfrak m_i),\qquad
 \alpha\longmapsto
       (\partial^{\mathrm M}_{V_{\mathfrak m_i}}\alpha)_i,
\]
is surjective. Here $\partial^{\mathrm M}_{V_{\mathfrak m_i}}$
is the Milnor tame-symbol map for the discrete valuation ring
$V_{\mathfrak m_i}$.
\end{sublemma}

\begin{proof}
Fix $i\in\{1,\ldots,t\}$ and a symbol
$\{\bar u_1,\ldots,\bar u_{j-1}\}$ in the $i$-th summand.
The ideals $\mathfrak m_1,\ldots,\mathfrak m_t$ are pairwise
comaximal. By the Chinese remainder theorem
\cite[\stackstag{00DT}]{Stacks}, choose $u_a\in V$ with
\[
 u_a\bmod\mathfrak m_i=\bar u_a,
 \qquad u_a\bmod\mathfrak m_k=1\quad(k\ne i).
\]
All these residues are nonzero, so $u_a\in V^\times$.
The localization map
\[
 \mathfrak m_i/\mathfrak m_i^2\longrightarrow
 \mathfrak m_iV_{\mathfrak m_i}/
                (\mathfrak m_iV_{\mathfrak m_i})^2
\]
is an isomorphism of $V/\mathfrak m_i$-vector spaces: the
source is annihilated by $\mathfrak m_i$, so every element
of $V\setminus\mathfrak m_i$ already acts invertibly on it.
The target has dimension one because $V_{\mathfrak m_i}$ is
a discrete valuation ring. Choose a nonzero element of
$\mathfrak m_i/\mathfrak m_i^2$. The ideals
$\mathfrak m_i^2$ and $\mathfrak m_k$, $k\ne i$, are again
pairwise comaximal: modulo $\mathfrak m_k$ the ideal
$\mathfrak m_i$ is the unit ideal, hence so is its square.
The same Chinese remainder theorem \cite[\stackstag{00DT}]{Stacks} gives $\varpi_i\in V$ with
that prescribed class modulo $\mathfrak m_i^2$ and with
$\varpi_i\equiv1\pmod{\mathfrak m_k}$ for $k\ne i$.
Consequently $\varpi_i$ has valuation one at $\mathfrak m_i$
and valuation zero at every other maximal ideal.

By \cref{not:milnor-tame-residues},
\[
 \partial^{\mathrm M}_{V_{\mathfrak m_i}}
       \{u_1,\ldots,u_{j-1},\varpi_i\}
       =\{\bar u_1,\ldots,\bar u_{j-1}\}.
\]
At every $\mathfrak m_k$ with $k\ne i$ all entries of this symbol
are units, so its tame symbol there is zero.
For $j=1$, the symbol $\{\varpi_i\}$ maps to $1$ in the $i$-th
summand $K_0^M(V/\mathfrak m_i)=\Z$ and to zero elsewhere.
Every target summand is generated by symbols, which proves
surjectivity.
\end{proof}

\begin{theorem}[The upper cohomological bound on semilocal sources]
\label{thm:semilocal-motivic-vanishing}
Let \(P\) be a semilocal essentially smooth \(\Z_{(\Sigma)}\)-algebra.
Then, for every \(j\geq0\),
\[
      H^q\bigl(\Z(j)^{\mathrm{cla}}(\operatorname{Spec}P)\bigr)=0
                                  \qquad(q>j).
\]
\end{theorem}

\begin{proof}
The assertion for the zero ring is immediate, since its cycle
complex is zero. Assume \(P\ne0\). The regular noetherian ring
\(P\) is a finite product of domains
\cite[\stackstag{0567}, \stackstag{030C}]{Stacks}, each of which is
a semilocal essentially smooth \(\Z_{(\Sigma)}\)-algebra, and
\(\Z(j)^{\mathrm{cla}}\) is additive (\cref{def:bloch-complex}). We may
therefore assume that \(P\) is a domain, with fraction field \(F\).
For \(j=0\), the complex \(\Z(0)^{\mathrm{cla}}(\Spec P)\) is zero in
cohomological degrees \(q>0\). Suppose that \(j\geq1\).

\smallskip\noindent
\emph{The Dedekind localization.}
Let \(\mathfrak v_1,\ldots,\mathfrak v_t\) be the minimal prime
ideals of the ideals \(pP\), for those \(p\in\Sigma\) with
\(pP\neq P\). Each \(\mathfrak v_i\) has height one
\cite[\stackstag{00KV}]{Stacks}, being minimal over a nonzero
principal ideal of the domain \(P\). If \(t\geq1\), let \(V\) be the
semilocalization of \(P\) at \(\mathfrak v_1,\ldots,\mathfrak v_t\)
(\cref{def:arithmetic-semilocal}). It is a semilocal Dedekind
domain with fraction field \(F\) and maximal ideals
\(\mathfrak v_iV\). If \(t=0\), put \(V=F\).

\smallskip\noindent
\emph{Injectivity of restriction to \(V\).}
We claim that the restriction
\[
 H^q\bigl(\Z(j)^{\mathrm{cla}}(\operatorname{Spec}P)\bigr)
   \longrightarrow
 H^q\bigl(\Z(j)^{\mathrm{cla}}(\operatorname{Spec}V)\bigr)
\]
is injective for every \(q\in\Z\). Let \(\alpha\) be a class in its
kernel. The ring \(V\) is the localization of \(P\) at the
multiplicative subset \(P\setminus\bigcup_{i=1}^t\mathfrak v_i\),
respectively \(P\setminus\{0\}\) if \(t=0\). By
\cref{sub:classical-cycle-spread}, applied to a presentation
\(P=S^{-1}C\) and to the inverse image in \(C\) of this subset,
there is an element \(f\) of this subset such that \(\alpha\)
vanishes in
\(H^q(\Z(j)^{\mathrm{cla}}(\operatorname{Spec}P[1/f]))\). For every
\(p\in\Sigma\), the ring \(P/pP\) is reduced, as a localization of
a smooth \(\mathbf F_p\)-algebra, and its minimal prime ideals are
among the \(\mathfrak v_i/pP\). Hence \(f\) is a nonzerodivisor on
\(P/pP\), and \(\Spec(P/fP)\) is a flat effective principal divisor
on \(\Spec P\) by \cref{sub:flat-principal-model}. Put \(m=2j-q\).
By the exact sequence \eqref{eq:cycle-localization-sequence} for
\(A=P\) and \(h=f\), the class \(\alpha\) lies in the image of
\[
 i_*:H_m\bigl(z^{j-1}(\Spec(P/fP),\bullet)\bigr)\longrightarrow
     H_m\bigl(z^j(\Spec P,\bullet)\bigr)
     =H^q\bigl(\Z(j)^{\mathrm{cla}}(\operatorname{Spec}P)\bigr).
\]
This map is zero: its target is zero if \(m<0\), and if \(m\geq0\)
this is \cref{prop:flat-principal-divisor-effacement}. Hence
\(\alpha=0\).

\smallskip\noindent
\emph{Vanishing for \(V\).}
For a field \(E\) and \(a,m\in\mathbf N\) with \(m<a\), we have
\(z^a(\Spec E,m)=0\), because \(\Delta^m_E\) has dimension \(m\).
Hence
\begin{equation}\label{eq:field-cycle-vanishing}
 H^q\bigl(\Z(a)^{\mathrm{cla}}(\Spec E)\bigr)=0\qquad(q>a).
\end{equation}
This proves the theorem if \(V=F\). Assume that \(t\geq1\). The
semilocal Dedekind domain \(V\) is a principal ideal domain, since
its nonzero ideals are invertible and hence free of rank one by
\cref{sub:semilocal-projective-free}. Let \(h\in V\) generate
\(\bigcap_{i=1}^t\mathfrak v_iV\). Then
\(V/hV\simeq\prod_{i=1}^t\kappa(\mathfrak v_i)\) by the Chinese
remainder theorem \cite[\stackstag{00DT}]{Stacks}, and
\(V[1/h]=F\). The exact sequence
\eqref{eq:cycle-localization-sequence} for \(A=V\) reads
\[
 \begin{split}
 H^{q-1}\bigl(\Z(j)^{\mathrm{cla}}(\operatorname{Spec}F)\bigr)
 &\longrightarrow \bigoplus_{i=1}^t
        H^{q-2}\bigl(\Z(j-1)^{\mathrm{cla}}(\operatorname{Spec}\kappa(\mathfrak v_i))\bigr)\\
 &\longrightarrow
       H^q\bigl(\Z(j)^{\mathrm{cla}}(\operatorname{Spec}V)\bigr)
 \longrightarrow
       H^q\bigl(\Z(j)^{\mathrm{cla}}(\operatorname{Spec}F)\bigr).
 \end{split}
\]
By \eqref{eq:field-cycle-vanishing}, the second and fourth terms
vanish for \(q>j+1\). For \(q=j+1\), the fourth term vanishes and
the first arrow is the simultaneous Milnor tame-symbol map
\begin{equation}\label{eq:simultaneous-tame-map}
                      K_j^M(F)\longrightarrow
                            \bigoplus_{i=1}^t K_{j-1}^M(\kappa(\mathfrak v_i)).
\end{equation}
By naturality of localization, the component of the boundary at
\(\mathfrak v_i\) is the boundary for the discrete valuation ring
\(V_{\mathfrak v_iV}\). Under the identification of diagonal
motivic cohomology of fields with Milnor \(K\)-theory, it is the tame
symbol by \cite[proof of Corollary~4.4, p.~786]{Geisser04}. By
\cref{sub:semilocal-simultaneous-tame-symbols}, the map
\eqref{eq:simultaneous-tame-map} is surjective. Hence
\(H^q(\Z(j)^{\mathrm{cla}}(\operatorname{Spec}V))=0\) for \(q>j\), and
the injectivity of restriction from \(P\) to \(V\) proved in the first part gives the same vanishing for \(P\).
\end{proof}

\subsection{The semilocal comparison and its applications}

\begin{theorem}[Finite connective \'etale comparison]
\label{thm:semilocal-finite-connective-comparison}
Let \(A\) be a nonzero semilocal ring and let \(N\geq 2\)
be an integer invertible in \(A\). For every \(n\in\mathbf N\), the
finite comparison map \eqref{eq:finite-comparison-map} induces an
isomorphism
\[
 K_n^{\mathrm{cn}}(A;\Z/N)\xrightarrow{\ \sim\ }
 \pi_n R\Gamma_{\mathrm{\acute et}}\bigl(A,
        \tau_{\leq 2n}\mathcal K^{\mathrm{\acute et}}_{A,N}\bigr),
\]
natural in \(A\). If \(A\) is regular and noetherian, then
\(K_n^{\mathrm{cn}}(A;\Z/N)=K_n(A;\Z/N)\).
\end{theorem}

\begin{proof}
For the choice of \(\Sigma\) in \cref{lem:finite-comparison-from-source-bound},
the required vanishing on semilocal essentially smooth sources is
\cref{thm:semilocal-motivic-vanishing}. Apply \cref{lem:finite-comparison-from-source-bound}.
\end{proof}

\begin{theorem}[Injectivity over semilocal domains essentially of finite type]
\label{thm:finite-type-semilocal-gersten-injectivity}
Let $\ell$ be a prime number, and let $R$ be a regular semilocal
domain essentially of finite type over a field or over $\Z$
\textup{(\cref{def:arithmetic-semilocal})}, with fraction field $F$,
such that $\ell\in R^\times$. Suppose that the closed residue fields
of $R$ have at most one positive characteristic. For every $\nu\geq1$
and every $n\in\Z$, the restriction map
\[
         K_n(R;\Z/\ell^\nu)\longrightarrow K_n(F;\Z/\ell^\nu)
\]
is injective.
\end{theorem}

\begin{proof}
For \(n<0\) the coefficient groups vanish by regularity
\cite[Proposition~6.8(b)]{TT} and \eqref{eq:coefficient-sequence}.
For \(n\geq0\), apply the natural comparison of
\cref{thm:semilocal-finite-connective-comparison}. As
\(\mathcal K^{\mathrm{\acute et}}_{R,\ell^\nu}\) and
\(\mathcal K^{\mathrm{\acute et}}_{F,\ell^\nu}\) are connective
(\cref{sub:finite-k-etale-coefficients}\textup{(i)}), the restriction
map identifies with \(\Phi^{[0,2n]}_n(R)\to\Phi^{[0,2n]}_n(F)\)
for modulus \(\ell^\nu\), which is injective by
\cref{prop:finite-type-band-injection}.
\end{proof}

\begin{corollary}[Finite \'etale algebras]
\label{cor:finite-etale-algebra-k-injectivity}
Let $R$ be a noetherian regular local ring with fraction field $F$,
and let $A$ be a finite \'etale $R$-algebra. For every natural
number $N\geq2$ invertible in $R$ and every $n\in\Z$, the
restriction map
\[
         K_n(A;\Z/N)\longrightarrow K_n(A\otimes_R F;\Z/N)
\]
is injective.
\end{corollary}

\begin{proof}
\emph{Descent.}
As in the proof of
\cref{thm:gersten-injectivity}, we
may assume that $N=\ell^\nu$ for a prime number $\ell\in R^\times$.
Let $x\in K_n(A;\Z/N)$ have zero image in $K_n(A\otimes_RF;\Z/N)$,
and let $k$ be as in \cref{lem:finite-type-witness}. The ring $R$ is
the union of the directed set of its finitely generated
$k$-subalgebras. By \cite[\stackstag{0BTV}]{Stacks}, there are a
finitely generated $k$-subalgebra $R_1\subseteq R$ and a finite
\'etale $R_1$-algebra $A_1$ with an isomorphism of $R$-algebras
$A_1\otimes_{R_1}R\simeq A$. By \cref{sub:k-finite-continuity}, the
functor
\[
 \Phi:B\longmapsto K_n(A_1\otimes_{R_1}B;\Z/N)
\]
from $R_1$-algebras to abelian groups commutes with filtered
colimits, and $x\in\Phi(R)$ has zero image in $\Phi(F)$. Let $R_0$,
$F_0$, $R_1\to R_0\to R$ and $x_0\in\Phi(R_0)$ be as in
\cref{lem:finite-type-witness}, and put $A_0=A_1\otimes_{R_1}R_0$.
Then $x_0\in K_n(A_0;\Z/N)$ has image $x$ in $K_n(A;\Z/N)$ and zero
image in $K_n(A_0\otimes_{R_0}F_0;\Z/N)$, and it suffices to prove
that $x_0=0$. This holds if $A_0=0$.

\smallskip\noindent
\emph{The components of $A_0$.}
Suppose that $A_0\neq0$. Since $R_0\to R$ is local,
$\ell\in R_0^\times$. By
\cref{sub:arithmetic-finite-etale}\textup{(ii)}, there is an
isomorphism of $R_0$-algebras $A_0\simeq\prod_{j=1}^rB_j$ in which
every $B_j$ is a regular semilocal domain essentially of finite type
over $k$, and the closed residue fields of $B_j$ have the
characteristic of the residue field of $R_0$. The ring $B_j$ is
nonzero and flat over the domain $R_0$, so $R_0\to B_j$ is injective.
Hence $B_j\otimes_{R_0}F_0$ is the localization of the domain $B_j$
at the image of $R_0\setminus\{0\}$, which does not contain $0$; it
is a domain, finite over the field $F_0$, and therefore the field
$\Frac B_j$. Consequently
\[
 K_n(A_0;\Z/N)=\bigoplus_{j=1}^rK_n(B_j;\Z/N),
 \qquad
 K_n(A_0\otimes_{R_0}F_0;\Z/N)
    =\bigoplus_{j=1}^rK_n(\Frac B_j;\Z/N),
\]
and the $j$-th component of $x_0$ has zero image in
$K_n(\Frac B_j;\Z/N)$. By
\cref{thm:finite-type-semilocal-gersten-injectivity} for the rings
$B_j$, every component of $x_0$ is zero.
\end{proof}

\begin{corollary}[The second Moore group on semilocal rings]
\label{cor:semilocal-K2-injectivity}
Let $R$ be a connected noetherian regular semilocal ring
with fraction field $F$. For every natural number $N\geq2$ invertible in $R$,
the map $K_2(R;\Z/N)\to K_2(F;\Z/N)$ is injective.
\end{corollary}

\begin{proof}
By \cref{sub:moore-primary-decomposition} we may assume that
\(N\) is a prime power. Apply the natural comparison of
\cref{thm:semilocal-finite-connective-comparison} at $n=2$.
For $A=R,F$, \cref{sub:finite-k-etale-coefficients} gives the
fibre sequences of spectral sheaves
\[
\begin{gathered}
 (\Z/N(2))[4]\longrightarrow
 \tau_{\leq4}\mathcal K^{\mathrm{\acute et}}_{A,N}\longrightarrow
 \tau_{\leq2}\mathcal K^{\mathrm{\acute et}}_{A,N},\\
 (\Z/N(1))[2]\longrightarrow
 \tau_{\leq2}\mathcal K^{\mathrm{\acute et}}_{A,N}\longrightarrow \Z/N.
\end{gathered}
\]
Here the coefficient sheaves denote their Eilenberg--Mac Lane sheaves
as in \cref{con:etale-k-sheaves}. Taking derived global sections of the
second sequence shows that the middle term has $\pi_3=0$ and
$\pi_2=H^0_{\mathrm{\acute et}}(A,\Z/N(1))$, since for an abelian sheaf
$M$, $\pi_iR\Gamma_{\mathrm{\acute et}}(A,M)=H^{-i}_{\mathrm{\acute et}}(A,M)$
vanishes for $i>0$. Taking
homotopy of derived global sections of the first sequence and using
the natural comparison gives a commutative diagram with exact rows:
\[
\begin{tikzcd}[column sep=small]
 0\arrow[r] & H^2_{\mathrm{\acute et}}(R,\Z/N(2))\arrow[r]\arrow[d,hook]
 & K_2(R;\Z/N)\arrow[r]\arrow[d]
 & H^0_{\mathrm{\acute et}}(R,\Z/N(1))\arrow[d,hook] \\
 0\arrow[r] & H^2_{\mathrm{\acute et}}(F,\Z/N(2))\arrow[r]
 & K_2(F;\Z/N)\arrow[r]
 & H^0_{\mathrm{\acute et}}(F,\Z/N(1)).
\end{tikzcd}
\]
Restriction on the left is injective by
\cref{prop:semilocal-degree-two-injectivity}.
Restriction on the right is injective because $\Spec R$ is connected:
a section of a finite locally constant sheaf is determined by its
value at a geometric generic point. Exactness gives the conclusion.
\end{proof}

\subsection{Semilocal rings with a regular model}
\label{subsec:semilocal-regular-model}

\begin{corollary}[Semilocal rings with a regular model: \(K\)-theory]
\label{cor:semilocal-regular-model-k-injection}
Let \(R\) be a noetherian semilocal domain with fraction field \(F\),
whose closed residue fields have at most one positive characteristic,
and let \(N\geq2\) be a natural number invertible in \(R\).
\begin{enumerate}
\item Suppose that there is a regular homomorphism \(A\to R\) from a
regular ring \(A\) essentially of finite type over a field or over
\(\Z\). Then, for every \(n\in\Z\), the restriction map
\(K_n(R;\Z/N)\to K_n(F;\Z/N)\) is injective.
\item Let \(p\) be a prime number. Suppose that \(R\) is regular of
characteristic zero, that every closed residue field of \(R\) has
characteristic \(0\) or \(p\), and that
\(p\notin\mathfrak m^2R_{\mathfrak m}\) for every maximal ideal
\(\mathfrak m\) of \(R\) containing \(p\). Then the conclusion of
\textup{(i)} holds.
\end{enumerate}
\end{corollary}

\begin{proof}
\textup{(i)} By \cref{sub:moore-primary-decomposition}, we may assume
that \(N=\ell^\nu\) for a prime number \(\ell\in R^\times\) and
\(\nu\geq1\). Let \(x\in K_n(R;\Z/N)\) have zero image in
\(K_n(F;\Z/N)\), and let \((B_i)_{i\in I}\) be as in
\cref{lem:semilocal-regular-approximation}. By
\cref{sub:k-finite-continuity} and
\cref{sub:principal-open-witness}\textup{(ii)} for \(B=\Z\), there
are \(i\in I\) and \(x_i\in K_n(B_i;\Z/N)\) with image \(x\) in
\(K_n(R;\Z/N)\) and zero image in \(K_n(\Frac B_i;\Z/N)\). By
\cref{lem:semilocal-regular-approximation}, the ring
\(B_i\) is a regular semilocal domain essentially of finite type over
a field or over \(\Z\), and its closed residue fields have at most
one positive characteristic. In the notation of
\cref{sub:localized-filtered-presentation}, \(B_i=S_i^{-1}P_i\), and
the image of \(\ell\) in \(P_i\) lies in no prime ideal
\(\mathfrak p_{i,j}\), because its image in \(R\) lies in no maximal
ideal; hence \(\ell\in B_i^\times\). By
\cref{thm:finite-type-semilocal-gersten-injectivity}, \(x_i=0\).
Hence \(x=0\).

\smallskip\noindent
\textup{(ii)} By \cref{sub:unramified-regular-over-zp}, the ring map
\(\Z_{(p)}\to R\) is regular, and \(\Z_{(p)}\) is a regular ring
essentially of finite type over \(\Z\). Apply \textup{(i)} to
\(A=\Z_{(p)}\).
\end{proof}

\begin{remark}[Semilocal rings containing a field]
\label{rem:semilocal-equicharacteristic-classical}
Let \(R\) be a noetherian regular semilocal domain containing a
field, with prime field \(k\). The ring map \(k\to R\) is regular, as
in the proof of \cref{lem:finite-type-witness}, so
\cref{cor:semilocal-regular-model-k-injection}\textup{(i)} applies to
\(A=k\). In this case the conclusion is classical, and it holds for
every \(N\geq2\): the descent in the proof of
\cref{cor:semilocal-regular-model-k-injection}\textup{(i)} reduces it
to the rings \(B_i\) of
\cref{lem:semilocal-regular-approximation}, which are
regular semilocal domains essentially of finite type over \(k\). For
these rings the Gersten sequence is pure exact
\cite[Chapter~V, Corollary~9.6.4]{Weibel}, and the first step of the
proof of \cref{lem:equicharacteristic-injection} gives the
injectivity of \(K_n(B_i;\Z/N)\to K_n(\Frac B_i;\Z/N)\).
\end{remark}


\begin{thebibliography}{EHKSY20}
\raggedright

\bibitem[AOV08]{AOV08}
Dan Abramovich, Martin Olsson, and Angelo Vistoli, \emph{Tame stacks
in positive characteristic}, Annales de l'Institut Fourier
\textbf{58} (2008), no.~4, 1057--1091.
\href{https://doi.org/10.5802/aif.2378}{doi:10.5802/aif.2378}.
Corrigendum, ibid.\ \textbf{64} (2014), no.~3, 945--946.
\href{https://doi.org/10.5802/aif.2869}{doi:10.5802/aif.2869}.
\href{https://arxiv.org/html/math/0703310v1}{arXiv:math/0703310v1}
(Section~2 is numbered differently there; we use the numbering of the
published version).

\bibitem[Ane25]{AnelPiFinite}
Mathieu Anel, \emph{The category of $\pi$-finite spaces},
Journal of Pure and Applied Algebra \textbf{229} (2025), no.~6,
Article~107925.
\href{https://doi.org/10.1016/j.jpaa.2025.107925}{doi:10.1016/j.jpaa.2025.107925}.
\href{https://arxiv.org/html/2107.02082v4}{arXiv:2107.02082v4}
(we use the numbering of the published version).

\bibitem[Bac25]{BachmannEffective}
Tom Bachmann, \emph{The very effective covers of $KO$ and $KGL$
over Dedekind schemes}, Journal of the European Mathematical
Society \textbf{27} (2025), no.~11, 4705--4712.
\href{https://doi.org/10.4171/JEMS/1458}{doi:10.4171/JEMS/1458}.
\href{https://arxiv.org/html/2201.02786v2}{arXiv:2201.02786v2}.

\bibitem[BEM25]{BEMMotivic25}
Tom Bachmann, Elden Elmanto, and Matthew Morrow,
\emph{$\mathbb A^1$-invariant motivic cohomology of schemes},
preprint, 2025.
\href{https://arxiv.org/html/2508.09915v1}{arXiv:2508.09915v1}
(statement numbers as in this version).

\bibitem[Ber17]{Bergh}
Daniel Bergh, \emph{Functorial destackification of tame stacks with
abelian stabilisers}, Compos. Math. \textbf{153} (2017), no.~6, 1257--1315.
\href{https://doi.org/10.1112/S0010437X17007084}{doi:10.1112/S0010437X17007084}.
\href{https://arxiv.org/html/1409.5713v1}{arXiv:1409.5713v1}
(statement numbers and Algorithm~E as in this version).

\bibitem[BR19]{BerghRydh}
Daniel Bergh and David Rydh, \emph{Functorial destackification and
weak factorization of orbifolds}, preprint, 2019.
\href{https://arxiv.org/html/1905.00872v1}{arXiv:1905.00872v1}.

\bibitem[BM21]{BhattMathew21}
Bhargav Bhatt and Akhil Mathew, \emph{The arc-topology},
Duke Mathematical Journal \textbf{170} (2021), 1899--1988.
\href{https://doi.org/10.1215/00127094-2020-0088}{doi:10.1215/00127094-2020-0088}.
\href{https://arxiv.org/html/1807.04725v3}{arXiv:1807.04725v3}.

\bibitem[BS15]{BhattScholze}
Bhargav Bhatt and Peter Scholze, \emph{The pro-\'etale topology for
schemes}, Ast\'erisque \textbf{369} (2015), 99--201.
\href{https://doi.org/10.24033/ast.960}{doi:10.24033/ast.960}.
\href{https://www.math.uni-bonn.de/people/scholze/proetale.pdf}{Author's copy}.

\bibitem[BO74]{BlochOgus}
Spencer Bloch and Arthur Ogus, \emph{Gersten's conjecture and the
homology of schemes}, Annales scientifiques de l'\'Ecole Normale
Sup\'erieure (4) \textbf{7} (1974), no.~2, 181--201.
\href{https://doi.org/10.24033/asens.1266}{doi:10.24033/asens.1266}.
\href{http://www.numdam.org/item/ASENS_1974_4_7_2_181_0/}{Numdam}.

\bibitem[Bou24]{Bouis}
Tess Bouis, \emph{Motivic cohomology of mixed characteristic
schemes}, preprint, 2024, revised 2026.
\href{https://arxiv.org/html/2412.06635v3}{arXiv:2412.06635v3}.

\bibitem[Bou25]{BouisWeibel25}
Tess Bouis, \emph{Weibel vanishing and the projective bundle formula
for mixed characteristic motivic cohomology}, preprint, 2025.
\href{https://arxiv.org/html/2507.16501v1}{arXiv:2507.16501v1}.

\bibitem[BK25]{BouisKundu25}
Tess Bouis and Arnab Kundu, \emph{Beilinson--Lichtenbaum phenomenon
for motivic cohomology}, preprint, 2025, revised 2026.
\href{https://arxiv.org/html/2506.09910v2}{arXiv:2506.09910v2}.

\bibitem[CMM21]{CMM21}
Dustin Clausen, Akhil Mathew, and Matthew Morrow,
\emph{$K$-theory and topological cyclic homology of henselian pairs},
Journal of the American Mathematical Society \textbf{34} (2021),
no.~2, 411--473.
\href{https://doi.org/10.1090/jams/961}{doi:10.1090/jams/961}.
\href{https://arxiv.org/html/1803.10897v2}{arXiv:1803.10897v2}.

\bibitem[Dah18]{Dahlhausen18}
Christian Dahlhausen, \emph{Milnor $K$-theory of complete discrete
valuation rings with finite residue fields}, Journal of Pure and
Applied Algebra \textbf{222} (2018), no.~6, 1355--1371.
\href{https://doi.org/10.1016/j.jpaa.2017.07.002}{doi:10.1016/j.jpaa.2017.07.002}.
\href{https://arxiv.org/html/1509.01087v3}{arXiv:1509.01087v3}
(statement numbers as in this version).

\bibitem[DF21]{DemarcheFlorence}
Cyril Demarche and Mathieu Florence,
\emph{Splitting families in Galois cohomology}, Annales
scientifiques de l'\'Ecole Normale Sup\'erieure (4) \textbf{54}
(2021), no.~3, 779--792.
\href{https://doi.org/10.24033/asens.2470}{doi:10.24033/asens.2470}.
\href{https://arxiv.org/html/1711.06585v2}{arXiv:1711.06585v2}
(statement numbers as in this version).

\bibitem[ER21]{EdidinRydh}
Dan Edidin and David Rydh, \emph{Canonical reduction of stabilizers
for Artin stacks with good moduli spaces}, Duke Math. J.
\textbf{170} (2021), no.~5, 827--880.
\href{https://doi.org/10.1215/00127094-2020-0050}{doi:10.1215/00127094-2020-0050}.
\href{https://arxiv.org/abs/1710.03220v3}{arXiv:1710.03220v3}
(statement numbers as in this version).

\bibitem[EHIK21]{EHIK}
Elden Elmanto, Marc Hoyois, Ryomei Iwasa, and Shane Kelly,
\emph{Cdh descent, cdarc descent, and Milnor excision},
Mathematische Annalen \textbf{379} (2021), 1011--1045.
\href{https://doi.org/10.1007/s00208-020-02083-5}{doi:10.1007/s00208-020-02083-5}.
\href{https://arxiv.org/html/2002.11647v3}{arXiv:2002.11647v3}.

\bibitem[EHKSY20]{EHK}
Elden Elmanto, Marc Hoyois, Adeel A. Khan, Vladimir Sosnilo, and
Maria Yakerson, \emph{Modules over algebraic cobordism},
Forum of Mathematics, Pi \textbf{8} (2020), e14.
\href{https://doi.org/10.1017/fmp.2020.13}{doi:10.1017/fmp.2020.13}.
\href{https://arxiv.org/html/1908.02162v2}{arXiv:1908.02162v2}.

\bibitem[EM23]{ElmantoMorrow}
Elden Elmanto and Matthew Morrow, \emph{Motivic cohomology of
equicharacteristic schemes}, preprint, 2023, revised 2026.
\href{https://arxiv.org/html/2309.08463v2}{arXiv:2309.08463v2}.

\bibitem[Fel26a]{FeldCounterexample}
Niels Feld,
\emph{Finite-coefficient Gersten injectivity fails in ramified
mixed characteristic}, preprint, 2026.
\href{https://arxiv.org/html/2608.05005v2}{arXiv:2608.05005v2}.

\bibitem[Fel26b]{FeldHenselian}
Niels Feld,
\emph{Finite-coefficient $K$-theory of henselian valued fields and
Gersten injectivity}, preprint, 2026.
\href{https://arxiv.org/html/2608.18789v2}{arXiv:2608.18789v2}.

\bibitem[Fel26c]{FeldSurface}
Niels Feld, \emph{Gersten exactness with invertible coefficients
for excellent regular local surfaces}, preprint, 19 November 2026.
\href{https://hexagonmath.org/2611.00625v1}{Hexagon:2611.00625v1}.

\bibitem[Gab92]{Gabber92}
Ofer Gabber, \emph{$K$-theory of Henselian local rings and Henselian
pairs}, in \emph{Algebraic $K$-theory, commutative algebra, and
algebraic geometry (Santa Margherita Ligure, 1989)}, Contemporary
Mathematics \textbf{126}, American Mathematical Society, Providence,
RI, 1992, 59--70.
\href{https://doi.org/10.1090/conm/126/00509}{doi:10.1090/conm/126/00509}.

\bibitem[GT25]{GaidaiTurlov}
Ivan Gaidai-Turlov, \emph{The Bloch--Ogus theorem over DVR},
preprint, 2025, revised 2026.
\href{https://arxiv.org/abs/2506.06892v3}{arXiv:2506.06892v3}
(statement numbers as in this version).

\bibitem[GLL13]{GLLIndex13}
Ofer Gabber, Qing Liu, and Dino Lorenzini,
\emph{The index of an algebraic variety},
Inventiones Mathematicae \textbf{192} (2013), 567--626.
\href{https://doi.org/10.1007/s00222-012-0418-z}{doi:10.1007/s00222-012-0418-z}.
\href{https://arxiv.org/html/1209.2828v1}{arXiv:1209.2828v1}.

\bibitem[GLL15]{GLLHypersurfaces15}
Ofer Gabber, Qing Liu, and Dino Lorenzini,
\emph{Hypersurfaces in projective schemes and a moving lemma},
Duke Mathematical Journal \textbf{164} (2015), no.~7, 1187--1270.
\href{https://doi.org/10.1215/00127094-2877293}{doi:10.1215/00127094-2877293}.
\href{https://arxiv.org/html/1404.5366v1}{arXiv:1404.5366v1}.

\bibitem[Gei04]{Geisser04}
Thomas Geisser, \emph{Motivic cohomology over Dedekind rings},
Mathematische Zeitschrift \textbf{248} (2004), 773--794.
\href{https://doi.org/10.1007/s00209-004-0680-x}{doi:10.1007/s00209-004-0680-x}.
\href{https://www2.rikkyo.ac.jp/web/geisser/Dedekind.pdf}{Author's copy}.

\bibitem[GL00]{GeisserLevine}
Thomas Geisser and Marc Levine,
\emph{The $K$-theory of fields in characteristic $p$},
Inventiones Mathematicae \textbf{139} (2000), no.~3, 459--493.
\href{https://doi.org/10.1007/s002220050014}{doi:10.1007/s002220050014}.
\href{https://www2.rikkyo.ac.jp/web/geisser/p-part.pdf}{Author's copy}.

\bibitem[Gil86]{Gillet}
Henri Gillet,
\emph{Gersten's conjecture for the $K$-theory with torsion
coefficients of a discrete valuation ring},
Journal of Algebra \textbf{103} (1986), no.~1, 377--380.
\href{https://doi.org/10.1016/0021-8693(86)90193-6}
{doi:10.1016/0021-8693(86)90193-6}.

\bibitem[GL87]{GilletLevine}
Henri Gillet and Marc Levine,
\emph{The relative form of Gersten's conjecture over a discrete
valuation ring: the smooth case},
Journal of Pure and Applied Algebra \textbf{46} (1987),
no.~1, 59--71.
\href{https://doi.org/10.1016/0022-4049(87)90043-0}
{doi:10.1016/0022-4049(87)90043-0}.

\bibitem[HW19]{HaesemeyerWeibel19}
Christian Haesemeyer and Charles A. Weibel,
\emph{The norm residue theorem in motivic cohomology},
Annals of Mathematics Studies \textbf{200}, Princeton University
Press, Princeton, NJ, 2019.
\href{https://doi.org/10.1515/9780691189635}{doi:10.1515/9780691189635}.

\bibitem[ILO14]{ILO14}
Luc Illusie, Yves Laszlo, and Fabrice Orgogozo (eds.),
\emph{Travaux de Gabber sur l'uniformisation locale et la
cohomologie \'etale des sch\'emas quasi-excellents},
Ast\'erisque \textbf{363--364} (2014).
\href{https://doi.org/10.24033/ast.935}{doi:10.24033/ast.935}.
\href{https://arxiv.org/abs/1207.3648v1}{arXiv:1207.3648v1}
(Expos\'e~X, pp.~148--152 in the preprint; the page references in the
text use the \href{https://www.numdam.org/item/AST_2014__363-364__167_0.pdf}{printed pagination}).

\bibitem[Jan88]{Jannsen}
Uwe Jannsen, \emph{Continuous \'etale cohomology}, Mathematische
Annalen \textbf{280} (1988), 207--245.
\href{https://doi.org/10.1007/BF01456052}{doi:10.1007/BF01456052}.
\href{https://epub.uni-regensburg.de/26684/1/jannsen12.pdf}{Author's copy}.

\bibitem[Ker10]{Kerz10}
Moritz Kerz, \emph{Milnor $K$-theory of local rings with finite
residue fields}, Journal of Algebraic Geometry \textbf{19} (2010),
173--191.
\href{https://doi.org/10.1090/S1056-3911-09-00514-1}{doi:10.1090/S1056-3911-09-00514-1}.

\bibitem[Lev01]{LevineLocalization}
Marc Levine, \emph{Techniques of localization in the theory of
algebraic cycles}, Journal of Algebraic Geometry \textbf{10}
(2001), no.~2, 299--363.
\href{https://citeseerx.ist.psu.edu/document?doi=9622ccdf2df87d21c760a444b7b04cd40b1662b0\&repid=rep1\&type=pdf}{Author manuscript}.

\bibitem[L\"ud24]{LudersMilnor24}
Morten L\"uders, \emph{On the relative Gersten conjecture for
Milnor $K$-theory in the smooth case}, Journal of Pure and Applied
Algebra \textbf{228} (2024), no.~11, Article 107718.
\href{https://doi.org/10.1016/j.jpaa.2024.107718}{doi:10.1016/j.jpaa.2024.107718}.
\href{https://arxiv.org/html/2010.02622v4}{arXiv:2010.02622v4}
(statement numbers as in this version).

\bibitem[L\"ud26]{Luders22}
Morten L\"uders, \emph{Bloch--Ogus theory for smooth and semi-stable
schemes in mixed characteristic}, \'Epijournal de G\'eom\'etrie
Alg\'ebrique \textbf{10} (2026), Article~19.
\href{https://doi.org/10.46298/epiga.2026.10021}{doi:10.46298/epiga.2026.10021}.
\href{https://arxiv.org/html/2201.11553v4}{arXiv:2201.11553v4}.

\bibitem[Lur09]{HTT}
Jacob Lurie, \emph{Higher Topos Theory}, Annals of Mathematics Studies,
vol.~170, Princeton University Press, 2009;
\href{https://arxiv.org/abs/math/0608040}{arXiv:math/0608040};
\href{https://www.math.ias.edu/~lurie/papers/HTT.pdf}{author's updated text}.

\bibitem[Lur17]{HA}
Jacob Lurie, \emph{Higher algebra}, version of 18 September 2017.
\url{https://www.math.ias.edu/~lurie/papers/HA.pdf}.

\bibitem[Lur]{Kerodon}
Jacob Lurie, \emph{Kerodon}. \url{https://kerodon.net/}.

\bibitem[NSS15]{NikolausSchreiberStevenson}
Thomas Nikolaus, Urs Schreiber, and Danny Stevenson,
\emph{Principal $\infty$-bundles---general theory},
Journal of Homotopy and Related Structures \textbf{10} (2015),
no.~4, 749--801.
\href{https://doi.org/10.1007/s40062-014-0083-6}{doi:10.1007/s40062-014-0083-6}.
\href{https://arxiv.org/html/1207.0248v3}{arXiv:1207.0248v3}.

\bibitem[OAI26a]{OpenAIGersten3}
OpenAI, \emph{An integral degree-three Gersten counterexample},
OpenAI Math Release preprint, 26 September 2026.
\href{https://github.com/openai/math/tree/main/preprints/An-Integral-Degree-Three-Gersten-Counterexample-September-26-2026}{Preprint and source}.

\bibitem[OAI26b]{OpenAIGersten5}
OpenAI, \emph{An integral counterexample to Gersten's conjecture},
OpenAI Math Release preprint, 25 September 2026.
\href{https://github.com/openai/math/tree/main/preprints/An-Integral-Counterexample-to-Gerstens-Conjecture-September-25-2026}{Preprint and source}.

\bibitem[Pan03]{Panin}
Ivan A. Panin,
\emph{The equicharacteristic case of the Gersten conjecture},
Proceedings of the Steklov Institute of Mathematics
\textbf{241} (2003), 154--163.
\href{https://www.mathnet.ru/eng/tm394}{Journal entry}.

\bibitem[PS25]{PaninStavrova25}
Ivan Panin and Anastasia Stavrova,
\emph{On the Gille theorem for the relative projective line},
preprint, 2023, revised 2025.
\href{https://arxiv.org/html/2305.16627v3}{arXiv:2305.16627v3}.

\bibitem[PZ03]{PaninZainoulline}
Ivan Panin and Kirill Zainoulline, \emph{Variations on the
Bloch--Ogus theorem}, Documenta Mathematica \textbf{8} (2003),
51--67.
\href{https://doi.org/10.4171/DM/137}{doi:10.4171/DM/137}.
\href{https://arxiv.org/abs/math/0203128}{arXiv:math/0203128}.

\bibitem[Pop18]{Popescu19arXiv}
Dorin Popescu, \emph{On a question of Swan}, with an appendix by
K\k{e}stutis \v{C}esnavi\v{c}ius, preprint version.
\href{https://arxiv.org/pdf/1803.06956v6}{arXiv:1803.06956v6}.

\bibitem[Pop19]{Popescu19}
Dorin Popescu, \emph{On a question of Swan}, with an appendix by
K\k{e}stutis \v{C}esnavi\v{c}ius,
Algebraic Geometry \textbf{6} (2019), no.~6, 716--729.
\href{https://doi.org/10.14231/AG-2019-030}{doi:10.14231/AG-2019-030}.

\bibitem[Qui73]{Quillen}
Daniel Quillen,
\emph{Higher algebraic $K$-theory. I},
in \emph{Algebraic $K$-theory I: Higher $K$-theories},
Lecture Notes in Mathematics \textbf{341}, Springer, Berlin, 1973,
85--147.
\href{https://doi.org/10.1007/BFb0067053}{doi:10.1007/BFb0067053}.

\bibitem[Ryd13]{RydhGeometricQuotients}
David Rydh, \emph{Existence and properties of geometric quotients},
J. Algebraic Geom. \textbf{22} (2013), no.~4, 629--669.
\href{https://doi.org/10.1090/S1056-3911-2013-00615-3}{doi:10.1090/S1056-3911-2013-00615-3}.
\href{https://arxiv.org/html/0708.3333v2}{arXiv:0708.3333v2}.

\bibitem[Sak20]{Sakagaito20}
Makoto Sakagaito, \emph{A note on Gersten's conjecture for \'etale
cohomology over two-dimensional henselian regular local rings},
Comptes Rendus. Math\'ematique \textbf{358} (2020), no.~1, 33--39.
\href{https://doi.org/10.5802/crmath.9}{doi:10.5802/crmath.9}.
\href{https://arxiv.org/abs/1903.01677v2}{arXiv:1903.01677v2}
(statement numbers in the text refer to the published version).

\bibitem[SS23]{SchmidtStrunk23}
Johannes Schmidt and Florian Strunk, \emph{A Bloch--Ogus theorem
for henselian local rings in mixed characteristic}, Mathematische
Zeitschrift \textbf{303} (2023), Article~82.
\href{https://doi.org/10.1007/s00209-023-03223-8}{doi:10.1007/s00209-023-03223-8}.

\bibitem[Ser97]{SerreGC}
Jean-Pierre Serre, \emph{Galois cohomology}, translated from the
French by Patrick Ion, Springer-Verlag, Berlin, 1997.
\href{https://doi.org/10.1007/978-3-642-59141-9}{doi:10.1007/978-3-642-59141-9}.

\bibitem[Ska17]{Skalit}
Christopher Skalit,
\emph{Regular morphisms and Gersten's conjecture}, preprint, 2017.
\href{https://arxiv.org/html/1710.00303v1}{arXiv:1710.00303v1}.

\bibitem[Spi18]{Spitzweck}
Markus Spitzweck, \emph{A commutative $\mathbf P^1$-spectrum representing
motivic cohomology over Dedekind domains}, M\'emoires de la
Soci\'et\'e Math\'ematique de France \textbf{157} (2018).
\href{https://arxiv.org/html/1207.4078v3}{arXiv:1207.4078v3}
(statement numbers as in this version).

\bibitem[Sta26]{Stacks}
The Stacks Project Authors,
\emph{The Stacks Project}, 2026.
\url{https://stacks.math.columbia.edu}.

\bibitem[TT90]{TT}
Robert W. Thomason and Thomas Trobaugh,
\emph{Higher algebraic $K$-theory of schemes and of derived
categories}, in \emph{The Grothendieck Festschrift}, vol.~III,
Progress in Mathematics \textbf{88}, Birkh\"auser, Boston, 1990,
247--435.
DOI of the reprint:
\href{https://doi.org/10.1007/978-0-8176-4576-2_10}
{doi:10.1007/978-0-8176-4576-2\_10}.

\bibitem[Voe11]{Voevodsky11}
Vladimir Voevodsky, \emph{On motivic cohomology with
$\Z/l$-coefficients}, Annals of Mathematics (2) \textbf{174} (2011),
no.~1, 401--438.
\href{https://doi.org/10.4007/annals.2011.174.1.11}{doi:10.4007/annals.2011.174.1.11}.

\bibitem[Wei13]{Weibel}
Charles A. Weibel,
\emph{The $K$-book: An introduction to algebraic $K$-theory},
Graduate Studies in Mathematics \textbf{145},
American Mathematical Society, Providence, RI, 2013.
\href{https://doi.org/10.1090/gsm/145}{doi:10.1090/gsm/145}.
\href{https://sites.math.rutgers.edu/~weibel/Kbook/Kbook.pdf}{Author's copy}.

\end{thebibliography}
\end{document}